%% file: AAmain.tex
\documentclass[12pt,fleqn]{report}
\pdfoutput=1
\usepackage{amssymb,amsmath,uwthesis}
\usepackage{graphicx}
\usepackage{bbm}
\usepackage{bm} 
\usepackage[usenames]{color}
\usepackage{multicol}
\usepackage{calc} 
\usepackage{amscd} 
\newcommand{\bfit}[1]{\textbf{\textit{#1}}}

\def\definingA#1{\textbf{#1}\index{#1|bb}} 
\def\definingB#1#2{\textbf{#2}\index{#1|bb}} 
\newcommand{\defining}[2][]{\ifthenelse{\equal{#1}{}}{\definingA{#2}}{\definingB{#1}{#2}}}
\def\first#1{\textit{#1}\index{#1}} 
\def\mention#1{\textit{#1}\index{#1}} 
\def\gdefining#1{\glslink{#1}{\textbf{#1}}} 
\usepackage{makeidx} 
\usepackage[nottoc]{tocbibind}
\makeindex
\usepackage{hyperref}
\hypersetup{
    unicode=false,          
    pdftoolbar=true,        
    pdfmenubar=true,        
    pdffitwindow=false,     
    pdfstartview={FitH},    
    pdftitle={My title},    
    pdfauthor={E.A. Johnson},  
    pdfsubject={dissertation}, 
    pdfcreator={E.A. Johnson},   
    pdfproducer={E.A. Johnson}, 
    pdfkeywords={reconnection} {heat} {entropy}, 
    pdfnewwindow=true,      
    colorlinks=true,        
    linkcolor=blue,         
    citecolor=green,        
    filecolor=magenta,      
    urlcolor=red            
}
\usepackage[toc]{glossaries} 
\usepackage{glossaries}
\newglossary{vocab}{v.gls}{v.glo}{Vocabulary Glossary}
\newglossary{symbol}{s.gls}{s.glo}{Symbol Glossary}
\makeglossaries
\input{macros}

\input{glossary.tex}
\begin{document}
%
\title{Gaussian-Moment Relaxation Closures for Verifiable Numerical
  Simulation of Fast Magnetic Reconnection in Plasma}
\author{Evan Alexander Johnson}
\degree{Doctor of Philosophy}
\dept{Mathematics}
\thesistype{dissertation}
\beforepreface
\renewcommand{\baselinestretch}{1.5}  
\prefacesection{Abstract}
\input abs

\renewcommand{\baselinestretch}{1}  
\prefacesection{Acknowledgements}
\input ack
\listoffigures
\listoftables
%

\pagebreak
\afterpreface 
\setlength{\parskip}{1ex plus 0.5ex minus 0.2ex} 
 \include{chap1} 
 \include{chap2} 
 \include{chap3} 
 \include{chap4} 
 \include{chap5} 
 \include{chap6} 
 \include{chap7} 
 \include{chap8} 
\appendix       
\include{appendix2}

\include{appendix3}

\include{appendix4}
%
%
\bibliographystyle{siam}
\bibliography{thesis}
\newglossarystyle{twocolumn}{
    \glossarystyle{list}
    \renewcommand*{\glossarypreamble}{\begin{multicols*}{2}} 
    \renewcommand*{\glossarypostamble}{\end{multicols*}}
}
\printglossary[type=vocab, style=twocolumn]
\printglossary[type=symbol, style=twocolumn]
\phantomsection
\printindex
%

\end{document}

%% file: macros.tex
\def\magenta{\color{magenta}}

\def\black{\color{black}}
\definecolor{MyDarkBlue}{rgb}{0,0.10,0.60}
\def\dblue{\color{MyDarkBlue}}
\newcommand{\inmag}[1] {{\color{magenta} #1}}
\newcommand{\inred}[1] {{\color{red} #1}}

\newcommand{\indarkblue}[1] {{\color{MyDarkBlue} #1}}

\def\fhrule{\vspace{12pt}\hrule}
\def\hfrule{\hrule\vspace{8pt}}
\def\hfrule{\hrule\vspace{12pt}}
\usepackage{accents}
\usepackage{mathdots}
\DeclareMathOperator{\sgn}{sgn}
\DeclareMathOperator{\minmod}{minmod} 
\DeclareMathOperator{\sech}{sech}
\DeclareMathOperator{\Span}{span} 
\DeclareMathOperator{\adj}{adj} 
\DeclareMathOperator{\Ant}{Ant} 
\DeclareMathOperator{\Sym}{Sym} 
\DeclareMathOperator{\SymB}{Sym2} 
\DeclareMathOperator{\SymC}{Sym3} 
\DeclareMathOperator{\id}{{\bm{\mathbbm{1}}}} 
\def\intv{\int_\v}
\def\intc{\int_\c}
\usepackage{stmaryrd}
\def\mypara{\sslash} 
\newcommand{\deviator}[1]{#1^\circ}
\def\topdeviator{\deviator}
\newcommand\tensorb[1]{\underline{\underline{{#1}}}}
\newcommand\tensorc[1]{\underline{\underline{\underline{{#1}}}}}

\renewcommand\sup[1]{\ensuremath{{}^{(#1)}}}
\def\Zpara{\parallel}
\def\Maxwell{\mathcal{M}}
\def\fMaxwell{f_\theta}
\def\Gauss{\mathcal{G}}
\def\fGauss{f_\Theta}
\def\curl{\nabla\times}
\def\CCCC{\underline{\underline{\underline{\underline{C}}}}}

\def\Zdelta{\delta}
\def\Mcoefinv{\tviscosity}
\def\Mcoef{\Mcoefinv\inv}
\def\Gcoef{C^\ZRT}
\def\dTTinv{\deviator{(\TT^{-1})}}
\def\dPshapeinv{\deviator{(\Pshape^{-1})}}
\def\Zdotp{{\,\boldsymbol\cdot\,}}
\def\Zddotp{{\,\boldsymbol:\,}}
\newcommand*{\Zdddotp}{%
  \mathrel{\vcenter{\offinterlineskip 
  \hbox{$\boldsymbol\cdot$}\vskip-.35ex\hbox{$\boldsymbol\cdot$}\vskip-.35ex\hbox{$\boldsymbol\cdot$}}}}
\def\fd{\bar f}
\def\Zstimes{\widetilde\otimes}
\def\Zsvee{\widetilde{\vee}}
\def\Zsveebar{\widetilde{\veebar}}
\def\ZsSym{\widetilde{\Sym}}
\def\inv{^{-1}}
\def\mean#1{\langle #1 \rangle}
\DeclareMathOperator*{\Prandtl}{\ensuremath{Pr}} 
\DeclareMathOperator*{\error}{\ensuremath{error}} 
\newcommand\note[1]{\emph{[#1]}}
\def\QED{Q.E.D.}
\def\e{\mathrm{e}}
\def\i{\mathrm{i}}
\def\ZvA{v_A}
\def\ZC{C}
\def\ZD{\nabla}
\def\ZDiv{\nabla\Zdotp}
\def\ZET{\mathbb{E}}
\def\ZE{\mathbf{E}}
\def\ZJ{\mathbf{J}}
\def\ZM{\mathbf{M}}
\def\ZNrg{\mathcal{E}}
\def\ZNrgTot{\mathcal{\widetilde E}}
\def\ZO{\mathcal{O}}
\def\ZPTd{\ZPT^d}
\def\ZPT{\mathbb{P}}
\def\ZPshape{{\bm{\pi}}}
\def\ZQt{Q^t}
\def\ZQTt{\mathbb{Q}^t}
\def\ZQf{Q^f}
\def\ZQTf{\mathbb{Q}^f}
\def\ZQ{Q}
\def\ZQT{\mathbb{Q}}
\def\ZRT{{\mathbb{R}}}
\def\ZR{\mathbf{R}}
\def\ZTTinv{\ZTT^{-1}}
\def\ZTT{\mathbb{T}}
\def\ZThetaT{\Theta}
\def\ZTinv{T^{-1}}
\def\Za{\mathbf{a}}
\def\ZdPT{\deviator{\mathbb{P}}}
\def\ZdPshape{{\topdeviator{\ZPshape}}}
\def\ZdTT{\deviator{\ZTT}}
\def\Zdmt{\mathrm{d}\widetilde m}
\def\Zdm{\mathrm{d}m}
\def\ZdqT{\deviator{\ZqT}}
\def\Zdstress{\deviator{\Zstress}}

\def\ZeT{\mathbbm{e}}
\def\Zebas{\mathbf{e}}

\def\Zgyrofreq{\omega_c}

\def\ZheatConductivity{k}
\def\ZtheatConductivity{\widetilde{\bm{k}}}
\def\ZTheatConductivity{\bm{k}}
\def\ZtHeatConductivity{\widetilde{\bm{K}}}
\def\ZTHeatConductivity{\bm{K}}
\def\Zidfour{\id\diamond\id}

\def\Zlaplacian{\nabla^2}
\def\Zme{m_\e}
\def\Zmi{m_\i}
\def\Zmred{\check{\mu}}
\def\Zmte{\mt_\e}
\def\Zmti{\mt_\i}

\def\Zmt{\widetilde m}
\def\Zmut{\widetilde{\mred}}
\def\Znrg{\mathrm{e}}
\def\ZpTot{\widetilde p}
\def\Zpv{{\widetilde\v}}
\def\p{\mathrm{p}}
\def\ZqT{\bm{\mathbbm{q}}}
\def\Zq{\mathbf{q}}
\def\Zreals{\mathbb{R}}
\def\Zresistivity{\eta}
\def\Ztresistivity{\widetilde{\bm{\eta}}}
\def\ZTresistivity{\bm{\eta}}
\def\Zstrain{\bm{e}}
\def\Zstress{\bm{\sigma}}
\def\s{\mathrm{s}}

\def\Zu{\mathbf{u}}
\def\Zviscosity{\mu}
\def\Ztviscosity{\bm{\widetilde\mu}}
\def\ZTviscosity{\bm{\mu}}
\def\Zv{\mathbf{v}}
\def\Zw{\mathbf{w}}
\def\Zxb{\mathbf{x}}
\def\ZB{\mathbf{B}}
\def\Zbhat{\mathbf{b}}
\def\Zc{\mathbf{c}}
\def\Zdstrain{\topdeviator{\Zstrain}}
\def\Zptime{\tau}
\def\Zhtime{\widetilde\tau}
\def\Zpgf{\varpi}
\def\Zhgf{\widetilde\varpi}
\def\Zqdens{\sigma}
\def\Zmdens{\rho}
\def\Zndens{n}
\def\ZdebyeLength{\lambda_D}
\def\Zdcs{\overline{\delta}}
\def\ZtC{C^T}
\def\ZC{C}
\def\Zf{f}
\def\ZK{K}
\def\b{\bhat}
\def\B        {\glshyperlink[{\black \ZB          }]{B}}

\def\D        {\glshyperlink[{\black \ZD          }]{D}}
\def\Div      {\glshyperlink[{\black \ZDiv        }]{Div}}
\def\E        {\glshyperlink[{\black \ZE          }]{E}}
\def\ET       {\glshyperlink[{\black \ZET         }]{ET}}
\def\J        {\glshyperlink[{\black \ZJ          }]{J}}
\def\M        {\glshyperlink[{\black \ZM          }]{M}}
\def\Nrg      {\glshyperlink[{\black \ZNrg        }]{Nrg}}
\def\NrgTot   {\glshyperlink[{\black \ZNrgTot     }]{NrgTot}}
\def\O{\ZO}
\def\PT       {\glshyperlink[{\black \ZPT         }]{PT}}
\def\PTd      {\glshyperlink[{\black \ZPTd        }]{PTd}}
\def\Pshape   {\glshyperlink[{\black \ZPshape     }]{Pshape}}
\def\Q        {\glshyperlink[{\black \ZQ          }]{Q}}
\def\Qf       {\glshyperlink[{\black \ZQf         }]{Qf}}
\def\Qt       {\glshyperlink[{\black \ZQt         }]{Qt}}
\def\QT       {\glshyperlink[{\black \ZQT         }]{QT}}
\def\QTf      {\glshyperlink[{\black \ZQTf        }]{QTf}}
\def\QTt      {\glshyperlink[{\black \ZQTt        }]{QTt}}
\def\R        {\glshyperlink[{\black \ZR          }]{R}}
\def\RT       {\glshyperlink[{\black \ZRT         }]{RT}}
\def\TT       {\glshyperlink[{\black \ZTT         }]{TT}}
\def\TTinv    {\glshyperlink[{\black \ZTTinv      }]{TTinv}}
\def\ThetaT   {\glshyperlink[{\black \ZThetaT     }]{ThetaT}}
\def\Tinv     {\glshyperlink[{\black \ZTinv       }]{Tinv}}
\def\a        {\glshyperlink[{\black \Za          }]{a}}
\def\bhat     {\glshyperlink[{\black \Zbhat       }]{bhat}}
\def\c        {\glshyperlink[{\black \Zc          }]{c}}
\def\dPT      {\glshyperlink[{\black \ZdPT        }]{dPT}}
\def\dPshape  {\glshyperlink[{\black \ZdPshape    }]{dPshape}}
\def\dTT      {\glshyperlink[{\black \ZdTT        }]{dTT}}
\def\dm       {\glshyperlink[{\black \Zdm         }]{dm}}
\def\dmt      {\glshyperlink[{\black \Zdmt        }]{dmt}}
\def\dqT      {\glshyperlink[{\black \ZdqT        }]{dqT}}
\def\dstress  {\glshyperlink[{\black \Zdstress    }]{dstress}}
\def\dt       {\glshyperlink[{\black \Zdt         }]{dt}}
\def\eT       {\glshyperlink[{\black \ZeT         }]{eT}}
\def\ebas     {\glshyperlink[{\black \Zebas       }]{ebas}}

\def\gyrofreq {\glshyperlink[{\black \Zgyrofreq   }]{gyrofreq}}
\def\idfour   {\glshyperlink[{\black \Zidfour     }]{idfour}}
\def\idpara   {\glshyperlink[{\black \Zidpara     }]{idpara}}
\def\idperp   {\glshyperlink[{\black \Zidperp     }]{idperp}}
\def\idskew   {\glshyperlink[{\black \Zidskew     }]{idskew}}
\def\idtens   {\glshyperlink[{\black \Zidtens     }]{idtens}}
\def\laplacian{\glshyperlink[{\black \Zlaplacian  }]{laplacian}}
\def\me       {\glshyperlink[{\black \Zme         }]{me}}
\def\mi       {\glshyperlink[{\black \Zmi         }]{mi}}
\def\mred     {\glshyperlink[{\black \Zmred       }]{mred}}
\def\mt       {\glshyperlink[{\black \Zmt         }]{mt}}
\def\mte      {\glshyperlink[{\black \Zmte        }]{mte}}
\def\mti      {\glshyperlink[{\black \Zmti        }]{mti}}

\def\mut      {\glshyperlink[{\black \Zmut        }]{mut}}
\def\nrg      {\glshyperlink[{\black \Znrg        }]{nrg}}
\def\para    {{\glshyperlink[{\black \Zpara       }]{para}}}
\def\pTot     {\glshyperlink[{\black \ZpTot       }]{pTot}}
\def\pv      {{\glshyperlink[{\black \Zpv         }]{pv}}}
\def\q        {\glshyperlink[{\black \Zq          }]{q}}
\def\qT       {\glshyperlink[{\black \ZqT         }]{qT}}
\def\reals    {\glshyperlink[{\black \Zreals      }]{reals}}
\def\strain   {\glshyperlink[{\black \Zstrain     }]{strain}}
\def\stress   {\glshyperlink[{\black \Zstress     }]{stress}}
\def\tr       {\glshyperlink[{\black \Ztr         }]{tr}}
\def\u       {{\glshyperlink[{\black \Zu          }]{u}}}
\def\uc      {{\glshyperlink[{\black \Zuc         }]{uc}}}
\def\v       {{\glshyperlink[{\black \Zv          }]{v}}}

\def\w       {{\glshyperlink[{\black \Zw          }]{w}}}
\def\xb      {{\glshyperlink[{\black \Zxb         }]{xb}}}
\def\dstrain{\glshyperlink[{\black \Zdstrain         }]{dstrain}}
\def\viscosity        {\glshyperlink[{\black \Zviscosity        }]{viscosity}}
\def\resistivity      {\glshyperlink[{\black \Zresistivity      }]{resistivity}}
\def\heatConductivity {\glshyperlink[{\black \ZheatConductivity }]{heatConductivity}}
\def\theatConductivity{\glshyperlink[{\black \ZtheatConductivity }]{theatConductivity}}
\def\TheatConductivity{\glshyperlink[{\black \ZTheatConductivity }]{TheatConductivity}}
\def\tHeatConductivity{\glshyperlink[{\black \ZtHeatConductivity }]{tHeatConductivity}}
\def\THeatConductivity{\glshyperlink[{\black \ZTHeatConductivity }]{THeatConductivity}}
\def\tviscosity       {\glshyperlink[{\black \Ztviscosity        }]{tviscosity}}
\def\tresistivity     {\glshyperlink[{\black \Ztresistivity      }]{tresistivity}}
\def\theatConductivity{\glshyperlink[{\black \ZtheatConductivity }]{theatConductivity}}
\def\Tviscosity       {\glshyperlink[{\black \ZTviscosity        }]{Tviscosity}}
\def\Tresistivity     {\glshyperlink[{\black \ZTresistivity      }]{Tresistivity}}
\def\TheatConductivity{\glshyperlink[{\black \ZTheatConductivity }]{TheatConductivity}}
\def\stimes           {\glshyperlink[{\black \Zstimes            }]{stimes}}

\def\sveebar          {\glshyperlink[{\black \Zsveebar           }]{sveebar}}
\def\sSym             {\glshyperlink[{\black \ZsSym              }]{sSym}}
\def\dotp             {\glshyperlink[{\black \Zdotp              }]{dotp}}
\def\ddotp            {\glshyperlink[{\black \Zddotp             }]{ddotp}}
\def\dddotp           {\glshyperlink[{\black \Zdddotp            }]{dddotp}}
\def\pgf              {\glshyperlink[{\black \Zpgf               }]{pgf}}
\def\hgf              {\glshyperlink[{\black \Zhgf               }]{hgf}}
\def\ptime            {\glshyperlink[{\black \Zptime             }]{ptime}}
\def\htime            {\glshyperlink[{\black \Zhtime             }]{htime}}
\def\qdens            {\glshyperlink[{\black \Zqdens             }]{qdens}}
\def\mdens            {\glshyperlink[{\black \Zmdens             }]{mdens}}
\def\ndens            {\glshyperlink[{\black \Zndens             }]{ndens}}
\def\vA               {\glshyperlink[{\black \ZvA                }]{vA}}

\def\dcs              {\glshyperlink[{\black \Zdcs               }]{dcs}}
\def\tC               {\glshyperlink[{\black \ZtC                }]{tC}}
\def\C                {\glshyperlink[{\black \ZC                 }]{C}}
\def\f                {\glshyperlink[{\black \Zf                 }]{f}}
\def\Pr               {\glshyperlink[{\black \Prandtl            }]{Pr}}
\def\K                {\glshyperlink[{\black \ZK                 }]{K}}

\def\LCi              {\glshyperlink[{\dblue \ZC_\i              }]{C}}
\def\LCe              {\glshyperlink[{\dblue \ZC_\e              }]{C}}
\def\LCie             {\glshyperlink[{\magenta \ZC_{\i\e}        }]{C}}
\def\LCei             {\glshyperlink[{\magenta \ZC_{\e\i}        }]{C}}
\def\LDivdPT#1{\glslink{dPT}{\inred{\ZDiv\ZdPT_#1}}}
\def\LDivqT#1{\glslink{qT}{\inred{\ZDiv\ZqT_#1}}}
\def\LDivq#1{\glslink{q}{\inred{\ZDiv\Zq_#1}}}
\def\LdPT#1{\glslink{dPT}{\inred{\ZdPT_#1}}}
\def\LqT#1{\glslink{qT}{\inred{\ZqT_#1}}}
\def\LR#1{\glslink{R}{\inmag{\ZR_#1}}}
\def\LTresistivity{\glslink{Tresistivity}{\inmag{\ZTresistivity}}}
\def\LQ#1{\glslink{Q}{\inmag{\ZQ_#1}}}
\def\LQT#1{\glslink{Q}{\inmag{\ZQT_#1}}}
\def\LRT#1{\glslink{RT}{\indarkblue{\ZRT_#1}}}
\def\LQt#1{\glslink{Qt}{\inmag{\ZQt_#1}}}
\def\LQf#1{\glslink{Qf}{\inmag{\ZQf_#1}}}
\def\LQTt#1{\glslink{QTt}{\inmag{\ZQTt_#1}}}
\def\LQTf#1{\glslink{QTf}{\inmag{\ZQTf_#1}}}

\def\Ldelta   {\glshyperlink[{\black \Zdelta      }]{delta}}
\def\Ldiamond {\glshyperlink[{\black \diamond     }]{diamond}}

%% file: glossary.tex
\def\s{\mathrm{s}}

\newglossaryentry{Boltzmann equation}{type=vocab, name={Boltzmann equation},
  description={\index{Boltzmann equation}
    conservation of particle density in phase space}}
\newglossaryentry{kinetic equation}{type=vocab, name={kinetic equation},
  description={\index{kinetic equation} the \glshyperlink{Boltzmann equation}}}
\newglossaryentry{Gaussian-moment model}{type=vocab, name={Gaussian-moment model},
  description={\index{Gaussian-moment model} a model of a gas which evolves
    all quadratic velocity moments.  In three dimensions of space the
    Gaussian-moment model evolves ten monomial moments and is therefore
    also referred to as the \glshyperlink{ten-moment model}}}
\newglossaryentry{ten-moment model}{type=vocab, name={ten-moment model},
  description={\index{ten-moment model} the \glshyperlink{Gaussian-moment model}}}
\newglossaryentry{Maxwellian-moment model}{type=vocab, name={Maxwellian-moment model},
  description={\index{Maxwellian-moment model} a model of a gas which evolves
  the subquadratic monomial moments (mass and momentum)
  and the energy (a quadratic moment).
  In three dimensions of space the
  Maxwellian-moment model evolves five monomial moments and is therefore
  also referred to as the \glshyperlink{five-moment model}}}
\newglossaryentry{five-moment model}{type=vocab, name={five-moment model},
  description={\index{five-moment model} the \glshyperlink{Maxwellian-moment model}}}
\newglossaryentry{two-fluid Maxwell}{type=vocab, name={two-fluid Maxwell model},
  description={\index{two-fluid Maxwell model} A fluid model of a two-species
  plasma which uses separate fluid models for the positive and negative species
  and which solves the full system of Maxwell's equations for the electromagnetic field}}
\newglossaryentry{quasineutrality}{type=vocab, name={quasineutrality},
  description={\index{quasineutrality}
    the assumption that the net charge density is approximately zero}}
\newglossaryentry{collision operator}{type=vocab, name={collision operator},
  description={\index{collision operator}
    the error when the particle density is substituted into the Vlasov equation.
    That is, it is assumed that particle density satisfies an evolution equation
    of the form
    \begin{gather*}
     \partial_t f + \v\dotp\D_{\xb}f+\a \dotp\D_{\v} f = C,
    \end{gather*}
    where $C$ is the collision operator.
    See section \ref{ConventionsUsed}
    }}

\newglossaryentry{0}{type=symbol, name={$0$}, sort={!0},
  description={origin (used in the context of reconnection, usually the X-point)}}
\newglossaryentry{para}{type=symbol, name={$\para$}, sort={!para},
 description={as a subscript on a vector, denotes the component parallel
   to the magnetic field direction $\bhat$;
   as a superscript denotes the component parallel to the out-of-plane
   direction $\ebas^\para=\ebas_z$}}
\newglossaryentry{perp}{type=symbol, name={$\perp$}, sort={!perp},
 description={as a subscript on a vector, denotes the component perpendicular
   to the magnetic field direction $\bhat$;
   as a superscript denotes the component perpendicular to the out-of-plane
   direction}}
\newglossaryentry{D}{type=symbol, name={$\ZD$}, sort={!del},
 description={del operator}}
\newglossaryentry{Div}{type=symbol, name={$\ZDiv$}, sort={!div},
 description={divergence operator}}
\newglossaryentry{laplacian}{type=symbol, name={$\Zlaplacian$}, sort={!laplacian},
 description={Laplacian operator, $\sum_i \partial_{x_i}^2$}}

\newglossaryentry{ebas}{type=symbol, name={$\ebas_i$}, sort={e^bas},
 description={elementary basis vector aligned with axis $i$}}
\newglossaryentry{delta}{type=symbol, name={$\Zdelta_\s$}, sort={@delta.s},
 description={inertial length}}
\newglossaryentry{diamond}{type=symbol, name={$\diamond$}, sort={!diamond},
 description={diamond product of tensors, defined by
 $(A\diamond B)_{ijkl} = A_{ik}B_{lj}$, i.e.,
 $A\diamond B = A\stimes B^T$}}
\newglossaryentry{stimes}{type=symbol, name={$\Zstimes$}, sort={!splicetimes},
 description={splice product of tensors, defined 
 in section \ref{SpliceTensorOperators}}}
\newglossaryentry{svee}{type=symbol, name={$\Zsvee$}, sort={!splicevee},
 description={splice vee product of tensors, defined 
 in section \ref{SpliceTensorOperators}}}
\newglossaryentry{sveebar}{type=symbol, name={$\sveebar$}, sort={!spliceveebar},
 description={splice symmetric product of tensors, defined 
 in section \ref{SpliceTensorOperators}}}
\newglossaryentry{veebar}{type=symbol, name={$\veebar$}, sort={!veebar},
 description={symmetric product of tensors, defined 
 in section \ref{TensorProducts}}}
\newglossaryentry{vee}{type=symbol, name={$\vee$}, sort={!vee},
 description={vee product; rescales the symmetric product $\veebar$ so that
 the vee product of symmetric tensors is the sum over all distinguishable
 permutations of the indices of their tensor product; see
 section \ref{TensorProducts}}}
\newglossaryentry{sSym}{type=symbol, name={$\ZsSym$}, sort={Symsplice},
 description={splice symmetrization of a second-order tensor, defined 
 in section \ref{SpliceTensorOperators}}}
\newglossaryentry{ddotp}{type=symbol, name={$\ddotp$}, sort={!ddotp},
  description={double dot product of tensors,
  defined by contracting two adjacent indices, e.g.\
  $A\Zddotp B = A_{ij}B_{ij}$}}
\newglossaryentry{dddotp}{type=symbol, name={$\dddotp$}, sort={!dddotp},
  description={triple dot product of tensors, e.g.\
  $A\Zdddotp B = A_{ijk}B_{ijk}$}}

\newglossaryentry{DebyeLength}{type=symbol, name={$\ZdebyeLength$}, sort={@lambda.D},
  description={Debye length, $\sqrt{\frac{\epsilon_0 T_0}{n_0 e^2}}$}}
\newglossaryentry{plasmafreq}{type=symbol, name={$\omega_p$}, sort={@zomega},
  description={plasma frequency; $\omega_{ps}^2 = \frac{n_0 e^2}{\epsilon_0 m_s}$}}

\newglossaryentry{a}{type=symbol, name={$\Za_\s$}, sort={a},
  description={particle acceleration, $(q_\s/m_\s)(\E+\v\times\B)$}}
\newglossaryentry{B}{type=symbol, name={$\ZB$}, sort={B}, description={magnetic field}}
\newglossaryentry{b}{type=symbol, name={$\Zbhat$}, sort={b},
  description={magnetic field direction, $\B/|\B|$}}
\newglossaryentry{C}{type=symbol, name={$C$}, sort={C},
  description={collision operator.  See note at \gls{collision operator}}}
\newglossaryentry{tC}{type=symbol, name={$\ZtC$}, sort={C},
  description={total collision operator, e.g.\ $\tC_\i = \C_\i + \C_{\i\e}$.
    See note at \gls{collision operator}}}
\newglossaryentry{c }{type=symbol, name={$c$}, sort={c}, description={speed of light}}
\newglossaryentry{c}{type=symbol, name={$\c_\s$}, sort={c},
  description={velocity in the reference frame of the fluid velocity
    of species $\s$ (an independent variable or the thermal velocity of a particle)}}
\newglossaryentry{dm}{type=symbol, name={$\Zdm$}, sort={dm},
  description={species mass difference, $\mi-\me$}}
\newglossaryentry{dmt}{type=symbol, name={$\Zdmt$}, sort={dm~}, description={$=\dm/e$}}
\newglossaryentry{e }{type=symbol, name={$e$}, sort={e}, description={charge on a proton}}
\newglossaryentry{e}{type=symbol, name={$\e$}, sort={e}, description={electron species index}}
\newglossaryentry{eT}{type=symbol, name={$\ZeT$}, sort={e},
  description={energy tensor per mass}}
\newglossaryentry{strain}{type=symbol, name={$\Zstrain$}, sort={e},
  description={strain rate, $\Sym(\D\u)$}}
\newglossaryentry{dstrain}{type=symbol, name={$\Zdstrain$}, sort={e^c},
  description={deviatoric strain rate $\Zstrain-\id\Div\u/3$}}
\newglossaryentry{dcs}{type=symbol, name={$\Zdcs_t$}, sort={@delta},
  description={bulk derivative; $\Zdcs_t^s\alpha:=\partial_t\alpha + \Div(\u_\s \alpha)$}}
\newglossaryentry{E}{type=symbol, name={$\ZE$}, sort={E}, description={electric field}}
\newglossaryentry{ET}{type=symbol, name={$\ZET$}, sort={E},
  description={energy tensor}}
\newglossaryentry{epsilon_0}{type=symbol, name={$\epsilon_0$}, sort={@epsilon.0},
  description={permittivity of space}}
\newglossaryentry{resistivity}{type=symbol, name={$\Zresistivity$}, sort={@heta},
  description={resistivity}}
\newglossaryentry{Tresistivity}{type=symbol, name={$\ZTresistivity$}, sort={@heta.c},
  description={resistivity tensor}}
\newglossaryentry{tresistivity}{type=symbol, name={$\Ztresistivity$}, sort={@heta.t},
  description={resistivity shape tensor,
  defined by $\Tresistivity=\resistivity\tresistivity$}}
\newglossaryentry{f}{type=symbol, name={$f_\s$}, sort={f},
  description={mass density of species $\s$ in phase space $(\xb,\v)$}}
\newglossaryentry{gamma}{type=symbol, name={$\gamma$}, sort={@c},
  description={Lorentz factor (derivative of time with respect to proper time)}}
\newglossaryentry{NrgTot}{type=symbol, name={$\ZNrgTot$}, sort={E_Tot},
  description={total energy per volume including both gas-dynamic energy
  and the magnetic energy $\mu_0\inv |\B|^2/2$}}
\newglossaryentry{Nrg}{type=symbol, name={$\ZNrg_\s$}, sort={E},
  description={energy per volume of species $\s$}}
\newglossaryentry{nrg}{type=symbol, name={$\Znrg_\s$}, sort={e},
  description={energy per mass of species $\s$}}
\newglossaryentry{i}{type=symbol, name={$\i$}, sort={i}, description={ion species index (usually proton or positron)}}
\newglossaryentry{J}{type=symbol, name={$\ZJ$}, sort={J}, description={current per volume}}
\newglossaryentry{K}{type=symbol, name={$\ZK_\s$}, sort={K}, 
  description={thermal equilibration coefficient}}
\newglossaryentry{heatConductivity}{type=symbol, name={$\heatConductivity$}, sort={k},
  description={heat conductivity}}
\newglossaryentry{TheatConductivity}{type=symbol, name={$\TheatConductivity$}, sort={k},
  description={heat conductivity tensor}}
\newglossaryentry{theatConductivity}{type=symbol, name={$\theatConductivity$}, sort={k},
  description={heat conductivity tensor shape,
  defined by $\TheatConductivity=\heatConductivity\theatConductivity$}}
\newglossaryentry{THeatConductivity}{type=symbol, name={$\THeatConductivity$}, sort={k},
  description={heat tensor conductivity tensor}}
\newglossaryentry{tHeatConductivity}{type=symbol, name={$\tHeatConductivity$}, sort={k},
  description={heat tensor conductivity tensor shape,
  defined by $\THeatConductivity=\THeatConductivity\tHeatConductivity$}}
\newglossaryentry{M}{type=symbol, name={$\ZM$}, sort={M}, description={momentum per volume}}
\newglossaryentry{m }{type=symbol, name={$m_\s$}, sort={m}, description={particle mass of species $\s$}}
\newglossaryentry{mean}{type=symbol, name={$\mean{\chi}$}, sort={!mean},
  description={statistical average over velocity space of $\chi$}}
\newglossaryentry{mte}{type=symbol, name={$\Zmte$}, sort={me~}, description={$=\me/e$}}
\newglossaryentry{mti}{type=symbol, name={$\Zmti$}, sort={mi~}, description={$=\mi/e$}}
\newglossaryentry{me}{type=symbol, name={$\Zme$}, sort={me}, description={electron mass}}
\newglossaryentry{mi}{type=symbol, name={$\Zmi$}, sort={mi}, description={ion mass}}
\newglossaryentry{mred}{type=symbol, name={$\Zmred$}, sort={@mu.red},
  description={the reduced mass, $\frac{m_\i m_\e}{m_\e+m_\e}$}}
\newglossaryentry{viscosity}{type=symbol, name={$\Zviscosity$}, sort={@mu},
  description={the viscosity, $p\tau$}}
\newglossaryentry{tviscosity}{type=symbol, name={$\Ztviscosity$}, sort={@mu.t},
  description={the viscosity shape tensor, defined by $\Tviscosity=\viscosity\tviscosity$}}
\newglossaryentry{Tviscosity}{type=symbol, name={$\ZTviscosity$}, sort={@mu.T},
  description={the viscosity tensor}}
\newglossaryentry{mut}{type=symbol, name={$\Zmut$}, sort={@mu.te},
  description={$=\Zmred/e$}}
\newglossaryentry{ndens}{type=symbol, name={$\Zndens$}, sort={n},
  description={particle number per volume
    (of either species in a neutral two-species plasma with equal
    charge on each species)}}
\newglossaryentry{ne}{type=symbol, name={$n_\e$}, sort={ne}, description={electron number density}}
\newglossaryentry{ni}{type=symbol, name={$n_\i$}, sort={ni}, description={ion number density}}
\newglossaryentry{gyrofreq}{type=symbol, name={$\Zgyrofreq$}, sort={@zomega.c},
  description={gyrofrequency (alias cyclotron frequency), $q|\B|/m$}}
\newglossaryentry{PT}{type=symbol, name={$\ZPT$}, sort={P},
  description={pressure tensor}}
\newglossaryentry{dPT}{type=symbol, name={$\ZdPT$}, sort={P^c},
  description={deviatoric part of pressure tensor}}
\newglossaryentry{PTd}{type=symbol, name={$\ZPTd$}, sort={P^d}, description={drift pressure tensor}}
\newglossaryentry{pTot}{type=symbol, name={$\ZpTot$}, sort={p_Tot},
  description={total pressure per volume including both gas-dynamic pressure
  and the magnetic pressure $\mu_0\inv |\B|^2/2$}}
\newglossaryentry{p }{type=symbol, name={$p$}, sort={p}, description={scalar pressure}}
\newglossaryentry{p}{type=symbol, name={$\p$}, sort={p}, description={particle index}}
\newglossaryentry{Pr}{type=symbol, name={$\Prandtl$}, sort={Pr}, description={Prandtl number}}
\newglossaryentry{Pshape}{type=symbol, name={$\ZPshape$}, sort={@pi.b},
  description={pressure tensor shape, $\ZPT/p$}}
\newglossaryentry{pgf}{type=symbol, name={$\Zpgf$}, sort={@pi.var},
  description={gyrofrequency per pressure isotropization rate, $\ptime\gyrofreq$}}
\newglossaryentry{hgf}{type=symbol, name={$\Zhgf$}, sort={@pi.var.heat},
  description={gyrofrequency per heat flux relaxation rate, $\htime\gyrofreq$}}
\newglossaryentry{dPshape}{type=symbol, name={$\ZdPshape$}, sort={@pi.d},
  description={deviatoric part of $\ZPshape$, i.e.\ $\dPT/p$}}
\newglossaryentry{q}{type=symbol, name={$\q$}, sort={q}, description={heat flux}}
\newglossaryentry{q }{type=symbol, name={$q$}, sort={q}, description={particle charge}}
\newglossaryentry{Q}{type=symbol, name={$Q_\s$}, sort={Q}, 
  description={heat source due to collisions with other species;
    $\Q_\s = \Qf_\s+\Qt_\s $}}
\newglossaryentry{Qf}{type=symbol, name={$\ZQf_\s$}, sort={Qf}, 
  description={heat source due to resistive drag}}
\newglossaryentry{QTf}{type=symbol, name={$\ZQTf_\s$}, sort={Qf},
  description={tensor heating due to resistive drag}}
\newglossaryentry{Qt}{type=symbol, name={$\ZQt_\s$}, sort={Qt}, 
  description={heat source due to thermal equilibration}}
\newglossaryentry{QTt}{type=symbol, name={$\ZQTt_\s$}, sort={Qt},
  description={tensor heating due to thermal equilibration}}
\newglossaryentry{QT}{type=symbol, name={$\ZQT_\s$}, sort={Q},
  description={tensor heating due to collisions with other species;
    $\QT_\s = \QTf_\s+\QTt_\s $}}
\newglossaryentry{qT}{type=symbol, name={$\ZqT_\s$}, sort={q}, description={heat flux tensor}}
\newglossaryentry{dqT}{type=symbol, name={$\ZdqT$}, sort={q^c},
  description={deviatoric heat flux tensor, i.e.\ the traceless part,
    $\qT-\SymC(\id\q)/5$}}
\newglossaryentry{mdens}{type=symbol, name={$\Zmdens$}, sort={@rho},
  description={mass per volume}}
\newglossaryentry{reals}{type=symbol, name={$\Zreals$}, sort={R_e}, description={set of all real numbers}}
\newglossaryentry{R}{type=symbol, name={$\ZR_\s$}, sort={R},
  description={resistive drag force on species $\s$
  due to collisions with other species}}
\newglossaryentry{RT}{type=symbol, name={$\ZRT_\s$}, sort={R_s},
  description={relaxation (isotropization) tensor}}
\newglossaryentry{S}{type=symbol, name={$S$}, sort={S}, description={entropy per volume}}
\newglossaryentry{sMaxwell}{type=symbol, name={$s_\Maxwell$}, sort={sMaxwell},
  description={gas-dynamic entropy per volume}}
\newglossaryentry{s}{type=symbol, name={$s$}, sort={s}, description={entropy per mass}}
\newglossaryentry{stress}{type=symbol, name={$\Zstress$}, sort={@sigma},
  description={stress tensor, $-\ZPT$}}
\newglossaryentry{dstress}{type=symbol, name={$\Zdstress$}, sort={@sigma.d},
  description={deviatoric stress tensor, $-\ZdPT$}}
\newglossaryentry{qdens}{type=symbol, name={$\Zqdens$}, sort={@sigma.q},
  description={charge per volume}}
\newglossaryentry{Sym}{type=symbol, name={$\Sym$}, sort={Sym},
  description={symmetric part of its tensor argument}}
\newglossaryentry{SymB}{type=symbol, name={$\SymB$}, sort={Sym2},
  description={twice the symmetric part of a tensor with two indices}}
\newglossaryentry{SymC}{type=symbol, name={$\SymC$}, sort={Sym3},
  description={thrice the symmetric part of a tensor with three indices}}
\newglossaryentry{ptime}{type=symbol, name={$\Zptime$}, sort={@tau},
  description={relaxation/collision period
  for deviatoric pressure and for the BGK collision operator}}
\newglossaryentry{htime}{type=symbol, name={$\Zhtime$}, sort={@tau.t},
  description={relaxation period for heat flux
  and for the Gaussian-BGK collision operator; $\Zhtime=\ptime/\Pr$}}
\newglossaryentry{TT}{type=symbol, name={$\ZTT$}, sort={T},
  description={temperature tensor}}
\newglossaryentry{dTT}{type=symbol, name={$\ZdTT$}, sort={T'},
  description={deviatoric part of temperature tensor}}
\newglossaryentry{T}{type=symbol, name={$T$}, sort={T},
  description={temperature}}
\newglossaryentry{ThetaT}{type=symbol, name={$\ZThetaT_\s$}, sort={@htheTa},
  description={pseudo temperature tensor, $\ZTT_\s/m_\s$}}
\newglossaryentry{theta}{type=symbol, name={$\theta_\s$}, sort={@htheta},
  description={pseudo temperature, $T_\s/m_\s$}}
\newglossaryentry{u}{type=symbol, name={$\u_\s$}, sort={u},
  description={fluid velocity of species $\s$}}
\newglossaryentry{v}{type=symbol, name={$\Zv$}, sort={v},
  description={velocity (of a particle or as an independent variable)}}
\newglossaryentry{pv}{type=symbol, name={$\Zpv$}, sort={v},
  description={proper velocity ($\gamma \v$)}}
\newglossaryentry{vt}{type=symbol, name={$v_{ts}$}, sort={v_t},
  description={thermal velocity}}
\newglossaryentry{vA}{type=symbol, name={$\ZvA$}, sort={v_A},
  description={Alfv\'en speed}}
\newglossaryentry{w}{type=symbol, name={$\Zw_\s$}, sort={w},
  description={drift velocity of species $\s$ relative to the bulk fluid velocity $\u$}}
\newglossaryentry{xb}{type=symbol, name={$\Zxb$}, sort={x},
  description={position in space (e.g.\ of a particle or as an independent variable)}}

%% file: abs.tex
%

The motivating question for this dissertation was to identify
the minimal requirements for fluid models of plasma to allow
converged simulations that agree well with converged kinetic
simulations of fast magnetic reconnection. We show
that truncation closure for the deviatoric pressure or for the
heat flux results in singularities. Due to the strong pressure
anisotropies that arise in magnetic reconnection we propose
Gaussian-moment two-fluid MHD with isotropization of the pressure
tensor and a Gaussian-BGK closure for the heat flux tensor as
the simplest model that is likely to agree reasonably well in
the diffusion region with kinetic simulations of fast magnetic
reconnection.

For two-dimensional problems invariant under 180-degree rotation
about the origin, we show that if the entropy production, heat
flux, diffusive entropy flux, or deviatoric pressure vanishes
in a neighborhood of the origin then any steady state solution
with nonzero reconnection rate must be singular. In particular,
models which simulate any species using a Vlasov equation or
an adiabatic five-moment or ten-moment model cannot support
converged steady nonsingular magnetic reconnection.
Therefore, for such problems, converged simulation of steady
magnetic reconnection requires that a nonzero collision operator
be explicitly specified.

To study dynamic nonlinear magnetic reconnection we simulate the
GEM magnetic reconnection challenge problem with an adiabatic
two-fluid-Maxwell model with pressure isotropization. Our
deviatoric pressure tensor agrees well with published kinetic
simulations at the time of peak reconnection, but sometime
thereafter the numerical solution becomes unpredictable and
develops near-singularities that crash the simulation unless
positivity limiters are applied. To explain these difficulties we
show that steady reconnection requires heat flux and argue that
sustained reconnection approximates steadily driven reconnection.

This prompts the need for a 10-moment gyrotropic heat flux
closure. Using a Chapman-Enskog expansion with a Gaussian-BGK
collision operator yields a heat flux closure for a magnetized
10-moment charged gas which generalizes the closure of
McDonald and Groth. We argue for this closure against an
entropy-respecting closure.

%% file: ack.tex
I want to express appreciation for some of the people
who have been important to me.

\textbf{Professional acknowledgements}

Firstly, my advisor, \textbf{James Rossmanith}. He has approximated
the ideal in an advisor. I love working with him. He is patient,
generous with his time, and treats students with respect. He has
encouraged me and has underscored what I need to do without ever
chastising or disparaging me for my ignorance or lack of progress.

I have greatly valued his openness. He knows his field and is
transparent about the extent of his knowledge. He structures
things cleanly. He is a model of clarity in how he presents
something to an audience and a superb teacher. He is clear and
transparent in his own thought process and is adept in engaging
the thought process of another. Consequently, he has imparted to
me not just what he knows but how he thinks.

I admire his dedication, reliability, and availability, in
his work, in his involvement with students, and in his
relationship with his family. 

Secondly, I would also like to express appreciation to \textbf{Jerry Brackbill}. 
I was struck from the moment I met him by his enthusiasm and openness.
Over the past year he has given me extensive help and counsel.
I specifically want to thank Jerry and his wife, Isabel, for coming all
the way from Portland so that he could serve on my committee.

In addition, I would like to thank a number of people in the plasma
physics community who have been helpful to me, especially when I was getting started.

\textbf{Nick Murphy} was the person I went to for help when
learning the basics of magnetic reconnection.  He pointed me
toward important recent developments, such as the role of secondary instabilities.

When I have become stuck on critical gritty details,
\textbf{Ping Zhu} has been the person I can go to who will work
through the details with me.

\textbf{Ammar Hakim} was the one who originally inspired me
to work on the ten-moment model.  His work (with Shumlak and Loverich)
on two-fluid simulation of the GEM problem laid the foundation
on which James and I built.  Ammar has pointed me to interesting problems
and has challenged me with critical questions about my numerical approach.

I would like to thank \textbf{Carl Sovinec}, both for serving
on my committee and for his hospitality to me and James
in including and involving us in seminars and arranging for
us to meet with visiting speakers such as Uri Shumlak
and Jim Drake.  Carl embodies the culture of openness and
thoughtful consideration that characterizes the plasma physics
group at UW-Madison.  I specifically wish to thank him for his
careful reading of my dissertation and his thoughtful responses.

I would also like to thank \textbf{Daniel den Hartog}
for personal counsel and professional advice.  I admire the
role that he serves both in the plasma community, which he
has helped me to understand, and in the
Christian community of which I am a part.

Thanks to David Levermore for pointing me to the role
of entropy evolution and suggesting that I use
the Gaussian-BGK collision operator.



\def\housemate{\textit}
\def\friend{\textit}
\textbf{Personal acknowledgements}

I would like to express appreciation to the people whose
presence has given meaning to my life and the work that
I have devoted myself to, and it is to them that I would
like to dedicate this work.

To my housemates who have been not just enduring friends but
men who have cared about what I care about and with whom I have
joined in developing our vision of the world and our role in it:

To \housemate{Angelo Scherer}: Your thoughtful questions have
elicited many of my best thoughts and ideas. By walking into the
light you have laid the foundation for building the community
life that you envision. May the Lord continue to use you to bring
disconnected people into spiritually transforming communities.

To \housemate{Jon Shea}: Your flood of passion and heart of
compassion have drawn me out of myself and helped me to see into
the heart and experience of others. When I lent you books and
cultural resources you responded with such enthusiasm that it was
really you who converted me and got me excited about the writings
of Wendell Berry and Lao Tzu. May the Lord bless and guide you
and Rebecca in all your ways.

To \housemate{Miles Kirby}: Your understanding spirit has lifted
me. Your respectfulness, positivity, and genuine honesty has
shown me by example better ways to handle interpersonal conflict.
Thanks for plugging me into the creative, communitarian social
networks that you are so active in. Your passion for population
health is based not in abstraction but in ground-level personal
relationships and shared experience. The woman you married
complements beautifully your gift of sincere personal connection.
I am thrilled that you and Katie can go to Kenya for two years,
and I hope to meet you again in Africa!

To \housemate{Ryan Doucette}: Thanks for praying and gardening
with me! The books you have introduced me to have expanded my
vision of life. Your appreciation for children, old books, gardening,
animals, and the north woods of Wisconsin is unique and of great
value; may the Lord give you opportunity to serve your community
with your love for particular places and particular things. Your
prayerful and peaceable spirit has calmed me.

To \housemate{Michael Peterson}: I have appreciated living with
someone who also looks for a pattern of community life that is attached
to and stewards the land. You are an example of simplicity, cheerfulness,
steadfastness, and gratitude. Your cooking skill and musicianship
have spiced and enriched my life.  May your life's work and dreams
come to a full fruition.

To those who have profoundly influenced my life:

To \friend{Paul Meyer}: You opened my world.
You introduced me to the Great Books tradition,
the Catholic magisterial intellectual tradition,
and indirectly to the wisdom of Eastern Orthodox spirituality.
Your respectfulness, forthrightness, and humility
have built my trust and brought down my barriers.

To \friend{John Vogel} and the faculty at
\emph{Trinity School River Ridge}:
Among you I experienced an intellectual community
patterned after the body of Christ,
and it is from you that I carry my vision of 
what a community of learners is to be and to do.
John, I am indebted to you for what you were
to me as my mentor teacher.  Your fatherly love
for your students, your exuberant love for
mathematics, and your genuine interest in the full range of
pursuits of the other faculty embody my ideal of a teacher.

To \friend{Mark Whitters} and the \emph{Servants of the Word}:
Through you I experienced a brotherhood that shares a common life
of purposeful manhood. It seems evident that the gift of the
Sword of the Spirit communities that you serve is the quality of
their family life, and by forgoing having families of your own
you have established a foundation of family life that promises
to have a transforming effect on communities for generations
to come. Your common life of prayer and service embodies the
pattern of life that I desire to carry with me and share with
others. Mark, you tell it like it is. You built my trust and I
appreciated the respect you conveyed and your good humor. DOC, your
involvement with Cornerstone (both in Detroit and Uganda) has
been an inspiration to me. Thank you to Nico, Stan, Ed, and many
others who helped me discern my life calling.

To \friend{Wafik Lotfallah}: Thank you for praying with me, for
sharing your wisdom, for helping me to understand the historical
experience of Christian communities in the Islamic world, and for
introducing me to Saint Mary and Saint Rewais Coptic Orthodox
Church in Madison. You are a model to me professionally as a
professor and mathematician and personally by your respectful
candor, perceptive understanding, and peaceful spirit which
reaches across barriers of mistrust. We pray for the continued
renewal of Egyptian society.

To those who profoundly deepened my life by sharing with me from
their own experience what it means to set aside the Lord's day,
particularly Peter Kim and Phil Johnson.

To those who helped me as a young man to grow into maturity:

To \friend{Walter Schultz}: You love truth. 
Your dictum --- ``wisdom, when it encounters truth, bows to it''
--- expresses the pattern of your life. You are truly a professor
of philosophy. You have been a model to me of a man who cares
for his students in a way that is not superficial but deep and
life-changing. Thank you for praying for me each month.
I appreciate the personal interest that you and Mary have 
always expressed in my life.


To \friend{Eric Thomsen}: You have been my model of cheer,
personal regard, humility, and intellectual seriousness.
With you I could push the boundaries and still belong and be respected.
You live out the truth that ``the fear of the Lord is the beginning of wisdom.''
When I was seeking my life direction I remember you exclaiming,
``You're a researcher!''

To \friend{Carl Fischer}: Thank you for your open home and
for treating my friends as your friends.
The love of Jesus shines in you.

To \bfit{Geneva Campus Church}: Geneva has been the center of
my life for the last eight years. In you I have found a group of
Christians who care about the world, who are engaged with the
life and calling of the university, and who are thoughtful and
wise. You have provided a context in which I and other
students have been nurtured in maturity; by
humility and truth, you have modeled family relationship and community life.

To \textit{Sylvia Boomsma}: You and Bob have been to me and to
many other students like second parents. I was moved by your love
and commitment to one another and by your faith, hope, and love
that came out so clearly as Bob was dying. You have both been
earthy, real, and thoughtful. Sylvia, I appreciate the prayer
notices that you send out, in particular for those who suffer
persecution. In the days before Bob died he urged me to do plasma
physics ``to the $n$th degree,'' and it has been a great pleasure
to attempt it.

To Marcia Bosscher: Thank you for being such a hospitable neighbor.
Your open home has been a blessing to many.

To Mike and Beth Winnowski: You personally represent the culture
of forthrightness, thoughtful response, and transparency that I
have valued at Geneva.  

To others with whom I have been engaged in prayer groups:

To Terry Morrison: I appreciate your steadfast interest in
the course of my life.

To Gayle Reed: Your personal and physical trials have deepened and refined you.
The books by Brother Yun and Andrew Murray that you gave me changed my life.



To friends I have made here in Madison:

To \friend{Ron and Margaret Miller}:
Ron, you are a man who cares about others 
and senses how they feel.  Thanks to you and to your mom,
Margaret, for the hospitality you have shown to me.

To friends I have made over the past year:

To Micah Behr and Li Ke: thanks for your encouragement
and understanding.

To two visiting political science professors from China:

To \friend{Bai Wengang}: I admire your commitment to truth
regardless of the cost and your wisdom in tracing modern
political patterns to their ancient precedents. Ancient Chinese
political wisdom --- that legitimate government is derived from
the mandate of Heaven, that good government consists in the
exaltation of the virtuous, that the foundation of government is
a community of standards and the cultivation of virtue (Micius),
and that an orderly society turns on the formation of a moral
atmosphere and an understanding of duties (Confucius) ---
provides needed correctives to modern American and Chinese public
life.


To \friend{Zheng Xiaohua}: I appreciate your concern for good
local government and caring communities. You have provoked me
to seek a deeper understanding of the American experiment with
democracy.

To \bfit{my family}:

To my brother, Michael
and my sister, Laura:
I appreciate the questions you ask,
the carefully considered judgments you offer,
your hospitality,
the love you have for your families,
and the passions and convictions that compel you.
Michael, you have been a force of stability and reason.
Laura, you have been a discerning observer who understands
my dreams.  You were the one who told me I should study physics.

To my parents: You have invested immeasurably in me. You have
imparted to me a love of truth, the value of kindness, and a
way of reasoning through an issue. Thank you for shepherding
my mathematical development. Mom, I appreciate your wisdom and
sympathy. Dad, thank you for teaching me to ask questions about
the world.

To two men who have died who were my closest friends,
who were there for me at my lowest point
and who I wish I had been there for:

To my grandfather, \emph{Aldridge Johnson}:
You were simple, sincere, and affectionate.

To \emph{Matthew Beise}: You were my first
adult friend and the person who most shared
what I loved -- physics, Greek, and Hebrew.
I miss you.

And to the one who will ultimately judge this work,
to Jesus Christ, the true and faithful martyr,
the firstborn from the dead and the ruler of the kings of the earth.

%% file: chap1.tex
%
\chapter{Introduction}
\label{introduction}   

For convenience, most mathematical symbols are hyperlinked to
the symbol glossary, which in turn references points in the text
where each term is defined or discussed.

\section{Conventions}
\label{ConventionsUsed}

For clarity and to avoid misunderstand we specify that 
\begin{itemize}
 \item \textit{The term ``\defining{collisions}'' is used 
   in the broadest possible sense 
   to include all microscale interactions of particles
   with the electromagnetic field.}  The language of
   collisions is thus used in an axiomatic sense,
   rather than in the typically more restricted sense which
   refers specifically to a detailed description of
   particle-particle (Coulomb) interactions.
   In particular, it is assumed that the evolution equation
   of the density functions
   $\f_\i(t,\xb,\v)$ for the ions and
   $\f_\e(t,\xb,\v)$ for the electrons
   satisfy kinetic equations of the form
   \begin{alignat*}{9}
     \partial_t \f_\i &+ \v\dotp\D_{\xb}\f_\i&&+\a_\i \dotp\D_{\v} \f_\i &&= \tC_\i &&= \C_\i &&+ \C_{\i\e}, \\
     \partial_t \f_\e &+ \v\dotp\D_{\xb}\f_\e&&+\a_\e \dotp\D_{\v} \f_\e &&= \tC_\e &&= \C_\e &&+ \C_{\e\i},
   \end{alignat*}
   where we call $\glslink{tC}{\tC_\s}$, $\C_\s$ and $\C_{\s\p}$
   \glslink{collision operator}{collision operators};
   $\C_\i$, $\C_\e$, and $\C_{\i\e}+\C_{\e\i}$ conserve mass,
   momentum, and energy, and the appropriate entropy inequalities
   are assumed to be satisfied.
 \item The terms \gls{kinetic equation} and \gls{Boltzmann equation}
   are used as synonyms.
 \item The terms \gls{Gaussian-moment model} and \gls{ten-moment model}
   are used as synonyms.
 \item The terms \gls{Maxwellian-moment model} and \gls{five-moment model}
   are used as synonyms.
\end{itemize}

%

%
%
%
\section{Plasma}  

A gas is a fluid composed of freely moving particles.
\defining{Plasma} is a gas of charged particles,
which interact with the electromagnetic field.
Most of the universe consists of plasma threaded
by magnetic field lines.  The negatively charged particles
are typically electrons but may also be negatively charged
dust grains.  The positively charged particles 
are ions, typically protons or positrons.  In this document
we are primarily concerned with \defining{two-species} plasmas,
which consist of a species of positively charged particles and a species
of negatively charged particles.
In the case of \defining{pair plasmas},
the negatively charged particles are electrons
and the positively charged particles are positrons.
Positrons are the antiparticles of electrons; they
have the same mass but opposite charge.
In the case of \defining{hydrogen plasmas},
the positively charged particles are protons and
the negatively charged particles are electrons.
The ratio of proton mass to electron mass is large (approximately 1836).

Plasma is a type of magnetohydrodynamic fluid\footnote{Technically
a magnetohydrodynamic fluid should be quasi-neutral
and have non-relativistic flow speeds.}. A
\defining{magnetohydrodynamic (MHD) fluid} is a fluid that
conducts electricity and has a magnetic field. In the presence
of a magnetic field electrical current results in a force on
the fluid perpendicular to the direction of the magnetic field
and the direction of the current. An MHD description of a fluid
becomes necessary when the magnetic field is strong enough to
modify fluid flow (enough to affect phenomena and quantities of
interest).

Each species in a plasma can be regarded as a separate fluid.
One can define bulk fluid variables by appropriately summing or
averaging the fluid variables of each species.
Since these fluids occupy the same space, they can interact
directly through frictional drag force and thermal heat exchange.
The motion of a species fluid relative to the velocity of the
bulk fluid is called the \defining{drift velocity} of the species.
A \defining{two-fluid} model is used for a two-species plasma
and represents each species as a distinct fluid.
\defining[two-fluid MHD]{Two-fluid MHD} regards the (bulk) MHD 
fluid as composed of positively and negatively charged fluids
whose charge densities are assumed to cancel.

Magnetic field evolution is determined by the electric field.
The electric field is given by \defining{Ohm's law}, which
specifies the electric field in terms of electrical current,
magnetic field, fluid velocity, and the pressures of the
positive and negative particles.

Written in full, Ohm's law is the evolution equation for
electrical current solved for the electric field. The full Ohm's
law is complicated, and Ohm's law can only be used by assuming
that fluid quantities have been specified by approximate closure
relationships.

Simplified models of MHD assume a simplified Ohm's law. In the
classical description of electromagnetism, magnetic field is
independent of reference frame but electrical field is not.
\defining{Ideal MHD} assumes that the electric field is zero
in the reference frame of the fluid. \defining{Resistive MHD}
assumes that in the reference frame of the fluid the electric
field equals the resistivity times the electrical current.

For a two-species plasma consisting of equal densities
of positive and negative charge, \defining{resistive Hall MHD}
is a more accurate model (in comparison to resistive MHD) which
defines the electric field in the reference frame of total drift,
which we define to be the fluid velocity plus the sum of the
drift velocities of the positively and negatively charged gases;
ideal Hall MHD assumes that the electric field is zero in this
frame, i.e., that the resistivity is zero \cite{HornigSchindler96}.

In general, ideal models of MHD imply the existence of a
\first{flux-transporting velocity} such that the electric field is
zero in the reference frame of this velocity. In ideal models of
MHD the component of the electric field parallel to the magnetic
field is always zero (\cite{sturrock94}, page 189).

\section{Magnetic reconnection}

A \defining{magnetic field line} is a curve through space
that is everywhere parallel to the local magnetic field. The
\defining{magnetic flux} through an infinitesimal surface
element is the surface area times the component of the
magnetic field perpendicular to the surface. The magnetic field is
\defining{divergence-free}; this means that the total flux of
magnetic field out of a closed region equals the total flux
into the region. As a consequence, by choosing representative
magnetic field lines appropriately, it is possible to think of
the strength of the magnetic field as proportional to the density
of magnetic field lines. We can then say that the magnetic flux
through a surface is proportional to the number of magnetic field
lines passing through the surface.

The magnetic field evolution equation says that the rate of
change of magnetic flux through a surface equals the \emph{circulation integral}
of the electric field around the boundary of the surface. In ideal MHD
the electric field in the reference frame of the fluid is zero.
Therefore, the magnetic flux through a surface element that is
carried with the fluid cannot change (because any circulation integral
of a uniformly zero electric field is zero). The result is that in ideal
MHD magnetic field lines are carried with the fluid and cannot
change their topology, i.e.\ how they are connected
(\cite{sturrock94}, p186).
This is known as the \emph{frozen-in flux condition}. In regions
where ideal MHD does not hold, magnetic field lines can slip through the fluid
and break and \defining{reconnect}, i.e.\ change their topology
(how they are connected).

\begin{figure}
     \begin{tabular}{c c}
    time = 1
    &
    time = 4
    \\
    \includegraphics[width=.4\linewidth,height=.2\linewidth]{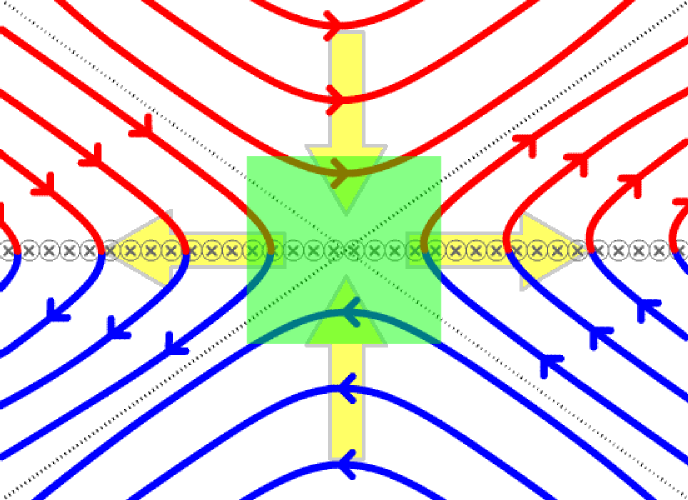}
    & 
    \includegraphics[width=.4\linewidth,height=.2\linewidth]{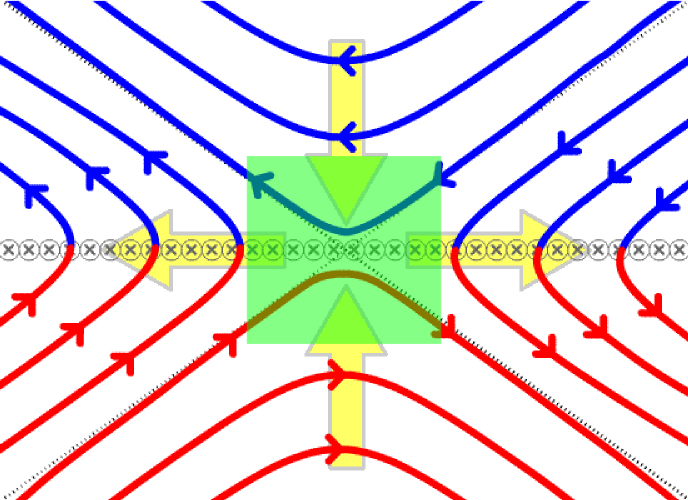}
    \\
    time = 5
    &
    time = 6
    \\
    \includegraphics[width=.4\linewidth,height=.2\linewidth]{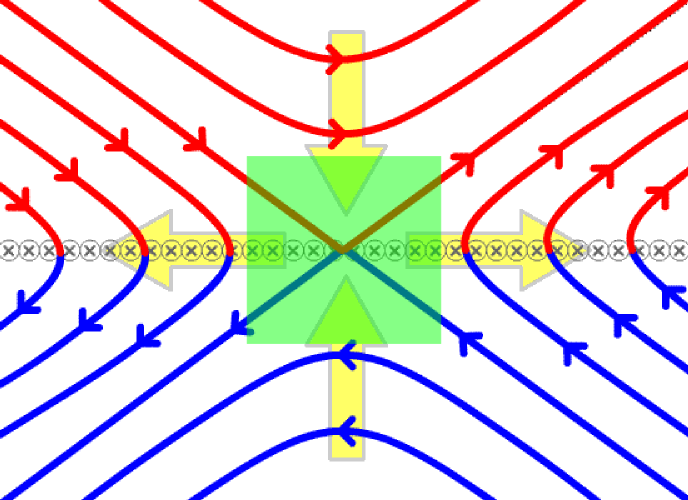}
    & 
    \includegraphics[width=.4\linewidth,height=.2\linewidth]{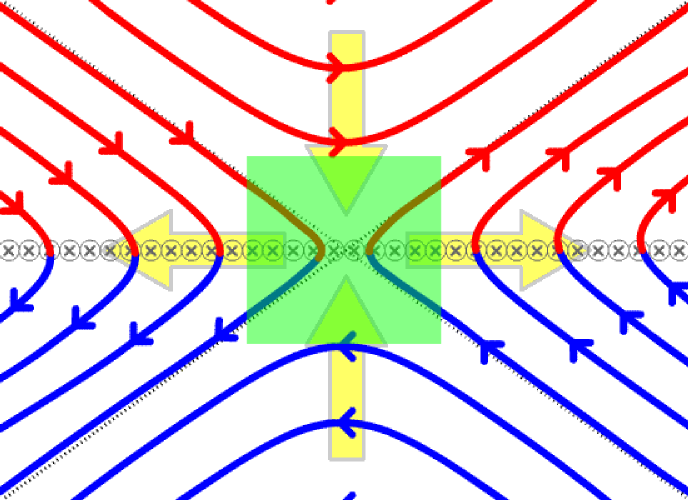}
    \end{tabular}
   \caption{Cartoon of symmetric antiparallel 2D separator reconnection (wikipedia).}
     At the origin the magnetic field is zero.
     Symmetry forces physical velocities to be zero at the origin,
     including the bulk fluid velocity, the fluid velocity of each species,
     and any flux-transporting flow.
     The magnetic field lines which intersect the origin are
     called \defining{separatrices} and partition the domain
     into four regions. 
     In an ideal plasma model the separatrices would be frozen
     in the fluid and therefore it would be impossible to transfer flux
     across the separatrices.
     The reconnected flux is defined to be the amount of
     flux transferred from one region to another.
     The green rectangle at the center signifies the diffusion region.
     The ratio of its long side to its short side is called its
     aspect ratio.
  \fhrule
  \label{wikiRecon}
 \end{figure}

The amount of reconnection is difficult to define in a precise
and general way. For the example of 2D separator reconnection
as depicted in figure \ref{wikiRecon},
a set of magnetic field lines called separatrices partitions the
spatial domain into four regions such that in an ideal model flux
could not be transferred from one domain to another.
In this case we can define the amount of reconnection to be
the amount of flux transferred from one of these regions to a
neighboring region.

A definition of magnetic reconnection in general geometry
requires a covariant (relativistic) description of the
electromagnetic field \cite{HornigSchindler96}.
Since this dissertation only studies 2D
problems symmetric under 180-degree rotation in the plane about
the origin, it defines magnetic reconnection only for this case;
the rate of reconnection turns out simply to be the out-of-plane
electric field component at the origin, as argued in section \ref{2DmagneticReconnection}.


\section{Fast magnetic reconnection}

\subsection{Background}

Magnetic reconnection is ultimately controlled by Faraday's
evolution equation for the magnetic field,
\begin{gather*}
  \partial_t \B + \curl\E = 0;
\end{gather*}
so the evolution of magnetic field $\B$ is determined by the electric field $\E$.
Plasma physics studies the evolution of plasmas on time scales larger than
the plasma period $\glslink{plasmafreq}{\omega_p}^{-1}$ (the time of oscillation
in response to charge separation) and larger than
the Debye length \gls{DebyeLength} (the distance traveled by an electron
in a plasma period moving at the thermal velocity \gls{vt},
which is the distance over which charge shielding occurs) \cite{HW04}.
On these scales the assumption of \gls{quasineutrality} holds and
\mention{Ohm's law} \eqref{OhmsLawFullTerms} for the electric field applies.
For a plasma consisting of electrons and ions, Ohm's law is
\begin{equation} 
 \begin{aligned}
  \E =& \B\times\u  & \hbox{(ideal term)}
    \\ &+ \Tresistivity\dotp\J  & \hbox{(resistive term)}
    \\ &+ \tfrac{\mi-\me}{e\mdens}\J\times\B & \hbox{(Hall term)}
    \\ &+ \tfrac{1}{e \mdens}
          \Div\left(\me\PT_\i-\mi\PT_\e\right)
         & \hbox{(pressure term)}
    \\ &+ \tfrac{\mi\me}{e^2 \mdens}\left[
          \partial_t\J + \Div\left(\u\J+\J\u-\tfrac{\mi-\me}{e\mdens}\J\J\right)
        \right] & \hbox{(inertial term),}
 \end{aligned}
 \label{OhmsLawStandard}
\end{equation}
where $\u$ is fluid velocity, $\J$ is electrical current density,
$\mdens$ is mass density, and $\PT_\i$ and $\PT_\e$ are the ion and
electron pressure tensors; the constants are the ion mass $\mi$,
the electron mass $\me$, and the magnitude of the charge on an
electron, $e$.  The resistivity $\Tresistivity$ requires a closure,
typically a function of (electron) temperature unless an anomalous
resistivity is defined.

This work restricts consideration to problems symmetric under 180-degree
rotation, because it allows a simple analysis of the X-point
which identifies constraints and requirements for reconnection.
As argued in section \ref{2DmagneticReconnection},
at the X-point Ohm's law reduces to
\begin{gather} 
 \begin{aligned}
  \E =& \Tresistivity\dotp\J  & \hbox{(resistive term)}
    \\ &+ \tfrac{1}{e \mdens}
          \Div\left(\me\PT_\i-\mi\PT_\e\right)
         & \hbox{(pressure term)}
    \\ &+ \tfrac{\mi\me}{e^2 \mdens}\left[
          \partial_t\J + \Div\left(\u\J+\J\u-\tfrac{\mi-\me}{e\mdens}\J\J\right)
        \right] & \hbox{(inertial term)}
 \end{aligned}
 \label{ohmsLawAtTheXpoint}
\end{gather}
where only the out-of-plane component survives. Since the rate of
reconnection is the out-of-plane component of the electric field at
the origin, this confirms that one of these three terms must be
nonzero to support reconnection.

Magnetohydrodynamic (MHD) models of plasma
explictly assume a form of Ohm's law and evolve a system
sufficient to determine the quantities it involves.
Ideal MHD discards all terms in Ohm's law except the ideal term
and is the simplest model of plasma. In ideal MHD magnetic
reconnection is not possible. The next simplest model, resistive
MHD, also retains the resistive term. Resistive MHD allows
magnetic reconnection to occur, because magnetic field lines can
diffuse through the plasma; as we will discuss, steady reconnection
rates are slow for resistive MHD unless an anomalous resistivity
is used.  Resistive Hall MHD includes the Hall
term as well, allowing much faster rates of steady reconnection. 
Models which include terms of Ohm's law beyond resistive MHD
are collectively known as \mention{extended MHD} models.

\subsection{Historical background for fast magnetic reconnection}

The notion of magnetic reconnection was introduced by Dungey
(\cite{dungey53}, 1953). In the subsequent decades people
attempted to explain reconnection in terms of resistive MHD.
Sweet (\cite{Sweet58}, 1958) and Parker (\cite{Parker57},
1957) developed a model of steady two-dimensional magnetic
reconnection. The Sweet-Parker model assumes that magnetic
reconnection occurs in a thin, rectangular \mention{diffusion
region} containing nearly antiparallel magnetic field lines.
The \defining{aspect ratio} of a rectangle denotes the ratio
consisting of the length of the long side divided by the length
of the short side. In the Sweet-Parker model of reconnection
the aspect ratio of the diffusion region is assumed to be
large, implying that the magnetic field lines have a ``Y-type''
configuration. Plasma is assumed to flow slowly into the broad
sides of the rectangle and rapidly out of the narrow sides of
the rectangle. In the diffusion region, gradients in the
magnetic field result in diffusion of the magnetic field
(and hence slipping of magnetic field lines through the
plasma) and a sheet of electrical current. Using this
model, Parker (\cite{Parker63}, 1963) estimated a typical rate
of magnetic reconnection in solar flares based on Spitzer's
formula for resistivity (which is based on Coulomb collisions
and asserts that resistivity is a function of temperature,
see \cite{Spitzer62}) and found that the predicted rate of
reconnection was at least one hundred times too small to account
for observed solar flare time scales.

The essential impediment to fast reconnection in the Sweet-Parker
model is the tension between the need for a narrow diffusion
region (so that magnetic field gradients can be sufficiently
strong) and a wide diffusion region (so that plasma can flow
out of the narrow sides of the rectangle rapidly enough).
Subsequent models of steady reconnection obtained faster rates of
reconnection by assuming anomalously high values of resistivity
in the diffusion region and/or by assuming an ``X-type'' magnetic
field configuration (in particular a diffusion region
with a small aspect ratio) which allows strong
magnetic field gradients at the X-point while opening up the
outflow so that it is not throttled by being confined to a narrow
rectangle \cite{priestForbes00}. The first such model was devised
by Petschek (\cite{petschek64}, 1964) and consisted of an X-type
magnetic field geometry with a miniature Sweet-Parker region at
the center. The Petschek model allowed much faster reconnection
rates than the Sweet-Parker model, and thus Sweet-Parker
reconnection became designated as \mention{slow reconnection}
while reconnection rates on the order given by the Petschek model
were identified as \mention{fast reconnection}.

When numerical simulation of magnetic reconnection became
feasible, it was found that magnetic reconnection resulted
in Y-type configurations and slow reconnection if a uniform
resistivity was assumed, whereas X-type configurations and fast
reconnection occurred if anomalously high values of resistivity
were assumed near the X-point \cite{birnPriest07}. There are physical reasons
to expect anomalously high values of resistivity near the
X-point. First, resistive drag may depend nonlinearly on electric
current. Electrical currents are strong in the reconnection
region. Electrical currents represent relative drift of ions
and electrons. If this relative drift becomes strong enough,
a streaming instability develops, limiting interspecies drift
and greatly increasing resistivity. Second, resistivity may be
\emph{spatially} dependent in a weakly collisional plasma where
fluid closures cannot be rigorously justified. Spitzer's formula
for resistivity assumes collisional transport theory, which
is applicable when the mean free path of a particle is small
relative to the length scale of variations in the magnetic field
and gas-dynamic quantities.  Particle mean free paths are much larger than 
the width of current sheets or diffusion layers where reconnection
occurs, and so collisional transport theory is not applicable
even when, as in the solar corona, it is applicable to
large-scale structures (see \cite{priestForbes00}, page 45); thus,
the reconnection region is governed by collisionless physics
in essentially all space and laboratory plasmas
where magnetic reconnection is important \cite{brackbill:priv11}.

Thus, for resistive MHD the game of modeling reconnection naturally
became to determine anomalous values of resistivity that account
for fast reconnection. By assuming a \emph{spatially} dependent
anomalous resistivity one can essentially prescribe a desired
rate of reconnection. Some space modelers have taken the approach
of prescribing spatially anomalous resistivities which give
results that seem to agree with observational and statistical
data. Such an approach can be effective in a specific problem
domain such as space weather modeling, but in general we prefer
the simplest models with the greatest explanatory and predictive success,
which are based in physical principles, and which give physical
insight.

One can obtain fast rates of reconnection using an anomalous
resistivity that is dependent on current but spatially
independent. The project to formulate a spatially independent
anomalous resistivity that accurately models reconnection
has fallen short, however, on two accounts. First, there are many
formulas for anomalous resisivity, and a simple
basis for these formulas has been elusive
\cite{satoHayashi79, article:Vasyliunas75, brackbill:priv11}.
Second, in weakly collisional regimes steady reconnection is supported primarily
by the divergence of the species pressure tensors rather than
by resistive drag, and therefore attempting to attribute the
reconnection electric field to a resistive term is artificial
and does not promise to give physical insight; the end road
of insisting on an anomalous resistivity framework is that an
appropriate formula for the anomalous resistivity will probably
be found to be in terms of the divergence of the electron pressure
tensor!

Although a uniform resistivity results in a long, thin current
sheet and slow reconnection, this configuration is often unstable.
Furth, Killeen, and Rosenbluth (\cite{furth63}, 1963) found
that a current sheet with an aspect ratio of about $2\pi$ or
greater is unstable to spontaneous reconnection which
forms magnetic islands.  This process is called the \mention{tearing mode
instability}, and the magnetic islands are referred to as
\mention{plasmoids}.  Bulanov \emph{et al.} (\cite{bulanov78}, 1978)
repeated the tearing mode calculation of \cite{furth63}
assuming linear outflow along the current sheet and found that
the outflow had a stabilizing effect.

Jumping forward to the past decade,
Loureiro, Schekochihin, and Cowley
(2007, \cite{article:LouSchCow07}) performed an asympototic
analysis in the inverse of the aspect ratio of the current sheet
and showed that for a current sheet of sufficiently large length $L$
(for which the Lundquist number $\mathrm{Lu}:=L \vA/\resistivity$
is greater than a critical value of about $4\times 10^4$,
corresponding to an aspect ratio of about 200;
see \cite{huangBh10}), a chain of plasmoids rapidly forms. 
Subsequent simulations have confirmed that the ejection of these
plasmoids allows fast reconnection rates even in resistive MHD.
The consequence is that although resistive MHD with uniform
resistivity does not admit fast reconnection, for sufficiently
large (e.g.\ astrophysical) domains statistically steady 
fast reconnection can be expected via the cascading formation and
ejection of plasmoids \cite{huangBh10}.  Nevertheless, one can
still make the categorical assertion that resistive MHD with
uniform resistivity does not support steady-state fast reconnection.

In laboratory plasmas Lundquist numbers are typically on the order
of $10^3$ (page 44 in \cite{priestForbes00}),
which is too small to give rise to plasmoid-mediated reconnection;
instead one expects (slow) Sweet-Parker reconnection if the ion inertial length
$\glslink{delta}{\Zdelta_s}$ is larger than the Sweet-Parker layer
thickness.
For ion skin depth smaller than the Sweet-Parker layer thickness
one expects fast, \first{Hall-mediated reconnection} instead,
as discussed in the next section.
(See p315 in \cite{zweibel09} and Figure 1 in \cite{huang10}.)

\subsection{The GEM problem}

The inadequacy of resistive MHD to account for reconnection electric
fields in the diffusion region lead to studies of reconnection
using using extended MHD.
The historical development of these studies is traced in \cite{shay01},
and lead to the following observations regarding 2D separator reconnection.
As shown for collisional tearing in \cite{Terasawa83}, the Hall
current effect becomes important for a current sheet whose width
is comparable to the ion inertial length. Outside of the current
sheet gradients are small, and the ideal term dominates in Ohm's
law. Within the current sheet the Hall term becomes significant
and the ions decouple from the magnetic field lines, defining the
\mention{ion diffusion region}. In a smaller region the electron
pressure term (or inertial term) becomes significant and the
electrons decouple from the magnetic field as well, defining the
\mention{electron diffusion region}.

The culmination of these studies and observations was the
formulation and simulation of the \mention{Geospace Environmental
Modeling magnetic reconnection challenge problem} (\mention{GEM
problem}) in 2001 \cite{article:GEM}. The GEM problem was
formulated to study the ability of plasma models to resolve
fast magnetic reconnection \cite{article:GEM}. The GEM problem
identified two-fluid/Hall effects as critical to allow
fast reconnection. Although ideal Hall MHD does not admit
reconnection, resistive Hall MHD (even with small resistivity)
was found to admit fast reconnection \cite{shay01}, as if the Hall term were
a catalyst accelerating the must slower rate of reconnection that
occurs in resistive MHD without the Hall term.

I remark that ideal Hall MHD (Ohm's law using only the ideal
term and the Hall term) does not allow reconnection because
a flux-transporting flow exists; magnetic field lines are
essentially frozen to the electrons. The ideal Hall MHD
simulations in \cite{shay01} were able to get fast reconnection
rates because of the presence of numerical resistivity; the
results therefore cannot be converged. Finding 3 at the end of
section 1 suggests that in their simulations reconnection is
supported by numerical resistivity with an anomalously high value
near the X-point. It appears that it is still an open question
whether converged fast reconnection is possible in resistive
Hall MHD with uniform resistivity. If the answer is no, then one
may conclude that short of anomalous closures nonzero
divergence of the pressure tensor is necessary for converged
steady fast magnetic reconnection.
This suggests a study of of reconnection with a two-fluid model
with resistivity but without viscosity, with and without the inertial term.

\subsubsection{Pair plasma GEM simulations}

Since the seminal GEM problem studies had identified the Hall
effect as the essential physics to admit fast reconnection,
it was natural to investigate reconnection rates in pair
plasmas, for which $\me=\mi$ and the Hall term of Ohm's law
is absent. Particle simulations by Bessho and Bhattacharjee
of antiparallel reconnection
\cite{article:BeBh05, article:BeBh07, bessho10} have demonstrated
that fast reconnection rates occur even in the case of pair
plasmas; they find that the divergence of the pressure tensor is
the term of equation \eqref{ohmsLawAtTheXpoint} that primarily
supports the reconnection electric field.
For the guide-field case, in which the out-of-plane component of
the magnetic field at the origin is nonzero,
Chac\'on \emph{et al.} \cite{article:ChSiLuZo08}
subsequently demonstrated that steady fast reconnection
is possible in a viscous incompressible model of pair plasma
if viscosity dominates but not if resistivity dominates.

In this work we simulate the pair plasma version of the
GEM problem with a two-fluid adiabatic model without resistivity
and show that rates of reconnection are still fast (although our
rate of reconnection is only 60\% of the rate in the PIC simulations
reported e.g.\ in \cite{article:BeBh07}).

\subsubsection{Two-fluid GEM simulations}

The fluid models used in the seminal GEM problem studies
did not include the pressure and inertial terms of
Ohm's law.  It is therefore natural to ask whether
the inclusion of these terms would allow
significantly improved agreement of
fluid simulation of magnetic reconnection 
with kinetic simulations.

In 2005 Hakim, Loverich, and Shumlak simulated
fast magnetic reconnection with an adiabatic
inviscid \glshyperlink[five-moment]{five-moment model}
\glshyperlink{two-fluid Maxwell} model which
implies an Ohm's law that includes
the Hall term, the inertial term, and a pressure term with
scalar pressures \cite{article:HaLoSh06}\footnote{
This paper used the finite volume wave propagation method
described in \cite{book:Le02}.  Loverich, Hakim, and Shumlak
also performed a complementary study using the Discontinuous
Galerkin method at about the same time; this study was
(finally!) accepted for publication in \cite{LoHaSh11}.}.
Their figure 10 shows their reconnected flux values
superimposed on the reconnection rates reported in the
seminal GEM problem papers and arguably shows improved
agreement with particle simulations in comparison
to the Hall MHD simulations.
In 2007 Hakim submitted simulations of the GEM problem
with a two-fluid Maxwell model which uses a hyperbolic
(adiabatic inviscid) five-moment model for the electrons
and a hyperbolic \glshyperlink{ten-moment} model for the ions and
again obtained reconnection rates that agree
well with kinetic simulations \cite{article:Hakim08}.\footnote{
In particular, Hakim attains one nondimensionalized
unit of reconnected flux at a nondimensionalized
time of about 17.6, in comparison to the values of
15.7 obtained using a PIC simulation \cite{pritchett01}
and 17.7 using a Vlasov simulation \cite{article:SmGr06};
see table \ref{table:recon}.
}

The models used by Loverich, Hakim, and Shumlak are hyperbolic.
In particular, for the electrons they use hyperbolic
five-moment gas dynamics, which uses truncation closures for
the deviatoric pressure and heat flux. As a proxy for Ohm's law
\eqref{OhmsLawStandard} one may consider the electron momentum equation
\eqref{momentumEvolution} solved for the electric field:
\begin{gather}
 \begin{aligned}
  \E &= \B\times\u_\e  & \hbox{(ideal term)}
  \\ &+ \frac{-\R_\e}{q_\e \ndens_\e}  & \hbox{(resistive term)}
  \\ &+ \frac{\Div\PT_\e}{q_\e \ndens_\e}  & \hbox{(pressure term)}
  \\ &+ \frac{m_\e}{q_\e} d_t \u_\e  & \hbox{(inertial term)}.
 \end{aligned}
\end{gather}
At the X-point symmetry the ideal term disappears, simplifying this
to the equivalent of equation \eqref{ohmsLawAtTheXpoint}.
In the simulations of Loverich, Hakim, and Shumlak the resistive
term is zero and the pressure term vanishes at the X-point
because the deviatoric pressure is zero. This would force the
electron velocity at the X-point to ramp with reconnected flux.
As our simulations indicate (see figure \ref{fig:0}), this almost
certainly is not realized in their simulations for later times,
and therefore their solutions are presumably not converged for later times.
Furthermore, based on kinetic simulations (see e.g.\ figure
\ref{fig:6}), we expect the pressure term to dominate at the
X-point.

We were therefore motivated to consider two-fluid model with
viscosity in both species fluids. An easy way to implement
viscosity in a gas is to represent it using a
\glshyperlink{ten-moment model}
(also known as a \glshyperlink{Gaussian-moment model}
with relaxation to isotropy.  A ten-moment model with relaxation
to isotropy agrees with a viscous five-moment model for fast
isotropization rates and small pressure anisotropies. For slow
isotropization rates, large pressure anisotropies can develop.
In kinetic simulations, strong pressure anisotropy is in fact
observed at the X-point (see figure 6 of \cite{article:SmGr06}),
adding to our motivation for use of a ten-moment isotropizing model
for both electrons and ions rather than a five-moment viscous model.




Using the ten-moment two-fluid model to simulate the GEM problem,
we obtain qualitatively good agreement with kinetic simulations
for plots of the pressure tensor at the time of peak reconnection
rate (roughly 16--18, where the unit of time is a typical ion
angular gyroperiod as defined in the GEM problem); at this time
the pressure is highly agyrotropic in the immediate vicinity of
reconnection point. In the approximate time interval from 18 to
28 the rate of reconnection remains approximately constant while
the electron temperature tensor becomes increasingly singular
near the X-point. For a coarse mesh this temperature singularity
does not interfere with the normal progression and saturation
of reconnection, probably due to numerical thermal diffusion,
but when the mesh is refined the singularity becomes sharp and
ultimately prevents normal progression of reconnection. The specific
terminating behavior is erratic and hard to predict, but in
general at the X-point the temperature parallel to the outflow
axis will become very cold while the temperature parallel to the
inflow axis will become very hot. Unless positivity limiters
are applied the simulation crashes when a non-positive-definite
temperature tensor develops. (Typically the anisotropically hot spot
at the origin splits into two hot spots along the outflow axis
before this crash occurs.) If positivity limiters are applied then a
secondary island at the origin is likely to form. If symmetry
is enforced then this stops reconnection, but if not then
spontaneous symmetry-breaking can eject the island and allow
reconnection to proceed.

\subsection{Entropy production and heat flux requirements
  for steady magnetic reconnection}

The characteristic behavior of dynamical systems turns critically
on the character of their equilibria. To understand the
difficulties encountered with the GEM problem when an adiabatic
model is used, we consider the entropy production and heat flux
requirements for steady magnetic reconnection. We show that for
models which lack a mechanism for heat flux, in reconnection
problems that are symmetric under 180-degree rotation about the
origin (as holds for the GEM problem), solutions which exhibit
steady reconnection must be singular.

\subsection{Gyrotropic ten-moment heat flux closure}

We therefore consider what an appropriate 10-moment heat
flux closure would be.  In the presence of a magnetic
field appropriate diffusive closures need to be
\emph{gyrotropic} rather than isotropic.  We therefore
generalize the isotropic ten-moment heat flux closure 
recently formulated by McDonald and Groth to the gyrotropic case.

\section{Model equations}

For subsequent reference, in this section we list the systems of equations that define
the models discussed in this work.
All simulations have been performed using the 10-moment two-fluid Maxwell
model or the 5-moment two-fluid Maxwell model.

\begin{figure}[t]
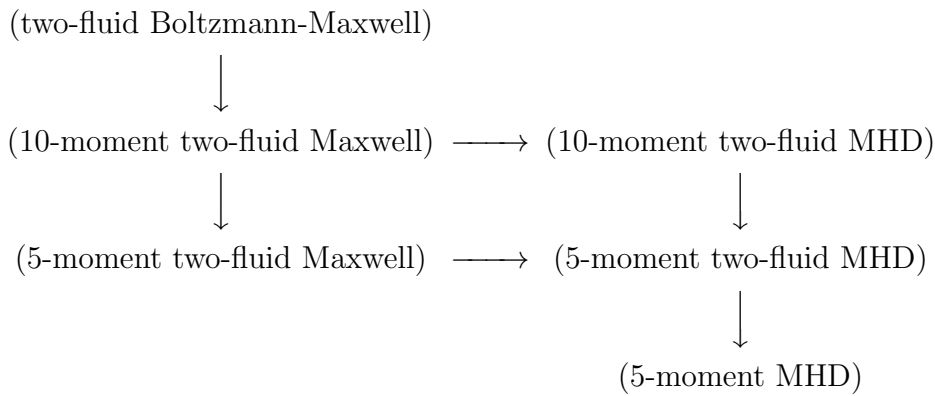

\hfrule
\textbf{\large Model hierarchy}
\[
    \begin{CD}
      \hbox{(two-fluid Boltzmann-Maxwell)} \\
        @VVV \\
      \hbox{(10-moment two-fluid Maxwell)} @>>> \hbox{(10-moment two-fluid MHD)} \\
        @VVV                     @VVV \\
      \hbox{(5-moment two-fluid Maxwell)} @>>> \hbox{(5-moment two-fluid MHD)} \\
         @.                       @VVV \\
                          @.   \hbox{(5-moment MHD)} \\
    \end{CD}
\]
\caption{Hierarchy of models considered in this work.}
  Motion down and to the right indicates use of simplifying limits.
  Motion down reduces the number of moments evolved.
  Motion to the right takes light speed to infinity.
 \fhrule
 \label{fig:modelHierarchy}
\end{figure}

The relationship among models is laid out in figure \ref{fig:modelHierarchy}.
As the standard of truth we take the Boltzmann-Maxwell model
displayed in figure \ref{fig:boltzmannMaxwell}.
The 10-moment two-fluid Maxwell equations that we solve
are displayed in figure \ref{fig:10mom2fluidMaxwell}.
The 5-moment two-fluid Maxwell equations that we solve
are displayed in figure \ref{fig:5mom2fluidMaxwell}.
Taking the light speed to infinity converts the two-fluid Maxwell
models to MHD models.
In solving the two-fluid Maxwell equations our goal is
to approximate the MHD systems, which we attempt to do by
using a sufficiently high light speed.
We therefore display the equations of 10-moment 2-fluid MHD
in figure \ref{fig:10mom2fluidMHD}
and display the equations of 5-moment 2-fluid MHD
in figure \ref{fig:5mom2fluidMHD}.
This work calculates intraspecies collisional closure coefficients
using a Chapman-Enskog expansion with a Gaussian-BGK collision operator.
The resulting formulas for closure coefficients are displayed
for the ten-moment model in figure \ref{fig:10momCoef}
and for the five-moment model in figure \ref{fig:5momCoef}.

\begin{figure}[t]
\hfrule
 \textbf{\large Boltzmann-Maxwell model}
\begin{itemize}
\item \textbf{Kinetic/Boltzmann equations:}
\begin{alignat*}{7}
  \hspace{-10mm}
\partial_t \f_\i &+ \v\dotp\D_{\xb}\f_\i&&+\a_\i \dotp\nabla_{\v} \f_\i &&= \LCi &&+ \LCie \\
  \hspace{-10mm}
\partial_t \f_\e &+ \v\dotp\D_{\xb}\f_\e&&+\a_\e \dotp\nabla_{\v} \f_\e &&= \LCe &&+ \LCei
\end{alignat*}
\item \textbf{Maxwell's equations:}
\begin{alignat*}{5}
  \hspace{-10mm}
  &\partial_t \B + \curl\E = 0, &\quad& \D\dotp\B=0 \\
  \hspace{-10mm}
  &\partial_t \E - c^2\curl\B = -\J/\epsilon_0, &\quad& \D\dotp\E=\qdens/\epsilon_0
  %
\end{alignat*}
\item \textbf{Definitions:}
\begin{gather*}
  \hspace{-10mm}
  \qdens := \sum_\s \frac{q_\s}{m_\s} \int \f_\s \, d\v, \\
  \hspace{-10mm}
  \J := \sum_\s \frac{q_\s}{m_\s} \int \v \f_\s \, d\v
  \\
  \begin{alignedat}{7}
   \hspace{-10mm}
   \a_\i &:= \tfrac{e}{m_i}&&\left(\E+\v\times\B\right), \\
   \hspace{-10mm}
   \a_\e &:= \tfrac{-e}{m_e}&&\left(\E+\v\times\B\right)
  \end{alignedat}
\end{gather*}
\end{itemize}
\caption{Equations of the two-species Boltzmann-Maxwell model (standard of truth).}

 The interspecies collision operators $ \inmag{C_{ie}}$ and $ \inmag{C_{ei}}$
 are generally ignored in this work, but the intraspecies collision operators
 $\indarkblue{C_i}$ and $\indarkblue{C_e}$ play a critical role.
 \label{fig:boltzmannMaxwell}
  \fhrule
\end{figure}

\begin{figure}[t]
\hfrule
 \textbf{\large Ten-moment two-fluid Maxwell model}

 \textbf{Gas dynamics equations}
 \vspace{-1ex}
 {
 \begin{alignat*}{10}
    \mspace{-30mu}
    \dcs_t^\i
    \begin{bmatrix}
      \mdens_\i \\
      \mdens_\i\u_\i \\
      \mdens_\i\eT_\i
    \end{bmatrix}
    &+
    \begin{bmatrix}
      0 \\
      \Div\PT_\i \\
      \SymB(\Div(\PT_\i\u_\i))+\LDivqT{\i}
    \end{bmatrix}
    &&= 
    \qdens_\i
    \begin{bmatrix}
      0 \\
      \E + \u_\i\times\B \\
      {\SymB(\u_\i\E+\eT_\i\times\B)}
    \end{bmatrix}
    &&+
    \begin{bmatrix}
      0 \\
      \LR{\i} \\
      \LRT{\i} + \LQT{\i}
    \end{bmatrix}
   \\
    \mspace{-30mu}
    \dcs_t^\e
    \begin{bmatrix}
      \mdens_\e \\
      \mdens_\e\u_\e \\
      \mdens_\e\eT_\e
    \end{bmatrix}
    &+
    \begin{bmatrix}
      0 \\
      \Div\PT_\e \\
      \SymB(\Div(\PT_\e\u_\e))+\LDivqT{\e}
    \end{bmatrix}
    &&= 
    \qdens_\e
    \begin{bmatrix}
      0 \\
      \E + \u_\e\times\B \\
      {\SymB(\u_\e\E+\eT_\e\times\B)}
    \end{bmatrix}
    &&+
    \begin{bmatrix}
      0 \\
      \LR{\e} \\
      \LRT{\e} + \LQT{\e}
    \end{bmatrix}
 \end{alignat*}
 }

 \begin{multicols}{2}
\textbf{Maxwell's equations:}
\begin{alignat*}{5}
    \mspace{-30mu}
  &\partial_t \B + \curl\E = 0, &\quad& \D\dotp\B=0 \\
    \mspace{-30mu}
  &\partial_t \E - c^2\curl\B = -\J/\epsilon_0, &\quad& \D\dotp\E=\qdens/\epsilon_0
\end{alignat*}
 \textbf{Definitions:}
 \begin{gather*}
    \mspace{-30mu}
     \qdens_\i = en_\i, \quad
     \qdens_\e = -en_\e, \quad
     \qdens = \qdens_\i + \qdens_\e
     \\
    \mspace{-30mu}
     \J = \qdens_\i\u_\i + \qdens_\e\u_\e
     \\
    \mspace{-30mu}
     \glslink{dcs}{\Zdcs_t^s}\alpha:=\partial_t\alpha + \Div(\u_\s \alpha)
   \end{gather*}

 \columnbreak
 \textbf{Closures:}
 \begin{gather*}
    \vspace{-1ex}
    \hspace{-1.5em}
   \glslink{RT}{\indarkblue{\ZRT_\s}} = -\ptime_\s\inv \dPT_\s
   \\ \hspace{-1.5em}
     \LqT{\s} = -\tfrac{2}{5}\THeatConductivity_\s
         \dddotp\SymC\left(\Pshape_\s\dotp\D\TT_\s\right) 
   \\ \hspace{-1.5em}
     \LR{\i} = -\LR{\e}
    = \qdens_\i \qdens_\e \glslink{Tresistivity}{\inmag{\ZTresistivity}}\dotp(\u_\e-\u_\i)
   \\ \hspace{-1.5em}
     \LQT{\s} = ?
       \LQTf{\s} + \LQTt{\s} ?
   \\ \hspace{-1.5em}
     \LQTf{\e} = \frac{m_\i}{m_\i+m_\e}\SymB
      \left(
       (\alpha^\para_\e-\alpha^\perp_\e)\LR{\i}\w_\e
       +\alpha^\perp_\e\LR{\i}\dotp\w_\e\id
      \right)\, ? \\
   \\ \hspace{-1.5em}
     \LQTt{\e} = \tfrac{2}{3} \inmag{K_{\e\i}}n_\e n_\i\left(\TT_\i-\TT_\e\right) ?
\end{gather*}
 \end{multicols}
  \caption{Equations of the 10-moment 2-fluid Maxwell model.}
  The interspecies collisional terms
      $\inmag{\ZR_\s}$ and
      $\inmag{\ZQT_\s}$
  are generally ignored in this work.
  The intraspecies collisional terms
    $\indarkblue{\ZRT_\s}$
    are used in the simulations and play a critical role;
    we study reconnection as the isotropization rates
    $\ptime_\s$ are dialed between $0$ and $\infty$.
  The simulations neglect $\inred{\ZqT_\s}$, evidently causing
    late-time singularities.
 \label{fig:10mom2fluidMaxwell}
  \fhrule
\end{figure}

\begin{figure}[t]
\hfrule
 \textbf{\large Five-moment two-fluid Maxwell model}

 \textbf{Gas dynamics equations}
 \vspace{-1ex}
 {
 \begin{alignat*}{10}
    \mspace{-30mu}
    \dcs_t^\i
    \begin{bmatrix}
      \mdens_\i \\
      \mdens_\i\u_\i \\
      \mdens_\i\nrg_\i
    \end{bmatrix}
    &+
    \begin{bmatrix}
      0 \\
      \D p_\i + \LDivdPT{\i} \\
      {\Div(\u_\i p_\i) + \inred{\ZDiv(\Zu_\i\Zdotp\ZdPT_\i)}+\LDivqT{\i}}
    \end{bmatrix}
    &&= 
    \qdens_\i
    \begin{bmatrix}
      0 \\
      \E + \u_\i\times\B \\
      {\u_\i\dotp\E}
    \end{bmatrix}
    &&+
    \begin{bmatrix}
      0 \\
      \inmag{\ZR_\i} \\
      \inmag{Q_\i}
    \end{bmatrix}
   \\
    \mspace{-30mu}
    \dcs_t^\e
    \begin{bmatrix}
      \mdens_\e \\
      \mdens_\e\u_\e \\
      \mdens_\e\nrg_\e
    \end{bmatrix}
    &+
    \begin{bmatrix}
      0 \\
      \D p_\e + \LDivdPT{\e} \\
      {\Div(\u_\e p_\e) + \inred{\ZDiv(\Zu_\e\Zdotp\ZdPT_\e)}+\LDivqT{\e}}
    \end{bmatrix}
    &&= 
    \qdens_\e
    \begin{bmatrix}
      0 \\
      \E + \u_\e\times\B \\
      {\u_\i\dotp\E}
    \end{bmatrix}
    &&+
    \begin{bmatrix}
      0 \\
      \inmag{\ZR_\e} \\
      \inmag{Q_\e}
    \end{bmatrix}
 \end{alignat*}
 }

 \begin{multicols}{2}
\textbf{Maxwell's equations:}
\begin{alignat*}{5}
    \mspace{-30mu}
  &\partial_t \B + \curl\E = 0, &\quad& \D\dotp\B=0 \\
    \mspace{-30mu}
  &\partial_t \E - c^2\curl\B = -\J/\epsilon_0, &\quad& \D\dotp\E=\qdens/\epsilon_0
\end{alignat*}
 \textbf{Definitions:}
 \begin{gather*}
    \mspace{-30mu}
     \qdens_\i = en_\i, \quad
     \qdens_\e = -en_\e, \quad
     \qdens = \qdens_\i + \qdens_\e
     \\
    \mspace{-30mu}
     \J = \qdens_\i\u_\i + \qdens_\e\u_\e
     \\
    \mspace{-30mu}
     \glslink{dcs}{\Zdcs_t^s}\alpha:=\partial_t\alpha + \Div(\u_\s \alpha)
   \end{gather*}

 \columnbreak
 \textbf{Closures:}
 \begin{gather*}
    \vspace{-1ex}
    \hspace{-1.5em}
   \LdPT{\s} = -\Tviscosity\ddotp\SymB(\D\u)^\circ
   \\ \hspace{-1.5em}
     \LqT{\s} = -\TheatConductivity\dotp\D T
   \\ \hspace{-1.5em}
     \LR{\i}
    = -\LR{\e}
    = \qdens_\i \qdens_\e \glslink{Tresistivity}{\inmag{\ZTresistivity}}\dotp(\u_\e-\u_\i)
   \\ \hspace{-1.5em}
     \LQ{\s} = ?
     \LQf{\s} + \LQt{\s} (?) \\
   \\ \hspace{-1.5em}
     \LQf{\i} = \R_\i\dotp(\u_\e-\u_\i)\frac{m_\e}{m_\i+m_\e} (?) \\
   \\ \hspace{-1.5em}
     \LQt{\i} = \K_{\i\e}n_\i n_\e(T_\e-T_\i) (?)
 \end{gather*}
 \end{multicols}
  \caption{Equations of the 5-moment 2-fluid Maxwell model.}
  The interspecies collisional terms
      $\LR{\s}$ and
      $\LQ{\s}$
  are generally ignored in this work.
  The simulations also neglect $\LdPT{\s}$ and $\LqT{\s}$,
    evidently causing late-time singularities.
 \label{fig:5mom2fluidMaxwell}
  \fhrule
\end{figure}

\begin{figure}[t]
\hfrule
 \textbf{\large Ten-moment two-fluid MHD}
 \\
 \textbf{Pressure evolution:}
 \begin{alignat*}{10}
    n_\i d_t \TT_\i &+ \SymB(\PT_\i\dotp\D\u_\i) &&+ \LDivqT{\i}
      &&= \tfrac{q_\i}{m_\i}\SymB(\PT_\i\times\B) &&+ \LRT{\i} &&+ \LQT{\i}
 \\ 
    n_\e d_t \TT_\e &+ \SymB(\PT_\e\dotp\D\u_\e) &&+ \LDivqT{\e}
      &&= \tfrac{q_\e}{m_\e}\SymB(\PT_\e\times\B) &&+ \LRT{\e} &&+ \LQT{\e}
 \end{alignat*}

 \begin{multicols}{2}
 \textbf{mass and momentum:}
 \begin{gather*}
    \partial_t \mdens + \Div(\u\mdens) = 0
 \\ \mdens d_t \u + \Div(\PT_\i+\PT_\e+\PTd) = \J\times\B
 \end{gather*}
 \textbf{Electromagnetism}
   \begin{gather*}
      \partial_t \B + \curl\E = 0, \quad \Div\B = 0
     \\ \J = \mu_0\inv \curl\B
   \end{gather*}
 \textbf{Ohm's law}
   \begin{equation*} 
     \begin{aligned}
     \hspace{0.5em}
     \E &= \Tresistivity\dotp\J  + \B\times\u  
           + \tfrac{m_\i-m_\e}{e\mdens}\J\times\B 
       \\ &+ \tfrac{1}{e \mdens}
             \Div\left(\me\PT_\i-\mi\PT_\e\right)
       \\ &+ \tfrac{\mi\me}{e^2 \mdens}\left[
             \partial_t\J + \Div\left(\u\J+\J\u-\tfrac{\mi-\me}{e\mdens}\J\J\right)
           \right]
     \end{aligned}
   \end{equation*}

 \columnbreak
 \textbf{Definitions:}
 \begin{gather*}
   \begin{aligned}
   d_t &:=\partial_t + \u_\s\dotp\D
\\ \PTd&:=\mdens_i\w_\i\w_\i+\mdens_e\w_\e\w_e
   \end{aligned}
\\ \begin{aligned}
   \w_\i &= \frac{\me\J}{e\mdens}, \quad & \u_\i &= \u + \w_\i
\\ \w_\e &= \frac{\mi\J}{-e\mdens},\quad & \u_\e &= \u + \w_\e
   \end{aligned}
 \end{gather*}
 \textbf{Closures:}
 \begin{gather*}
    \vspace{-1ex}
    \hspace{-1em}
   \LRT{\s} = -\ptime_\s\inv \dPT_\s
   \\ \hspace{-1em}
     \LqT{\s} = -\tfrac{2}{5}\THeatConductivity_\s
         \dddotp\SymC\left(\Pshape\dotp\D\TT_\s\right) 
   \\ \hspace{-1em}
     -\LR{\i} = \LR{\e} = \ndens e\LTresistivity\dotp\J
   \\ \hspace{-1em}
     \LQT{\s} = \LQTf{\s} + \LQTt{\s} ?
 \end{gather*}
 \end{multicols}
 \caption{Equations of 10-moment 2-fluid MHD.}
 This system is derived but not simulated in this work.
 Simulations instead solve the two-fluid Maxwell equations
 with light speed intended to be sufficiently high to
 approximate this system.
  \fhrule
  \label{fig:10mom2fluidMHD}
\end{figure}

\begin{figure}[t]
\hfrule
 \textbf{\large Five-moment two-fluid MHD}
 \\
 \textbf{Pressure evolution:}
  \begin{alignat*}{8}
   \tfrac{3}{2}n d_t T_\i &+ p_\i \Div \u_\i &&+ \LdPT{\i}\ddotp\D\u_\i &&+ \LDivq{\i} &&= \LQ{\i}\\
   \tfrac{3}{2}n d_t T_\e &+ p_\e \Div \u_\e &&+ \LdPT{\e}\ddotp\D\u_\e &&+ \LDivq{\e} &&= \LQ{\e}
  \end{alignat*}

 \begin{multicols}{2}
 \textbf{mass and momentum:}
 \begin{gather*}
    \partial_t \mdens + \Div(\u\mdens) = 0
 \\ \mdens d_t \u + \Div(\PT_\i+\PT_\e+\PTd) = \J\times\B
 \end{gather*}
 \textbf{Electromagnetism}
   \begin{gather*}
      \partial_t \B + \curl\E = 0, \quad \Div\B = 0
     \\ \J = \mu_0\inv \curl\B
   \end{gather*}
 \textbf{Ohm's law}
   \begin{equation*} 
     \begin{aligned}
     \hspace{0.5em}
     \E &= \Tresistivity\dotp\J  + \B\times\u  
           + \tfrac{m_\i-m_\e}{e\mdens}\J\times\B 
       \\ &+ \tfrac{1}{e \mdens}
             \Div\left(\me(p_\i\idtens+\LdPT{\i})-\mi(p_\e\idtens+\LdPT{\e})\right)
       \\ &+ \tfrac{\mi\me}{e^2 \mdens}\left[
             \partial_t\J + \Div\left(\u\J+\J\u-\tfrac{\mi-\me}{e\mdens}\J\J\right)
           \right]
     \end{aligned}
   \end{equation*}

 \columnbreak
 \textbf{Definitions:}
 \begin{gather*}
   \begin{aligned}
   d_t &:=\partial_t + \u_\s\dotp\D
\\ \PTd&:=\mdens_i\w_\i\w_\i+\mdens_e\w_\e\w_e
   \end{aligned}
\\ \begin{aligned}
   \w_\i &= \frac{\me\J}{e\mdens}, \quad & \u_\i &= \u + \w_\i
\\ \w_\e &= \frac{\mi\J}{-e\mdens},\quad & \u_\e &= \u + \w_\e
   \end{aligned}
 \end{gather*}
 \textbf{Closures:}
 \begin{gather*}
    \vspace{-1ex}
    \hspace{-1em}
   \LdPT{\s} = -\Tviscosity\ddotp\SymB(\D\u)^\circ
   \\ \hspace{-1em}
     \LqT{\s} = -\TheatConductivity\dotp\D T
   \\ \hspace{-1em}
     -\LR{\i} = \LR{\e} = \ndens e\LTresistivity\dotp\J
   \\ \hspace{-1em}
     \LQ{\s} = \LQf{\s} + \LQt{\s} ?
 \end{gather*}
 \end{multicols}
 \caption{Equations of 5-moment 2-fluid MHD.}
 This system is derived but not simulated in this work.
 Simulations instead solve the two-fluid Maxwell equations
 with light speed intended to be sufficiently high to
 approximate this system.
  \label{fig:5mom2fluidMHD}
  \fhrule
\end{figure}

{
\def\w{\id_{\wedge}}
\def\l{\id_{\Vert}}
\def\p{\id_{\perp}}
\def\v{\hgf}
\begin{figure}[t]
 \hfrule
 \textbf{\large Ten-moment relaxation closure}
 \begin{multicols}{2}
 \textbf{Implicit closure 
     \eqref{implicitHeatFluxTensorClosureCopy}}
   \begin{gather*}
     \hspace{-2em}
     {\small
     \qT + \SymC(\hgf\b\times\qT) =
       -\tfrac{2}{5}\heatConductivity \SymC\left(\Pshape\dotp\D\TT\right),
     }
   \end{gather*}
  \textbf{Definitions:}
   \begin{gather*}
     \begin{aligned}
     \b&:=\B/|\B| \\
     \gyrofreq&:=\tfrac{q_\s}{m_\s}|\B| \\
     \hgf&:=\htime\gyrofreq \\
     \pgf&:=\ptime\gyrofreq \\
     \htime &:= \ptime/\Pr
     \end{aligned}
   \end{gather*}
 \textbf{Explicit closure 
   \eqref{qTclosureForm}}
  \begin{gather*}
    \hspace{-2em}
    {\small
    \begin{aligned}
    \qT&= -\tfrac{2}{5}\heatConductivity\tHeatConductivity
          \dddotp\SymC\left(\Pshape\dotp\D\TT\right)
    \\ &= -\tfrac{2}{5}\heatConductivity\SymC\left(\tHeatConductivity'
          \dddotp\SymC\left(\Pshape\dotp\D\TT\right)\right)
    \end{aligned}
    }
  \end{gather*}
  where equation \eqref{3MexpansionCopy2} gives
   \begin{equation*}
    {\small
    \begin{aligned}
     \tHeatConductivity
        &=\left(\l^3 + \tfrac{3}{2}\l(\p^2+\w^2)\right)
       \\& + \frac{3}{1+\v^2}\left(\p\l^2 - \v \w\l^2\right)
       \\& + \frac{3}{1+4\v^2}\left(\frac{\p^2-\w^2}{2}\l - 2\v \w\p\l\right)
      \\  &+ (k_0 \p^3 + k_1 \w\p^2 + k_2 \w^2\p + k_3 \w^3),
    \end{aligned}
    }
   \end{equation*}
   \textbf{Definitions:}
   \begin{gather*}
     \begin{aligned}
      \l&:=\b\b \\
      \p&:=\id-\l \\
      \w&=\id\times\b
     \end{aligned}
   \end{gather*}
   \textbf{Relations:}
   \begin{gather*}
     \tfrac{2}{5}\heatConductivity = \frac{p\ptime}{m\Pr}
   \end{gather*}
   The remaining coefficients are \eqref{heatCoefsM}
  \begin{equation*}
   {\small
   \begin{aligned}
     k_3 &:= \frac{-6\v^3}{1+10\v^2+9\v^4} = -\tfrac{2}{3}\v\inv + \O(\v^{-3}), \\
     k_2 &:= \frac{6\v^2 + 3\v(1+3\v^2) k_3}{1+7\v^2} = \O(\v^{-2}) , \\
     k_1 &:= \frac{-3\v + 2\v k_2}{1+3\v^2} = -\v\inv + \O(\v^{-3}), \\
     k_0 &:= 1+\v k_1 = \O(\v^{-2}).
   \end{aligned}
   }
  \end{equation*}
 \end{multicols}
 \caption{Ten-moment Gaussian-BGK closure.}
    All tensor products in the formula 
    for $\tHeatConductivity$ are splice symmetric products.
    The formula for 
    $\tHeatConductivity'$ is exactly the same as for $\tHeatConductivity$
    except that the products may be taken simply to be splice products
    as in \eqref{qTclosure}.
  \label{fig:10momCoef}
  \fhrule
\end{figure}

\begin{figure}[t]
 \hfrule
 \textbf{\large Five-moment relaxation closure}
 \begin{multicols}{2}
 \textbf{Implicit closure
     \eqref{implicitHeatFluxClosureCopy},
     \eqref{implicitDeviatoricStressClosureCopy}}
   \begin{gather*}
     \hspace{-1em}
     \begin{aligned}
       &\q + \hgf\b\times\q = -\heatConductivity \D T, \\
      &\dPT + \SymB(\pgf\b\times\dPT) = -\viscosity 2 \dstrain
     \end{aligned}
   \end{gather*}
  \textbf{Definitions:}
   \begin{gather*}
     \begin{aligned}
     \b&:=\B/|\B| \\
     \gyrofreq_\s&:=\tfrac{q_\s}{m_\s}|\B| \\
     \hgf&:=\htime_\s\gyrofreq_\s \\
     \pgf&:=\ptime_\s\gyrofreq_\s \\
     \htime &:= \ptime/\Pr
     \end{aligned}
   \end{gather*}
 \textbf{Explicit closure
   \eqref{expHeatFluxInternal}:}
   \begin{gather*}
     \begin{aligned}
       \q &= -\heatConductivity\theatConductivity\cdot\D T, \\
     \dPT &= -2\viscosity \tviscosity\ddotp \dstrain
       \\ &= -2\viscosity\Sym\left(\tviscosity'\ddotp \dstrain\right)
     \end{aligned}
   \end{gather*}
   where by
     \eqref{theatConductivityClosure} and
     \eqref{tviscosityClosure},
   \begin{gather*}
     \hspace{-1em}
    {\small
    \begin{aligned}
    \theatConductivity :=& \l + \frac{1}{1+\v^2}\Big(\p - \v \w\Big), \\
    \tviscosity :=& \left(\l\l+\frac{\p\p+\w\w}{2}\right) \\
      +& \frac{2}{1+\v^2}(\p\l - \v \w\l) \\
      +& \frac{1}{1+4\v^2}\left(\frac{\p\p-\w\w}{2} -  2\v \w\p\right).
    \end{aligned}
    }
    \label{5momClosureCoef}
   \end{gather*}
   \textbf{Definitions:}
   \begin{gather*}
     \begin{aligned}
      \l&:=\b\b \\
      \p&:=\id-\l \\
      \w&=\id\times\b
     \end{aligned}
   \end{gather*}
   \textbf{Relations:}
   \begin{gather*}
     \begin{aligned}
       \viscosity_\s &= p_\s\ptime_\s \\
       \tfrac{2}{5}\heatConductivity_\s &= \frac{\viscosity_\s}{m_\s\Pr_\s}
     \end{aligned}
   \end{gather*}
 \end{multicols}
 \caption{Five-moment Gaussian-BGK closure.}
    The tensor products in the formula for $\tviscosity$ in \eqref{5momClosureCoef}
    are splice symmetric products.  The formula for
    $\tviscosity'$ is exactly the same as for $\tviscosity$
    except that the products may be taken simply to be splice products.
  \label{fig:5momCoef}
  \fhrule
\end{figure}

\begin{figure}[t]
 \hfrule
 \textbf{\large Relaxation coefficients}
 \begin{multicols}{2}
 \textbf{Relaxation periods}
   \begin{gather*}
     \begin{aligned}
      \ptime_0 &:= \frac{12 \pi^{3/2}}{\ln\Lambda}\left(\frac{\epsilon_0}{e^2}\right)^2,
      \qquad \eqref{baseIsoPeriod}
   \\ \ptime'_{\i\i} &= \ptime_0 \sqrt{m_\i}\frac{T_\i^{3/2}}{n_\i}\ \ \hbox{ or}
         \\          &= \ptime_0 \sqrt{m_\i}\frac{\sqrt{\det(\TT_\i)}}{n_\i},
   \\ \ptime'_{\e\e} &= \ptime_0 \sqrt{m_\e}\frac{T_\e^{3/2}}{n_\e}\ \ \hbox{ or}
         \\          &= \ptime_0 \sqrt{m_\e}\frac{\sqrt{\det(\TT_\e)}}{n_\e},
   \\ \ptime_\i &= .96 \ptime'_{\i\i},
   \\ \ptime_\e &= .52 \ptime'_{\e\e},
     \end{aligned}
   \end{gather*}

 \textbf{Temperature-determined coefficients}
   \begin{gather*}
     \begin{aligned}
      \viscosity_\i =& \ptime_\i p_\i, 
   \\ \viscosity_\e =& \ptime_\e p_\e, 
   \\ \tfrac{2}{5}\heatConductivity_\i =& \frac{\viscosity_\i}{m_\i \Pr_\i},
   \\ \tfrac{2}{5}\heatConductivity_\e =& \frac{\viscosity_\e}{m_\e \Pr_\e},
   \\ {\Pr}_\i =& .61 \approx \tfrac{2}{3},
   \\ {\Pr}_\e =& .58 \approx \tfrac{2}{3},
   \\ \resistivity_0 :=& \frac{m_\e}{e^2 n_\e}\frac{\sqrt{2}}{\ptime'_{\e\e}},
   \\ \resistivity_\para :=& .51 \resistivity_0,
   \\ \lim_{\pgf\to\infty}\resistivity_\perp =& \resistivity_0
     \end{aligned}
   \end{gather*}

  These closures for a hydrogen plasma have been derived by Braginskii
  and others (see section \ref{BraginskiiClosureCoefficients})
  using a Coulomb collision operator
  and assuming a strongly collisional plasma.
  Fast magnetic reconnection is a collisionless phenomenon
  and therefore substantial deviation is expected from
  these closures.  Therefore, we are content with 
  rough and simple approximations of Braginskii's
  coefficients.  In particular, 
  there is likely no great loss in 
  assuming a scalar resistivity with value $\resistivity_0$
  or likewise assuming isotropic relaxation for other
  nondiffusive closure coefficients (in comparison
  to the Braginskii closure).
 \end{multicols}
 \caption{Relaxation coefficients}
 \fhrule
\end{figure}
}
 \clearpage 

%
%


%% file: chap2.tex
\def\i{\mathrm{i}}
\def\v{\mathbf{v}}
\def\c{\mathbf{c}}
\chapter{Models of Plasma}\label{models}

The purpose of this chapter is primarily to develop the equations
listed at the end of the previous chapter.


\glsadd{Div}
\glsadd{D}
\glsadd{ebas}

\section{Overview}


We consider plasma to consist of charged particles of two
types of particles (called \first{species}) subject only to
electromagnetic forces. We ignore gravitational effects, quantum
effects, and nuclear forces.

Therefore, as our fundamental standard of truth we assume that the
electromagnetic field is governed by Maxwell's equations
of electromagnetism and that particle motion is governed by
the Lorentz force law and the relativistic version of
Newton's second law.

Our practical standard of truth will be the nonrelativistic
Boltzmann equation.  This involves two simplifications:
(1) we represent particle distributions with a continuum
distribution and (2) we neglect relativity.
Each of these simplifications entails issues and problems.

Regarding the first simplification, the Boltzmann equation is
incomplete as a standard of truth until a collision operator is
specified. One of the main points of this dissertation is to
study the dependence of reconnection on the choice of collision
operator.

Regarding the second simplification, we expect fundamental
physical laws to be invariant under change of reference frame.
Physical laws are \defining{Lorentz-invariant} if they remain
unchanged under Lorentz transformations. A \defining{Lorentz
transformation} is a linear transformation of space-time that
leaves light speed invariant and respects the direction of time
and the scale and orientation of space-time. Physical laws are
\defining{Galilean-invariant} if they remain unchanged under
Galilean transformations. A \defining{Galilean transformation}
is a linear transformation of space-time that leaves time
and distance invariant and respects the direction of time
and the scale and orientation of space. We say that physical
laws \index{relativity principle} \emph{satisfy a relativity
principle} if they are Galilean-invariant or Lorentz-invariant.
We say that physical laws are \first{relativistic} (in the
sense of Einstein's theory of special relativity) if they are
Lorentz-invariant.

The non-relativistic Boltzmann-Maxwell system is an intermediate
system which is neither fully Lorentz-invariant nor fully
Galilean-invariant. Maxwell's equations (with prescribed
Lorentz-invariantly-defined current and charge density) are
Lorentz-invariant. The nonrelativistic Boltzmann equation (with
prescribed Galilean-invariantly-defined electromagnetic field)
are Galilean-invariant. One of these systems must be modified for
the system to satisfy a relativity principle.

The simplifying assumptions of MHD provide a fully
Galilean-invariant system, essentially by taking the light speed
to infinity. MHD models that admit fast magnetic reconnection
admit fast waves and require implicit numerical methods, which
are not easy to implement.

The non-relativistic Boltzmann-Maxwell system that I have
chosen to implement is something of a toy model. It leads
to fluid models that are conceptually simple and can be
implemented with explicit numerical methods. Nonrelativistic
two-fluid-Maxwell models arise from taking moments of the
non-relativistic Boltzmann-Maxwell system. My simulations study
the ability of non-relativistic two-fluid-Maxwell models to
agree with the two-species non-relativistic Boltzmann-Maxwell system.
In the relativistic regime one would instead prefer a fully relativistic
two-fluid-Maxwell system. In the low-speed limit one would
prefer an MHD system.  I therefore discuss the equations
and properties of MHD.
Simulating with them would be a natural extension of this dissertation.

\subsection{Kinetic models}

We generally refer to models based on the evolving of particle
states according to fundamental laws as 
\defining{particle models}.
\first{Particle-in-cell (PIC)} methods simulate plasma by
a rescaling of the fundamental laws which leaves macroscopic
physical quantities basically unchanged but reduces the number
of particles to a computationally feasible number.
At any given time the state of a particle is specified
by its position in \defining{phase space},
defined to be the pair $(\xb,\v)$ consisting of the
particle's position \gls{xb} and velocity \gls{v}.

Particle models are a type of kinetic model. \defining{Kinetic
models} evolve a represention of particle positions and
velocites. Continuum kinetic models represent particles via
\first{particle density functions} which specify the density of
particles with a particular location and velocity. The simplest
general continuum kinetic model is the \gls{Boltzmann equation}.
The Boltzmann equation represents each species $\s$
with a particle mass density $\glslink{f}{\f_\s}(t,\xb,\v)$
which is a function of time ($t$), spatial location (\gls{xb})
and particle velocity (\gls{v}). The Boltzmann equation depends
on the choice of a \defining{collision operator}
$\glslink{C}{C_\s}$, which specifies a model for how particles collide.
If the collision operator is entirely neglected, then the
Boltzmann equation is instead referred to as the \defining{Vlasov
equation}.

Simulation with kinetic models is highly expensive, because of
the need to evolve a representation of the detailed distribution
of particle velocities. Particle simulations typically track
millions of particles, and a typical Vlasov simulation might be
two orders of magnitude more expensive than a comparable particle
simulation.  The computational expensiveness of kinetic models
is one of the motivations for developing \first{fluid models}.

\subsection{Fluid models}

\subsubsection{Moments}

Fluid models of plasma or gas allow simulation with greatly
reduced expense. \defining[fluid models]{Fluid models} assume
that the distribution of particle velocities is characterized by
a small set of parameters which are typically \first{moments}
of the distribution. An \defining{$n$th order moment} specifies
at each time $t$ and location $\xb$ the average value (averaged
over all particles momentarily near $\xb$) of a sum of products
of $n$ components of $\v$. The most important examples of
moments are the (physically) \defining{conserved moments}:
mass density, momentum density, and energy density. The
conserved moments are conserved in collisions. Collisions
cause the velocity distribution to trend toward an equilibrium
distribution characterized by its conserved moments. An
equilibrium velocity distribution is bell-shaped and is known
as a \first{Maxwellian distribution}. We refer to the conserved
moments as \first{Maxwellian moments}.

\subsubsection{Five-moment model.} The simplest and most well-justified
model of a gas is the \defining{five-moment model}, which
represents the state of the gas in terms of the conserved
moments. In the three-dimensional space of our physical world
the conserved moments are the ``first five moments'': the mass
density $\glslink{mdens}{\Zmdens_\s}:=\int_{\v\in\reals^3}
\f_\s$ (which is the zeroth-order moment), the momentum density
$\glslink{M}{\ZM_\s}:=\mdens_\s\u_\s:=\int_{\v\in\reals^3}
\f_\s \v$ (three first-order moments, one for each component
of the velocity/momentum), and the energy density
$\glslink{Nrg}{\ZNrg}_\s:=\mdens_\s \glslink{nrg}{\Znrg}_\s:=\int_{\v\in\reals^3}
\f_\s (v_1^2+v_2^2+v_3^2)/2$ (the conserved second-order moment).
Here \gls{reals} denotes the real numbers.

\subsubsection{Primitive and conservation variables.} Note that the
fluid velocity $\glslink{u}{\Zu_\s}$ is defined to be the average velocity
of particles of species $\s$. We define the \defining{thermal velocity}
$\glslink{c}{\Zc_\s}:=\v-\u_\s$ to be the velocity of a
particle in the reference frame defined by the local fluid
velocity. The first moments of $\c_\s$ are by definition
zero; other moments of the particle velocity $\v$ are
equivalent to moments of the thermal velocity $\c$. For
example, $\mdens_\s e_\s = \mdens_s|\u_\s|^2/2 + (3/2)p_\s$,
where $p_\s:=\int_{\c\in\reals^3} \f_\s |\c|^2/3$ is the
pressure. The fluid velocity $\u_\s$ and the non-first-order
moments of $\c$ (including density and pressure) are called
\defining{primitive variables}. In contrast, we refer to moments
of $\v$ as \defining{conservation variables}, not because they
are all conserved, but because they would be conserved except for
nondifferentiated source terms that represent particle collisions
and interaction with the electromagnetic field. I remark that
fluid closures are naturally defined in terms of moments of $\c$;
closures thus defined in the reference frame of the fluid are by
definition independent of reference frame.

\subsubsection{Ten-moment model.}
\first{Extended-moment fluid models} use additional moments
of order two or higher to represent the particle velocity
distribution more accurately. For plasmas with nonzero
magnetic field, the simplest extended-moment fluid model is the
\emph{\defining{gyrotropic} (six-moment) model.} Instead of
evolving a single energy moment (or scalar pressure), the
gyrotropic model evolves a pressure component parallel to the
magnetic field and a pressure component perpendicular to the
magnetic field. This representation is natural in a strongly
magnetized plasma due to the rapid rotation of particles around
field lines and the relatively slow trend to equilibrium of
velocity components parallel to the magnetic field.

In unmagnetized regions of low-collisionality velocity
distributions do not have a magnetic field to align with and
trend slowly to equilibrium. In such regions pressure can become
highly anisotropic or agyrotropic. The \defining{ten-moment
model} evolves all six independent second-order monomial moments,
which define a full pressure tensor, allowing it to represent highly
anisotropic and agyrotropic pressures.

\subsubsection{Extended-moment models of plasma for reconnection}

The use of ten-moment gas dynamics is particularly relevant for
nearly collisionless magnetic reconnection near a magnetic null
point or null line. In the general case of steady two-dimensional
non-resistive magnetic reconnection symmetric
under 180-degree rotations about an X-point,
pressure cannot be gyrotropic in the vicinity
of the X-point (see section \ref{agyrotropy}).
Moreover, in the special case of antiparallel
reconnection with reflectional symmetry across both axes,
simulations of the GEM problem with a near-collisionless
ten-moment model develop strong agyrotropy in the vicinity of the
X-point. A diffusive closure for the viscous pressure in terms of
the fluid velocity gradient allows the five-moment (isotropic)
model to accurately simulate small perturbations from isotropy
and allows the six-moment (gyrotropic) model to accurately
simulate small perturbations from gyrotropy. But these diffusive
closures are based on the assumption of high collisionality.

In contrast, the use of a ten-moment model with relaxation
of the pressure tensor toward isotropy or gyrotropy not only
agrees with the viscous isotropic and gyrotropic models in the
high-collisionality limit, but also allows the isotropization
rate, equivalent to the collision rate, to be dialed over the
full range between zero and infinity, making it appropriate for
study and simulation of low-collisionality plasmas.
In particular, \textit{\textbf{the ten-moment model seems to be
the simplest fluid model
for which one can study the dependence of magnetic reconnection
on the collision operator as the collision rate is taken to zero}}.
This may give qualitiative insight into the behavior
of kinetic models as the collision operator is taken to zero.

In chapter \ref{MagneticReconnection} we argue that to simulate
sustained magnetic reconnection it is necessary
to admit heat flux at least in some cases. Therefore we derive
a ten-moment heat flux closure in this chapter.

Chapter \ref{MagneticReconnection} also includes a ten-moment
adiabatic simulation of a pair plasma version of the
\mention{GEM magnetic reconnection challenge problem}
as the collision rate is dialed between zero and infinity.
The absence of heat flux in the model results in instabilities
when simulating the GEM problem past the time of the peak
reconnection rate.

I propose doing a similar study of this problem for the
Boltzmann equation using a \first{Gaussian-BGK}
collision operator.

A study of vanishing collisionality using the ten-moment model is
admittedly artificial from a physical perspective. In particular,
as the collision rate goes to zero the heat flux must go to
infinity unless one artificially takes the Prandtl number to
infinity.  In fact, to study the effects of vanishing entropy production
in the low-collisionality limit it would be necessary to take
the heat flux to \emph{zero} as the collision rate goes to zero.
This is highly unphysical.

In the zero-collisionality limit of the Boltzmann equation
the rate of isotropization (i.e.\ the rate of 
decay of the deviatoric pressure)
and the rate of decay of the heat flux both go to zero.
The ten-moment model is able to study vanishing
isotropization because it \emph{evolves deviatoric pressure}.
But it is not able to study vanishing heat flux decay
because it uses a diffusive closure for the heat flux
(in terms of spatially differentiated state variables).

To study vanishing heat flux decay requires a model which
\emph{evolves heat flux}.  Heat flux is a third-order moment.
Models which evolve such superquadratic moments are called
\first{higher-moment models}.  

The challenge in formulating a fluid moment model which
evolves heat flux is to define a \emph{hyperbolic
closure} for the unevolved moments in the equations.
For the ten-moment model a hyperbolic closure is
to set the heat flux to zero.  This closure
guarantees that pressure and density will not go
negative and that the equations remain hyperbolic
(that is, well-posed, meaning that the solution depends
continuously on the data) as the solution evolves. 
For higher-moment models, however,
it is not clear how to define hyperbolic closures.
In particular, in order to approach Galilean-invariance,
hyperbolic higher-moment closures have to admit
wave speeds that approach infinity for states
arbitrarily close to equilbrium \cite{mcdonald09, mcdonald11}.

In contrast, a \mention{Lorentz-invariant hyperbolic
higher-moment closure} would presumably bound wave speeds by the
speed of light. A natural way to formulate such a closure is
to choose an assumed form for the family of possible particle
velocity distribution functions (just as one assumes a Maxwellian
distribution for the five-moment hyperbolic closure and a
Gaussian distribution (equation \eqref{GaussianDistributionCopy})
for the ten-moment hyperbolic closure). It is not known how to
choose such a distribution so that solutions remain hyperbolic
and \first{realizable} (that is, the evolved moments are those of
a distribution of the assumed form). To get a closure in explicit
form one needs to be able to compute in closed form the integrals
giving the unevolved moments in the equations. Torrilhon has used
a Pearson-IV velocity distribution family to obtain closed-form
closures for the Galilean-invariant Boltzmann equation
\cite{torrilhon09}, but for the Lorentz-invariant Boltzmann
equation, except in the case of a Maxwellian distribution of
proper velocities, which has the unique property that the
distribution may be represented as the product of one-dimensional
distributions, it seems highly unlikely that one can formulate
an assumed family of distributions for which the integrals
can be computed in closed form and for which solutions remain
hyperbolic and the moments remain realizable by a distribution.
One therefore would have to resort to numerical quadrature for
these integrals over three-dimensional velocity space, which
would be very expensive unless some simplifying technique can be
identified.

The ten-moment model not only allows a fluid study of
vanishing pressure isotropization (which generally implies that the pressure
is \emph{not} isotropic) but has the additional general benefit
that one can study nonvanishing viscosity (i.e.\ nonvanishing
deviatoric pressure) efficiently with an explicit numerical method.

Likewise, a higher-moment fluid model would not only allow a fluid study of
vanishing heat flux decay (which generally implies that the heat flux
is \emph{nonvanishing}) but would have the additional general benefit
that one could study nonvanishing heat conductivity (i.e.\ nonvanishing
heat flux) efficiently with an explicit numerical method.

In summary, what I seek is a Lorentz-invariant higher-moment
plasma model; it would serve as a fluid analog of the
relativistic Boltzmann equation with Gaussian-BGK closure.
The ultimate extension of this dissertation would be to simulate
fast magnetic reconnection in the low-collisionality limit with
a Lorentz-invariant higher-moment plasma model using an
explicit numerical method.  This is my dream.

\section{Equations of kinetic models}

\def\F{\mathbf{F}}
\subsection{Particle models (PIC)}

Particle models of plasma are based on the fundamental laws of
classical electrodynamics:\\
\textbf{Newton's second law},
\begin{align*}
  m_\p d_t \pv_\p = \F_\p, && \v_\p := d_t \xb_\p,
\end{align*}
the \textbf{Lorentz force law},
\begin{gather*}
  \F_\p = q_\p\left(\E + \v\times\B\right),
\end{gather*}
and \textbf{Maxwell's equations},
\begin{align}
   \label{MaxwellsEqns}
   \partial_t \B + \curl\E &= 0,                &\Div\B &= 0,
\\ \partial_t \E - c^2\curl\B &= -\J/\epsilon_0,  &\Div\E &= \qdens/\epsilon_0,
\end{align}
where \gls{B} is magnetic field, \gls{E} is electric field,
and current \gls{J} and charge density
\footnote{Electromagnetism books often use $\rho$ for charge density,
but I'm already using this for mass density, so I took the
next letter of the alphabet, as does, for example, Balescu \cite{Balescu88}.}
\gls{qdens}
are given by
\begin{align*}
   \J& = \sum_\p q_\p \v_\p \delta_{\xb_\p},
    &\qdens&= \sum_\p q_\p \delta_{\xb_\p};
\end{align*}
here
\gls{c } is the speed of light,
\gls{epsilon_0} is permittivity of space,
\gls{p} is particle index,
$\glslink{m }{m}_p$ is particle mass,
$\glslink{q }{q}_p$ is particle charge,
$\glslink{xb}{\Zxb}_p$ is particle position,
$\glslink{v}{\Zv}_p$ is particle velocity,
and $\glslink{pv}{\Zpv}_\p$ (i.e.\ $\gamma_\p \v_\p$) is the proper velocity,
where $\glslink{gamma}{\gamma} := \left(1-(v/c)^2\right)^{-1/2}$, the Lorentz factor,
is the rate of elapse of time with respect to proper time,
where the elapse of proper time is defined in the reference
frame of a particle moving with velocity $\v$.

In case $|v|\ll c$, $\gamma\approx 1$
and we can use the nonrelativistic
approximation $\pv\approx\v$.

\subsection{Boltzmann/Vlasov equations}


The \gdefining{Boltzmann equation}
or \gdefining{kinetic equation} is written
\def\s{\mathrm{s}}
\def\fs{\f_\s}
\def\a{\mathbf{a}}
\def\as{\a_\s}
\begin{gather}
  \partial_t\fs+\D_{\xb}\dotp(\v \fs)+\D_{\pv}\dotp(\as \fs) = C_\s/\gamma;
  \label{BoltzmannEquation}
\end{gather}
here
$\glslink{f}{\fs}(t,\xb,\pv)$ is the particle density function of species $\s$
as a function of time and \defining{phase space} $(\xb,\pv)$,
$\as = d_t \pv = (q_\s/m_\s)(\E+\v\times\B)$ is particle acceleration,
and we recall that
$\v=d_t\xb$ is velocity and
$\pv=\gamma\v$ is proper velocity.

The Maxwell source terms are
provided by the relations
\begin{align*}
   \J& = \sum_\s q_\s \int \fs \v d\pv,
  &\qdens&:= \sum_\s q_\s \int \fs d\pv.
\end{align*}

The collision operator $C_\s$ specifies the rate of change of
the particle distribution due to particle collisions
and depends on the choice of collision model.

If we set $C_\s=0$ then we get the \defining{Vlasov equation}.
If we regard $\fs$ as a linear superposition of spike functions
then the Vlasov system with Maxwell's equations is equivalent to
the fundamental equations of electrodynamics written in terms
of individual particles. This is because at the microscale
level particles do not actually collide; they interact with
an electromagnetic field that becomes highly irregular in the
vicinity of another particle.

Use of the Vlasov equation to resolve interaction of individual
particles would require impossibly high resolution and would
be enormously expensive computationally. (Actually, one would
need a mechanism to maintain the singular electromagnetic field
near individual particles, which would force it to be a particle
method.) For this reason, we partition the electromagnetic field
into a relatively large but smoothly varying \first{macro-scale} field
that governs the interaction of particles with the plasma as a
whole and a relatively small but highly irregular \first{microscale
field} that mediates localized particle interactions. Rather than
evolve the microscale electromagnetic field and microscale particle
distribution function, we assume a smooth electromagnetic field and smooth
particle distribution functions and model particle interactions via
a collision operator. The \first{Coulomb collision operator} assumes that
particle interactions occur pairwise and randomly (independently
and in proportion to particle density) and are governed by the
inverse square force law of electrostatics.  Assuming that
the particle density function is smooth (reflecting a smooth
estimate of the probability of the number of particles in a local
region), the value of the collision operator will also be smooth.

Thus, while in (microscale) theory the Boltzmann equation is an
approximation to the Vlasov equation, in (macro-scale) practice
the Vlasov equation is an approximation to the Boltzmann equation
under the assumption that particle collisions may be neglected
for the phenomena or scales of interest.


\section{Fluid balance equations of Galilean-invariant gas dynamics}

Fluid equations are obtained from the Boltzmann equation by evaluating
velocity moments and assuming closures for the unevolved moments
that arise in terms of the evolved moments.
In this section we evaluate velocity moments of the Boltzmann
equation to obtain generic balance laws for mass, momentum, and
various forms and consequences of energy evolution,
including pressure, temperature, and entropy evolution.
The fluid equations that we derive are Galilean-invariant rather
than Lorentz-invariant.
We do not consider closures in this section.

\subsection{Notational conventions for fluid moments}

\subsubsection{Velocity-space integrals}
Fluid quantities are functions of time and spatial
location that are defined in terms of integrals of
a particle distribution function over velocity space. 
By default integrals over velocity
are taken over the full domain of possible velocities:
$\intv :=\int_{\v\in\reals^3}$.

Let $\f(t,\xb,\v)$ be the particle mass distribution of a gas.
The primary fluid quantities that describe the gas are
are the density $\mdens(t,\xb):=\intv \f$
and the velocity defined by $\mdens\u(t,\xb):=\intv \v \f$.
\first{Primitive variables} are most naturally defined
in terms of moments of the thermal velocity
$\c(t,\xb,\v):=\v-\u(t,\xb)$.  When integrating
moments of $\c$, I typically write
$\intc$ for $\intv$, since 
$\int_{\c\in\reals^3}=\int_{\v\in\reals^3}$.

\def\vint#1{\langle #1 \rangle}
\def\mean#1{\langle #1 \rangle}
\def\dens#1{\langle #1 \rangle}
Authors often use angle brackets to denote velocity-space
integrals, following one of three conventions:
\begin{enumerate}
  \item[1.] $\vint{\f\chi} := \vint{\f\chi}_{\footnotesize \int} := \intv \f\chi$.
    That is, $\vint{\dotp}$ simply denotes integration
    over velocities.  In this notation the distribution
    is explicit and may thus be replaced with another distribution
    or the collision operator.
  \item[2.] $\dens{\chi} := \dens{\chi}_f := \intv \f\chi$.
    That is, $\vint{\chi}(t,\xb)$ denotes the density per volume
    of the quantity $\chi$.  This definition is intermediate between
    the other two.
  \item[3.] $\mdens\glslink{mean}{\mean{\chi}} := \mdens\overline{\mean{\chi}} := \intv \f\chi$.
    That is, $\vint{\chi}(t,\xb)$ denotes the statistical average
    value of $\chi$ at time $t$ over all particles at point $\xb$.
    It is a density per mass.
\end{enumerate}
In this work we adopt convention (3) and simply write integrals
for the other two cases.

\subsubsection{Tensor products}
\label{TensorProducts}

In this document products of tensors are by default tensor
products and powers of tensors are by default (unsymmetrized
and uncontracted) tensor powers.

Let $A$ be an $n$th order tensor and let $B$ be an $m$th order tensor.
Then by $AB$ we mean the tensor product $A\otimes B$.
In contrast, by $AB$ some authors mean the \emph{symmetric}
tensor product $\Sym(AB)$, where $\Sym$ denotes the average
over all permutations of subscripts.  Except in section
\ref{ExplicitClosures} I always explicitly indicate symmetrization.

Some authors define the vee product by $A\vee B:=\Sym(AB)$,
just as some authors define the wedge product by
$A\wedge B:=\Ant(A\otimes B)$, where $\Ant$ denotes the average over
all permutations of subscripts each multiplied by the sign
of its permutation.
Most authors, however, such as John Lee \cite{book:Lee03},
define the wedge product by
$A\wedge B:=\frac{(n+m)!}{n!m!}\Ant(A\otimes B)$;
in case $A$ and $B$ are antisymmetric tensors this
defines their wedge product to be the sum over all
``distinguishable'' signed permutations: specifically,
there exists a partition of the $(n+m)!$ signed
permutations of $A\otimes B$ into $\frac{(n+m)!}{n!m!}$
classes each containing $n!m!$ identical terms.
For the antisymmetric product Lee instead uses the notation
$A\barwedge B:=\Ant(A\otimes B)$.

Analogously, I define the vee product by
\begin{gather}
  A\glslink{vee}{\vee} B:=\frac{(n+m)!}{n!m!}\Sym(A\otimes B);
  \label{veeDefinition}
\end{gather}
in case $A$ and $B$ are symmetric tensors this
defines their vee product to be the sum over all
``distinguishable'' permutations: specifically,
there exists a partition of the $(n+m)!$ 
permutations of $A\otimes B$ into $\frac{(n+m)!}{n!m!}$
classes each containing $n!m!$ identical terms.
For the symmetric product I instead use the notation
$A\glslink{veebar}{\veebar} B:=\Sym(A\otimes B)$.

\subsubsection{Convective derivatives}
It is natural to describe the evolution a fluid quantity
$a(t,\xb)$ in terms of its rate of change as an observer moves
with the fluid.
Given a velocity field $\u(t,\xb)$,
we define the \defining{convective derivative}
of $a$ relative to $\u$ to be
$d_t^\u a := \partial_t a  + \u\dotp\D_\xb a$,
the rate of change of $a$ along a trajectory
moving with the fluid.  We write $d_t:=d_t^\u$
when the relevant velocity field $\u$ is clear;
$\u$ is generally understood to be the velocity
of the fluid under consideration.  Thus,
when describing quantities defined in terms of
an individual species, $\u$ is by default taken to be
the velocity of the species fluid, and 
when describing quantities summed over all species
$\u$ is by default taken to be the overall bulk
fluid velocity including both species.

\def\dcs{\overline{\delta}}
For convenience define the \defining{``bulk derivative''} by
$\dcs_t \alpha := \partial_t\alpha + \Div(\u\alpha)$ for all $\alpha$.
Suppose that mass density $\mdens(t,\xb)$
is conserved:
$\dcs_t\mdens = \partial_t +\Div(\mdens\u) = 0$.
Then the convective derivative of $a$ is
related to the bulk derivative
of $\mdens a$:
\begin{gather*}
  \mdens d_t a = \dcs_t(\mdens a) = \partial_t(\mdens a) + \Div(\mdens\u a).
\end{gather*}
Here $a$ would be a density per mass and
$\mdens a$ would be a density per volume.
While the convective derivative is conventional
and widely used, the bulk derivative does not enjoy
a common conventional name or notation.
For this reason densities per mass (or per particle
number) are often preferred.

\subsection{Conservation moments}

Our starting point for deriving Galilean-invariant
fluid equations is the Galilean-invariant Boltzmann equation
(equation \eqref{BoltzmannEquation} with $\gamma$ taken as $1$),
\begin{gather}
  \label{GalBoltzmannEqn}
  \partial_t\fs+\D_\xb\dotp(\v \fs)+\D_\v\dotp(\as \fs) = C_\s,
\end{gather}
where $\as = (q_\s/m_\s)(\E + \v\times\B)$ is acceleration
due to the Lorentz force.  In deriving fluid balance laws we
take the electric field $\E$ and the magnetic field $\B$ as prescribed.
To ensure that the Lorentz force is invariant under a Galilean
transformation of inertial reference frame we assume
(that is, pretend) that the electromagnetic field transforms according to
\begin{gather*}
  \B' = \B, \\
  \E' = \E + d\v \times\B,
\end{gather*}
where $d\v:=\v-\v'$ is the velocity of the primed
reference frame measured in the unprimed frame;
this approximates the actual transformation
of electromagnetic field and velocities when $d\v$
is much smaller than the speed of light.

Let $\chi(\v)$ be a generic conservation velocity moment.
Multiply both sides by $\chi$ and integrate by parts.
Get the generic moment evolution equation
 \begin{gather*}
   \partial_t\intv\fs \chi + \Div\intv\fs\v\chi
 \end{gather*}
Taking $\chi=1$ gives an evolution equation for \defining{density}
$\glslink{mdens}{\Zmdens}_\s := \intv \fs$,
\begin{gather}
  \label{densityEvolution}
  \partial_t \mdens_\s + \Div(\mdens_\s\u_\s) = 0
\end{gather}
and taking $\chi=1/m_\s$ gives an evolution equation for
\defining{particle number density} $\glslink{ndens}{n}_\s := \mdens_\s/m_\s$,
\begin{gather}
  \label{spcNumDensEvolution}
  \partial_t n_\s + \Div(n_\s\u_\s) = 0.
\end{gather}
Taking $\chi=\v$ gives an evolution equation for
\defining{momentum density}\footnote{with apologies to the
magnetization $\vec M$ and to those who prefer to mind their p's
and q's.  The letter $P$ has gotten out of hand. 
This dissertation is \emph{not} brought to you by the letter
``\textbf{P}'' (nor ``e'', nor ``p'', nor ``s'', nor ``i'').}
$\glslink{M}{\ZM_\s}:=\mdens_\s\u_\s := \intv \v\fs$,
\begin{gather}
  \label{momentumEvolution}
  \partial_t (\mdens_\s\u_\s) + \Div(\mdens_s\u_\s\u_\s + \PT_\s)
  = q_\s n_\s(\E+\u_\s\times\B) + \R_\s,
\end{gather}
where $\glslink{PT}{\ZPT_\s}:=\intc \f_\s\c\c = \mdens_\s(\eT_\s - \u_\s\u_s)$ is
the \defining{pressure tensor}
and $\glslink{R}{\ZR_\s}:=\intv C_\s \v$ is the (resistive) \defining{drag force}
due to collisions with other species.
Taking $\chi=\v\v$ gives an evolution equation for the \defining{energy tensor}
$\glslink{ET}{\ZET_\s}:=\mdens_\s\glslink{eT}{\ZeT}_\s:= \intv \v\v\fs$,
\begin{gather}\label{ETevolution}
  \partial_t (\mdens_\s\eT_\s)
  + \Div(\mdens_s\u_\s\eT_\s) + \SymB(\Div(\PT_\s\u_\s)) + \Div \qT_\s \\
  = q_\s n_\s \SymB(\u_\s\E+\eT_\s\times\B)
    + \SymB(\u_\s\R_\s)+\RT_\s + \QT_\s,
\end{gather}
where $\glslink{qT}{\ZqT_\s}:=\intc \f_\s\c\c\c$ is the heat flux tensor,
\gls{Sym} denotes the symmetric part of its tensor argument,
\gls{SymB} denotes twice the symmetric part of a second-order tensor,
$\RT_\s+\QT_\s:=\intv \v\v C_\s$ is collisional source of tensor
energy, where the isotropization tensor
$\glslink{RT}{\ZRT_\s} :=\intv \v\v C_{\s\s}$ is the result
of intraspecies collisions and the heating tensor
$\glslink{QT}{\ZQT_\s}:=\sum_{\p\ne\s} \intv \v\v C_{\s\p}$ is the result
of collisions with other species.
Finally, taking half the trace of the energy tensor evolution equation
(or taking $\chi=|\v|^2/2$) gives an evolution equation for the
energy $\glslink{Nrg}{\ZNrg}_\s:=\mdens_\s\nrg_\s:=\intv |\v|^2\fs$,
\begin{gather}\label{scalarEnergyEvolution}
  \partial_t (\mdens_\s\nrg_\s) + \Div(\mdens_s\u_\s\nrg_\s + \u_\s\dotp\PT_\s + \q_\s)
  = q_\s n_\s \u_\s\dotp\E + Q_\s,
\end{gather}
where
$\glslink{nrg}{\Znrg_\s}$ is the energy per mass,
$\glslink{q}{\Zq_\s}:=\intc \f_\s \c |\c|^2/2 = \tr \qT_\s/2$
is the heat flux, and
$\glslink{Q}{Q_\s}:=\sum_{\p\ne\s}\intv C_{\s\p} |v|^2/2$
is the heat source due to collisions with other species.

\subsection{Pressure evolution}

Multiplying the momentum evolution equation \eqref{momentumEvolution}
by $2\u_\s$, taking the symmetric part,
and assuming a smooth solution gives the kinetic energy
tensor evolution equation
\begin{gather*}
  \partial_t(\mdens_\s\u_\s\u_\s) + \Div(\mdens_\s\u_\s\u_\s\u_\s)
    + \SymB(\u_\s\Div\PT_\s) \\
     = ({q_\s}/{m_\s})\SymB(\mdens_\s\u_\s\E+\mdens_\s\u_\s\u_\s\times\B) + \SymB(\u_\s\R_\s).
\end{gather*}
Subtracting this kinetic energy tensor evolution equation
from the energy tensor evolution equation \eqref{ETevolution}
gives the pressure tensor evolution equation
\begin{gather}
   \label{PTevolution}
   \partial_t \PT_\s + \Div(\u_\s\PT_\s) 
     + \SymB(\PT_\s\dotp\D\u_\s) + \Div\qT_\s
     = ({q_\s/m_\s})\SymB(\PT_\s\times\B) + \RT_\s + \QT_\s
\end{gather}
(where we have used that
$\SymB(\Div(\PT_\s\u_\s)-\u_\s\Div\PT_\s)
 = \SymB(\PT_\s\dotp\D\u_\s)$).
Alternatively, to obtain the pressure tensor evolution equation
directly, multiply both
sides of the Boltzmann equation by $\c\c$ and integrate
by parts.

Scalar pressure evolution is half the trace of pressure tensor
evolution:
\begin{gather}
   \label{pEvolution}
   (3/2)(\partial_t p_\s + \Div(\u_\s p_\s))
     + p_\s \Div \u_\s + \dPT\ddotp\D\u_\s + \Div\q_\s
     = Q_\s,
\end{gather}
where $\glslink{p }{p_\s}:=\tr\PT_\s/3 = \int_\c \f_\s|c|^2/3$
is the \defining{scalar pressure},
$\dPT_\s\ddotp\PT_\s-\id p_\s$ is the deviatoric pressure
tensor, and \gls{ddotp} denotes contraction over
two adjacent indices.

\subsection{Temperature evolution}

Conserved variables such as energy are
volume densities which depend on the choice of reference frame.
Pressure represents volume density of energy
in the reference frame of the fluid.
\emph{Temperature} represents mass density
of energy in the reference frame of the fluid.
It is the natural quantity with which to
represent energy evolution (1) in terms of
convective derivatives and (2) independent
of reference frame.

For a five-moment gas we define the
\defining{temperature} $T_\s$ to be the average
thermal energy per particle,
$\glslink{T}{T_\s} := p_\s/n_\s = m_\s\mean{|\c|^2}/3$,
and we define the \defining{``quasi'' temperature}
to be the thermal energy density per mass,
$\glslink{theta}{\theta_\s} := p_\s/\mdens_\s = \mean{|\c|^2}/3$.
Note that $T_\s = m_\s \theta_\s$.
 
For a ten-moment gas we define the
the \defining{temperature tensor} by
$\glslink{TT}{\ZTT_\s} := \PT_\s/n_\s = m_\s\mean{\c\c}$
and we define the \defining{``quasi'' temperature tensor}
by $\glslink{ThetaT}{\ZThetaT_\s} := \PT_\s/\mdens_\s = \mean{\c\c}$.

Using $n_\s d_t\TT_\s = \dcs_t(n_\s \TT_\s) = \dcs_t(\PT_\s)$,
the pressure evolution equation \eqref{PTevolution} may
also be regarded as a temperature evolution equation,
\begin{gather}\label{nTTevolution}
   n_\s d_t \TT_\s
     + \SymB(\PT_\s\dotp\D\u_\s) + \Div\qT_\s
     = ({q_\s/m_\s})\SymB(\PT_\s\times\B) + \RT_\s + \QT_\s.
\end{gather}
Dividing by number density,
\begin{gather} \label{TTevolution}
   d_t \TT_\s 
     + \SymB(\TT_\s\dotp\D\u_\s) + n_\s^{-1}\Div\qT_\s
     = ({q_\s/m_\s})\SymB(\TT_\s\times\B) + n_\s^{-1}\RT_\s + n_\s^{-1}\QT_\s.
\end{gather}
Scalar temperature evolution is half the trace of temperature
tensor evolution.  Equivalently,
scalar pressure evolution \eqref{pEvolution}
may be regarded as scalar temperature evolution
using $n_\s d_t T_\s = \dcs_t(n_\s T_\s) = \dcs_t(p_\s)$:
\begin{gather}
   \label{nTevolution}
   (3/2)n d_t T_\s
     + n T_\s \Div \u_\s + n \dTT\ddotp\D\u_\s + \Div\q_\s
     = Q_\s;
\end{gather}
divided by number density,
\begin{gather} \label{Tevolution}
   (3/2)d_t T_\s 
     + T_\s \Div \u_\s 
     + \dTT_\s\ddotp\D \u_\s
     + n_\s^{-1}\Div\q_\s
     = n_\s^{-1}Q_\s;
\end{gather}
here $\glslink{dTT}{\ZdTT_\s}:= \TT_\s - T_\s\id$
(where $\id$ is the identity tensor)
is the \defining{deviatoric part} of the temperature tensor;
we have separated it out in anticipation of the entropy-respecting
five-moment closure.

\section{Entropy evolution}

\def\fs{\f_\s}
\subsection{Kinetic entropy evolution}

The \defining[entropy density!statistical]{(statistical) entropy density}
of the particle density distribution $\f_\s(t,\xb,\v)$
is defined to be
\begin{align*}
  S_\s &:= \intv \eta_\s, &&\hbox{where} &
       \eta_\s &:= -(\beta_\s\fs\ln\fs+\alpha_\s \fs);
\end{align*}
$S_\s$ increases with randomness and
measures the expected ``surprise''
or unlikelihood of the actual distribution of particles. 
The total gas-dynamic entropy is the sum of the 
entropies of the distributions.
Collision operators are expected to cause the
total statistical entropy to increase.

The constant $\alpha_\s$ is arbitrary and may be
freely chosen so that the entropy of equilibrium
distributions has a simple formula.
The constant $\beta_\s$ must be positive.

To ensure that interaction of particles of
different species increases the total entropy,
the scaling of $\beta_\s$ should be consistent
among species.
If $\f_\s$ is taken to be a \emph{number} density
of particles in phase space then
$\beta_\s$ should be equal (e.g.\ to 1) for all particle
species;
if $\f_\s$ is taken to be a \emph{mass} density
then $\beta_\s$ should be proportional to 
$m_\s^{-1}$.

%
One way to justify the definition of statistical entropy
is to discretize state space and approximate
the evolution of $\f_\s$ by a Markov chain.

Recall the Boltzmann equation,
which says that particles are
conserved under flow in phase space.
Since $\D_\xb\dotp\v=0$ and since $\D_\v\dotp\a_\s=0$
(because $\a_\s=(q_\s/m_\s)(\E+\v\times\B)$),
flow in phase space is incompressible and can
instead be written in terms of a material derivative
in phase space:
\begin{gather*}
  \partial_t\fs+\v\dotp\D_\xb\fs+\as\dotp\D_\v \fs = C_\s;
\end{gather*}
this is not the natural form to use when deriving fluid conservation laws,
but it is relevant in obtaining entropy evolution.

Multiply the Boltzmann equation by
$\eta_\s':=\frac{d\eta_\s}{d \f_\s}$.
Using the chain rule,
\begin{gather*}
  \partial_t\eta_\s+\v\dotp\D_\xb\eta_\s+\as\dotp\D_\v \eta_\s
  = \eta'_\s C_\s.
\end{gather*}
Using phase space incompressibility to put it back into conservation form,
\begin{gather*}
  \partial_t\eta_\s+\D_\xb\dotp(\v \eta_\s)+\D_\v\dotp(\as \eta_\s) = \eta' C_\s.
\end{gather*}
Integrate over velocity space.
Get the entropy evolution equation
\begin{gather*}
  \partial_t S_\s
  +\D_\xb\dotp\underbrace{\intv \v\eta_\s}_{\hbox{Call $\Phi_\s$}}
  = \underbrace{\intv \eta' C_\s}_{\hbox{Call $\Psi_\s$}};
\end{gather*}
$\Phi_\s$ is the \defining{entropy flux}
and $\Psi_\s$ is the rate of \defining{entropy production}.
The content of Maxwell's ``H'' theorem is that $\Psi_\s$
is zero precisely when $\f_\s$ is a \first{Maxwellian}
(equilibrium) velocity distribution; else it is strictly
positive. Therefore, collisions cause distributions to trend
toward equilibrium. In general we impose the H theorem as an
assumption and a modeling requirement for collision operators.

To study properties and closure, we want to put entropy
evolution in frame-invariant form.  Define the molar entropy
by $s_\s := S_\s/n_\s$.  Then
\begin{gather*}
  n_\s d_t s_\s
  +\D_\xb\dotp\underbrace{\intc \c\eta_\s}_{\hbox{Call $\Phi'_\s$}}
  = \Psi_\s.
\end{gather*}
We refer to $\Phi'_\s$ as the \defining{diffusive entropy flux}.
It is the portion of the entropy flux that is not accounted for
by fluid motion.  We will see that it corresponds roughly
to heat flux (divided by temperature).

\subsection{Maxwellian limit}

In the absence of other effects, collisions cause entropy to
trend to a distribution that maximizes entropy subject to the
constraint that mass, momentum, and energy must be conserved.
Solving this constrained maximum problem by variational calculus
shows (see equation \eqref{MaxwellianDistribution})
that the equilibrium distribution is a Maxwellian,
\begin{gather*}
  \f_\mathcal{M}
  := \frac{\mdens}{(2\pi\theta)^{3/2}}\exp\left(\frac{-|\v-\u|^2}{2\theta}\right).
\end{gather*}
The \defining{Gaussian distribution} is defined to be
the entropy-maximizing distribution subject to the
(not physically justified) constraint that all second-order velocity
moments are conserved, and equals (see \eqref{GaussianDistribution})
\begin{gather}
  \label{GaussianDistributionCopy}
  \f_\mathcal{G}
  := \frac{\mdens}{\sqrt{\det(2\pi\Theta)}}\exp\left(-\c\dotp\Theta^{-1}\dotp\c/2\right).
\end{gather}

We say that $\f$ is an \defining{even function of velocity} if $\f(-\v) = \f(\v)$.
Maxwellian and Gaussian distributions are even functions of thermal velocity.
If $\f$ is even then so is $\eta(\f)$.  Odd moments of even functions are
zero.  In particular, \emph{the diffusive entropy flux $\Phi' := \intc \c\eta$
and the heat flux $\q:=\intc \c|\c|^2 \f$ and tensor heat flux
$\qT:=\intc \c\c\c \f$ all vanish for Maxwellian and Gaussian distributions.}

\subsection{Gas-dynamic entropy evolution}

In a highly collisional gas,
the distribution function $\f$ is close to Maxwellian.
Assuming that the distribution is exactly Maxwellian
gives the \defining{hyperbolic five-moment model} of a gas;
it is equivalent to the Boltzmann equation for a collision
operator that instantaneously relaxes to equilibrium.

As obtained in equation \ref{MaxwellianEntropy},
the entropy density of a Maxwellian distribution is
$S_\Maxwell = n s_\Maxwell$, where
\begin{gather*}
  \glslink{sMaxwell}{s_\Maxwell} := \ln\left(\frac{T^{3/2}}{n}\right)
\end{gather*}
is the molar entropy.
We define the \defining{five-moment gas-dynamic entropy}
by this formula regardless of whether
the distribution is actually Maxwellian.

Assuming that the distribution is exactly Gaussian
gives the \defining{hyperbolic ten-moment model} of a gas;
it is equivalent to the Boltzmann equation for an
(artificial and unphysical) collision operator that
instantaneously relaxes to the entropy-maximizing
distribution subject to the constraint that all ten
moments are preserved.

As obtained in equation \ref{GaussianEntropy},
the entropy density of a Gaussian distribution is
$S_\Gauss = n s_\Gauss$, where
\begin{gather*}
  s_\Gauss := \ln\left(\frac{\sqrt{\det{\TT}}}{n}\right)
\end{gather*}
is the molar entropy.
We define the \defining{ten-moment gas-dynamic entropy}
by this formula regardless of whether
the distribution is actually Gaussian.

Note that
\begin{gather*}
  s \le s_\Gauss \le s_\Maxwell,
\end{gather*}
with equality only when the assumed distributions match.

\subsubsection{Five-moment entropy evolution}

We can obtain an evolution equation for each gas-dynamic entropy
analogous to our equation for the convective derivative
of kinetic entropy.
Taking the convective derivative of five-moment entropy 
and using the continuity equation $d_t n = -n\Div\u$
gives
$
  d_t s_\Maxwell = d_t((3/2)\ln T-\ln n)
  = (3/2)\Tinv d_t T + \Div\u;
$
using the temperature evolution equation \eqref{Tevolution}
(divided by $T$)
allows us to eliminate not only $d_t T$ but also the
$\Div\u$ term:
\begin{gather} \label{s5evolutionSimple}
  n\,d_t s_\Maxwell + \Tinv\Div\q + \Tinv\dPT\ddotp\D\u = \Tinv Q;
\end{gather}
that is,
\begin{gather} \label{s5evolution}
  n\,d_t \glslink{sMaxwell}{s_\Maxwell} + \Div(\Tinv\q)
  = \q\dotp\D \Tinv - \Tinv\dPT\ddotp\Sym(\D\u) + \Tinv Q;
\end{gather}
here we have separated out the \mention{thermal entropy production}
$\q\dotp\D \Tinv$ from the divergence of the
\mention{diffusive entropy flux} $\Phi'_\Maxwell := \Tinv\q$.
The term $\Tinv Q$ represents entropy source
due to collisional exchange with other species
(via thermal exchange and resistive drag force);
since I am interesed in non-resistive two-fluid plasmas
I will generally ignore this term.

The entropy-evolution equations
reveal that for smooth solutions
entropy is conserved in the absence of heat flux.

\def\AA{{\underline{\underline{A}}}}
\subsubsection{Ten-moment entropy evolution}
To get ten-moment entropy evolution we
take the convective derivative of ten-moment
entropy $s_\Gauss = \ln(\det\TT)/2 - \ln n$.
Recall the Jacobi formula for the differential of the determinant,
$d\det \AA = \tr(\adj(\AA)\dotp d\AA)$,
that is,
$d\ln \det \AA = \tr(\AA\inv\dotp d\AA)$.
Thus, $d_t\ln\det\TT = \TTinv:d_t\TT$.
Again using the continuity equation $d_t\ln n = -\Div\u$,
\begin{gather}
  \label{dtsG}
  2 d_t s_\Gauss = \TTinv:d_t\TT + 2\Div\u.
\end{gather}
Recall the temperature evolution equation \eqref{TTevolution},
\begin{gather*}
   n_\s d_t \TT_\s 
     + n_\s \SymB(\TT_\s\dotp\D\u_\s) + \Div\qT_\s
     = n_\s ({q_\s/m_\s})\SymB(\TT_\s\times\B) + \RT_\s + \QT_\s.
\end{gather*}
Applying $\TTinv:$ to this equation,
using the identities
\begin{gather*}
  \TTinv\ddotp d_t\TT =  d_t \ln \det\TT, \\
  \TTinv\ddotp\Sym(\TT\dotp\D\u) = \Div\u=-d_t\ln n, \hbox{ and} \\
  \TTinv\ddotp\Sym(\TT\times\B) = 0,
\end{gather*}
and substituting into equation \eqref{dtsG} yields
\begin{gather}
  \label{sGaussEvolutionSimple}
  2 n\,d_t s_\Gauss + \TTinv\ddotp\Div\qT_\s
     = \TTinv\ddotp\RT_\s + \TTinv\ddotp\QT_\s;
\end{gather}
that is,
\begin{gather*} \label{sGaussEvolution}
  2 n\,d_t s_\Gauss + \Div(\TTinv\ddotp\qT_\s)
     = \qT_\s\dddotp\D\TTinv + \TTinv\ddotp\RT_\s + \TTinv\ddotp\QT_\s,
\end{gather*}
where \gls{dddotp} denotes contraction over three adjacent indices.
Here we have separated out the \mention{thermal entropy production}
$\qT_\s\dddotp\D\TTinv$ from the divergence of the
\mention{diffusive entropy flux} $\Phi'_\Gauss := \TTinv\ddotp\qT_\s$.
The term $\TTinv\ddotp\RT_\s$ represents entropy production
due to intraspecies collisional exchange
(which effects relaxation toward a Maxwellian),
and the term $\TTinv\ddotp\QT_\s$ represents entropy production
due to collisions with other species
(via thermal exchange and resistive drag force);
again, since I am interesed in non-resistive two-fluid plasmas
I will generally ignore this term.

\section{Entropy-respecting forms for closure}
\label{entropyRespectingClosure}


While gas-dynamic entropy of a physical/kinetic gas can decrease,
we may impose the assumption that gas-dynamic entropy cannot
decrease as a requirement that closure relations must satisfy.
For the five-moment gas dynamic equations, when collisions are
sufficiently predominant to keep velocity distributions near
Maxwellian, closures can be well-justified, and in particular the
requirement that gas-dynamic entropy must increase
can be justified as follows. Since Maxwellian distributions
maximize entropy subject to \emph{physical} constraints, if
the deviation from the Maxwellian distribution is of order
$\epsilon$ the deviation of the entropy will merely be of order
$\epsilon^2$. Therefore, gas-dynamic entropy is an accurate
approximation to kinetic entropy, and for a collisional gas we
expect it to increase.

In contrast, closures of the ten-moment equations which do not
cause it to approximate the five-moment model are difficult to
justify unless one assumes an idiosynchratic collision operator
which kills deviations from a Gaussian much faster than it
relaxes a Gaussian distribution to a Maxwellian. For such a model
heat flux is negligible in comparison to viscosity, i.e., the
Prandtl number is much larger than 1. For a monatomic gas the
Prandtl number is approximately 2/3; for other gases and fluids
the Prandtl number is smaller, often much smaller.

Given a deviation of order $\epsilon$ from a Gaussian
distribution, to conclude that the deviation of the ten-moment
entropy is of order $\epsilon^2$ we would need to assume that
the deviation from the artificial constraint (that the tensor
pressure is invariant) is of order $\epsilon^2$. Regardless of
how artificially high a finite Prandtl number is, for sufficiently
small $\epsilon$ this artificial constraint does not hold.
The requirement that ten-moment gas-dynamic entropy should be
nondecreasing for ten-moment closures thus lacks adequate
physical justification except in the adiabatic case
where the Prandtl number is infinite.

Having stated these prefatory caveats, we now proceed to
derive the form of entropy-respecting 5-moment and
10-moment closures.

\subsection{Five-moment closure}

Recall the five-moment entropy evolution equation \eqref{s5evolution},
which we now rewrite as:
\begin{gather} \label{sMevolution}
  n\,d_t s_\Maxwell + \Div(\Tinv\q)
  = \q\dotp\D \Tinv - n\dPshape\ddotp\dstrain + \Tinv Q;
\end{gather}
here $\glslink{dPshape}{\ZdPshape}:=\Tinv\dTT$
is the deviatoric part of the ``shape'' of the temperature tensor and
$\glslink{dstrain}{\Zdstrain} := \deviator{\Sym(\D\u)} := \Sym(\D\u) - \id\Div\u/3$
is called the \defining{deviatoric strain rate}
and is the traceless part of the \defining{strain rate}
tensor $\glslink{strain}{\Zstrain} := \Sym(\D\u)$.
We now require the intraspecies entropy source terms
$\q\dotp\D \Tinv$ and $-\dPshape\ddotp\dstrain$
to be nonnegative as a closure requirement
and deduce the form of the closure.

\subsubsection{Isotropic linearized heat flux closure}

To ensure that $\q\dotp\D \Tinv$
is nonnegative we make $\q$ a function of $\D \Tinv$.
We will approximate the function as linear
(with state-dependent coefficients).
Since
$\q$ should be zero in equilibrium 
(for which $\D \Tinv=0$), and
assuming that the derivative of the map $\D \Tinv\mapsto \q$
is nonzero, in the near-Maxwellian limit
such a linearized closure is rigorously justifiable.

Physical laws should be invariant under rotation for closed systems.
In the absence of a magnetic field (or other external
influence) sufficiently strong to break this fundamental symmetry
we expect intraspecies closure relations to be isotropic.
The general form of a linear isotropic closure is 
\def\kn{\overline{\kappa}}
\begin{gather}\label{scalarHeatFluxClosure}
 \q_\s = \kn_\s \D \Tinv_\s  = -\heatConductivity_\s \D T_\s,
\end{gather}
where \gls{heatConductivity} is called the \defining{heat conductivity}
and we infer that
\begin{gather}
  \kn_\s = T_\s^2 \heatConductivity_\s.
  \label{heatConductivityVersusKappa}
\end{gather}
Nonnegativity of $\q_\s\dotp\D \TTinv_\s =\kn\|\D \TTinv_\s\|^2$
is ensured as long as $\kn$ is nonnegative.

\def\b{\bhat}
\def\idperp{\id_\perp}
\def\idskew{\id_\wedge}
\def\idpara{\id_\parallel}
\def\Kperp{\kappa_\perp}
\def\Kskew{\kappa_\wedge}
\def\Kpara{\kappa_\parallel}
\subsubsection{Gyrotropic heat flux closure.}
In the presence of a sufficiently strong magnetic field $\B=\|\B\|\bhat$
we cannot assume that intraspecies collisions are 
governed by isotropic physics, and 
instead we can merely assume a gyrotopic linear closure
$\q=-\TheatConductivity\dotp\D T$, where $\TheatConductivity$ is a \defining{gyrotropic tensor};
that is, $\TheatConductivity$ is invariant under rotations around an axis
aligned with the magnetic field.
In particular,
\begin{gather}
  \label{entropyRespectingHeatFluxForm}
  \q = -\left(\Kperp\idperp
            +\Kskew\idskew
            +\Kpara\idpara\right)\dotp\D T_\s ,
\end{gather}
where we have used that since
$\TheatConductivity$ is gyrotropic second-order tensor it is
therefore a linear combination of the perpendicular, skew, and
parallel gyrotropic tensors
\begin{align*}
             \idperp &:=\id-\b\b,
           & \idskew &:=\id\times\b,
           & \idpara &:=\b\b.
\end{align*}
For this closure, to ensure that
$\q\dotp\D\TTinv = T^{-2}\D T\dotp\TheatConductivity\dotp\D T\ge 0$,
$\TheatConductivity$ must be positive definite;
that is, the parallel and perpendicular
heat conductivities must be nonnegative,
\begin{align}
  \label{entropyRespectingHeatFluxClosureRequirements}
  \Kperp\ge 0, && \Kpara\ge 0.
\end{align}

\subsubsection{Isotropic viscous stress closure}
\label{isoViscStressClosure}

To ensure that the viscous entropy production
term $-\dPshape\ddotp\dstrain$ in equation
\eqref{sMevolution}
is nonnegative, we make $\dPshape$ a function of $\dstrain$.
The simplest such closure is to make $\dPshape$ a linear
isotropic function of $\dstrain$.
Then the facts that $\dstrain$ is symmetric and that
$\dPshape$ must be symmetric and traceless imply that
$\dPshape$ is proportional to $\dstrain$,
\begin{gather}\label{dPshapeClosure}
  -\dPshape = 2\ptime\dstrain,
\end{gather}
where we will see that $\glslink{ptime}{\ptime}$ is a relaxation period;
equivalently, the deviatoric stress $-\glslink{dPT}{\ZdPT}:=-n\dTT$
is proportional to the deviatoric strain,
\begin{gather}\label{dPTclosure}
  -\dPT = 2\viscosity\dstrain,
\end{gather}
where the \defining{viscosity} $\glslink{viscosity}{\Zviscosity}:=p\ptime$
must be positive to respect entropy.
In the absence of a symmetry-breaking magnetic field,
in the near-Maxwellian-limit, assuming that the function
$\dstrain\mapsto\dPT$ has nonzero derivative at $\dstrain=0$,
such a closure can be rigorously justified.

\subsubsection{Gyrotropic viscous stress closure}
\label{GyrotropicViscousStressClosure}

In the presence of a sufficiently strong magnetic field,
we might merely assume that $\dPshape$ is a \emph{gyrotropic}
function of $\dstrain$:
\def\AA{\underline{\underline{A}}}
\begin{align}\label{gyrotropicTempDeviator}
  -\dPshape &= 2\ptime\tviscosity\ddotp\dstrain, &&\hbox{that is,}& 
  -\dPT &= 2\Tviscosity\ddotp\dstrain,
\end{align}
where $\glslink{Tviscosity}{\ZTviscosity}:=\tviscosity\viscosity$
is the viscosity tensor and \gls{tviscosity} is its nondimensional shape;
$\tviscosity$ is a gyrotropic tensor which evidently
must be symmetric and traceless in its first two coefficients
and which is without loss of generality symmetric in its
last two coefficients.
This leads to five distinct coefficients of viscosity.
To ensure that entropy production $\Tinv\dstrain\ddotp2\Tviscosity\ddotp\dstrain$
is strictly positive for a non-Maxwellian distribution, we impose the
positive-definiteness (and invertibility) criterion
$\AA\ddotp\tviscosity\ddotp\AA > 0$ for any $\AA\ne 0$.
In the isotropic case $\tviscosity$ is 
the identity tensor $\id\Ldiamond\id$ for linear transformations
on the space of second-order tensors, whose components are
$\delta_{ijkl}=\delta_{ik}\delta_{jl}$.

\subsection{Ten-moment closure}

Recall the ten-moment entropy evolution equation \eqref{sGaussEvolution},
which we now write as
\begin{gather}
  \label{sGaussEvolution}
  2 n\,d_t s_\Gauss + \Div(\TTinv\ddotp\qT)
     = \qT\dddotp\Sym(\D\TTinv) + \dTTinv\ddotp\RT + \TTinv\ddotp\QT;
\end{gather}
here we have used that $\tr(\RT)=0$ (because
$\RT:=\intv \v\v C$ and $\intc |\v|^2 C = 0$ by conservation of energy)
to replace
$\TTinv$ with its deviatoric part $\dTTinv:=\TTinv-\id\,\tr(\TTinv)/3$.

We impose the closure requirement that the intraspecies entropy
source terms must be nonnegative and zero at equilibrium.

\subsubsection{Intraspecies collision closure}

By selecting an artificial collision operator that instantaneously
relaxes to a Gaussian distribution we may justify the requirement
that the entropy production of the intraspecies collision term
$\dTTinv\ddotp\RT$ should be positive.
Clearly we can ensure positivity by making $\RT$ a function of
$\dTTinv$.  However, to keep the closure simple
we instead prefer to make $\RT$ a function of $\dTT$;
this is equivalent for small deviations from isotropy
(i.e.\ equilibrium), which is when we can justify
linearized closures anyway.

\textbf{Isotropic intraspecies collision closure.}
Imposing that $\RT$ is a linear isotropic function
of $\dTT$, equivalently of $\dPT$, yields the closure
\begin{gather}
  \label{isoRTclosure}
  \RT = -\ptime\inv \dPT,
\end{gather}
where $\glslink{ptime}{\ptime}$ is the \first{relaxation period}
and $\ptime\inv$ is called the \first{relaxation rate}.
Then
\def\V{\mathbf{V}}
\def\W{\mathbf{W}}
\begin{align*}
  \dTTinv\ddotp\RT/n &= \TTinv\ddotp\RT/n 
          \\ &= \TTinv\ddotp\left(T\idtens-\TT\right)/\ptime
          \\ &= \left(\tr(\TTinv)\tr(\TT)/3-3\right)/\ptime
          \\ &= 3\tan^2(\theta)/\ptime,
\end{align*}
where $\theta$ is the angle between the vector
$ 
  \V :=
    ( 
      \sqrt{T_1},
      \sqrt{T_2},
      \sqrt{T_3} 
    )^T 
  $
  and the vector 
  $
  \W :=
    ( 
      \sqrt{T_1}^{-1}, 
      \sqrt{T_2}^{-1}, 
      \sqrt{T_3}^{-1} 
    )^T, 
$ 
where $T_1$, $T_2$, and $T_3$ are the eigenvalues
of the temperature tensor.
Indeed,
  $\V\dotp\V = \tr\TT,$
  $\W\dotp\W = \tr(\TTinv),$, and
  $\V\dotp\W = 3;$
since
$\cos^2\theta = \frac{(\V\dotp\W)^2}{\|\V\|^2\|\W\|^2}$,
$\tan^2\theta = \sec^2\theta-1
  = \frac{(\tr\TT)\tr(\TTinv)}{9}-1$,
which is strictly greater than $0$ unless
$\V\times\W=\mathbf{0}$ (which says that the temperature
tensor is isotropic).

\def\taucoef{\underline{\underline{\underline{\underline{C}}}}_\ptime}
\textbf{Gyrotropic intraspecies collision closure.}
More generally, we assume that
$\RT$ is a linear \emph{gyrotropic} function
of $\dTT$, equivalently, of $\dPT$:
\begin{gather}\label{RTclosure}
  \RT = -\ptime\inv \Gcoef\ddotp\dPT,
\end{gather}
where $\Gcoef$ is a nondimensional gyrotropic tensor of coefficients.
To respect entropy, based on equation \eqref{sGaussEvolution},
we require that
\begin{align}
  \label{entropyRespectingGyrotropicRelaxationClosureRequirement}
  \dTTinv\ddotp\RT &= -\dTTinv\ddotp\Gcoef\ddotp\dTT > 0
\end{align}
whenever $\dTT\ne 0$, that is, whenever the pressure
is not isotropic. 
Recall that $\glslink{dPshape}{\ZdPshape}:=\Tinv\dTT$.
In the limit $\|\dPshape\|\to 0$,
$\dPshapeinv\approx-\dPshape$.  Therefore, for small
$\|\dPshape\|$, equation
\eqref{entropyRespectingGyrotropicRelaxationClosureRequirement}
implies that for all $\dPshape$
\begin{align}
  \dPshape\ddotp\Gcoef\ddotp\dPshape > 0,
\end{align}
which implies that $\Gcoef$ has an inverse
when regarded as a linear transformation relating traceless tensors.
We impose invertibility as a minor auxiliary closure requirement,
which we justify by showing in section \ref{equivalenceOfClosures}
the asymptotic equivalence of the five-moment and ten-moment closures.


\subsubsection{Entropy-respecting closure for tensor heat flux}

The requirement that the entropy production of
the instraspecies heat flux tensor be positive is difficult
to justify.  One could try to justify it with a sequence of
collision operators which relax to a Gaussian increasingly
rapidly, but in such a limit there is no heat flux anyway.
Nevertheless, we work out the consequences of this requirement.

The local production of ten-moment gas-dynamic entropy
due to heat flux is seen in equation
\eqref{sGaussEvolution} to be
\begin{gather}
  \label{sGentropyProduction}
  \qT\dddotp\Sym(\D\TTinv).
\end{gather}
To ensure that this quantity is positive,
we posit that $\qT$ is a linear gyrotropic
function of its complement in this inner product,
\def\KKK{K_{[6]}}
\begin{gather}
     \qT=\KKK\dddotp\Sym(\D\TTinv),
\end{gather}
where $\KKK$ is a gyrotropic tensor which
must be symmetric in its first three indices
and which without loss of generality we require to
be symmetric in its last three indices;
the entropy production \eqref{sGentropyProduction} is then
\begin{gather}
  \label{sGentropyProduction}
    \Sym(\D\TTinv)\dddotp\KKK\dddotp\Sym(\D\TTinv),
\end{gather}
which is guaranteed to be positive if $\KKK$
satisfies the positive-definiteness criterion
\def\AAA{A}
\def\BBB{B}
\begin{gather}
  \label{sGentropyProductionVersion2}
    \AAA\dddotp\KKK\dddotp \AAA > 0
\end{gather}
for all $\AAA\ne 0$.

\def\qTL{\qT_\mathrm{L}}
\def\kLa{\overline{\kappa}_1}
\def\kLb{\overline{\kappa}_2}
\textbf{Isotropic case.}
In the absence of a magnetic field $\KKK$
should be isotropic.
Recall that in three-dimensional space
(or any odd-dimensional space) every even-order
isotropic tensor is a linear combination of tensor products
of the identity tensor \cite{Jeffreys72}.
So if $\AAA$ and $\BBB$ are symmetric
third-order tensors and $\AAA$
is a linear isotropic function of $\BBB$,
then $\AAA = \kLa \BBB + \kLb \id\vee\,\tr\,\BBB$
for some $\mu$ and $\lambda$,
where recall the definition of the vee product $\vee$
from \eqref{veeDefinition}.
Therefore, we can write the general form of
a linear isotropic entropy-respecting heat flux
closure,
\begin{gather}
  \label{qTlevermore}
  \qTL = \kLa \D\vee\TT\inv + \kLb \id\vee\tr(\D\vee\TT\inv),
\end{gather}
where $\qTL$ denotes the entropy-respecting
heat flux closure proposed by Levermore for
ten-moment gas dynamics (C.D. Levermore, presented
in talk slides e.g.\ at Kinetic FRG Young Researchers Workshop:
\emph{Kinetic Description of Multiscale Phenomena: Modeling, Theory,
and Computation}, University of Maryland, College Park, 5 March 2009).
To determine for what coefficients the entropy
production given by \eqref{sGentropyProduction}
is positive, we compute:
\begin{align*}
    \SymC(\D\TTinv)\dddotp\qTL
    &= (\D\veebar\TTinv)\dddotp \qTL
\\  &= (\D\veebar\TTinv)\dddotp
    \left(\kLa 3\D\veebar\TT\inv + \kLb 9\id\veebar\tr(\D\veebar\TT\inv)\right)
\\  &= 3 \kLa \lVert \D\veebar\TT\inv \rVert^2
        + 9 \kLb \lVert\tr(\D\veebar\TT\inv)\rVert^2;
\end{align*}
here we have used that
\begin{align*}
       &(\D\veebar\TTinv)\dddotp(\id\veebar\tr(\D\veebar\TT\inv))
  \\ = &(\D\veebar\TTinv)\dddotp(\id\otimes\tr(\D\veebar\TT\inv))
  \\ = &(\D\veebar\TTinv)_{ijk} \delta_{ij}(\D\veebar\TT\inv)_{mmk}
  \\ = &\tr(\D\veebar\TTinv)\dotp\,\tr(\D\veebar\TT\inv)
  \\ = &\lVert\tr(\D\veebar\TTinv) \rVert^2.
\end{align*}
So to ensure that entropy is respected we would require that
$\kLa\ge 0$ and $\kLb\ge 0$.

To compare with other closures, we express $\D\TT\inv$
in terms of $\D\TT$.
Taking the differential of the identity
$\TT\inv\dotp\TT = \id$ and solving yields
\begin{align*}
  d\TT\inv &= -\TT\inv\dotp d\TT\dotp\TT\inv
            = -T^{-2}\Pshape\inv\dotp d\TT\dotp\Pshape\inv
            = -T^{-2}d\TT\ddotp(\Pshape\inv\Ldiamond\Pshape\inv).
\end{align*}
In particular,
\begin{align*}
  \D\TT\inv &= -T^{-2}\D\TT\ddotp(\Pshape\inv\Ldiamond\Pshape\inv)
\end{align*}
  and
\begin{align}
  \D\vee\TT\inv &= -T^{-2}\SymC(\D\TT\ddotp(\Pshape\inv\Ldiamond\Pshape\inv)).
  \label{LevermoreComparisonForm}
\end{align}

In the isotropic case, where $\Pshape=\id$,
   $\TT = \id T$ and
\begin{gather}
 \begin{aligned}
   \D\vee\TT\inv &= \id\vee\D T\inv,&
   \tr(\D\vee\TT\inv) &= 5 \D T\inv,&
 \end{aligned}
 \label{isotropicSimplifications}
\end{gather}
and the expressions above simplify.
We then have
\begin{align}
 \label{isotropicLevermore}
 \begin{aligned}
  \qTL&= \kLa \D\vee\TT\inv + \kLb \id\vee\tr(\D\vee\TT\inv)
   \\& = \kLa \id\vee\D T\inv + 5 \kLb \id\vee\D T\inv
   \\& = (\kLa + 5 \kLb) \id\vee\D T\inv.
 \end{aligned}
\end{align}
So the entropy production is
\begin{align*}
    (\D\veebar\TTinv)\dddotp \qTL
    &=  \id\veebar\D T\inv \dddotp (\kLa + 5 \kLb) \id\vee\D T\inv
\\  &=  5(\kLa + 5 \kLb) \lVert \D T\inv\rVert^2,
\end{align*}
where we have used that
$
     \id\veebar\D T\inv\dddotp \id\vee\D T\inv
    =\id \D T\inv\dddotp \id\vee\D T\inv
    =(\delta_{ij}\partial_k T\inv)(
        \delta_{ij}\partial_k T\inv
      + 2\delta_{ik}\partial_j T\inv)
    =5 \partial_k T\inv
       \partial_k T\inv
    = 5 \lVert\D T\inv\rVert^2.
$

\subsection{Equivalence of ten-moment and five-moment stress closure
for near-isotropy}
\label{equivalenceOfClosures}


In the near-Maxwellian limit, that is, in case pressure
anisotropy is small (as is generally the case if isotropization
is rapid) we can show that the ten-moment viscous stress closure
for the relaxation tensor $\RT$ given by equation
\eqref{RTclosure}, equivalently,
\begin{gather*}
-\ptime\RT/p = \Gcoef\ddotp\dPshape,
\end{gather*}
is asymptotically equivalent to the five-moment viscous stress closure
for the deviatoric pressure $\dPT$ given by equation
\eqref{gyrotropicTempDeviator}, equivalently,
\begin{gather*}
  -2\ptime\dstrain = \Mcoef\ddotp\dPshape.
\end{gather*}
(So evidently we need $\RT/p=2\dstrain$.)

The strategy is to
match up the ten-moment pressure tensor evolution equation 
with the five-moment closure for the viscous stress.
Since the five-moment closure is expressed in terms
of the deviatoric pressure, it is natural to proceed by
writing an evolution equation for the deviatoric part of
the pressure.  We will instead write an evolution equation
for the \emph{shape} of the pressure tensor
where (recall that) the nondimensional shape of the pressure
(or temperature) tensor $\glslink{Pshape}{\ZPshape}$
and its deviatoric part $\glslink{dPshape}{\ZdPshape}$
are defined by
\begin{align*}
  \Pshape &:= \frac{\TT}{T} = \frac{\PT}{p}, &
  \Pshape &=: \idtens+\dPshape, &&\hbox{ so} &
  \dPshape &:= \frac{\dTT}{T} = \frac{\dPT}{p},
\end{align*}
where $\glslink{dPT}{\ZdPT}$ is of course the deviatoric part of the pressure tensor.
This formulation is convenient because we expect the ten-moment
closure to agree with the five-moment closure when the deviatoric
pressure is small, that is, when $\dPshape$ is much smaller than
the identity tensor $\id$.

Recall temperature tensor evolution \eqref{TTevolution}:
\begin{gather*}
   d_t \TT
   + \SymB(\TT\dotp\D\u)
     + n\inv\Div\qT 
     = (q/m)\SymB(\TT\times\B) + \RT/n,
\end{gather*}
where we assume a single species and
neglect heating due to interaction with other species.
Divide temperature tensor evolution by $T$ and get
\begin{gather*}
   T\inv d_t \TT
   + \SymB(\Pshape\dotp\D\u)
     + p\inv\Div\qT 
     = \SymB(\Pshape\times q\B/m) + \RT/p.
\end{gather*}
Using that
$\Pshape=\id+\dPshape$,
and that
$T\inv d_t\TT = T\inv d_t(T\Pshape) = d_t\dPshape+\Pshape d_t\ln T$,
we get an evolution equation for the deviatoric part of
the shape of the pressure tensor,
\begin{gather}
 d_t\dPshape+(\id+\dPshape) d_t\ln T + \SymB(\D\u)
  + \SymB(\dPshape\dotp\D\u)
  + p\inv\Div\qT \nonumber \\
   = \tfrac{q}{m}\SymB(\dPshape\times\B) + \RT/p.
   \label{dPshapeEvolutionA}
\end{gather}
Half the trace of this equation gives
\begin{gather}\label{halfTrace}
  (3/2)d_t\ln T + \Div\u +\dPshape\ddotp\D\u + p\inv\Div\q = 0,
\end{gather}
which incidently we can rewrite as evolution of the five-moment entropy,
$s_\Maxwell := \ln(T^{3/2}/n)$,
\begin{gather}
  \label{fiveMomentEntropyEvolution}
  d_t s_\Maxwell = (3/2)d_t\ln T + \Div\u = -\dPshape\ddotp\D\u - p\inv\Div\q
\end{gather}
(compare \eqref{s5evolutionSimple}).
We could simplify the analysis at this point by identifying
the assumptions which imply that entropy is conserved on
the relaxation time scale $\ptime$, but instead we delay
making approximating assumptions in order to reap the benefit
of writing an exact evolution equation for $\dPshape$
which allows us to more sharply match up the closures
and identify the domain of agreement.

Solving equation \eqref{halfTrace} for $d_t\ln T$
and substituting back into equation \eqref{dPshapeEvolutionA},
which we first rewrite as
\begin{align*}
 \RT/p =& \id d_t\ln T + \SymB(\D\u) + \SymB(\dPshape\dotp\D\u) + \p\inv\Div\qT
     \\ &+ \dPshape d_t\ln T + d_t\dPshape - \SymB(\dPshape\times q\B/m)
\end{align*}
gives a manifestly traceless equation (which as it happens can be
interpreted as an evolution equation for $\dPshape$):
\begin{align}
 \RT/p = &2\dstrain + \deviator{\SymB(\dPshape\dotp\D\u)} + \p\inv\Div\dqT
  \nonumber\\ &-(2/3)\dPshape\left(\Div\u+\dPshape\ddotp\D\u+\p\inv\Div\q\right)
  \nonumber\\ &+ d_t\dPshape - \SymB(\dPshape\times q\B/m),
  \label{dPshapeEvolution}
\end{align}
where $\dstrain:=\deviator{\Sym(\D\u)}$ is the strain deviator,
$\glslink{dqT}{\ZdqT}:=\qT-(2/5)\SymC(\id\q)$
is the \defining{deviatoric heat flux tensor}, and
$\deviator{\SymB(\dPshape\dotp\D\u)} = \SymB(\dPshape\dotp\D\u)-(2/3)\id\dPshape\ddotp\D\u$
denotes the deviatoric part of $\Sym(\dPshape\dotp\D\u)$.
Note that \gls{SymC} denotes thrice the symmetric part of a third-order tensor.

Recall the gyrotropic closure \eqref{RTclosure}, i.e.\
\begin{gather*}
  -\ptime\RT/p = \Gcoef\ddotp\dPshape.
\end{gather*}
Multiplying \eqref{dPshapeEvolution}
by $\ptime$ gives an equation in terms of the nondimensional
quantities $\Gcoef\ddotp\dPshape$ and
\begin{align*}
\ptime\D\u, && \ptime p\inv\Div\qT, &&
d_{(t/\ptime)}\dPshape, && \hbox{and} && \ptime q\B/m
\end{align*}
which reads
\begin{align}
 -\Gcoef\ddotp\dPshape = &2\ptime\dstrain + \deviator{\SymB(\dPshape\dotp\ptime\D\u)} + \ptime\p\inv\Div\dqT
  \nonumber\\ &-(2/3)\dPshape\left(\ptime\Div\u+\dPshape\ddotp\ptime\D\u+\ptime\p\inv\Div\q\right)
  \nonumber\\ &+ d_{(t/\ptime)}\dPshape - \SymB(\dPshape\times \ptime q\B/m).
  \label{nondimensional_dPshapeEvolution}
\end{align}
We want to show that this agrees with the five-moment closure
\eqref{gyrotropicTempDeviator},
\begin{gather}
  -\Mcoef\ddotp\dPshape = 2\ptime\dstrain;
  \label{fiveMomentClosure}
\end{gather}
since we do not want closure coefficients that depend on
differentiated quantities, we will 
discard all terms after the first plus sign
except for the magnetic field term.

We will need to assume that the deviatoric pressure is small, that is,
$\|\dPshape\|\ll 1$.

\def\tscale{\ptime_0}
Let $\tscale$ designate the time scale defined by $\D\u$.
That is, $\|\D\u\|=\O(\tscale\inv)$.
The five-moment closure \eqref{fiveMomentClosure}
then implies that $\dPshape = \O(\ptime/\tscale)$.
To discard terms such as $\dPshape\ptime\Div\u$
we need that $\ptime/\tscale\ll 1$.
To discard the heat flux terms we need that
the divergence of the deviatoric heat flux is very small,
\begin{align*}
  \|\ptime p\inv\Div\dqT\|\ll \|\dPshape\|,
  && \hbox{i.e.}
  && \|p\inv\Div\dqT\|\ll 1/\tscale,
  && \hbox{i.e.}
  && \|p\inv\Div\dqT\|\ll \|\D\u\|
\end{align*}
and that the divergence of the nondeviatoric part is small,
\begin{align*}
  \|\dPshape \ptime p\inv\Div\q\|\ll \|\dPshape\|,
  && \hbox{i.e.}
  && \|p\inv\Div\q\|\ll 1/\ptime.
\end{align*}

Of the unwanted terms in equation \eqref{nondimensional_dPshapeEvolution}
it remains to discard the term $ d_{(t/\ptime)}\dPshape.$
We need that
\begin{align*}
  \|d_{(t/\ptime)}\dPshape\| \ll \|\dPshape\| = \O(\ptime/\tscale).
\end{align*}
This says that the shape of the stress (or strain) deviator
changes little on the time scale of a relaxation period.

We now consider the magnetic field term
$\SymB(\dPshape\times \ptime q\B/m)$
in equation \eqref{nondimensional_dPshapeEvolution}.
Define the \defining{gyrofrequency}
$\glslink{gyrofreq}{\Zgyrofreq}:=q\|B\|/m$
and the nondimensionalized gyrofrequency
$\glslink{pgf}{\Zpgf}:=\ptime\gyrofreq$,
which is the rate of gyration divided by the
rate of relaxation.
We can neglect the magnetic field term and use the
same closure coefficients for both closures if
the magnetic field is sufficiently small so that
\begin{align*}
  \SymB(\dPshape\times \ptime q\B/m) \ll \dPshape,
  && \hbox{i.e.,}
  && \pgf \ll 1,
\end{align*}
which says that the gyrofrequency is much smaller than
the relaxation period.  If the effect of the magnetic
field is negligible then $\Mcoef$ and $\Gcoef$
should be the identity and
we might as well have restricted our study to isotropic closures.

What if the magnetic field is large?
Then we need the approximate closure identity
\begin{gather}
 \boxed{
 -\Gcoef\ddotp\dPshape = 2\ptime\dstrain - \SymB(\dPshape\times \pgf \b)
 }
 \label{approximateClosureIdentity}
\end{gather}
(from \eqref{nondimensional_dPshapeEvolution},
where $\glslink{b}{\Zbhat}:=\B/\|\B\|$ is the magnetic field direction vector)
to match up with the five-moment closure of equation \eqref{fiveMomentClosure},
\begin{gather*}
  -\Mcoef\ddotp\dPshape = 2\ptime\dstrain.
\end{gather*}
Evidently we need that
\begin{align}
  \Mcoef\ddotp\dPshape = \Gcoef\ddotp\dPshape - \SymB(\dPshape\times \pgf \b).
  \label{coefMatching}
\end{align}
We rewrite the last term as
\def\A{\underline{\underline{A}}}
\def\M{\underline{\underline{M}}}
\def\N{\underline{\underline{N}}}
\begin{align}
  \SymB(\dPshape\times \pgf \b)
  = \SymB(\id\dotp\dPshape\dotp \id\times\pgf\b)
  = \SymB(\id\dotp\dPshape\dotp\idskew\pgf),
\end{align}
where we define
$\idskew:=\id\times\b = \b\times\id$,
which projects onto the plane orthogonal to $\b$
and then rotates 90 degrees in this plane.
Define the \defining{diamond product} of two tensors with the convention that
$(\M\glslink{diamond}{\Ldiamond}\N)\ddotp\A = \M\dotp\A\dotp\N$.
Using that $\idskew$ is antisymmetric,
\begin{align}
  \SymB(\id\dotp\dPshape\dotp \idskew)
  = \SymB(\id\Ldiamond\idskew\ddotp\dPshape)
  = (\id\Ldiamond\idskew-\idskew\Ldiamond\id)\ddotp\dPshape.
\end{align}
So \eqref{coefMatching} becomes
\begin{align*}
  \Mcoef\ddotp\dPshape = (\Gcoef + \pgf(\idskew\Ldiamond\id - \id\Ldiamond\idskew))\ddotp\dPshape.
\end{align*}
We infer that
$\Mcoef = \Gcoef + \pgf(\idskew\Ldiamond\id - \id\Ldiamond\idskew),$
that is,
\begin{gather}
  \boxed{
  \Mcoef = \Gcoef + (\id\times\pgf\b)\Ldiamond\id - \id\Ldiamond(\pgf\b\times\id)
  },
 \label{coefRelation}
\end{gather}
where we have used that $\b\times\id=\id\times\b$.

As we will see in section \ref{PerturbativeStressClosure},
for $\Gcoef$ the closure coefficients are quite simple
(the identity four-tensor) whereas the closure coefficients
for $\tviscosity$ must take into account the rotational
effects of the magnetic field because it defines a closure
in terms of a differentiated quantity.
Equation \eqref{approximateClosureIdentity}
agrees with equation (5.122) in the book \cite{book:Woods04} by
Woods if one simply takes $\Gcoef=\idfour$.
He arrives at this form by assuming that the response of the
viscous stress to the deviatoric strain is delayed by
a collision period.  During this time the affect of the
magnetic field is to rotate the deviatoric strain.


I remark that to compare the five- and ten-moment
closures one must use equation \eqref{coefRelation}
to translate closure coefficients
(or just \eqref{approximateClosureIdentity}
if all one wants to do is test the agreement of the models).


\section{Perturbative closure}
\label{PerturbativeClosure}

In section \ref{entropyRespectingClosure}
we deduced forms of closure by requiring entropy to increase.
This gave the functional form of the closure, but it did not
give any estimate of the closure coefficients.
A method that determines closure coefficients 
is to use a Chapman-Enskog perturbative expansion.

\subsection{BGK and Gaussian-BGK collision operators}

To perform a Chapman-Enskog expansion we have to assume 
the form of the collision operator.  The simplest
choice of collision operator which satisfies the physical
constraints that it (1) conserves the conserved
moments (mass, momentum, and energy) and (2) respects entropy
is a collision operator which relaxes the distribution toward
an entropy-maximizing distribution.  Appendix \ref{distributions}
calculates that distribution which maximizes entropy subject
only to these physical constraints is a \defining{Maxwellian
distribution}, of the form \eqref{MaxwellianDistribution}
\begin{gather*}
  \fMaxwell = \frac{\mdens}{(2\pi\theta)^{3/2}}\exp\left(\frac{-|\v-\u|^2}{2\theta}\right),
\end{gather*}
where $\fMaxwell$ is mass density in velocity space
and $\glslink{theta}{\theta}:=\mean{|\c|^2/2} = T/m$
defines the \mention{pseudo-temperature}.
A distribution which maximizes entropy subject to the additional
constraint that all quadratic velocity moments are conserved is
a \defining{Gaussian distribution}, of the form \eqref{GaussianDistribution}
\begin{gather*}
  \f_\Gauss = \frac{\mdens}{\sqrt{\det(2\pi\Theta)}}
    \exp\left(-(\v-\u)\dotp\Theta^{-1}\dotp(\v-\u)/2\right),
\end{gather*}
where
$\glslink{ThetaT}{\ZThetaT} := \mean{\c\c}$
is the \defining{pseudo-temperature tensor}.
The distribution $\f_\Gauss$ is a normal distribution along
any axis through the origin and is the product of independent
normal distributions defined along three orthogonal principal axes.
It agrees with the Maxwellian
distribution in the isotropic case $\ThetaT=\idtens\theta$.

Recall the Galilean-invariant Boltzmann equation
(equation \eqref{GalBoltzmannEqn} assuming a single default species),
\begin{gather}
  \label{SoleBoltzmannEqn}
  \partial_t \f+\D_\xb\dotp(\v \f)+\D_\v\dotp(\a \f) = C.
\end{gather}
The BGK collision operator is
\def\Cmaxwell{C_\theta}
\begin{gather*}
  \Cmaxwell = \frac{\fMaxwell - \f}{\ptime},
\end{gather*}
where $\ptime$ is the \defining{relaxation period}.
For the BGK collision operator, the \mention{Prandtl number}
(i.e.\ the ratio of the rate of thermal diffusion to momentum diffusion,
see Appendix section \ref{PrandtlNumber})
is $\glslink{Pr}{\Pr}=1$ (for a monatomic gas), as we will see.

\def\tTheta{{\widetilde\Theta}}
\def\Cgauss{C_\tTheta}
In order to admit a tunable Prandtl number, the \mention{Gaussian-BGK}
collision operator $\Cgauss$,
also called the \mention{Ellipsoidal-Statistical BGK} collision operator,
was introduced by Holway \cite{holway66} and
is often used instead.  This model relaxes to a Gaussian distribution
for a temperature tensor which is an affine combination of the
isotropized temperature tensor (of the Maxwellian distribution)
and the temperature tensor (of the Guassian distribution with
the same quadratic moments).
The Gaussian-BGK collision operator is
\def\tautTheta{\tau_\tTheta}
\def\ftTheta{\f_\tTheta}
\begin{gather}
  \label{GaussBGKcol}
  \Cgauss = \frac{\ftTheta - \f}{\htime},
\end{gather}
where \gls{htime} is the Gaussian-BGK relaxation period,
$\ftTheta$ is the Gaussian distribution for the
temperature tensor
\begin{gather}
  \tTheta := (1-\nu)\theta\idtens+\nu\Theta,
  \label{tThetaDefinition}
\end{gather}
and $\nu$ is a tunable parameter which corresponds
to the Prandtl number, as we will see, according
to the relation $\Pr(1-\nu)=1$.
Define $\ptime$ by the requirement that
the viscosity satisfies $\viscosity=p\ptime$.
Then $\ptime = \Pr\htime$;
$\ptime$ is considered to be the relaxation period
of the standard BGK collision operator.
The collision operator $\Cgauss$ respects entropy
over the full range of $\nu$ values for which 
$\Theta\ge 0$
(is positive definite) implies
$\tTheta\ge0$,
that is, for $0\le(1-\nu)\le 3/2$,
corresponding to the range of
Prandtl numbers $2/3\le\Pr\le\infty$, which is surprising
since in the case $(1-\nu)>1$
(where $\theta\idtens$ is ``super-weighted'')
the affine combination \eqref{tThetaDefinition}
is nonconvex \cite{andries:GBGK99}.  The Gaussian-BGK model thus
gives an entropy-respecting collision operator
which allows the heat flux to be tuned from zero
($\Pr=\infty$) essentially up to the full heat flux for a monatomic
gas: for a monatomic gas with Maxwell molecules
$\Pr = 2/3$, and this is a good approximation for a broad range of
physical collision operators including hard spheres
and Coulomb collisions (in particular,
$\Pr\approx.61$ for ions and $\Pr\approx.58$
for electrons in the Braginskii closure \cite{article:Braginskii65}).

Henceforth in this section \ref{PerturbativeClosure}
we take as our standard of truth
the Galilean-invariant Boltzmann equation \eqref{SoleBoltzmannEqn}
with Gaussian-BGK collision operator \eqref{GaussBGKcol}
\begin{gather}
  \partial_t \f+\D_\xb\dotp(\v \f)+\D_\v\dotp(\a \f)
  = \Cgauss := \frac{\ftTheta - \f}{\htime},
  \label{StandardBoltzmannEqn}
\end{gather}
where $\ftTheta$ is given by the Guassian distribution equation
\label{GaussianDistribution} for $\tTheta$,
\begin{gather*}
  \ftTheta = \frac{\mdens}{\sqrt{\det(2\pi\tTheta)}}
    \exp\left(-(\v-\u)\dotp\tTheta^{-1}\dotp(\v-\u)/2\right),
\end{gather*}
and $\tTheta$ is given by equation \eqref{tThetaDefinition},
\begin{gather*}
  \tTheta := (1-\nu)\theta\idtens+\nu\Theta.
\end{gather*}

\subsection{Chapman-Enskog expansion}


Consider the Boltzmann equation
\begin{gather*}
  D \f = C[\f],
\end{gather*}
where $D=\partial_t \f+\D_\xb\dotp(\v \f) +\D_v\dotp(\a \f)$
is the convective derivative in phase space
and $C$ is the collision operator.
Assume the BGK collision operator
\begin{gather*}
  C = \frac{\fMaxwell - \f}{\ptime}.
\end{gather*}
This collision operator is linear and satisfies
$C[\fMaxwell]=0$.
(In contrast, the general Gaussian-BGK
collision operator causes $\tTheta$ to be modified
over time, resulting in nonlinear evolution of $\f$.)

The idea of the Chapman-Enskog expansion is to begin
with a guess $\f_0$ (which is the distribution about which
one is expanding, e.g.\ a Maxwellian), and solve by
successive approximations:
\begin{align*}
    D \f_0 &= C \f_1
 \\ D \f_1 &= C \f_2
 \\ D \f_2 &= C \f_3
 \\       &\vdots
\end{align*}
where it is assumed that one can solve e.g.\ for
$\f_1$ in terms of $\f_0$.  Hopefully the sequence
$\f_n$ converges rapidly to a solution $\f$.

A formalism which effects such a sequence is as follows.
Expand $\f$ as
\def\Cd{\bar C}
\begin{gather*}
  \f = \f\sup{0}
    + \epsilon \fd\sup{1}
    + \epsilon^2 \fd\sup{2}
    + \cdots,
\end{gather*}
where $\epsilon$ is a formal smallness parameter.
Let $C=\Cd/\epsilon$:
\begin{gather}\label{epsBoltzmann}
  D \f= \frac{\Cd[\fd]}{\epsilon}.
\end{gather}
Substituting the expansion of $\f$
and matching powers of $\epsilon$ yields the infinite sequence
\begin{align*}
    \epsilon^{-1}:&& 0 &= \Cd[\f\sup{0}],
 \\ \epsilon^0:&& D \f\sup{0} &= \Cd[\fd\sup{1}],
 \\ \epsilon^1:&& D \fd\sup{1} &= \Cd[\fd\sup{2}],
 \\ \epsilon^2:&&            &\vdots
\end{align*}
Multiplying each equation by its $\epsilon^n$ factor
and summing yields equation \eqref{epsBoltzmann}
(assuming convergence).
Choosing $\f\sup{0}=\fMaxwell$ means that
$C[\f\sup{0}]=0$ and we can get the chain
of approximations started.

To map onto the description in terms of successive approximations,
$\f_n = \sum_{k=0}^{n} \f\sup{k}$,
where $\f\sup{k}:=\epsilon^k \fd\sup{k}$.

The formal smallness parameter $\epsilon$ serves two purposes.
First, it specifies how to match up terms (to determine
which terms are used for the correction and
which for the previous estimate).
Second, it signifies that the corrections are expected
to decay; otherwise the series would not converge.
Once we have used $\epsilon$ to match up terms 
and get an approximate closure we can set $\epsilon=1$
as a shortcut to restoring the original collision
operator and making the replacements
$\f\sup{n}:=\epsilon^n \fd\sup{n}$.

\subsection{Perturbative stress closure}
\label{PerturbativeStressClosure}

To derive the five-moment closure for the deviatoric
pressure we take the deviatoric part of the pressure tensor evolution equation.
Neglecting collisions with other species, the pressure tensor evolution
equation \eqref{PTevolution} reads
\begin{gather}
   \partial_t \PT + \Div(\u\PT) 
     + \SymB(\PT\dotp\D\u) + \Div\qT
     = ({q/m})\SymB(\PT\times\B) + \RT,
  \label{solePTevolution}
\end{gather}
where recall that $\RT:=\intv \v\v C$.
Note that $\RT=\intc \c\c C$, by conservation
of mass ($\intv C=0$ and momentum ($\intc \c C=0$).
Also, $\RT$ is traceless by conservation of energy
($\intv |\c|^2 C=0$).
To prepare to take the deviatoric part 
we separate into isotropic and traceless parts:
\begin{align*}
  \SymB(\PT\dotp\D\u)
     &= \SymB((p\id+\dPT)\dotp\D\u)
  \\ &= p \SymB(\D\u) + \SymB(\dPT\dotp\D\u)
  \\ &= p \deviator{\SymB(\D\u)} + 2p\Div\u\id + \SymB(\dPT\dotp\D\u)
\end{align*}
Taking the deviatoric part of \eqref{PTevolution} thus gives
\begin{gather}
   \!\!\!\!\!\!
   \partial_t \dPT + \Div(\u\dPT) 
     +  p \deviator{\SymB(\D\u)} + \SymB(\dPT\dotp\D\u)
     + \Div\dqT
     = \tfrac{q}{m}\SymB(\dPT\times\B) + \RT,
  \label{dPTevolution}
\end{gather}
where recall that
$\glslink{dqT}{\ZdqT}:=\qT-(2/5)\SymC(\id\q)$
is the \mention{deviatoric heat flux tensor}.

If a BGK or Gaussian-BGK collision operator is used then
we can calculate the relaxation term:
\begin{align*}
  \RT &= \intc \c\c \Cgauss
       = \intc \c\c \frac{\ftTheta-\f}{\htime}
       = \htime\inv((1-\nu)p\id+\nu\PT - \PT)
       = \frac{1-\nu}{\htime}(p\id- \PT).
\end{align*}
That is, using $\ptime:=\htime/(1-\nu)$,
\begin{gather}
  \label{idRTclosure}
  \RT = -\dPT/\ptime.
\end{gather}
Note that this agrees with 
\eqref{RTclosure} in case $\Gcoef=\idfour$.

To perform a Chapman-Enskog expansion we
assume an initial guess (a Maxwellian)
for the velocity distribution
and substitute into the pressure tensor
evolution equation \eqref{solePTevolution}
The right side will be zero but the
left side will not.  So we modify the
solution used on the right hand side
so that it agrees with value assumed
by the left hand side for the original guess.
One could then substitute the modified
distribution into the left hand side
to obtain yet another correction from
the right hand side.  Hopefully this
process would converge to a solution
of the equation.

Formally, we assume a velocity distribution
expansion of the form
\begin{gather*}
  \f = \f\sup{0}
    + \epsilon \f\sup{1}
    + \epsilon^2 \f\sup{2}
    + \cdots,
\end{gather*}
where $\epsilon$ is a formal smallness
parameter which conveys the expectation
that the series is intended to converge
and which will indicate how to match up terms.
We seek a near-Maxwellian stress closure,
so we assume that $\f\sup{0}$ is Maxwellian
(with the same conserved moments as $\f$).
The expansion of $\f$ implies moment expansions such as
\begin{gather*}
  \PT = p\id
    + \epsilon \PT\sup{1}
    + \epsilon^2 \PT\sup{2}
    + \cdots,
\\\qT = \epsilon \qT\sup{1}
    + \epsilon^2 \qT\sup{2}
    + \cdots.
\end{gather*}
Using the closure equation \eqref{idRTclosure}
$\RT = -\dPT/\ptime$,
the right hand side of the pressure tensor evolution equation is
\begin{gather*}
   \partial_t \PT + \Div(\u\PT) 
     + \SymB(\PT\dotp\D\u) + \Div\qT
     = ({q/m})\SymB\left(\frac{\PT\times\B}{\epsilon}\right)
        - \frac{\dPT}{\epsilon\ptime},
\end{gather*}
where we have replaced $\B$ with $\B/\epsilon$
and $\ptime\inv$ with $\ptime\inv/\epsilon$
firstly as a way of indicating that $\ptime$ is
expected to be large and $\B$ could be large
and secondly in order to show how to match up terms
when we substitute the expansions in $\epsilon$.
The evolution equation for the deviatoric pressure
  \eqref{dPTevolution}
now reads
\begin{gather}
   \nonumber
   \partial_t \dPT + \Div(\u\dPT) 
     +  p \deviator{\SymB(\D\u)} + \SymB(\dPT\dotp\D\u)
     + \Div\dqT
     \\
     = ({q/m})\SymB\left(\frac{\dPT\times\B}{\epsilon}\right)
        - \frac{\dPT}{\epsilon\ptime}.
\end{gather}
Since the right hand side vanishes for $\f\sup{0}$
(a Maxwellian) we can match up powers of $\epsilon$.
For $\epsilon^0$ we obtain
\begin{gather}
  p \deviator{\SymB(\D\u)}
    = (q/m)\SymB(\dPT\sup{1}\times\B) - \dPT\sup{1}/\ptime.
  \label{eqn}
\end{gather}
At this point we can
``eliminate'' $\epsilon$ either by
replacing $\epsilon\dPT\sup{1}$ with $\dPT\sup{1}$,
$\B/\epsilon$ with the original expression $\B$, and
$\ptime\epsilon$ with the original expression $\ptime$,
all of which leaves equation \eqref{eqn} unchanged in appearance,
or (more simply) we can just set the formal smallness
parameter $\epsilon$ to 1.
Multiply \eqref{eqn} by $\ptime$.
Assume $\dPT\sup{1}\approx\dPT$.
Get the following implicit closure for the deviatoric stress:
\begin{gather}
  \boxed{
  p\ptime 2 \dstrain
    = \SymB(\dPT\times\pgf\b) - \dPT
  },
  \label{implicitDeviatoricStressClosure}
\end{gather}
where recall that $\pgf:=\ptime q|\B|/m = \ptime\gyrofreq$,
which agrees with the Stokes closure
\eqref{approximateClosureIdentity}
in case $\Gcoef=\idfour$.

\subsection{Perturbative heat flux closure}

McDonald and Groth have obtained a heat flux closure
using a Chapman-Enskog expansion.
This section follows their derivation
in section III.A in \cite{article:McDonaldGroth08},
generalizing to the case of nonzero magnetic field.

Multiplying the Boltzmann equation
\eqref{StandardBoltzmannEqn} by $\c\c\c$
and integrating by parts gives the evolution equation
\def\qTfour{\qT^{[4]}}
\begin{gather}
  \nonumber
  \partial_t\qT+\Div(\u\qT)
  + \SymC(\qT\dotp\D\u)
  - \SymC(\PT(\Div\PT)/\mdens)
  + \Div\qTfour
  \\
  = \SymC(\qT\times q\B/m) + \intc \c\c\c C,
  \label{heatFluxEvolution0}
\end{gather}
where $\qTfour := \intc \c\c\c\c \f$
is the fourth-order primitive moment.
For the Gaussian-BGK collision operator we compute that
\begin{gather*}
  \intc \c\c\c C
  = \intc \c\c\c \frac{\ftTheta-\f}{\htime}
  = \frac{0-\qT}{\htime}
  = \frac{-\qT}{\htime}
  = \frac{-\Pr}{\tau}\qT.
\end{gather*}
In a Chapman-Enskog expansion we assume 
a velocity distribution expansion
\begin{gather*}
  \f = \f\sup{0}
    + \epsilon \f\sup{1}
    + \epsilon^2 \f\sup{2}
    + \cdots,
\end{gather*}
where $\f\sup{0}$ is assumed to be the Gaussian
distribution $\f_\ZThetaT$ of $\f$.
This implies moment expansions such as
\begin{gather*}
  \qT = \epsilon \qT\sup{1}
    + \epsilon^2 \qT\sup{2}
    + \cdots,
 \\
  \qTfour = \qTfour\sup{0}
    + \epsilon \qTfour\sup{1}
    + \epsilon^2 \qTfour\sup{2}
    + \cdots.
\end{gather*}
For a Gaussian distribution we compute
(see equation \eqref{fourthMomentGauss}) that
$\qTfour= \SymC(\PT\PT)/\mdens$.
Following \cite{article:McDonaldGroth08},
Let
\def\Kfour{K^{[4]}}
\begin{gather*}
  \Kfour: = \intc \c\c\c\c (\f-\fGauss)
          = \qTfour - \SymC(\PT\PT)/\mdens
\end{gather*}
denote the deviation of the fourth moment
from the value for the Gaussian distribution.
Using
$ \SymC(\Div(\PT\PT/\mdens))
  - \SymC(\PT(\Div\PT)/\mdens)
 = \SymC\left(\PT\dotp\D\left(\frac{\PT}{\mdens}\right)\right),
$
equation \eqref{heatFluxEvolution0} becomes
\begin{gather}
  \nonumber
  \partial_t\qT+\Div(\u\qT)
  + \SymC(\qT\dotp\D\u)
 + \SymC\left(\PT\dotp\D\left(\frac{\PT}{\mdens}\right)\right)
  + \Div\Kfour
  \\
  = \frac{1}{\epsilon}
    \left(\SymC(\qT\times q\B/m) - \frac{\qT}{\htime}\right),
  \label{heatFluxEvolution}
\end{gather}
where we have rescaled the collision operator
(including the magnetic field) on the right hand side
by a factor of $\epsilon$ by replacing $\B$ with $\B/\epsilon$
and $\ptime$ with $\ptime\epsilon$.
Substituting the expansions
\begin{gather*}
  \qT = \epsilon \qT\sup{1}
    + \epsilon^2 \qT\sup{2}
    + \cdots,
 \\
  \Kfour = \epsilon \Kfour\sup{1}
    + \epsilon^2 \Kfour\sup{2}
    + \cdots
\end{gather*}
and matching powers of epsilon yields the
implicit heat flux closure equation
\begin{gather*}
   \SymC\left(\PT\dotp\D\left(\frac{\PT}{\mdens}\right)\right)
  = \SymC(\qT\times q\B/m) - \frac{\qT\sup{1}}{\htime}.
\end{gather*}
Set the formal smallness parameter $\epsilon$ to 1.
Multiply by $\htime=\ptime/\Pr$.  Take $\qT\approx\qT\sup{1}$.
Get the following implicit closure for the heat flux tensor:
\begin{gather}
   \frac{\ptime}{m\Pr}\SymC\left(\PT\dotp\D\left(\frac{\PT}{\mdens}\right)\right)
  = \SymC(\qT\times\hgf\b) - \qT,
  \label{McDonaldImplicitHeatFluxTensorClosure}
\end{gather}
which by setting $\hgf=0$
is observed to be a generalization of the heat
flux closure given by McDonald and Groth in equation (23) 
of \cite{article:McDonaldGroth08};
here we define
\begin{gather*}
   \glslink{hgf}{\Zhgf}
    := \frac{\pgf}{\Pr}
     = \frac{\ptime\gyrofreq}{\Pr}
     = \htime \gyrofreq
     = \htime q|\B|/m.
\end{gather*}
We rewrite \eqref{McDonaldImplicitHeatFluxTensorClosure} in the form
\begin{gather}
 \boxed{
     \qT + \SymC(\hgf\b\times\qT) =
    -\tfrac{2}{5}\heatConductivity \SymC\left(\Pshape\dotp\D\TT\right)
  },
  \label{implicitHeatFluxTensorClosure}
\end{gather}
where $\heatConductivity$ is the \mention{heat conductivity},
related to the Prandtl number by the formula \eqref{PrFormula}
\begin{align*}
  \glslink{Pr}{\Pr} = \frac{5}{2}\frac{p\ptime}{m\heatConductivity},
\end{align*}
and we have used that $\PT/\mdens=\TT/m$ and $\Pshape:=\PT/p$.
The explicit solution to equation 
\eqref{implicitHeatFluxTensorClosure}
is obtained in the appendix as equation \eqref{qTclosure}.

Half the trace of equation \eqref{implicitHeatFluxTensorClosure}
will give an implicit closure for the heat flux.
Using that $\tr \SymC(\qT\times\b)=2\q\times\b$
and that
$
\tr \SymC(\Pshape\dotp\D\TT)
= 3\Pshape\dotp\D T + 2\Pshape\ddotp\D\TT
= 5\Pshape\dotp\D T + 2\Pshape\dotp\D\dTT,
$
half the trace of \eqref{implicitHeatFluxTensorClosure} is
\begin{gather}
  \q + \hgf\b\times\q =
   -\heatConductivity \left(\Pshape\dotp\D T + \tfrac{2}{5}\Pshape\dotp\D\dTT\right).
\end{gather}
Neglecting deviatoric parts (which would be of order $\epsilon$
in a Chapman-Enskog expansion about a Maxwellian distribution)
and recalling that $\Pshape=\id+\dPshape$ yields
\begin{gather}
   \boxed{
   \q + \hgf\b\times\q =
   -\heatConductivity \D T
   },
   \label{implicitHeatFluxClosure}
\end{gather}
which agrees with the implicit heat flux closure equation
(5.107) in \cite{book:Woods04} in case $\Pr=1$.

To get an explicit heat flux closure we need to solve 
the implicit closures for the heat flux.
In section \ref{HeatFlux} we solve this equation for $\q$
to obtain the heat flux closure \eqref{expHeatFlux},
{
\def\w{{\id_\wedge}}
\def\l{{\id_\Vert}}
\def\p{{\id_\perp}}
\def\v{{\hgf}}
\def\o{\stimes}
\begin{gather}
     \q = -\heatConductivity
     \left(\b\b + \frac{1}{1+\v^2}\Big((\id-\b\b) - \hgf \b\times\id\Big)\right)\cdot\D T,
     \label{expHeatFlux}
\end{gather}
I remark that this closure satisfies the form of
an entropy-respecting closure specified in equations
\eqref{entropyRespectingHeatFluxForm} and
\eqref{entropyRespectingHeatFluxClosureRequirements}.
}

\subsection{Comparison of heat flux closures}

In the absence of a magnetic field, the McDonald-Groth
(Chapman-Enskog Gaussian-BGK) heat flux closure
\eqref{implicitHeatFluxTensorClosure} is
\def\qTM{\qT_\mathrm{M}}
\begin{gather*}
     \qTM = -\tfrac{2}{5}\heatConductivity \SymC\left(\Pshape\dotp\D\TT\right),
\end{gather*}
whereas the entropy-respecting closure \eqref{qTlevermore} is
\begin{gather}
  \qTL = \kLa \D\vee\TT\inv + \kLb \id\vee\tr(\D\vee\TT\inv).
\end{gather}
In the isotropic case (see \eqref{isotropicSimplifications}),
\begin{align*}
   \D\vee\TT\inv = T^{-2} \id\vee\D T
   \hbox{and}
   \id\vee\tr(\D\vee\TT\inv) &= 5 T^{-2}\id\vee\D T,
\end{align*}
and these closures simplify to (see \eqref{isotropicLevermore})
\begin{align*}
  \qTL &= -T^{-2}(\kLa + 5 \kLb) \id\vee\D T&
     &\hbox{and}&
  \qTM &= -\tfrac{2}{5}\heatConductivity \id\vee\D T.
\end{align*}
So for consistent closure coefficients we require that
\begin{gather*}
   \kLa + 5 \kLb = \tfrac{2}{5}T^2\heatConductivity,
\end{gather*}
which is analogous to equation \eqref{heatConductivityVersusKappa}.
In general, however, the two closures disagree.
Recall from \eqref{LevermoreComparisonForm} that
\begin{align*}
  \D\vee\TT\inv = -T^{-2}\SymC(\D\TT\ddotp(\Pshape\inv\Ldiamond\Pshape\inv)).
\end{align*}
The entropy-respecting closures with consistent closure coefficients are thus
\begin{align*}
  \qT_1 &= -\tfrac{2}{5}\heatConductivity
    \underbrace{\SymC\left(\D\TT\ddotp\Pshape\inv\Ldiamond\Pshape\inv\right)}
      _{\hbox{$\id\vee\D T$ if $\Pshape=\id$}}
    && \hbox{(if $\kLb=0$) and}
  \\
  \qT_2 &= -\tfrac{2}{25}\heatConductivity
    \id\vee\underbrace{\tr\left(\SymC\left(\D\TT\ddotp\Pshape\inv\Ldiamond\Pshape\inv\right)\right)}
      _{\hbox{$5 \D T$ if $\Pshape=\id$}}
    && \hbox{(if $\kLa=0$)}
\end{align*}
and convex combinations thereof, in contrast to the McDonald-Groth closure
\begin{gather*}
     \qTM = -\tfrac{2}{5}\heatConductivity \SymC\left(\Pshape\dotp\D\TT\right).
\end{gather*}
In the isotropic limit $\Pshape\approx\id$ and the two
closures agree.
If $\Pshape\approx\id$ and
if $\D\dTT$ is small (i.e.\ $|\D\dPshape|\ll|\D \ln T|$)
then the closures simplify to
\def\qTI{\qT_{\id}}
\begin{gather*}
     \qTM \approx \qTL \approx 
     \qTI := -\tfrac{2}{5}\heatConductivity \id\vee\D T.
\end{gather*}
I remark that the entropy-respecting
constraint assumed by the Levermore closure is not physically
justified except in the near-Maxwellian limit,
whereas the McDonald-Groth closure is physically justifiable
in the near-Gaussian limit.
One wonders whether the closures $\qTL$ and $\qTI$
are well-posed.

\section{Intraspecies closure coefficients}

The five-moment and ten-moment gas-dynamic closures
are determined by the relaxation period
$\glslink{ptime}{\ptime}$
of the deviatoric stress and the relaxation period
$\htime$ of the heat flux, related to $\ptime$
by $\htime=\ptime/\glslink{Pr}{\Pr}$.
The viscosity is given by $\glslink{viscosity}{\Zviscosity}=p\ptime$
and the thermal conductivity by
$\glslink{heatConductivity}{\ZheatConductivity}=\frac{5}{2}\frac{\viscosity}{m\Pr}$.

To obtain full closure we need to specify
$\ptime$ in terms of the evolved moments.
We therefore assume that $\ptime$ is a function
of the evolved moments, that is, of the
Maxwellian or Gaussian distribution
which is in bijective correspondence with the evolved moments.
This is a simplifying assumption, since
for physical collision operators $\ptime$ and $\htime$
actually depend on the details of the distribution.
So in general we write, e.g., $\ptime(\f)$.

\subsection{Viscosity and heat conductivity are functions of temperature
 independent of density.}

By Galilean relativity $\ptime$ is a function
of the density and the primitive moments 
of order greater than 1 (e.g.\ the pressure).
In the five-moment case we can write
$\ptime(n,T)$ and in the ten-moment case
$\ptime(n,\TT)$.
We want to characterize these functions.

Most classical collision operators satisfy
\begin{gather}
  C(\lambda \f) = \lambda^2 C(\f).
  \label{secondOrderHomogeneousProperty}
\end{gather}
We will argue that, for such collision operators, closures
for $\viscosity$ and $\heatConductivity$ depend only on temperature.

Property \ref{secondOrderHomogeneousProperty}
essentially says that doubling the number of particles
doubles the rate of collisions experienced by a particle
and thus doubles the rate of evolution of the shape
of the distribution without otherwise altering its shape.
Formally, if $\f$ satisfies the spatially homogeneous Boltzmann equation
\def\ft{\tilde f}
\begin{gather*}
  \partial_t \f = C(\f)
\end{gather*}
then so does $\ft(t) = \lambda \f(\lambda t)$.

Therefore, the relaxation periods scale according to
$\ptime(\lambda \f) = \lambda\inv\ptime(\f)$.
The moments are a function of $\f$ and scale as:
$n(\lambda \f) = \lambda n(\f)$,
$ T(\lambda \f)=T(\f)$,
and $\TT(\lambda \f)=\TT(\f)$.

Therefore, the viscosity $\viscosity=n T\ptime$
and the thermal conductivity scale as
$\viscosity(\lambda \f) = \viscosity(\f)$ and
$\heatConductivity(\lambda \f) = \heatConductivity(\f)$.

We now invoke the simplifying assumption e.g.\ that
$\viscosity(n,\TT)$.  This says that $\viscosity$ is the same
for any two distributions which share the same moments.
We show that $\viscosity$ is actually independent of $n$
by differentiating the constant
map $\lambda\mapsto \viscosity(\lambda \f) = \viscosity(\f)$
with respect to $\lambda$:
$0 = d_\lambda \viscosity(\lambda \f)
   = d_\lambda \viscosity(n(\lambda \f),\TT(\lambda \f))
   = d_\lambda \viscosity(\lambda n(\f),\TT(\f))
   = n \partial_n \viscosity$,
as needed.

In conclusion, for the ten-moment model
\begin{align*}
  \ptime(n,\TT) &= \frac{\viscosity(\TT)}{nT}, & \htime &= \ptime/\Pr.
   & \heatConductivity(\TT) &= \frac{5}{2}\frac{\viscosity(\TT)}{m\Pr}.
\end{align*}


\subsection{Positivity-preserving heat conductivity closure}
\label{PositivityPreservingClosure}

Contrary to the case for five-moment entropy
evolution, in the ten-moment case a local minimum
can decrease if a nonzero heat flux closure is used. 
The closure for $\heatConductivity$
is critically important in maintaining positivity.

To ensure that the heat flux closure maintains positivity
of the temperature, the heat conductivity
$\heatConductivity$ should go to zero as 
$\det(\TT)$ goes to zero.  This effectively shuts down heat
flow to prevent positivity violations.
Five-moment closures for $\heatConductivity$ 
naturally are specified in terms of the scalar temperature $T$.
When using such a closure formula in the ten-moment model,
$T$ (i.e.\ $\tr\TT/3$, the arithmetic average of the
eigenvalues of $\TT$) should be replaced with
$(\det\TT)^{1/3}$ (i.e.\ the geometric average of the
eigenvalues of $\TT$)
\cite{levermore:priv11}.

For the pressure tensor relaxation period $\ptime$ one
can choose to define $\ptime$ in terms of $\tr\TT$
or $\det\TT$.


In particular, the Braginskii five-moment closure
(discussed in the next subsection) is of the form
\begin{gather}
  \label{isoPeriod}
  \ptime = \ptime_0 \sqrt{m}\frac{T^{3/2}}{n}.
\end{gather}
In the ten-moment model one could use
\begin{gather*}
  \ptime = \ptime_0 \sqrt{m}\frac{(\det\TT)^{1/2}}{n}.
\end{gather*}
The viscosity is then
\begin{gather*}
  \viscosity = p\ptime = \ptime_0 \sqrt{m}T(\det\TT)^{1/2}
\end{gather*}

\def\tauBr{\ptime^\mathrm{Br}}
\subsection{Braginskii closure coefficients}
\label{BraginskiiClosureCoefficients}

In the absence of a magnetic field the viscosities
given by Braginskii \cite{article:Braginskii65} are\footnote{
  I have used the symbols $\tauBr_\i$ and $\tauBr_\e$
  to distinguish them from $\ptime_\i$ and $\ptime_\e$,
  which are defined in this document to be the isotropization
  period of the pressure tensor, equivalently 
  the collision period taking all species into account.

  My understanding is that $\tauBr_\i$ is intended to
  be the ion-ion collision period (which approximates
  the ion collision period), whereas
  $\tauBr_\e$ is intended to be the electron-ion collision
  period in a Lorentzian plasma; the factor of $\sqrt{2}$
  evidently arises from the fact that the reduced mass in an ion-ion
  collision is half the ion mass.  See \cite{article:Braginskii65}
  pages 220 and 277.

  Balescu objects to the artificial dissimilarity between
  ions and electrons that this convention introduces into
  the formulas and instead redefines $\ptime_\i$ so that its
  formula agrees with the formula for $\ptime_\e$
  (see his footnote on page 274 of \cite{Balescu88}),
  but Braginskii's definitions seem to have become fairly standard.
}
\begin{equation}
  \begin{matrix}
    \viscosity_\i = .96 \tauBr_\i p_\i = .96 \ptime'_{\i\i}p_\i,
 \\ \viscosity_\e = .73 \tauBr_\e p_\e = .52 \ptime'_{\e\e}p_\e,
  \end{matrix}
  \qquad
  \begin{matrix}
    \tauBr_\i := \ptime'_{\i\i},
 \\ \tauBr_\e := \frac{\ptime'_{\e\e}}{\sqrt{2}},
  \end{matrix}
  \qquad
  \begin{matrix}
    \ptime'_{\i\i} := \ptime_0 \sqrt{m_\i}\frac{T_\i^{3/2}}{n_\i},
 \\ \ptime'_{\e\e} := \ptime_0 \sqrt{m_\e}\frac{T_\e^{3/2}}{n_\e},
  \end{matrix}
  \label{ionElcViscosity}
\end{equation}
where the base isotropization period is
\begin{align}
 \ptime_0 &= \frac{12 \pi^{3/2}}{\ln\Lambda}\left(\frac{\epsilon_0}{e^2}\right)^2,
 \label{baseIsoPeriod}
\end{align}
where $\ln\Lambda$ is the Coulomb logarithm, discussed below.

In a one-species charged gas the ion viscosities would be
\begin{align*}
    \viscosity_\i = \ptime'_{\i\i}p_\i,
 \\ \viscosity_\e = \ptime'_{\e\e}p_\e.
\end{align*}
In the two-species gas, however, isotropization of each species
is accelerated by interspecies collisions.
Therefore, the viscosities are given by
\begin{align*}
    \viscosity_\i = \ptime_{\i}p_\i,
 \\ \viscosity_\e = \ptime_{\e}p_\e,
\end{align*}
where the overall isotropization rates are given by
\begin{gather*}
  \ptime_\i\inv \approx \ptime_{\i\i}\inv+\ptime_{\i\e}\inv,
\\\ptime_\e\inv \approx \ptime_{\e\e}\inv+\ptime_{\e\i}\inv.
\end{gather*}
Since the mass ratio is large,
\begin{alignat*}{5}
  \ptime_\i\inv &\approx \ptime_{\i\i}\inv &&\gg \ptime_{\i\e}\inv,
\\ \tfrac{1}{2}\ptime_\e\inv &\approx \ptime_{\e\e}\inv &&\approx \ptime_{\e\i}\inv
\end{alignat*}
explaining the differing coefficients in equation
\eqref{ionElcViscosity}.
In a pair plasma we expect that $\ptime_{\i}\inv = 2\ptime_{\i\i}\inv$.



The Coulomb logarithm is the logarithm of the
plasma parameter $\Lambda$ which 
is on the order of the number of particles in a
Debye sphere, roughly
\begin{align*}
 \Lambda &= n\lambda_D^3,
\end{align*}
where the Debye length
(see section \ref{CollisionlessNondimensionalization}) is given by
\begin{align*}
 \lambda_D &= \sqrt{\frac{\epsilon_0 T_e}{n_e e^2}};
\end{align*}
typically $10\lesssim\ln\Lambda\lesssim 20$.

Note that these collision periods are not identical to
the isotropization periods $\ptime_\s$ appearing elsewhere in this
dissertation, which I define to be the viscosity divided by the pressure:
$\viscosity=p\ptime$.

Comparing Braginskii's viscosities and heat conductivities,
\begin{align*}
  \viscosity_\e &= 0.73 p_\e\tauBr_\e, & \viscosity_\i &= 0.96 p_\i\tauBr_\i,
  \\
  \heatConductivity_\e &= 3.16 p_\e\tauBr_\e/m_\e, &
  \heatConductivity_\i &= 3.91 p_\i\tauBr_\i/m_\i,
\end{align*}
the Prandtl numbers for the unmagnetized Braginskii closure are
\begin{align*}
  {\Pr}_\e &= .58, & {\Pr}_\i &= .61,
\end{align*}
where we have used the definition \eqref{PrFormula}
\begin{gather*}
  {\Pr}_\s := \frac{5}{2}\frac{\viscosity_\s}{m \heatConductivity_\s}.
\end{gather*}
Recall that for a monatomic gas the Prandtl number should be
close to $2/3 = .\overline{66}$.

\section{Interspecies collisional closures}

Interspecies collisions are generally ignored
in this work, on the assumption that in a weakly
collisional plasma we can lump the large majority
of the microscale effects into the intraspecies
collision operators.  That is, we are assuming
that momentum and energy are largely conserved
within each species.  As an illustration of this
principle, note that even in the case of
Coulomb collisions one of the chief
effects of electron-ion collisions is to increase
the rate of isotropization of the
electrons.  Our Gaussian-BGK ``intraspecies'' collision
operator really includes the interspecies particle
interactions to the degree that these interactions
conserve momentum and energy within each species.

Nevertheless, it would be of interest to incorporate
an interspecies collision operator.  To obtain
a form for the interspecies collision operator one
can require interspecies collisions not to decrease
entropy.

\def\Cie{C_{\i\e}}
\def\Cei{C_{\e\i}}
For an interspecies BGK collision operator $C_{\i\e}(\f_\i,\f_\e)$,
one could make the loss term proportional to the distribution
function and define the gain term so as
to satisfy the following requirements:
\begin{itemize}
\item conservation of total mass, momentum, and energy:
  \begin{gather*}
    \int_\v(\C_{\i\e})  = 0, \\
    m_\i \int_\v(\v C_{\i\e}) + m_\e \int_\v(\v C_{\e\i}) = 0, \\
    m_\i \int_\v(|\v|^2 C_{\i\e}) + m_\e \int_\v(|v|^2 C_{\e\i}) = 0,
  \end{gather*}
\item bilinearity of
  $\Cie(f_\i,f_\e)$
  and
  $\Cei(f_\e,f_\i)$,
\item total entropy is nondecreasing in the near-Maxwellian limit,
  \begin{gather*}
  \int_\v(C_{\i\e} \log f_\i) + \int_\v(C_{\e\i} \log f_\e) \gtrapprox 0,
  \end{gather*}
\item agreement with Gaussian-BGK for agreeing distributions,
  e.g., if $f_\i=f_\e$ and $m_\i=m_\e$ 
  then $C_{\i\i}(f_\i,f_\i) = C_{\i\e}(f_\i,f_\e)$, and
\item Galilean invariance.
\end{itemize}

A heuristic to design such an interspecies collision operator is
first to specify a ``driftless interspecies collision operator''
for the case that the interspecies drift velocity is zero.
For the simple BGK case (relaxation to a Maxwellian) one could
calculate the Maxwellian distributions that the two distributions
would have in equilibrium and relax the distributions of the two
species toward their respective equilibrium Maxwellians at the
same rate.

To handle drifting distributions one can independently
relax the drift velocity to zero and relax the
distribution shapes toward an average.  For the
ten-moment model this leads to a fluid closure which relaxes 
toward a common temperature tensor.

%
So for the Maxwellian case the driftless collision
operators would be:
\def\Maxwell{\mathcal{M}}
\begin{gather*}
     C'_{\i\e} = (\Maxwell(n_\i,\u_\i,T) - f_\i)/\ptime,
  \\ C'_{\e\i} = (\Maxwell(n_\e,\u_\e,T) - f_\e)/\ptime,
\end{gather*}
where $\Maxwell$ denotes a Maxwellian distribution
and (for monatomic species)
\begin{gather*}
  T = (n_\i T_\i + n_\e T_\e)/(n_\i+n_\e)
\end{gather*}
to conserve total thermal energy.

To generalize to the Gaussian-BGK case, one would instead use
\def\Gauss{\mathcal{G}}
\def\TTg{\widetilde\TT}
\begin{gather*}
   \Cie = (\Gauss(n_\i,u_\i,\TTg) - f_\i)/\tau,
\\ \Cei = (\Gauss(n_\e,u_\e,\TTg) - f_\e)/\tau,
\end{gather*}
where
\def\TTgi{{\widetilde\TT_\i}}
\def\TTge{{\widetilde\TT_\e}}
\begin{align*}
      \TTg &= (n_\i\TTgi + n_\e\TTge)/(n_\i+n_\e),
  \\ \TTgi &= \nu\TT_\i + (1-\nu)T_\i\id,
  \\ \TTge &= \nu\TT_\e + (1-\nu)T_\e\id.
\end{align*}
This should respect entropy for near-Maxwellian distributions.

To generalize to the case $\u_\i \ne 0 \ne \u_\e$, we can use
the same collision operators specified above for the thermal
equilibration part and handle the resistive drag separately. The
constraints above on the collision operator are satisfied if
the resistive drag force on the species is equal and opposite,
so there is essentially complete freedom in specifying the
magnitude and direction of resistive drag and the allocation of
resistive heating among species and spatial directions. Closures
such as Braginskii's would imply the value of the drag force,
and one would expect heating to be allocated among species in
inverse proportion to particle mass. I am inclined to allocate
resistive heating primarily perpendicular to the direction of
drift velocity, although for some reason Miura and Groth \cite{MiuraGroth07}
allocate it primarily parallel to the drift velocity in their ten-moment
closure.



\section{MHD equations}

A \defining[magnetohydrodynamics (MHD)]{magnetohydrodynamic (MHD) fluid}
is a fluid that conducts electricity.  In this document we use MHD
to refer to models of plasma that evolve a single-fluid description
of mass density and momentum density.

MHD allows us to write a description of plasma evolution that is fully
Galilean-invariant.  Although Lorentz-invariant formulations of MHD
have also been formulated, in this discussion we take MHD to mean
Galilean-invariant MHD.  Consider Maxwell's equations \eqref{MaxwellsEqns}
in the form
\begin{align*}
   \partial_t \B + \curl\E &= 0,                &\Div\B &= 0,
\\ -c^{-2}\partial_t \E + \curl\B &= \mu_0\J,  &c^{-2}\Div\E &= \mu_0\qdens,
\end{align*}
where $\mu_0:=1/(c^2\epsilon_0)$.
In the limit $c\to\infty$ the displacement current
$\partial_t \E$ and the net charge $\qdens$ go to zero and we get
the Galilean-invariant system
\begin{equation}
  \label{MaxwellMHD}
  \begin{aligned}
   \partial_t \B + \curl\E &= 0, &\Div\B &= 0,
\\ \curl\B &= \mu_0\J,           &0 &= \mu_0\qdens.
  \end{aligned}
\end{equation}
Note that this does \emph{not} imply that $\Div\E=0$.
(In fact, $\Div\E$ is not a Galilean-invariant quantity
and in the Galilean limit transforms according to
$\Div\E' = \Div\E+d\v\dotp\J/\mu_0$, where $d\v=\v-\v'$.)
For this reason, in the MHD limit people often speak of
the assumption of charge \gls{quasineutrality} rather
than charge neutrality.  This should not be misunderstood.
The MHD model constitutes a self-consistent Galilean-invariant
model that assumes exact charge neutrality.

The MHD limit fundamentally alters the ``causal'' relationship
of electromagnetic quantities.  In the Lorentz-invariant Maxwell
equations a prescribed $\J$ and $\qdens$ (and initial conditions)
determine $\E$ and $\B$.  In the Galilean-invariant limit
$\J$ is determined from $\B$, and $\E$ must be externally supplied
from a fluid equation called \first{Ohm's law}.

\subsection{Ohm's law}

\defining{Ohm's law} is the evolution equation for current
density $\J$, solved for the electric field $\E$.
The evolution equation for $\J$ is obtained by summing
current evolution for each species over all species.
Current evolution for a single species $\s$
is momentum evolution \eqref{momentumEvolution}
times its charge-to-mass ratio $q_\s/m_\s$:
\begin{gather*}
  \partial_t \J_\s + \Div(\u_\s\J_\s + (q_\s/m_\s)\PT_\s)
  = (q_\s^2/m_\s) n_\s(\E+\u_\s\times\B) + (q_\s/m_\s)\R_\s,
\end{gather*}
where $\PT^q_\s:=(q_\s/m_\s)\PT_\s$
is the \first{electrokinetic pressure tensor}.
Summing over all species gives net current evolution,
\begin{gather*}
  \partial_t \J + \Div(\u\J + \J\u -\qdens\u\u + \sum_\s \qdens_\s \w_\s\w_\s)
    + \sum_\s (q_\s/m_\s)\Div \PT_\s
  \\ =
    \sum_\s (q_\s^2/m_\s) n_\s(\E+(\u+\w_\s)\times\B)
    + \sum_\s (q_\s/m_\s)\R_\s,
\end{gather*}
where the total momentum density $\mdens\u := \sum_\s \mdens_\s\u_\s$
defines the net fluid velocity $\u$
and the \first{species drift velocity} is defined by
$\glslink{w}{\Zw_\s}:=\u_\s-\u$.
To infer the species drift velocities $\w_\s$ from
current, for a multispecies fluid one must impose
constitutive assumptions, but for a two-species fluid
e.g.\ of ions $\i$ and electrons $\e$
the assumption of charge neutrality
\begin{gather*}
  \glslink{ni}{n_\i}=\glslink{ne}{n_\e}=:\glslink{ndens}{n}
\end{gather*}
allows one to infer species drift velocity from current.
Charge neutrality says that net current is independent of
reference frame.  In the reference frame of the fluid we can
solve the definitions of charge density and momentum density
\begin{align*}
  \J &= \J_\i + \J_\e, \\
     0 &= \frac{\glslink{mi}{\mi}}{q_\i} \J_\i + \frac{\glslink{me}{\me}}{q_\e} \J_\e
\end{align*}
for the drift velocities in terms of the current
and the ratios of mass to charge of the species.
For simplicity we take $\glslink{q }{q_\i}=e$ and $\glslink{q }{q_\e}=-e$,
where \gls{e } is the charge on a proton.  Then
\begin{align}\label{driftvel}
    \J_\s &= \frac{\mred}{m_\s}\J &\hbox{and}&&
    \w_\s &= \frac{\J_\s}{n q_\s}.
\end{align}
Current evolution simplifies to
\begin{gather*}
  \!\!\!\!\!\!\!\!\!\!\!\!
  \partial_t \J + \Div\left(\u\J+\J\u-\frac{\dmt}{\mdens}\J\J\right)
    + \Div\left(\frac{\PT_\i}{\mti}-\frac{\PT_\e}{\mte}\right)
  = \frac{en}{\mut}\left(\E+\left(\u-\frac{\dmt}{\mdens}\J\right)\times\B
    -\Tresistivity\dotp\J\right);
\end{gather*}
here the reduced mass \gls{mred} is defined by $\mred^{-1}:= \mi ^{-1} + \me ^{-1}$,
the mass difference is defined by $\glslink{dm}{\Zdm} = \mi -\me $,
and we use tilde to indicate division by $e$, so
$\glslink{mut}{\Zmut}:=\mred/e$, 
$\glslink{dmt}{\Zdmt}:=\dm/e$, 
$\glslink{mti}{\Zmti}:=\mi/e$, 
and
$\glslink{mte}{\Zmte}:=m_e/e$.
The resistivity \gls{resistivity} is related to the interspecies drag force by
$-\R_\i = \R_\e = n e\Tresistivity\dotp\J$.
The quantity $-\frac{\dmt}{\mdens}\J = \w_\i+\w_e$
is twice the velocity of the charges relative to the fluid.

Solving for $\E$ gives Ohm's law,
\begin{gather}\label{OhmsLaw}
  \!\!\!\!\!\!\!\!
  \E = \Tresistivity\dotp\J + \B\times\left(\u-\frac{\dmt}{\mdens}\J\right)
       + \frac{\mut}{e n}\left[
       \Div\left(\frac{\PT_\i}{\mti}-\frac{\PT_\e}{\mte}\right)
       + \partial_t\J + \Div\left(\u\J+\J\u-\frac{\dmt}{\mdens}\J\J\right)
       \right].
\end{gather}
We write Ohm's law in the form
\begin{gather}\label{OhmsLawSep}
  \E = \B\times\uc + \E',
\end{gather}
where
\begin{gather}\label{chargeVelocityDef}
  \uc := \u-\frac{\dmt}{\mdens}\J
\end{gather}
is the \defining{charge velocity} (defined to be the fluid
velocity plus the sum of the drift velocities of
both species) and where
\begin{gather*}
  \E' := \Tresistivity\dotp\J
       + \frac{\mut}{e n}\left[
       \Div\left(\frac{\PT_\i}{\mti}-\frac{\PT_\e}{\mte}\right)
       + \partial_t\J + \Div\left(\u\J+\J\u-\frac{\dmt}{\mdens}\J\J\right)
       \right]
\end{gather*}
is the nonideal component of the electric field.
The expression $\B\times\uc$ is the ideal electric field
of Hall MHD.

\def\muinv{\mu_0^{-1}}
In the MHD model the Ohm's law expression for electric field
is used in the evolution equation $\partial_t\B+\curl\E=0$
for the magnetic field, and Ampere's law (from \eqref{MaxwellMHD})
\begin{gather}
  \label{AmperesLaw}
  \J=\muinv\curl\B
\end{gather}
is used to define the current.
Thus the full evolution equation for the magnetic field is
{ 
  \small
\begin{align*}
  \!\!\!\!\!\!\!
  \partial_t \B + \curl\bigg(&
  \Tresistivity\dotp(\muinv\curl\B) + \B\times\left(\u-\frac{\dmt}{\mdens}\muinv\curl\B\right)
   + \frac{\mut}{e n} \Div\left(\frac{\PT_\i}{\mti}-\frac{\PT_\e}{\mte}\right)
\\ + &\frac{\mut}{e n}\left[
     \muinv\curl\partial_t\B
   + \muinv\Div\left(\u\curl\B+(\curl\B)\u-\muinv\frac{\dmt}{\mdens}(\curl\B)\curl\B\right)
       \right]\bigg) = 0.
\end{align*}
}
This is an implicit differential equation for $\B$
and requires an implicit numerical method.  Use of an
explicit method requires some sort of simplification.
Ideal MHD simplifies magnetic field evolution to
\begin{align*}
  \partial_t \B + \curl(\B\times\u) = 0,
\end{align*}
discarding all other terms; this is a hyperbolic system
and is naturally suited to an explicit method.
The Hall term
\begin{align*}
  \curl\left(\B\times\left(\frac{\dmt}{\mdens}\muinv\curl\B\right)\right)
\end{align*}
is strongly dispersive, especially if the $\partial_t\B$
terms is left out, and calls for an implicit method.

\subsection{Mass density evolution}

Since MHD assumes charge neutrality,
its representation of species densities should enforce this constraint.
In a charge-neutral two-species fluid
the density of each species can be inferred from the total
mass density.  Therefore we evolve the total mass density.
The evolution equation of (total) mass density
is the sum of the evolution equations \eqref{densityEvolution}
for the mass density of the individual species.
It reads
\begin{gather}
  \label{totalMassEvolution}
  \partial_t \mdens+\Div(\mdens\u)=0,
\end{gather}
where $\mdens=\mdens_\i+\mdens_\e$ and the total fluid velocity
is defined by the conservation requirement that the total momentum
density be the sum of the momentum densities of the individual species,
$\mdens\u:=\mdens_i\u_\i+\mdens_\e\u_\e$.

We can infer the mass density and number density of each species
from the total mass density using the relations
\begin{align*}
  \mdens&=(\mi+\me)n, &\mdens_\i&=n \mi, &\mdens_\e&=n \me.
\end{align*}

Recall from equation \eqref{driftvel}
that we can infer the species drift velocity
$\glslink{w}{\Zw_\s}:=\u_\s-\u$.
from the current:
\begin{align}\label{wsJ}
    \w_\s &= \frac{\mred}{m_\s}\frac{\J}{n q_\s}.
\end{align}
With this we may check that
mass evolution, equivalently
number density evolution \eqref{spcNumDensEvolution}
\begin{gather}
  \label{spcNumDensEvolutionCopy}
  \partial_t n_\s + \Div(n_\s \u_\s),
\end{gather}
is satisfied for each species.  Indeed,
$n_\s=n$ and
$
  \Div(n \u_\s)
  = \Div(n \u)
  + \Div(n \w_\s)
$
and $\Div(n \w_\s) = 0$ by equation \eqref{wsJ} because
by the Galilean-invariant Ampere's law $\Div\J=\mu_0\inv\Div\curl\B=0$.
So number density evolution for each species reduces to the same
assertion,
\begin{gather}\label{numberDensityEvolution}
  \partial_t n + \Div(\u n)=0,
\end{gather}
which is equivalent to the evolution equation
\eqref{totalMassEvolution}
for total mass density.

\subsection{Momentum density evolution}

MHD evolves an evolution equation for the net momentum density.
It is derived by summing the density evolution equations of the
individual species.
Knowledge of number density $n$, net momentum density $\mdens\u$,
and current density $\J$
is sufficient to infer the momentum density $\mdens_\s\u_\s$
of each species.  Indeed,
$\mdens_\s\u_\s = m_\s(\u+\w_\s)$
where $\w_\s = \frac{\mred}{m_\s}\frac{\J}{n q_\s}$
(equation \eqref{wsJ}).
Therefore, there is no need to evolve separate momentum evolution
equations for each species (and doing so would result in inconsistency
due to numerical error or use of an approximate Ohm's law).

From another viewpoint,
recall that Ohm's law is current evolution solved for electric field.
Exact current evolution is a linear combination
(with weights $e/m_\i$ and $-e/m_\e$)
of the momentum evolution equations of the two species.
Total momentum evolution is the sum of the momentum evolution
equations.  So we again see that momentum evolution plus Ohm's law
is equivalent to specifying momentum evolution for each species.

Summing the species momentum evolution equations \eqref{momentumEvolution},
\begin{gather*}
  \partial_t (\mdens_\s\u_\s) + \Div(\mdens_s\u_\s\u_\s + \PT_\s)
  = q_\s n_\s(\E+\u_\s\times\B) + \R_\s,
\end{gather*}
over ions $\i$ and electrons $\e$ gives net momentum evolution
\begin{gather}
  \label{MHDmomentumEvolution}
  \partial_t (\mdens\u) + \Div(\mdens\u\u + \PTd + \PT)
  = \J\times\B,
\end{gather}
where we will call 
\begin{gather*}
  \glslink{PTd}{\ZPTd}:=\mdens_\i\w_\i\w_\i+\mdens_\e\w_\e\w_\e
\end{gather*}
the \defining{drift pressure} tensor;
here $\E$ has disappeared because of charge neutrality and
we have used that $\R_\i+\R_\e=0$ by conservation of momentum.

We can compute the diffusion pressure in terms of the current
using equation \eqref{wsJ},
\begin{align*}
    \w_\s &= \frac{\mred}{m_\s}\frac{\J}{n q_\s}.
\end{align*}
We get
\begin{gather}
  \label{PTdformula}
  \PTd = \mte\mti\J\J/\mdens
\end{gather}

In addition to assuming \gls{quasineutrality}, MHD models typically
neglect second-order terms in $\w_\s$ such as $\PTd$.
A physical justification for doing so is the assumption that
interspecies drift velocity is dominated by the fluid velocity
$\u$ and the thermal velocity $\c$ of particles.
We wish to drop such second-order terms in $\w_\s$
because they are numerically difficult: since $\w_\s$
is defined in terms of $\curl\B$, second-order terms give rise
to nonlinear higher-order differential operators.

\subsection{Energy density evolution}

In contrast to the situation for mass and momentum,
the neutrality assumption of MHD does not require that
one replace the evolution equations for energy density of
the individual species with a net energy density.

Models which evolve separate energy evolution equations
for two species are called \defining{two-fluid MHD} models.
Models which evolve a total energy equation are
\defining{one-fluid MHD models}.

In this document we are chiefly interested in two-fluid MHD.
If we evolve a single energy equation then we must use
a constitutive assumption to infer the energy of the individual
species in order to define the pressure tensors that appear
in Ohm's law \eqref{OhmsLaw}.  A single energy evolution equation
is most often used when the pressure terms are neglected in Ohm's law.
A primary goal of this dissertation is to determine minimal modeling
requirements to simulate fast magnetic reconnection with a fluid model
without invoking resistive drag force (which is found to be insufficient
for fast reconnection unless resistivity is defined anomalously).
In the absence of resistivity the pressure term is necessary to
support steady magnetic reconnection.
We seek a model which supports fast reconnection even for pair plasma.
But in the case of pair plasma, if we evolve a combined pressure tensor
and assume that pressure is equally distributed among both species
then the pressure term disappears from Ohm's law and steady reconnection
cannot be supported.  If we evolve separate evolution equations for the
pressure of each species, however, even in the case of
\first{symmetric pair plasma}
we will see that the contributions of the pressure tensors
of the two species to the reconnection electric field will
add instead of cancel.



\subsection{Incompressible MHD}

Incompressible MHD assumes that the number density
$n$ is conserved along particle paths.
Recall number density evolution \eqref{numberDensityEvolution},
which we here write in the form
\begin{gather*}
  d_t n = -n\Div\u.
\end{gather*}
So incompressibility means that $\Div\u=0$.

In the compressible case we need to evolve energy
in order to infer the scalar pressure in the momentum evolution
equation.  In the incompressible case one instead infers
the scalar pressure from the incompressibility equation
as a constraint.  The deviatoric pressure can be obtained
either from deviatoric strain (for the five-moment closure)
or by evolving a deviatoric pressure tensor.

Again, we may check that the evolution equation
\eqref{spcNumDensEvolutionCopy} for number
density of each species is exactly satisfied when the
assumption of incompressibility is imposed.
As before, we need that $\Div(n_\s \w_\s)=0$.
The proof goes through without change, since it
is simply based on the definition of $\w_\s$ in terms of $\J$.

Incompressibility is justified in the presence of a strong,
slowly varying magnetic field if quantities vary slowly
in the direction of the magnetic field (see equation (6.157)
in \cite{hazeltine92}. 

\subsection{Entropy evolution for two-fluid MHD}

While neglect of the drift pressure means that the
momentum equation of two-fluid MHD usually differs from
the momentum equation of two-fluid Maxwell,
two-fluid MHD, whether compressible or not, satisfies
the exact same equations as two-fluid Maxwell for
density evolution and pressure evolution. 
Since entropy evolution is derived based only on
density evolution and pressure evolution and
not on momentum evolution, the entropy evolution equations
hold unchanged for both incompressible and compressible
two-fluid MHD.

\section{Summary}

For a summary of the results of this chapter, we refer the reader back
to the systems of equations immediately preceeding this chapter.

%% file: chap3.tex
\chapter{Steady Magnetic Reconnection}
\label{MagneticReconnection}

\def\A{\underline{\underline{A}}}
\def\b{\bhat}
\def\u{\mathbf{u}}
\def\i{\mathrm{i}}
\section{Definition of rotationally symmetric 2D magnetic reconnection}

\subsection{Translational symmetry}

We define a problem to be \defining{two-dimensional (2D)}
if the problem is invariant under translation parallel
to an axis which we momentarily call the translational axis.
Choose a plane perpendicular to the translational axis.
We will call it \defining{the plane}.
The solution is fully represented by its value on the plane.
So we can solve the problem on the plane, that is,
on a two-dimensional spatial domain rather than on
three-dimensional space.

\subsection{Rotational symmetry}

Suppose that the problem is also invariant under 180-degree
rotation around an axis parallel to the translational axis,
which we will call the \defining{out-of-plane axis}.
Then we say that the problem is \defining{rotationally symmetric}.
This chapter is concerned with 2D rotationally symmetric
problems.

We refer to the intersection of the plane and the out-of-plane axis
as the \defining{origin}, which we designate as \gls{0};
we choose a cartesian coordinate system
whose origin is at this point and one of whose coordinate
axes coincides with the out-of-plane axis. 
We will take the out-of-plane axis to be the \defining{$z$ axis}.
The other two axes lie in the plane.

The alignment of the other two axes is based on the symmetries of
the initial conditions and of the imposed boundary conditions.
We conform to the general convention that
the $x$ axis should be the ``outflow'' axis.
The $y$ axis will then be the ``inflow'' axis.\footnote{
It is more common in the reconnection literature to
use \defining{geocentric coordinates}, in which 
$y$ is the out-of-plane axis and $z$ is the inflow
axis.  The convention of geocentric coordinates
is that the $x$ axis connects Earth and Sun and
$z$ is perpendicular to the ecliptic
\cite{brackbill:priv11} 
}
The inflow axis is so-named because fluid flow
approaches the origin along the inflow axis
and diverges from the origin along the outflow axis.

Rotational symmetry implies that the in-plane component
of any vector (including the magnetic field) must be zero
at the origin.  Therefore, the origin is
a 2D null point of the in-plane component of the magnetic field.
In the remainder of this paragraph ignore the out-of-plane
component of the magnetic field and
take the phrase \emph{magnetic field} to refer to the in-plane
component of the magnetic field $\B^\perp$.
By linearizing the magnetic field near the origin
we can classify the origin as an \emph{X-point} or an \emph{O-point}.
Let $\A:=\D\B|_0$ denote the derivative of the in-plane
component of the magnetic field at the origin.
Then $\Div\B=0$ says that $\tr\A=0$.
We classify the point based on the eigenvalues
and eigenvectors of $\A$.
\begin{itemize}
\item \emph{X-point case.}
If $\A$ has real eigenvalues then they are equal and opposite.
If the real eigenvalues are nonzero then the eigenvectors
define a pair of \defining{separatrices}.
The separatrices are magnetic field lines which intersect
the origin.  The \emph{in-coming separatrix} is tangent to
the eigenvector with negative eigenvalue, and the
\emph{out-going separatrix} is tangent to the eigenvector
with positive eigenvalue.
In this case we refer to the origin as an \defining{X-point}.
\item \emph{Antiparallel case.}
If the eigenvalues are both zero, then $\A$ is either zero or nilpotent.
If $\A$ is nilpotent then near the origin the magnetic field lines
are antiparallel and are aligned with the eigenvector of $\A$.
(If $\A$ is zero then linearization is insufficient to classify
the topology of the magnetic field near the origin.)
\item \emph{O-point case.}
If the eigenvalues are complex then their real part is zero and
their imaginary parts are opposite,
so magnetic field lines near the origin are approximately ellipses.
In this case we refer to the origin as an \defining{O-point}.
\end{itemize}
\emph{In this document we are chiefly concerned with
rotationally symmetric 2D problems where the
origin is an X-point.}


Symmetry about the origin means that the in-plane component of
the fluid velocity vectors $\u^\perp_\s$ is also zero at the origin.
We can use the same type of eigenstructure analysis
that we used to classify magnetic field structure at the origin
in order to classify fluid flow near the origin.
For steady solutions conservation of particles
(i.e.\ $d_t n_\s := -n_\s\Div\u_\s$,
where recall that $d_t:=\partial_t+\u_\s\dotp\D$)
implies at the origin (where $d_t=0$) that
$\Div\u^\perp_\s=0$, so the classification
of fluid flow lines is exactly like the classification
of magnetic field lines.
Regardless of whether incompressibility holds,
by transforming into a rotating frame of reference
we can assume that $\Div\u^\perp_\s=\Sym(\Div\u^\perp_\s)$;
in this case eigenvalues are real and eigenvectors are
orthogonal. In the incompressible case the sum of the eigenvalues
is zero. So we can define inflow and outflow separatrices for
each species fluid (for irrotational flow or in a frame of
reference rotating with the fluid).
(Note, however, that when transforming into a rotating
reference frame steady solutions become periodic solutions.)

\subsection{Consequences of symmetries for tensor components}

We also consider problems that are symmetric under
reflection across a plane containing the out-of-plane axis.
Reflection across the $y$-$z$ plane is effected by negation
of the $x$ coordinate, so we refer to this reflection as
a \emph{reflection in $x$}.

We say that a tensor is a \defining{proper tensor}
if it is invariant under reflections.
A tensor that is negated under reflections
is called a \defining{pseudo-tensor}.
Pseudo-tensors are negated under reflections
because reflections reverse the orientation of space.
Vectors which are pseudo-tensors are called pseudo-vectors and
vectors which are proper tensors are called proper vectors.
The magnetic field $\B$ is a pseudo-vector.
So is the curl of any proper vector (for example,
the vorticity).  The curl of a pseudo-vector
is a proper vector.

Under a reflection in $x$, a component of a proper tensor
with an odd number of $x$ indices is negated; other components remain unchanged.
(Therefore a component of a pseudo-tensor
with an odd number of $x$ indices remains unchanged
and other components are negated.)
Reflection in $x$ followed by reflection in $y$ effects
180-degree rotation in the $x$-$y$ plane, i.e.\
around the $z$ axis.  Therefore, under rotation about
the $z$ axis, a component of a tensor with an odd number of
non-$z$ indices is negated; other components remain unchanged.
Rotations do not reverse the orientation of space, so pseudo-tensors
and proper tensors transform in the same way under a rotation.

Invariance under rotation (or reflection) means that negated
components must be zero.

The assumption of 180-degree rotational symmetry implies 
at the origin that tensor components with an odd number of
out-of-plane indices must be zero.  For example, vectors
at the origin must be parallel to the out-of-plane axis.
In particular, the magnetic field at the origin must be
out-of-plane and is called the \defining{guiding magnetic
field} or \defining{guide field}.

\def\reflectionalSymmetry{reflectional symmetry }
\subsection{Reflectional symmetry}
\label{Reflectional symmetry}

We will also consider problems where there is
symmetry under reflection across the in-plane
axes (the $x$ axis and the $y$ axis).  This implies
180-degree rotational symmetry about the out-of-plane axis.
More generally, for reflection in $x$,
reflection in $y$, and 180-degree rotation around $z$,
any two of these transformations is sufficient to
generate the third.  We refer to symmetry under this
set of transformations as
\defining{\reflectionalSymmetry}.

\emph{Assume in the remainder of this subsection
\ref{Reflectional symmetry}
that \reflectionalSymmetry holds.}

Symmetry under reflection in the $y$ axis means that on the $x$ axis the
only nonzero component of the magnetic field is $B_y$; likewise,
on the $y$ axis the only nonzero component of the magnetic field
is $B_x$. So at the origin the magnetic field must be zero, i.e.,
there is no guide field.

The fluid velocity is a proper vector.
Symmetry under reflection in the $y$ axis means that
on the $x$ axis $u_y=0$, and
symmetry under reflection in the $x$ axis means that
on the $y$ axis $u_x=0$.
Therefore the fluid velocity must be irrotational
when linearized about the origin, i.e.\ $\D\u^\perp_\s|_0$
must be symmetric, and the in-plane standard basis vectors
must be eigenvectors.  Assuming that $\D\u^\perp_\s|_0\ne 0,$
one in-plane axis (by convention the $x$ axis) must be
an outflow axis and the other (the $y$ axis) must be the inflow
axis.

If the origin is an X-point for the magnetic field, then the
separatrices must of course be symmetric with respect to the
axes. But note that because $\B$ is a pseudovector, in the case
of \reflectionalSymmetry it is still possible for the
origin to be a magnetic O-point.


If there is symmetry under reflection across e.g.\ the $y$ axis,
then the $\B_y$ component of the magnetic field must
be zero on the $y$ axis.  If there is also symmetry
under reflection across the $x$ axis then
the $\B_z$ component of the magnetic field must be zero
on the $y$ axis.  Under 180-degree rotation the magnetic
field must be invariant, and so we can conclude that
on the $y$ axis the magnetic field satisfies
$\B(y\ebas_y) = \B_x(y) \ebas_x$ and $\B(-y\ebas_y) = -\B_x(y) \ebas_x$.
(Of course the same sort of statements can be made regarding the $x$ axis.)
In other words, in the case of reflectional symmetry the magnetic
field lines are antiparallel on the $y$ axis.
By smoothness and symmetry the magnetic field is parallel to
the $x$ axis near the $y$ axis.
It is thus common to refer to reflectionally symmetric
reconnection problems with X-point magnetic field geometry
at the origin as \defining{antiparallel reconnection}.
(Note, however, that in general antiparallel magnetic field
implies magnetic field uniformly aligned e.g.\ with the $x$ axis.)

\section{2D magnetic reconnection}
\label{2DmagneticReconnection}

Recall that the ideal MHD model of plasma assumes
the ideal Ohm's law, which says that the electric field is
$ 
  \E = \B\times\u,
$ 
that is, zero in the reference frame of the fluid.
Inserting this into Faraday's law
$\partial_t\B+\curl\E=0$ gives the evolution equation
\begin{gather*}
  \partial_t\B+\curl(\B\times\u)=0.
\end{gather*}
This equation says that $\u$ is a
\first{flux-transporting flow} for $\B$.
It implies that magnetic flux 
and magnetic field lines
are transported with the fluid and that
the topology of the magnetic field therefore cannot change.

More generally, ideal Hall MHD assumes that
the electric field is zero in the reference
frame of the \mention{charge velocity} $\uc$.
Recall the full Ohm's law,
\begin{gather*}
  \partial_t\B+\curl(\B\times\uc + \E')=0,
\end{gather*}
where $\E'$ is the nonideal electric field.
If $\E'$ is zero then $\uc$ is a flux transporting
flow for the magnetic field and the topology of magnetic 
field lines cannot change.  When $\E'\ne 0$ magnetic
field lines can change their topology or \mention{reconnect}.

This raises the question, ``How should one define the \mention{rate
of magnetic reconnection}?''  A full answer to this question
reconnection in three-dimensional space
would really entail a covariant (Lorentz-invariant) definition
in space-time, but for two-dimensional reconnection
we can give elementary definitions.  
In the case of 180-degree rotational symmetry
we will be able to define the rate of reconnection
to be the out-of-plane component of the electric field
at the origin.  A reasoned account of this definition follows.

Faraday's law of magnetic induction
$
   \partial_t \B + \curl\E = 0
$
implies that the rate of change of magnetic flux through any
line segment in the plane is the difference of the out-of-plane
electric field values at its endpoints.
For the \first{GEM magnetic reconnection challenge problem}
conducting wall boundaries exist, at which the out-of-plane
component of the electric field must be zero. 
Magnetic field lines cannot pass through conducting wall boundaries 
(see section \ref{GEMboundaryConditions}).
We therefore can define the reconnected flux to be the change in
magnetic flux through a line segment extending from the origin to
the conducting wall boundary.  Then by Faraday's law the out-of-plane component
of the electric field at the origin is the rate of magnetic
reconnection.

For 2D steady state, Faraday's law implies that the out-of-plane
component of the electric field must be constant (so the presence
of conducting wall boundaries would not allow steady 2D reconnection).
Thus, to given a definition of reconnection which incorporates
steady driven reconnection we seek a more general definition of
reconnection based on Ohm's law \eqref{OhmsLawSep}
\def\uc{{\u_\mathrm{c}}}
\begin{gather*}
  \E = \B\times\uc+\E'
\end{gather*}
and its out-of-plane component
$
  \E^{\glslink{para}{\para}} = \B^\perp\times\uc^\perp+\E'^\para,
$
where $\E'$ denotes non-ideal electric field,
and $\B^{\glslink{perp}{\perp}}$ and $\uc^\perp$ are the in-plane
components of the magnetic field and the charge velocity.
At the origin the ideal term disappears (because 
$\B$ and $\uc$ must be parallel).  Assume
that $\E'$ vanishes in the \emph{ideal region},
which includes the whole domain except for a region
containing the origin called the \emph{diffusion region}.
Outside the diffusion region, Faraday's law
$
  \partial_t\B + \curl\E = 0
$
and its projection into the plane
$
  \partial_t\B^\perp + \curl(\E^\para) = 0
$
imply that $\uc$
(or $\uc_\perp = \frac{\E\times\B}{\B\dotp\B}$,
where $\uc_{\glslink{perp}{\perp}} := \uc-\uc\dotp\b\b$, and $\b:=\B/B$)
is a \mention{flux-transporting flow} for $\B$ and that the in-plane
component $\uc^\perp$
(or $(\uc^\perp)_\perp = \frac{\E^\para\times\B^\perp}{\B^\perp\dotp\B^\perp}$)
is a flux-transporting flow for $\B^\perp$.
Therefore, \emph{the rate at which in-plane magnetic flux
is convected across a point in the ideal region is
$\|\E^\para\|$ (i.e.\ $\|(\uc^\perp)_\perp\|\cdot\|\B^\perp\|$)}.
In steady state $\E^\para$ is constant.  So in general we say that
\emph{the rate of reconnection is the magnitude of the
out-of-plane component of the electric field at the origin}.

\section{Ohm's law (current evolution) at the origin}


\subsection{MHD}

For MHD Ohm's law (equations \ref{OhmsLaw}, \ref{chargeVelocityDef})
is assumed,
\begin{align}
  \label{OhmsLawTerms}
  \E =& \Tresistivity\dotp\J  & \hbox{(resistive term)}
    \\ \nonumber &+ \B\times\uc  & \hbox{(ideal Hall term)}
    \\ \nonumber &+ \frac{\mut}{e n}\left[
          \Div\left(\frac{\PT_\i}{\mti}-\frac{\PT_\e}{\mte}\right)
        \right] & \hbox{(pressure term)}
    \\ \nonumber &+ \frac{\mut}{e n}\left[
          \partial_t\J + \Div\left(\u\J+\J\u-\frac{\dmt}{\mdens}\J\J\right)
        \right] & \hbox{(inertial term).}
\end{align}
Each of these terms represents a component of the electric field
and thus may also be referred to as e.g.\ the resistive electric field.
The ideal term of Hall MHD is usually decomposed as
\begin{align*}
   \B\times\uc = \underbrace{\B\times\u}_{\hbox{(ideal term)}}
      + \underbrace{\frac{\dmt}{\mdens}\J\times\B}_{\hbox{(Hall term)}}.
\end{align*}

Assuming symmetry under 180-degree rotation in the plane,
at the origin only the out-of-plane component of the electric
field $\E^\para:=\E\dotp\ebas^\para$ survives and the Hall
term $\B\times\uc$ disappears:
\begin{align*}
  \E^\para &= (\Tresistivity\dotp\J)^\para 
   + \underbrace{\frac{\mut}{e n}
       \left[\Div\left(\frac{\PT_\i}{\mti}-\frac{\PT_\e}{\mte}\right)\right]^\para
       }_{\hbox{pressure term}}
   + \underbrace{\frac{\mut}{e n}
         \left[
         \partial_t\J^\para + \Div\left(\u\J+\J\u-\frac{\dmt}{\mdens}\J\J\right)
         \right]^\para}_{\hbox{inertial term}}
        & \hbox{at $0$}.
\end{align*}
Since $\E^\para|_0$ is the rate of reconnection,
this says that at the origin reconnection must be supported by 
the resistive term, the pressure term, or the inertial term.

In steady state the inertial term disappears, since
$\Div\J=0$ and since in steady state at the origin $\Div\u=0$:
\begin{equation*}
  \E^\para = (\Tresistivity\dotp\J)^\para 
       + \frac{\mut}{e n}\left[
       \Div\left(\frac{\PT_\i}{\mti}-\frac{\PT_\e}{\mte}\right)
       \right]^\para \qquad \hbox{at $0$ for $\partial_t=0$}.
\end{equation*}
Therefore \emph{steady reconnection must be supported by the
resistive term or the pressure term} \cite{article:Vasyliunas75}.

\subsection{Two-fluid-Maxwell}

For a two-fluid-Maxwell model neutrality does not necessarily
hold and therefore the Ohm's law assumed by MHD (which assumes
neutrality) cannot be assumed to strictly hold. To analyze the
constraints on the reconnecting electric field we therefore
revert to the one-species current evolution equations from which
Ohm's law is derived. Equivalently, we consider one-species
momentum evolution.

Recall the momentum evolution equation \eqref{momentumEvolution}
\begin{gather*}
  \mdens_\s d_t \u_\s + \Div \PT_\s
  = q_\s n_\s(\E+\u_\s\times\B) + \R_\s,
\end{gather*}
Solving for the electric field gives a single-species proxy ``Ohm's law'',
\begin{align}
  \label{proxyOhm}
  \E =&
      - (q_\s n_\s)\inv\R_\s & \hbox{(resistive term)}
  \\ \nonumber & + \B\times\u_\s & \hbox{(ideal term)}
  \\ \nonumber & + (q_\s n_\s)\inv\Div\PT_\s & \hbox{(pressure term)}
  \\ \nonumber & + (m_\s/q_\s) d_t \u_\s & \hbox{(inertial term)}
\end{align}
Assuming symmetry under 180-degree rotation in the plane,
at the origin only the out-of-plane component of the electric
field $\E^\para$ survives:
\begin{equation}
  \label{proxyOhmXpoint}
  \E^\para = - (q_\s n_\s)\inv\R_\s^\para
      + (q_\s n_\s)\inv(\Div\PT_\s)^\para
      + (m_\s/q_\s) d_t \u_\s^\para \qquad \hbox{at $0$}.
\end{equation}
In steady state the inertial term disappears:
\begin{equation*}
  \E^\para = - (q_\s n_\s)\inv\R_\s^\para
      + (q_\s n_\s)\inv(\Div\PT_\s)^\para
       \qquad \hbox{at $0$ for $\partial_t=0$}.
\end{equation*}

\subsection{Implications of Ohm's law at the origin}

In the following discussion ``Ohm's law'' may be taken as
the MHD Ohm's law \eqref{OhmsLawTerms} or the
two-fluid-Maxwell quasi-Ohm's law \eqref{proxyOhm}.
Assume 2D 180-degree rotational symmetry.

Since the out-of-plane component of the electric field
at the origin is the rate of reconnection,
Ohm's law at the origin implies a set of constraints on magnetic 
reconnection.

Note that the resistive term represents a frictional drag force
and results in heating (and thus entropy production) in an
energy-conserving model. In a model that conserves energy and
respects entropy but lacks diffusive entropy flux this means that
steady state is possible only if the resistive electric field is
zero at the origin.

If the electric field is zero at the origin then the pressure term
must be nonzero to support steady reconnection \cite{article:Vasyliunas75}.
That is, the divergence of the pressure must be nonzero.

Note that if $z$ is the out-of-plane axis then for $\PT=\PT_\s$
\begin{equation*}
  (\Div\PT)^\para = \partial_x \PT_{xz} + \partial_y \PT_{yz} \qquad \hbox{at 0}.
\end{equation*}

Note that in deriving Ohm's law one of our implicit regularity assumptions
was that the fluid density is nonzero.

\subsubsection{Agyrotropy is necessary for $\Div\PT|_0\ne 0$.}
\label{agyrotropy}

We say that a pressure tensor is \defining[isotropic pressure]{isotropic}
if $\PT=p\id$.  If the pressure is isotropic near
the origin (that is, in a neighborhood containing the origin)
then $\Div\PT_\s = \Div(p_\s\id) = \D p_\s$, which by
symmetry is zero at the origin.

We say that a pressure tensor is \defining[gyrotropic pressure]{gyrotropic}
if it is invariant under rotation around the direction vector
$\b:=\B/|\B|$ aligned with the magnetic field.  Then we can write
\begin{align}
  \PT &=p_\para\b\b + p_\perp(\id-\b\b)
    \nonumber \\&=p_\perp\id + (p_\para-p_\perp)\b\b,
   \label{gyrotropy}
\end{align}
where $\p_\para$ is called the \defining{parallel pressure}
and $\p_\perp$ is called the \defining{perpendicular pressure}.
At a point where the magnetic field is nonzero,
\begin{align*}
  \Div\PT =\D p_\perp + \B\cdot\D
     \underbrace{\left[(p_\para-p_\perp)\b/|\B|\right]}_{\hbox{if differentiable}},
\end{align*}
where we have used that $\Div\B=0$.
But $\B\cdot\D=0$ at 0,
and since $p_\perp$ is a scalar, $\D p_\perp=0$ at 0.
So $\Div\PT=0$ at the origin in the case of a nonzero guide field.
(This argument was given in \cite{article:HeKuBi04}.)

We can extend the notion of gyrotropy to include points where 
the magnetic field is zero
if isotropy holds when the magnetic field is zero (i.e.\
$p_\para=p_\perp$ when $\B=0$) and equation \eqref{gyrotropy}
defines a smooth function.  This is perhaps not enough to imply
that the expression marked \emph{if differentiable} is differentiable.
We could simply impose this assumption as an additional condition
for a gyrotropic function.

\def\A{\tensorb{A}}
\def\nhat{\mathbf{n}}
\def\Pxx{\PT_{xx}}
\def\Pyy{\PT_{yy}}
\def\Pxz{\PT_{xz}}
\def\Pyz{\PT_{yz}}
\def\Pzz{\PT_{zz}}
But I claim that gyrotropy implies $\Div\PT|_0=0$
even without imposing this extra assumption.
Let $\A=\D^\perp\B^\perp|_0$, that is, the two-dimensional
matrix which is the gradient of the projection of the
magnetic field onto the plane.  Note that $\tr\A=0$.
There are two orthogonal directions $\nhat$
in the plane such that $\nhat\dotp\A\dotp\nhat=0$.
Call the axes aligned with these directions
$x$ and $y$.  Since $\Div\PT|_0$ is aligned with
the out-of-plane axis, it is invariant under reflections,
so by averaging the data over
reflections across $x$ and $y$ we can assume
without loss of generality that symmetry holds
under reflections across $x$ and $y$.
\note{There is a hole in the proof at this point
because the condition \eqref{gyrotropy}
that defines gyrotropy is nonlinear and therefore not invariant
under averaging of data (otherwise I could choose any orthogonal axes). 
If $\D\B|_0 \ne 0$ probably I can patch
this up by using an estimate on the violation of symmetry
and an estimate of violation of agyrotropy, but I don't think
that I care enough about this proof to do that; instead, I think
I'll just limit my results to the symmetric case.}

Assume that symmetry holds
under reflections across $x$ and $y$.

On the $x$ axis $\B$ is parallel to $y$ or zero
(zero can occur in an arbitrarily small neighborhood
of $0$ only in case $\D\B|_0$ is zero or nilpotent),
and likewise on the $y$ axis
$\B$ is parallel to $x$.
Thus, gyrotropy implies that
along the $x$-axis $\Pxx=\Pzz$ and $\Pxz=0$ and
along the $y$-axis $\Pyy=\Pzz$ and $\Pyz=0$
(using that $\PT$ is isotropic at any points where $\B=0$).
Recall that at the origin only the out-of-plane component
$\Div\PT$ survives and thus
\begin{align*}
  \Div\PT = \partial_x P_{xz} + \partial_y P_{yz} \qquad\hbox{at 0}.
\end{align*}
Therefore $\Div\PT=0$ at 0.

So we can say in general that if the pressure is gyrotropic near the origin
then $\Div\PT_\s=0$.
That is, \emph{pressure cannot support reconnection without
agyrotropy in the vicinity of the origin}.

\section{Steady rotationally symmetric 2D reconnection requires heat flux.}
\label{steadyReconNeedsHeatflux}

The purpose of this section is to justify the following assertion,
which is one of the main results of this dissertation:
\begin{quote}
  In a 2D problem invariant under 180-degree rotation
  nonsingular steady reconnection is impossible
  in a plasma model which conserves mass, momentum, and energy 
  if the model implies that entropy production is zero
  in the vicinity of the origin
  or if the model implies that diffusive entropy flux (or heat
  flux) is zero in the vicinity of the origin.
\end{quote}
We attempt to make this statement precise and justify it as fully as possible.
The idea is to show that friction is necessary.
Friction produces heat, so if energy is conserved then
the heat must have a way to diffuse away from
the stagnation point at the origin.

In any gas model physical solutions must satisfy a set of
\defining{positivity conditions}. In a Boltzmann gas model the
particle density $f_\s$ must be positive. In a five-moment gas
model the density $\mdens_\s$ and pressure $p_\s$ must be positive.
In a ten-moment gas model the density $\mdens_\s$ must be positive
and the pressure tensor $\PT_\s$ must be positive definite.

In general, we say that a set of moments \mention{satisfies
positivity} or is \mention{physically realizable} if
there exists a distribution function with the specified
moments.  For a given collection of moments the set of
physically realizable values of the moments is a convex set.

We define a solution to be singular if it at some point the
state is not in the interior of the set of states that satisfy
positivity. In particular, a distribution is singular if it is
zero anywhere, a five-moment solution is singular if the density
or pressure is not strictly positive anywhere, and a ten-moment
solution is singular if the pressure tensor fails to be strictly
positive definite anywhere or if the density fails to be strictly
postitive. In the argument we will assume that solutions are
smooth. By convolution with a smooth approximate identity we can
make this assumption without loss of generality.

We first lay out the high-level argument.

\subsection{Outline of the physical argument}
\label{overviewOfArgument}

Before analyzing models we study the physics itself, taking the
Boltzmann equation as our standard of truth. We consider the the
possibilities at the origin for the physical solution based on
the constraints implied by the entropy evolution equation and the
momentum evolution equation and by the problem symmetries. We
then consider implications of the analysis for models in light of
the relationships assumed in our reasoning.

The outline of the argument is as follows. The argument focuses
attention on a single species and analyzes its entropy evolution.

If the collision operator $C$ is zero then the physics is
governed by the Vlasov equation. We assume a steady state
solution of the Vlasov equation at the origin and conclude that
it must be singular.

If the collision operator is nonzero then we consider entropy
evolution in the vicinity of the origin. Either there is entropy
production at the origin or there is not. If so, then entropy
and heat are produced at the origin and so there must also be
diffusive entropy flux and heat flux to balance the production
by dispersing the entropy and heat. (Note that entropy flux
approximately equal heat flux divided by temperature for
near-Maxwellian distributions.)

If there is no entropy production at the origin, then since the
collision operator is nonzero, the Boltzmann H theorem implies
that the solution at the origin is Maxwellian. In this case
5-moment gas dynamics applies near the origin. Assuming 5-moment
gas dynamics, if the deviatoric strain rate is nonzero at the
origin then we show that there is entropy production, contrary
to hypothesis. \emph{As we will show, if the deviatoric strain
rate is zero and if the heat flux is zero in the vicinity of the
origin then so is the divergence of the pressure tensor.} But
in steady state at the origin the only terms in the momentum
equation that can support a reconnection electric field are the
divergence of the pressure tensor and the resistive drag force.
If there is resistive drag force then there is entropy production
at the origin, contrary to hypothesis. So we may conclude that if
there is steady reconnection then there must be heat flux in the
vicinity of the origin even if there is no entropy production at
the origin.

The previous paragraph begins with two assertions
that the reader might question.

The first questionable assertion is that the only type of
distribution for which there is no entropy production is a Maxwellian.
The Boltzmann H theorem asserts that entropy production is
nonnegative and is zero precisely when the velocity distribution
is Maxwellian. This holds for a large class of collision
operators and is generally taken as a requirement for nonzero
collision operators. Nevertheless, one may conceive of a
collision operator for which states other than Maxwellian are
equilibria. In particular, the hyperbolic ten-moment model with
no heat flux or isotropization can be regarded as assuming an
idiosyncratic collision operator which instantaneously relaxes
to a Gaussian distribution; for this model Gaussian distributions
are equilibria.

The second questionable assertion is that if the solution at
the origin is Maxwellian then 5-moment gas dynamics applies
near the origin. The justification for this assertion is as
follows. First, five-moment gas dynamics holds rigorously in
the asymptotic limit as the solution approaches a Maxwellian
distribution. Second, since the velocity at the origin is zero,
the fluid near the origin stays near the origin for a long time
and has time to equilibrate with conditions at the origin; in
particular, flow along an inflow separatrix takes forever to
approach the origin. These statements hold to even higher order
in case the deviatoric strain rate is zero at the origin.

\subsection{Consequences of the argument for models}

We now identify the modeling consequences of the argument in
section \ref{overviewOfArgument} by carefully identifying the
assumptions of the argument and considering in what models these
assumptions hold.

The essential assumptions of the argument relevant to
fluid models are that the entropy evolution equation
\eqref{s5evolution} holds for each species and that
(1) the momentum evolution equation holds for each species
(as in the two-fluid-Maxwell case)
or more generally, that (2) Ohm's law holds
(as needed in the anlysis of two-fluid MHD).

The following models cannot support nonsingular steady magnetic
reconnection for 2D problems with rotational symmetry:
\begin{itemize}
  \item Models which evolve a form of the Vlasov equation.
  \item Adiabatic two-fluid models:
    \begin{itemize}
      \item whether two-fluid-Maxwell or two-fluid MHD,
      \item whether 5-moment (viscid or inviscid)
        or 10-moment (isotropizing or not),
      \item whether compressible or incompressible
        (although in the incompressible case one typically
        does not evolve energy and therefore steady
        reconnection is possible).
    \end{itemize}
  \item Ideal MHD or ideal Hall MHD.
\end{itemize}
I remark that Chac\'on et al.\ in \cite{article:ChSiLuZo08}
simulate steady fast reconnection in magnetized pair plasma
using an incompressible five-moment two-fluid model
without any heat flux in their equations.
This does not contradict the results of this section.
They do not have heat flux in their equations because
they do not solve an energy evolution equation.
They do not need explicit energy evolution because their
model is incompressible.  The scalar pressure is inferred
from the incompressibility constraint, and the deviatoric
pressure is the viscosity times the deviatoric strain rate.
So one may pretend that an energy evolution
equation (e.g.\ with nonzero heat flux) is being evolved,
but the rest of the equations do not depend on it.

For the incompressible two-fluid \emph{ten-}moment model
assumed in \cite{article:Br11}, however, the pressure
tensor \emph{is} evolved adiabatically in accordance
with \eqref{PTevolution},
and we can conclude that steady reconnection should not be possible
for the antiparallel case at least.

\subsection{Vlasov model}

We first consider the case that the collision operator is zero.
We argue that kinetic models require entropy production
for nonsingular rotationally symmetric steady reconnection.
Heuristically this holds because without collisions
particles sitting at the origin will be accelerated
without bound by the electric field.

Recall the Boltzmann equation,
\def\a{\mathbf{a}}
\def\v{\mathbf{v}}
\def\s{\mathrm{s}}
\begin{gather*}
  \partial_t f_\s + \v\dotp\D_\xb f_\s + \a_\s\dotp\D_\v f_\s = C_\s,
\end{gather*}
where $f=f_\s(t,\xb,\v)$ is particle density of species $\s$,
$\v$ is particle velocity,
$\a_\s=(q_\s/m\s)(\E+\v\times\B)$ is particle acceleration,
and $C_\s$ is the collision operator.
Drop the species subscript $\s$.
We claim that without the collision operator
steady reconnection is singular.

Assume rotational symmetry about the $z$-axis.
Then at the origin the Boltzmann equation simplifies to
$
  \partial_t f + (q/m)E_z\partial_{v_z} f = C.
$
Assume steady state ($\partial_t = 0$) and no collisions ($C=0$).
If $E_z = 0$ then there is no magnetic reconnection.
Otherwise $\partial_{v_z} f = 0$,
which says that $f$ is independent of $v_z$ when $\v$ is parallel to
the $z$-axis. But nonsingular distributions should go to $0$ as $z$
goes to infinity. So $f$ is $0$ when $\v$ is parallel to the
$z$-axis, which is a type of singularity.

\def\scoef{\widetilde\CCCC}
\subsection{Quadratic moment models need heat flux}

To complete the argument of section \ref{overviewOfArgument}
we first prove the following theorem.
\def\dS{\deviator{\tensorb{p}}}
\def\dSi{\deviator{\tensorb{\widetilde p}}}
\def\A{\underline{\underline{A}}}
\begin{theorem}
\label{steadyReconNeedsHeat}
  Assume a rotationally symmetric 2D problem
  in which a fluid species satisfies a steady-state
  entropy evolution equation of the form
  \begin{gather}\label{gyroEntropyEvol}
    \u\dotp\D s = \dSi\ddotp\scoef\ddotp\dS - \A\ddotp\Div\qT + \error
  \end{gather}
  in the vicinity of the origin, where
  $\dS$ and $\dSi$ are deviatoric tensors and where
  the coefficients $\scoef$ satisfy the positive-definiteness criterion
  $\dSi\ddotp\scoef\ddotp\dS>0$ for all $\dS\ne 0$.

  Suppose that $\Div\qT$, $\D\Div\qT$, $\laplacian\Div\qT$,
  $\error$, $\D\error$, $\D\D\error$,
  and $\Sym(\D\u)$ are all zero at the origin.
  
  Then $\D\dS=0$ at the origin.
\end{theorem}
\textbf{\em Proof of theorem.}
  Note that $\glslink{laplacian}{\Zlaplacian}$ denotes the laplacian.
  At 0 the left hand side is zero.  Since $\Div\qT=0$ at 0, $\dS=0$ at 0.

  For convenience use $'$ to denote $\partial_x$.
  We apply $\partial_x^2$ to the entropy evolution equation and evaluate
  at the origin.  Using that $\dS=0$ at 0,
  \begin{align}\label{d2sEvol}
    (\u\dotp\D s)'' &= {\dSi}'\ddotp\scoef\ddotp{\dS}' - (\A\ddotp\Div\qT)''
    & \hbox{ at 0.}
  \end{align}
  For the remainder of this proof all expressions
  and statements are evaluated at the origin.
  I claim that the left hand side is zero.
  The left hand side expands to
  $\u''\dotp\D s + \u'\dotp\D s' + \u\dotp\D s''$.
  The outer expressions are easily seen to be zero.
  For the first, $\u''\dotp\D s=0$ because $\D s$ is zero.
  For the last, $\u\dotp\D s''=0$ because $\u\dotp\D$ is zero.
  It remains to show that $\u'\dotp\D s'=0$, that is,
  $\partial_x\u\dotp\D \partial_x s=0$, that is,
  $\partial_x\u_x \partial_x^2 s
  +\partial_x\u_y \partial_y \partial_x s = 0$.
  But 
  $\partial_x\u_x=0$ because $\Sym(\D\u)=0$, and
  $\partial_y \partial_x s = 0$ because by rotational
  symmetry $s$ is even in $x$ and even in $y$.

  For the right hand side, since $\qT=0$ and $\D\qT=0$,
  it follows that $(\A\ddotp\Div\qT)'' = \A\ddotp\Div\qT''$.
  So we have shown that ${\dSi}'\ddotp\scoef\ddotp\dS{}' = \A\ddotp\Div\q''$.
  That is,
  \begin{gather*}
    \partial_x\dSi\ddotp\scoef\ddotp\partial_x\dS = \A\ddotp\Div\partial_x^2\q
      \ \hbox{ and}
 \\ \partial_y\dSi\ddotp\scoef\ddotp\partial_y\dS = \A\ddotp\Div\partial_y^2\q,
  \end{gather*}
  where we have used that the choice of axes is arbitrary.
  By the positive-definiteness assumption the left hand sides
  (and thus the right hand sides) must be nonnegative.
  But the hypothesis $\laplacian\Div\qT=0$ says that
  $\Div\partial_x^2\qT = -\Div\partial_y^2\qT$.
  So the right hand sides are zero.
  So by the strict positive-definiteness assumption on $\scoef$
  we conclude that
  $\partial_x\dS=0=\partial_y\dS$, that is,
  $\D\dS=0$.
  \QED

\begin{corollary}
  \label{fiveMomentNeedsHeat}
  Assume a rotationally symmetric 2D problem
  in which a fluid species satisfies the steady-state
  five-moment entropy evolution equation \eqref{fiveMomentEntropyEvolution}
  \begin{gather*}\label{gyroEntropyEvol}
    \u\dotp\D s_\Maxwell = -\dstrain\ddotp\Pshape - p\inv\Div\q + \error
  \end{gather*}
  with stress closure \eqref{gyrotropicTempDeviator}
  \begin{gather}
     \label{gyrotropicTempDeviatorCopy}
     -\dPshape=2\ptime\Mcoefinv\ddotp\dstrain
  \end{gather}
  in the vicinity of the origin, where
  the (gyrotropic) coefficients $\scoef$ satisfy the positive-definiteness
  (entropy-respecting) criterion
  $\dstrain\ddotp\Mcoefinv\ddotp\dstrain>0$ for all $\dstrain\ne 0$.

  Suppose that $\Div\q$, $\laplacian\Div\q$,
  $\error$, $\D\error$, and $\D\D\error$
  are all zero at the origin.

  Then $\Div\PT=0$ at the origin.
\end{corollary}
\textbf{\em Proof.}
  To map onto the theorem, take
  $\dSi:=\dS:=\dstrain$ and take
  $\scoef:=2\ptime\Mcoefinv$.
  Take $\qT:=(2/5)\SymC(\q\id)$.
  Take $\A:=(2p)\inv\idtens.$

  $\Div\q$ is a scalar and therefore $\D\Div\q=0$ at 0.
  At the origin, $\u\dotp\D s_\Maxwell = 0$, so if $\Div\q=0$ at 0
  then $\dstrain=0$ at 0.  But steady state conservation of mass 
  says that $\Div\u=0$ at 0.
  So $\Sym(\D\u)=0$ at 0, which is the remaining needed hypothesis of
  the theorem.

  The conclusion of the theorem says that $\D\dstrain=0$ at 0.
  I claim that this implies that $\D\PT=0$.
  This follows from the stress closure \eqref{gyrotropicTempDeviatorCopy},
  i.e.\ $-\dPshape=2\Tviscosity\ddotp\dstrain,$
  since at 0 $\dstrain=0$ and $\D\dstrain=0$.
  Taking the trace, $\Div\PT=0$.
  \QED
\begin{corollary} 
  \label{tenMomentNeedsHeatFlux}
  Assume a rotationally symmetric 2D problem
  in which a fluid species satisfies the steady state
  ten-moment temperature tensor evolution equation
  (see \eqref{PTevolution})
  \begin{gather}
     \label{PTevolutionCopy}
     n \partial_t \TT 
       + \SymB(\PT\dotp\D\u) + \Div\qT
       = ({q/m})\SymB(\PT\times\B) + \RT
  \end{gather}
  and the continuity equation
  \begin{gather*}
    \partial_t n + \Div(n\u) = 0
  \end{gather*}
  with isotropization closure \eqref{RTclosure}
  \begin{gather}
    \RT/n = -\ptime\inv \Gcoef\ddotp\dTT,
    \label{RTclosureCopy}
  \end{gather}
  in the vicinity of the origin, where
  the (gyrotropic) coefficients $\Gcoef$ satisfy the positive-definiteness
  (entropy-respecting) criterion
  \eqref{entropyRespectingGyrotropicRelaxationClosureRequirement}
  \begin{align*}
    -\dTTinv\ddotp\Gcoef\ddotp\dTT > 0.
  \end{align*}

  Suppose that $n\ne 0$ at the origin.
  Suppose that $\Div\qT$, $\D\Div\qT$, and $\laplacian\Div\qT$
  are all zero at the origin.

  Then $\Div\PT=0$ at the origin.
\end{corollary}
\textbf{\em Proof.}

  Recall that assuming the continuity equation and the
  pressure tensor evolution equation we derived 
  equation \eqref{sGaussEvolutionSimple}
  for ten-moment entropy, which by \eqref{RTclosureCopy}
  is
  \begin{gather}
    2 n \u\dotp\D s_\Gauss = -\TTinv\ddotp(\ptime\inv\Gcoef)\ddotp\dTT - \PT\inv\ddotp\Div\qT.
    \label{entropyEvolution}
  \end{gather}

First assume that $\Sym(\D\u)=0$ at the origin.
To map onto the hypotheses of the theorem
take $\dS = \dTT$,
$\dSi=-\dTTinv$,
$\scoef=(2\ptime)\inv\Gcoef$,
$\qT=\qT$, and
$\A=\PT\inv/2$.
The conclusion of the theorem says that $\D\dTT=0$ at 0.
But $\dPT$ is a scalar multiple of $\dTT$.
So $\D\dPT=0$ at 0.  So $\Div\PT=0$ at 0.

Now suppose that $\Sym(\D\u)\ne0$ at the origin.
That is, $\dstrain\ne0$ at 0
(since steady state says $\Div\u=0$ at 0).
Then we will show using the hypotheses of
the theorem that entropy production is nonzero
at the origin (contradicting the assumption
of steady state).  By entropy evolution 
it will be enough to show that $\RT\ne 0$ at 0,
that is (by the closure hypotheses \eqref{RTclosureCopy}),
$\dTT\ne 0$ at 0.


Henceforth in this proof all statements are evaluated at zero.
In steady state $d_t=0$.
Suppose that $\dTT=0$.  That is, $\dPshape=0$.
$\Div\qT=0$ implies $\Div\q=0$ and $\Div\dqT=0$.
So at the origin equation \eqref{PTevolutionCopy} 
reduces to
\begin{gather}
 2p\dstrain = \RT.
\end{gather}
But $\dstrain\ne 0$.
So $\RT\ne 0$, contrary to hypothesis.
So we may generally conclude that $\Div\PT=0$ at the origin.
\QED

\subsection{Steady reconnection requires heat flux.}

Section \ref{overviewOfArgument} promised that we
would show that if the deviatoric strain rate is zero and if
the heat flux is zero in the vicinity of the origin then so
is the divergence of the pressure tensor.

To complete the verification of this claim for the five-moment
model it remains to consider the hyperbolic case,
where the closure $\dPT=0$ is used;
this closure is \emph{non}strictly positive-definite,
so the theorems in the previous section do not apply.
But in this case isotropy holds and $\Div\PT=\D p$,
which must be zero at the origin.
So steady reconnection is generally not possible in the five-moment
model without heat flux.

\note{I worry that the following argument is not completely valid.
For a Coulomb collision operator deviations from a Maxwellian
consisting e.g.\ of fast-moving electrons can decay arbitrarily
slowly, so the assumption that 5-moment gas-dynamics holds
in the Maxwellian limit seems dubious.
Maybe I have to limit my conclusions to collision operators for
which there is a uniform estimate on the rate of decay of
perturbations.  Are there other subtle assumptions?
At least this should go through for collision operators
which share the properties of the Gaussian-BGK collision operator
that I perhaps implicitly assume.
}

We now argue that the same claim holds for the Boltzmann equation
with a collision operator which is strictly postitive for non-Maxwellian
distributions and for which 5-moment gas dynamics with
linearized deviatoric stress closure applies in the Maxwellian limit.
This should hold for physically reasonable collision operators.

\def\remainder{{\mathrm{remainder}}}
\def\fremainder{f_{\mathrm{remainder}}}
In steady state the entropy production at the origin should be zero
and therefore the velocity distribution at the origin should be
a Maxwellian $f\sup{0}$.  Let $\epsilon$ be a smooth even function
equal to 0 at the origin such that $|f-f\sup{0}|\le \epsilon$.
Then $\epsilon = \O(|\xb|^2).$
Near the origin we can define 
a Chapman-Enskog expansion
$
  f = f\sup{0}
    + \epsilon \fd\sup{1}
    + \fremainder.
$
where $\fremainder = \O(\epsilon^2)$.
\note{Can I fully justify this assumption?}
Note that in a normal Chapman-Enskog expansion 
$\epsilon$ is assumed independent of space;
so for each position $\xb$ we perform a Chapman-Enskog
expansion using $\epsilon(\xb)$ and then evaluate it at $\xb$.
In the Gaussian-BGK model we used a Chapman-Enskog expansion to
show that (see equation \eqref{implicitDeviatoricStressClosure})
\begin{gather}
  \label{PTclosure}
  \PT = p\id
    + 2\mu\Gcoef\ddotp\dstrain
    + \PT_{\error}
\end{gather}
where $\PT_{\error} = \O(\epsilon^2)$ and is smooth assuming
the other quantities in this equation are smooth.
In section \ref{GyrotropicViscousStressClosure}
we obtained this closure by using
(1) that if the solution is within $\O(\epsilon)$ of a Maxwellian
then within $\O(\epsilon^2)$ five-moment entropy is
nondecreasing and (2) that five-moment entropy must be approximately
nondecreasing for any small $\dstrain$. 
For the five-moment closure
$\mu$ is determined by temperature and
the coefficients $\Gcoef$ are functions of the five-moment state,
but for a Boltzmann model the coefficients $\mu\Gcoef$ would
also depend on the \emph{shape} of the deviation of the distribution
from a Maxwellian.

Since $\epsilon$ is smooth and even at $\xb$, we can write
$\epsilon=\O(|\xb|^2).$  Therefore, near the origin
the entropy evolution equation \eqref{gyroEntropyEvol}
holds to within $\O(\epsilon^2) = \O(|\xb|^4)$.
This is true even when the stress closure
\eqref{gyrotropicTempDeviatorCopy} is used,
since the error in the deviatoric pressure closure
is $\PT_{\error} = \O(\epsilon^2) = \O(|\xb|^4)$.
So the assumptions of corollary \ref{fiveMomentNeedsHeat}
are satisfied.

We conclude that for the Boltzmann equation with a ``reasonable''
collision operator steady reconnection requires that $\Div\q\ne 0$
at the origin.  The same conclusion holds for fluid models which
evolve a temperature evolution equation that agrees with
the five-moment model sufficiently well near a Maxwellian distribution.


\subsection{The ten-moment hyperbolic model does not support steady reconnection
  (at least for \reflectionalSymmetry)}

We have shown that without heat flux the ten-moment model with
nonzero isotropization cannot support steady reconnection.
For the hyperbolic case the closure $\RT_\s=0$ is used,
which is \emph{non}strictly positive-definite, so 
Corollary \ref{tenMomentNeedsHeatFlux} does not directly apply.

We argue that in the case of \reflectionalSymmetry
the hyperbolic ten-moment model does not support regular steady
reconnection if $\partial_z^n (\u_s)_z \ne 0$ for some $n$.
By convolution with an approximate identity
we can assume smoothness at the origin.



\def\Txx{\TT_{xx}}
Consider the case of \reflectionalSymmetry.
Assume that the $y$ axis is the inflow axis.
On the $y$ axis $\B$ is aligned with the $x$ axis.
Along the $y$ axis the temperature
evolution equation \eqref{TTevolution} for the $\Txx$
component simplifies to
\def\uy{u_y}
\def\ux{u_x}
\def\uxx{u_{x,x}}
\def\uyy{u_{y,y}}
\def\Dy{\partial_y}
\def\Dx{\partial_x}
\def\Cy{C_y}
\def\Cx{C_x}
\begin{gather}
  \uy\Dy \Txx = -2 \Txx \Dx\ux.
  \label{TxxEvolution}
\end{gather}
Look at the continuity equation on the $y$ axis,
\begin{gather*}
  \uxx+\uyy = \uy\Dy\ln n.
\end{gather*}
Suppose that $\uy = \Omega(y^{n+1})$.
That is, $\uyy = \Omega(y^n)$.
Then $\uyy+\uxx = \Div\u = \uy\Dy\ln n = \Omega(y^{n+1})\O(y^2) = \O(y^{n+3}) = \O(y^{n+1})$.
So $\uyy$ and $\uxx$ have leading coefficients
of the same magnitude but opposite
sign of order $n$ in $y$.
So we can write
\begin{gather}
  \label{uxxExpansion}
  \uxx= C y^n + \O(y^{n+1})
\end{gather}
and
\begin{gather}
  \label{uyyExpansion}
  \uyy = -C y^n + \O(y^{n+1}).
\end{gather}
Note that $\Txx$ is even.
Assume that it is also smooth.
Evaluate the $n$th derivative of both sides
of \eqref{TxxEvolution} at the origin
and get $0 = -2\Txx C n!.$  This contradicts
the assumption that $\Txx$ is strictly positive.
So no such smooth steady solution exists.

We can gain insight into why no smooth steady solution
exists if we substitute equations
\eqref{uxxExpansion} and
\eqref{uyyExpansion} into
\eqref{TxxEvolution}.
We get
\begin{gather*}
  y \Dy \Txx = 2(n+1) \Txx + \O(y).
  \label{TxxEvolution}
\end{gather*}
Near the origin we can ignore the $\O(y)$ term;
separating and integrating gives
\begin{gather}
  \label{TxxOrigin}
  \Txx = C_1 y^{2(n+1)}.
\end{gather}
Since $n$ must be nonnegative, this says that
$\Txx$ must be zero at the X-point.
This is physically what we would expect,
since there is expansion in the $x$ axis
as the flow approaches the origin along
the $y$ axis.

To treat the general case of rotational symmetry
we could transform into a rotating reference frame.
In this case we would have an oscillatory solution
rather than a steady state solution.
I am not sufficiently interested to pursue the
hyperbolic case further.
%

\subsection{Singular solutions}

Hitherto we have shown that adiabatic models do not
support nonsingular steady reconnection.
One may also consider singular solutions.

Jerry Brackbill has studied symmetric 2D reconnection using
an incompressible adiabatic ten-moment model.
He uses \mention{geocentric coordinates}, so he
refers to the inflow axis as the $z$ axis and the out-of-plane
axis as the $y$ axis.
By neglecting partial derivatives of the pressure term
with respect to $x$ in the momentum evolution equations
the evolution equations reduce to an ODE along the inflow axis.
Specifically he neglects his $\partial_x \PT_{xz}$
(which is dominated by his $\partial_z \PT_{zz}$)
and his $\partial_x \PT_{xy}$
(which is dominated by his $\partial_z \PT_{zy}$)
for both species.  Based on this simplification he
has implemented a steady-state solver \cite{article:Br11}
which shows excellent agreement with kinetic simulations
when the isotropization rate is appropriately tuned.
The implication of this dissertation, however, is that a
steady state solution does not actually exist for the
equations he solves.  He solves the temperature evolution
equation without making any approximating assumption,
and his assumption of incompressible flow implies that
the continuity equation is also exactly satisfied.
His model therefore satisfies the hypotheses of
Corollary \ref{tenMomentNeedsHeatFlux} and cannot
support a solution that is smooth, nonsingular, and
steady at the origin.  In particular, his pressure tensors
are not isotropic at the origin (see Fig. 1 in \cite{article:Br11})
and therefore the entropy evolution
equation cannot be (consistently) in steady state.

It would be of interest to attempt to drive his model to a
singular steady state. In lieu of a full analysis, in section
\ref{AdiabaticStagnation} I assume a smooth prescribed velocity
field (and magnetic field, which turns out to be irrelevant) and
solve the temperature evolution equation along the inflow axis
for the linearize problem in the vicinity of the origin.
It would also be of interest to determine a reasonable heat flux
that would allow his model to give a converged steady regular solution.

We now turn to dynamic simulations of the GEM magnetic
reconnection challenge problem using a 10-moment adiabatic model;
this problem seems to develop singularities shortly after peak
reconnection evidently due to the absence of a heat flux to
regularize the solutions.
The analysis of this section was motivated by difficulties
encountered when simulating the GEM problem as well as
a consideration of modeling requirements for steady reconnection.



%% file: chap4.tex
\chapter{GEM challenge problem}
\label{GEMchapter}

\section{GEM magnetic reconnection challenge problem}

As mentioned in the introduction, simulations of the
\defining{Geospace Environmental Modeling magnetic reconnection reconnection
challenge problem} (\mention{GEM problem}) with a variety of plasma
models (resistive MHD, Hall MHD, and particle models) identified the Hall effect as
critical to resolving fast magnetic reconnection \cite{article:GEM}.
It has been solved with a variety of plasma models and thus serves as a benchmark
for comparison of plasma models.  In particular, we
study the ability of the ten-moment adiabatic model to match
kinetic simulations of low-collisionality fast magnetic reconnection.

Since the Hall effect was identified as critical to magnetic reconnection,
it was natural to study whether fast reconnection can occur
in pair plasma, for which the Hall term of Ohm's law is zero.
Bessho and Bhattacharjee have simulated a pair plasma version
of the GEM problem using a particle code
\cite{article:BeBh05, article:BeBh07, bessho10} and have obtained
fast rates of magnetic reconnection.  Chac\'on et al.\ \cite{article:ChSiLuZo08}
subsequently demonstrated that steady fast reconnection is possible in a viscous
two-fluid model of magnetized pair plasma.  In their simulations
the reconnecting electric field is supported by the divergence of
the pressure tensors.


In studying the modeling requirements for magnetic reconnection
we seek the simplest problems that distinguish models.
Reconnection in pair plasma is a singular
and simple case which brings out modeling
requirements for fast reconnection with exceptional clarity.
In particular, \first{symmetric pair plasma} problems are 2D
antiparallel reconnection problems where there is symmetry under exchange of
species coupled with reflection in the out-of-plane axis;
for this case the number of equations that must be solved is halved.
The pair plasma version of the GEM problem studied by Bessho and Bhattacharjee 
is a symmetric pair plasma in the case of equal plasma temperatures.

This dissertation reports simulations of the GEM problem with
an adiabatic ten-moment two-fluid-Maxwell plasma model
for which the only collisional term retained is isotropization of
the pressure tensor.  I report simulations of the original GEM problem
and of symmetric pair plasma and compare with kinetic simulations.  

\subsection{Specification of the GEM problem}

We recall the definition of the GEM magnetic reconnection
challenge problem  \cite{article:GEM}.
The problem is symmetric under 180 degree rotation around the origin
and in the case of zero guide field is also symmetric under
reflection across both the horizontal and vertical axes.

\subsubsection{Nondimensionalization.}
The original GEM study  \cite{article:GEM}
nondimensionalizes by the ion inertial length
$\delta_i := \sqrt{\frac{m_i}{\mu_0 n_0 e^2}}$
and by the ion cylotron frequency $\Omega_i:=\frac{e B_0}{m_i}$.
(This nondimensionalization is worked out in section
\ref{CollisionlessNondimensionalization}.)
It then states that in these units the velocities
are normalized to the Alfv\`en speed $v_A$.
In actuality $v_A = \delta_i\Omega_i\sqrt{\frac{m_i}{m_i+m_e}}$
(see \eqref{AlfvenSpeed}),
a discrepancy which becomes significant in the pair plasma case
$m_i=m_e$.
In keeping with \cite{article:BeBh05,article:BeBh07},
when plotting we use $\delta_i$ and $\Omega_i$ computed from $m_i$ to
nondimensionalize time and space values, but use 
$v_A$ to nondimensionalize when displaying the electric field
(as indicated in the plots).
In the subsequent account we assume this nondimensionalization.

\subsubsection{Domain.}
(We designate the vertical axis $y$ and the out-of-plane
axis $z$, opposite to the convention employed in
 \cite{article:GEM, article:BeBh05,article:BeBh07}.)
The computational domain is the rectangular domain
$[-L_x/2,L_x/2]\times[-L_y/2,L_y/2]$.
We set $L_x=r_s 8\pi$ and $L_y=r_s 4\pi$,
where $r_s$ is a rescaling factor.
In the original GEM problem $r_s=1$ and
this scaling is used in this document except for
some pair plasma simulations where
results labeled \mention{``rescaled problem''} use $r_s=0.5$.
We also refer to this as the \mention{``half-scale''} problem.

\subsubsection{Boundary conditions.}
\label{GEMboundaryConditions}
The domain is periodic along the $x$-axis.
The boundaries parallel to the $x$-axis
are thermally insulating conducting wall boundaries.
A conducting wall boundary is a solid wall boundary
(with slip boundary conditions in the case of ideal plasma)
for the fluid variables, and any electric field at
the boundary must be perpendicular to the boundary (see
section 2.5 in \cite{griffiths89}).
Assuming the ideal MHD Ohm's law $\E=\B\times\u$,
which essentially says that the
electric field and the magnetic field are always perpendicular,
this implies that at the conducting boundary the magnetic field must
be parallel to the boundary.

So at the conducting wall boundaries
the five-moment two-fluid Maxwell variables satisfy
\begin{align*}
    \partial_y u_{x,s} &= 0,
 &  \partial_y B_x &= 0,
 &  E_x &= 0,
 & \partial_y \mdens_s &= 0,
 \\ u_{y,s} &= 0,
 &  B_y &= 0,
 &  \partial_y E_y &= 0,
 & \partial_y p_s &= 0 = \partial_y \pTot
 \\ \partial_y u_{z,s} &= 0,
 &  \partial_y B_z &= 0,
 &  E_z &= 0,
 & \partial_y \nrg_s &= 0 = \partial_y \NrgTot.
\end{align*}
These boundary conditions imply that gas particles are reflected off
the conducting wall.  Thus, 
for ten-moment gas-dynamic variables we use the boundary conditions
\begin{align*}
 \begin{aligned}
   \partial_y (\ET_\s)_{xx} &= 0, \\
   \partial_y (\ET_\s)_{xz} &= 0, \\
   \partial_y (\ET_\s)_{yy} &= 0, \\
   \partial_y (\ET_\s)_{zz} &= 0, 
 \end{aligned}
 && \hbox{and} &&
 \begin{aligned}
 (\ET_\s)_{xy} &= 0, \\
 (\ET_\s)_{yz} &= 0.
 \end{aligned}
\end{align*}

\def\BH{\B_H}
\def\JH{\J_H}
\def\phiH{\phi_H}
\def\dB{\B_P}
\def\dJ{\J_P}
\subsubsection{Initial conditions.}
The initial conditions are a perturbed Harris sheet equilibrium.
The unperturbed (Harris sheet) equilibrium is given by
\begin{align*}
    \BH(y) & =B_0\tanh(y/\lambda)\e_x,
  & p(y) &= \frac{B_0^2}{2 n_0} n(y),
 \\ n(y) &= n_0(0.2+\sech^2(y/\lambda)),
  & p_i(y) &= \frac{T_i}{T_i+T_e}p(y),
 \\ \E & =0,
  & p_e(y) &= \frac{T_e}{T_i+T_e}p(y).
\end{align*}
On top of this the magnetic field is perturbed by
\begin{align}
   \label{deltaB}
   \dB &= \curl(\psi\e_z) = \D\psi\times\e_z, &\hbox{ where }&&
   \psi(x,y)&=\psi_0 \cos(2\pi x/L_x) \cos(\pi y/L_y).
\end{align}

In the GEM problem the initial condition constants are
\begin{align*}
    T_i/T_e &= 5,
  & \lambda&=0.5,
  & B_0&=1,
  & n_0&=1,
  & \psi_0&=B_0/10.
\end{align*}
For pair plasma we reset the initial temperature ratio to 1 to get symmetry between
the species.   We set $\psi_0$ to $r_s^2 B_0/10$, so that in the vicinity
of the X-point our initial conditions agree (up to first-order
Taylor expansion) with the initial conditions of the GEM problem.

\subsubsection{Model parameters.}
The GEM problem specifies that the ion-to-electron mass ratio
is $\mi/\me=25$.  (The initial temperature ratio is defined to be
the square root of the mass ratio.)
We used this mass ratio for hydrogen plasma
and of course $\mi/\me=1$ for pair plasma.

\subsubsection{Auxiliary model parameters.}
The GEM problem does not specify the speed of light.
It seems to have been formulated with Galilean-invariant models in mind.
Published simulations with models that have a light speed
typically set the light speed to 20.
This seems to be a sufficiently close approximation to infinity
so that results are insensitive to the light speed.
The fastest wave speeds in the electrons are less than half of this.
In the ten-moment model the maximum wave speed in the electron
gas is initially about 4 and increases to about 7.

In the pair plasma version of the problem
we set the mass of each species to $r_s$
(rather than the GEM values of $1$ for ions
and $1/25$ for electrons).
We set the speed of light to 10 
(rather than 20 as in  \cite{article:BeBh05});
this seems close enough to infinity, since
the maximum gas wave speed is about 3.

\subsubsection{Implied initial quantities}

Using the initial conditions we calculate the precise initial current.
For convenience define
\begin{align*}
  \begin{aligned}
    \omega_x &:= 2\pi/L_x, \\
    \omega_y &:= \pi/L_y, 
  \end{aligned}
  &&\hbox{and} &&
  \begin{aligned}
    \theta_x &:= x\omega_x, \\
    \theta_y &:= y\omega_y.
  \end{aligned}
\end{align*}
So equation \eqref{deltaB} says
\begin{gather}
   \label{deltaBcopy}
   \begin{aligned}
     \psi&=\psi_0 \cos\theta_x \cos\theta_y,
   \end{aligned}
   \\
   \nonumber
   \begin{aligned}
     \hbox{so} &&
      \D\psi &= 
      \begin{bmatrix}
        -\omega_x \psi_0 \sin \theta_x \cos\theta_y  \\
        -\omega_y \psi_0 \cos \theta_x \sin\theta_y  \\
        0
      \end{bmatrix}
     & \hbox{and} &
     & \laplacian\psi&=-(\omega_x^2+\omega_y^2)\psi.
   \end{aligned}
\end{gather}
and therefore
\begin{gather}
   \begin{aligned}
      \dB &= \D\psi\times\ebas_z = 
      \begin{bmatrix}
        -\omega_y \psi_0 \cos \theta_x \sin\theta_y  \\
         \omega_x \psi_0 \sin \theta_x \cos\theta_y  \\
        0
      \end{bmatrix}
      &\hbox{and} \\
      \mu_0\dJ &= \curl\dB = -\laplacian\psi\ebas_z
          = (\omega_x^2+\omega_y^2)\psi\ebas_z.
   \end{aligned}
\end{gather}
Analogously, define
\begin{gather*}
  \begin{aligned}
  \phiH &:= B_0\lambda\ln\cosh(y/\lambda), 
  \end{aligned}
  \\
  \begin{aligned}
  &\hbox{so} &
  \D\phiH &= B_0 \tanh(y/\lambda) \ebas_x
  &\hbox{and} &&
  \laplacian\phiH &= (B_0/\lambda) \sech^2(y/\lambda).
  \end{aligned}
\end{gather*}
Then
\begin{gather*}
 \begin{aligned}
  \BH &= \curl(\phiH\,\ebas_z)
    = \D\phiH\times\ebas_z = B_0\tanh(y/\lambda)\ebas_y
  &\hbox{and} \\
  \mu_0\JH &= \curl\BH
    = -\laplacian\phiH\ebas_z = -(B_0/\lambda) \sech^2(y/\lambda)\ebas_z.
 \end{aligned}
\end{gather*}
Then
\begin{align*}
  \J &= \JH + \dJ = 
        \mu_0\inv\left(-\frac{B_0}{\lambda}\sech^2(y/\lambda)
              +(\omega_x^2+\omega_y^2) \psi(x,y)\right) \ebas_z.
\end{align*}

\subsection{Discussion of the initial perturbation}

The GEM problem begins with a Harris sheet equilibrium
and perturbs it by adding $\dB$ from equation
\eqref{deltaBcopy} to the magnetic field, leaving unchanged
the density, the net gas-dynamic momentum (zero),
and the net gas-dynamic pressure (i.e.\ energy).
The current of course changes, since for MHD it is
given by the curl of the magnetic field, and 
therefore, for full consistency of the two-fluid
equations with MHD,
one may consider the drift velocities
\begin{align*}
    \w_\s &= \frac{\mred}{m_\s}\frac{\J}{n q_\s}
\end{align*}
to be correspondingly perturbed (although this
probably does not matter much)
in the pressure evolution equation
and (in the case of the two-fluid Maxwell equations)
in the species momentum equations.


The effect of the perturbation of the magnetic field
is to take the momentum equation and Faraday's law
out of equilibrium.  This gets the fluid and the
magnetic field lines moving.

To see the effect of the perturbation, we invoke
the decomposition of
the current and magnetic field into Harris sheet and
perturbation components:
$\B = \B_H + \dB$ and  $\J = \J_H + \dJ$.
Substituting into momentum evolution
equation \eqref{MHDmomentumEvolution} gives
\begin{gather*}
  \partial_t(\mdens\u) +
  \underbrace{\Div(\mdens\u\u+\PT+\PTd)
  = \J_H\times\B_H}_{\hbox{these balance}}
  + \underbrace{\J_H\times\dB
  + \dJ\times\B_H}_{\hbox{these push out}}
  + \underbrace{\dJ\times\dB}_{\hbox{small}}.
\end{gather*}
The term labeled ``pushes out'' generates
a force away from the origin both on the
$x$-axis and on the $y$-axis.  It is not
obvious whether this will initially result in forward
or reverse reconnection.
In some of my simulations there is an initial
reverse reconnection and in others there is not.
In either case, forward reconnection eventually takes off.

\subsubsection{Initial conditions for the two-fluid model}

How should the initial electric field be specified
for the GEM problem?  The GEM problem does not explicitly
specify the initial electric field because it is
formulated with MHD in mind.
In the two-fluid model one needs to specify initial values
of the electric field.  Simply setting the electric field
to zero results in massive oscillations that could
destabilize the results.
A guiding principle is to specify quantities
so that the Harris sheet is an equilibrium.
The Harris sheet is designed so that net momentum 
evolution is satisfied.  To ensure that the momentum
evolution equations of both species are satisfied
we also need Ohm's law to be satisfied.  We can
simply compute the electric field using the appropriate
version of Ohm's law \eqref{OhmsLawTerms},
\begin{equation}
  \label{OhmsLawFullTerms}
  \begin{aligned}
  \E =& \Tresistivity\dotp\J  & \hbox{(resistive term)}
    \\ &+ \B\times\u  & \hbox{(ideal term)}
    \\ &+ \frac{\dmt}{\mdens}\J\times\B & \hbox{(Hall term)}
    \\ &+ \frac{\mut}{e n}\left[
          \Div\left(\frac{\PT_\i}{\mti}-\frac{\PT_\e}{\mte}\right)
        \right] & \hbox{(pressure term)}
    \\ &+ \frac{\mut}{e n}\left[
          \partial_t\J + \Div\left(\u\J+\J\u-\frac{\dmt}{\mdens}\J\J\right)
        \right] & \hbox{(inertial term).}
  \end{aligned}
\end{equation}
The ideal term, Hall term, and inertial term are all initially zero;
to see that the inertial term is indeed zero,
the flux component of the inertial term is zero, because
$\Div\J=0$, $\Div\u=0$, and $\J\dotp\D=0$ (because
$\J$ is initially out-of-plane);
we assume for the moment and will verify that $\partial_t\J=0$. 
Assuming the absence of resistivity, the elecric field must balance the
pressure term.  For the ten-moment two-fluid model the pressure should
be calculated by solving the equation in the next section,
and the initial electric field could be calculated to balance the
pressure term.  For the five-moment two-fluid model the pressure term
essentially reduces to a linear combination of scalar pressure gradients.
Using the initial pressure values of the GEM problem,
the pressure term is
\begin{gather*}
   \D\left(\frac{p_\i}{\mti}-\frac{p_\e}{\mte}\right)
     =  \frac{\mti T_\e - \mte T_\i}{T_\i+T_\e}
        \frac{B_0^2}{\mdens}\sech^2(y/\lambda)\tanh(y/\lambda)\ebas_y,
\end{gather*}
which we can take as $\E|_{t=0}$.
In this case $\E|_{t=0}$ is of the form $g(y)\ebas_y$,
so at time 0 $\partial_t \B = - \curl\E = 0$,
so $\partial_t \J|_{t=0}=0$, as assumed.

\subsubsection{Initial conditions for the ten-moment model}

In the ten-moment model what are the appropriate
initial values of the pressure tensor components?
In the simulations reported here the pressure tensor is
assumed initially isotropic. This results in strong
oscillations between the pressure and inertial terms
which decay over a time interval determined by the
isotropization period.

A better way to set the initial conditions would have been
to set $\dPT$ so that the Harris sheet is a near-equilibrium
in the adiabatic case.  Recall the adiabatic pressure tensor
evolution equation 
(see \eqref{nTTevolution}),
\begin{gather*}
   n_\s d_t \TT_\s + \SymB(\PT_\s\dotp\D\u_\s)
   = ({q_\s/m_\s})\SymB(\dPT_\s\times\B) - \dPT_\s/\ptime_\s.
\end{gather*}
For a Harris sheet $d_t=0$.  Using the decomposition $\PT=p\idtens+\dPT$,
we need
\begin{gather}
   \label{neededBalance}
   p_\s\SymB(\D\u_\s)+\SymB(\dPT_\s\dotp\D\u_\s)
   = ({q_\s/m_\s})\SymB(\dPT_\s\times\B) - \dPT_\s/\ptime_\s.
\end{gather}
This is six equations in five unknowns.
Taking the deviatoric part will give us five equations in five
unknowns, which we can solve for $\dPT$.
The original equation would be satisfied if its trace
is satisfied.
Half the trace is steady state scalar pressure evolution (see \eqref{pEvolution}),
\begin{gather}
     p_\s \Div \u_\s + \dPT_\s\ddotp\deviator{(\D\u_\s)} = 0.
\end{gather}
The initial conditions imply that $\Div\u=0$,
but we do not expect the entropy production term
$ \dPT\ddotp\deviator{(\D\u)}$ to be exactly zero,
since the initial velocity is (see equation \eqref{wsJ})
\begin{gather*}
  \u_\s = \w_\s = \frac{\mred}{m_\s}\frac{\J}{n q_\s}.
\end{gather*}
So we satisfy ourselves with enforcing the deviatoric part
of equation \eqref{neededBalance},
which is equivalent to the equation
\begin{gather*}
   \deviator{\SymB\left(\dPT_\s\dotp(\D\u_\s
     - {q_\s/m_\s}\idtens\times\B + \idtens/\ptime_\s)\right)}
   = -p_\s\deviator{\SymB(\D\u_\s)}.
\end{gather*}
This is a matrix equation of the form
\def\XX{\tensorb{X}}
\def\BB{\tensorb{B}}
\def\CC{\tensorb{C}}
\begin{gather*}
   \deviator{\SymB\left(\deviator{\XX}\dotp\BB\right)} = \deviator{\CC}.
\end{gather*}
It is five equations in five independent unknowns.

In the ten-moment case, attempting to
satisfy pressure evolution also comes
at the price of sacrificing that the net momentum evolution
equation \eqref{MHDmomentumEvolution}
\begin{gather}
  \partial_t (\mdens\u) + \Div(\mdens\u\u + \PTd + \PT)
  = e n(\J\times\B),
\end{gather}
is satisfied; this is consistent with the fact that the Harris
sheet equilibrium exactly satisfies the net momentum
equation in the five-moment two-fluid model only if viscosity
is neglected.

I remark that if one assumes isotropic pressure
then the Harris sheet equilibrium really does satisfy
the net momentum equation exactly even if the
diffusion pressure $ \PTd = \mte\mti\J\J/\mdens $
(see equation \eqref{PTdformula})
is included, because for the Harris sheet
$
  \Div\PTd = \mte\mti\J\dotp\D(\J/\mdens)=0,
$
using that $\J=\mu_0\inv\curl\B$ is directly out of plane.

\section{System of equations used to solve the GEM problem}

To simulate the GEM challenge problem we used the adiabatic
ten-moment two-fluid Maxwell model with pressure isotropization
in each species. When the pressure isotropization period is
zero the ten-moment model agrees with the adiabatic inviscid
five-moment model, as we have verified numerically.

The hyperbolic two-fluid Maxwell equations have been simulated
by Shumlak, Loverich, and Hakim.  The hyperbolic
five-moment two-fluid Maxwell system was simulated using
the Finite Volume wave propagation method
(described in \cite{book:Le02})
in Hakim's thesis \cite{Hakim06}
and specifically for the GEM problem 
in \cite{article:HaLoSh06}; the same system was
used to solve the GEM problem
using the discontinuous Galerkin method in \cite{LoHaSh11}.
Hakim solved the GEM problem using a hyperbolic two-fluid Maxwell system
using ten moments for the ions and five moments for the
electrons in \cite{article:Hakim08}; at the conclusion
of this paper Hakim proposes incorparation of collisions
through relaxation to a Maxwellian and closure
via use of a Chapman-Enskog type
expansion about the Gaussian distribution, as has
been carried out in this work.  I have not been
able to find studies of reconnection with
an isotropizing ten-moment two-fluid model,
although Hesse and Winske in \cite{hesseWinske93} simulate
collisionless ion tearing using particles for the ions
and a ten-moment isotropizing fluid for the electrons.
Miura and Groth \cite{MiuraGroth07} analyze dispersion in
an adiabatic ten-moment two-fluid model with BGK collision source terms.

To simulate the GEM challenge problem we implemented adiabatic
two-fluid-Maxwell models with five or ten moments for each
species. These models solve Maxwell's equations and solve a
separate compressible gas dynamics system for each species. The
general model uses ten-moment gas dynamics for each species and
isotropizes the pressure tensor at a tunable rate. We did not
include terms that represent collisional exchange between the
species (that is, resistive drag force or thermal equilibration).
We used a conservative shock-capturing method and therefore
represented the system of equations in conservation form.

\def\diminished{\color{cyan}}
\def\emphasized{\color{blue}}
To solve Maxwell's equations we solved the system
\begin{gather}
 \label{MaxwellWithCP}
  \partial_{t}
     \begin{bmatrix}
       \B \\
       \E \\
     \end{bmatrix}
  + \begin{bmatrix}
       \curl\E{\diminished +\chi\nabla\psi} \\
       -c^2\curl \B\\
     \end{bmatrix}
  = \begin{bmatrix}
       0 \\
       -\J/{\epsilon_0} \\
     \end{bmatrix},
\end{gather}
where $\J=e(n_\i\u_\i-n_\e\u_\e)$ is the current.
The correction potential $\psi$ is initially zero.
These equations imply a wave equation that propagates
the divergence constraint error $\Div\B$ at the speed $c\chi$.
We chose $\chi=1.05$.

The divergence constraint on the electric field is
$\Div\E=\qdens/\epsilon_0$, where $\qdens = e(n_\i-n_\e)$. The
initial conditions of the GEM problem satisfy the divergence
constraint, and physical solutions of Maxwell's equations
maintain the divergence constraint, but we do not attempt to
enforce this constraint numerically. In our simulations the
divergence constraint on the electric field remains approximately
satisfied in the sense that the error remains centered on zero
and its growth tapers. Attempts to apply correction potentials
to the electric field were counterproductive, increasing the
magnitude of the error In contrast, when we did not properly
apply correction potentials to the magnetic field, the error on
the magnetic field drifted from being centered at zero until the
solution became grossly unphysical.

The adiabatic pressure-isotropizing ten-moment gas-dynamic system
in conservation form, which we solved, is
\begin{equation}
 \label{TenMomentIsoSystem}
 \begin{aligned}
   \partial_t \mdens_\s + \Div(\mdens_\s\u_\s) &= 0,
 \\\partial_t (\mdens_\s\u_\s) + \Div(\mdens_\s\u_\s\u_\s + \PT_\s)
     &= {q_\s/m_\s}\mdens_\s (\E + \u_\s\times\B), 
 \\ \partial_t \ET_\s + 3\Div\Sym(\u_\s \ET_\s) -2\Div(\mdens_\s\u_\s\u_\s\u_\s)
      &= {q_\s/m_\s}2\Sym(\mdens_\s\u_\s\E+\ET_\s\times\B) + \RT_\s,
 \end{aligned}
\end{equation}
where the isotropization tensor is given by the closure
\begin{gather*}
  \RT_\s = \ptime_\s\inv\left(\idtens\tr\PT_\s/3-\PT_\s\right).
\end{gather*}
For the isotropization period $\ptime_\s$ we generally used
\begin{gather*}
  \ptime_\s = \ptime_0 \frac{\sqrt{m_\s\det\TT_\s}}{n_\s},
\end{gather*}
where $\ptime_0$ is a tunable parameter.
For pair plasma simulations, however, we
simply chose $\ptime_\s$ to be a uniform constant.

In the limit $\ptime_\s\to 0$ the ten-moment system simplifies
to the adiabatic inviscid five-moment model:
\begin{equation}
 \begin{aligned}
 \label{FiveMomentIsoSystem}
   \partial_t \mdens_\s + \Div(\mdens_\s\u_\s) &= 0,
 \\\partial_t (\mdens_\s\u_\s) + \Div(\mdens_\s\u_\s\u_\s) + \D p_\s
     &= {q_\s/m_\s}\mdens_\s (\E + \u_\s\times\B), 
 \\ \partial_t \Nrg_\s + \Div(\u_\s(\Nrg_\s+p_\s))
      &= \J_\s\dotp\E.
 \end{aligned}
\end{equation}

%% file: chap5.tex
\chapter{GEM pair plasma simulations}

This chapter reports simulations of GEM problem
described in the previous chapter for the case
of pair plasma.
We assumed zero guide field and equal temperature for both species.
For equal temperature, symmetry between the two species halves the
number of equations needed. For zero guide field the GEM problem is
symmetric across both the horizontal and vertical axes. We enforced
all these symmetries.

All simulations were computed on a $128$ by $64$ mesh
unless otherwise indicated.  To verify convergence we
coarsened the grid by a factor of two in each direction.

As a proxy for Ohm's law we plotted the accumulation integrals
of the terms of \eqref{proxyOhmXpoint} (the positron momentum
equation solved for the electric field):
\begin{align*}
 \text{electric term: }& -\int_0^t \E_z, \\
 \text{pressure term: }& -\int_0^t \frac{(\nabla\cdot\PT_\i)_z}{e n_\i},\\
 \text{inertial term: }& -\int_0^t \frac{m_\i}{e} \partial_t (u_\i)_z
                         -\frac{m_i}{e}\big((u_\i)_z|_t - (u_\i)_z|_{t=0}\big),\\
 \text{residual term: }& -\int_0^t \text{residual},
\end{align*}
where $\text{residual}$ 
represents numerical resistance and is what 
must be added to the pressure term and inertial term to yield
the reconnection electric field.  We have plotted the loss of magnetic 
flux across the $y$-axis and have verified that it is indistinguishable
from the electric term.

\section{Rescaled pair-plasma GEM problem}

To avoid the formation of magnetic islands, we modified (Bessho
and Bhattacharjee's version of) the GEM problem, shrinking the
dimensions of the domain and the particle mass to half their original values.
The effect of these changes is to rescale the nondimensional
parameters of the GEM problem, increasing the width of the
current sheet from $1$ to $\sqrt{2}$ times the ion inertial
radius and decreasing the domain size from $8\pi$ by $4\pi$ to
$8\pi/\sqrt{2}$ by $4\pi/\sqrt{2}$.  We set the speed of light
to be $10$ times the Alfv\'en speed $v_A$.

\begin{figure}
    \begin{center}
    \includegraphics[width=\linewidth,height=.6\linewidth]{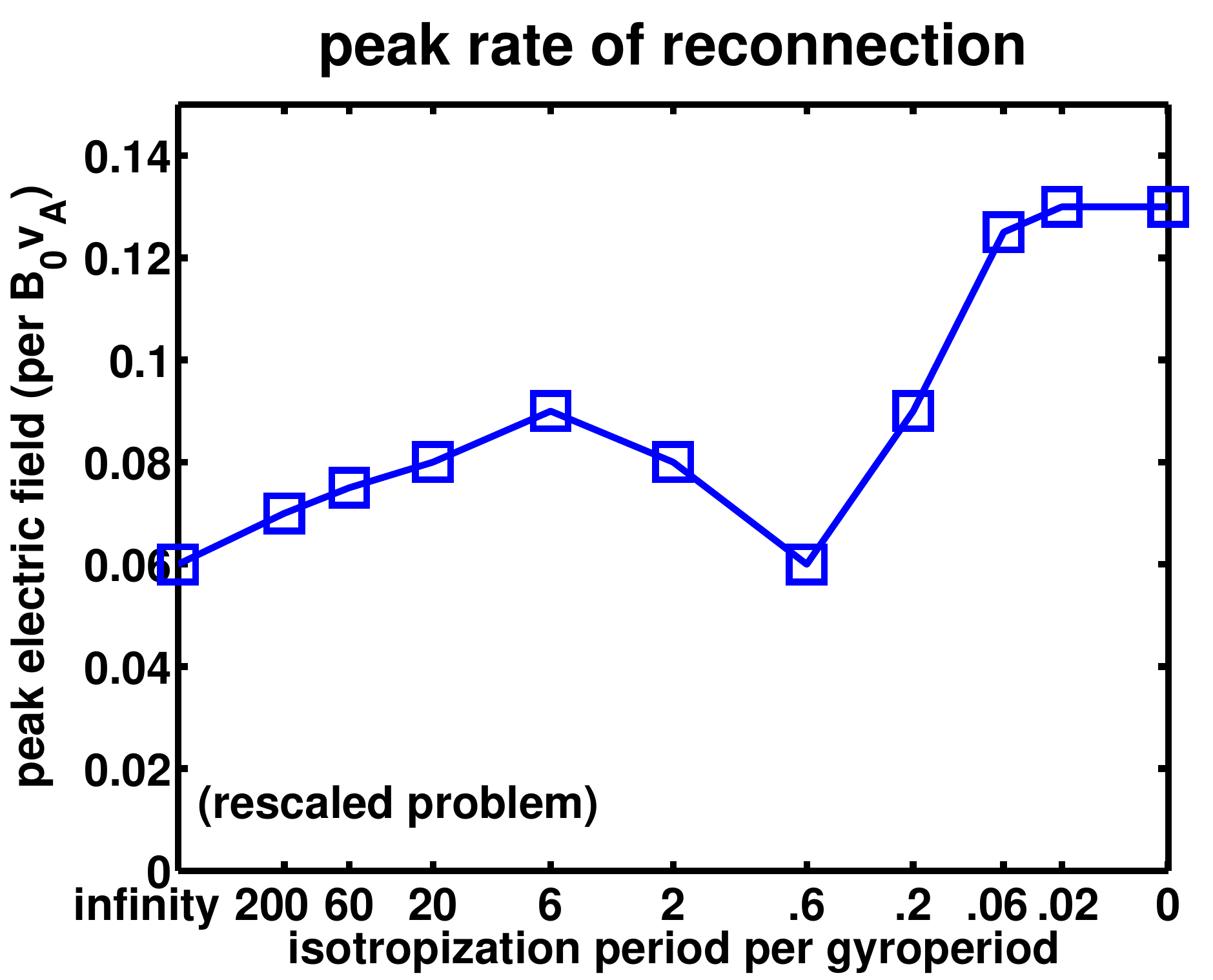}
    \end{center}
   \caption{Peak rate of reconnection versus isotropization period
     in the half-scale symmetric pair plasma GEM problem.}
     Observe that for intermediate isotropization
     rates increasing the rate of isotropization decreases the
     rate of reconnection, as reported in \cite{article:HeWi93},
     but for extreme rates of isotropization we observe the opposite trend.
  \fhrule
  \label{fig:recon_rate_per_isoperiod_halfscale}
 \end{figure}

We simulated this rescaled GEM problem for isotropization rates
ranging from 0 to instantaneous. For intermediate isotropization
rates increasing the rate of isotropization decreases the
rate of reconnection, which seems to agree qualitatively with
 \cite{article:HeWi93}, but for extreme rates of isotropization we
observed the opposite trend.
(See figure \ref{fig:recon_rate_per_isoperiod_halfscale}).
 
\begin{figure}
   \begin{center}
    \includegraphics[width=\linewidth,height=.6\linewidth]{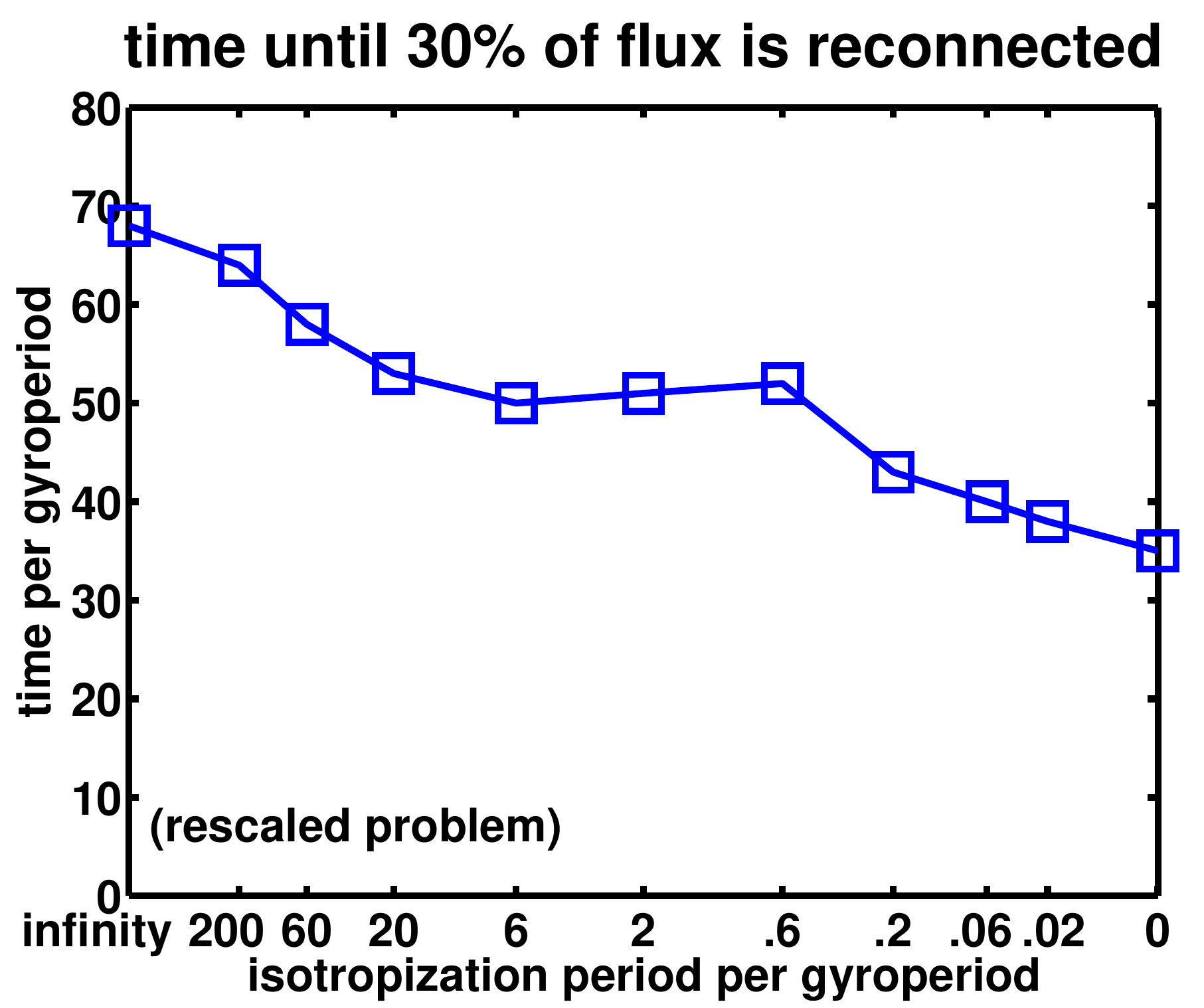}
   \end{center}
   \caption{Time until 30\% of flux is reconnected 
     in the half-scale symmetric pair plasma GEM problem.}
     For all simulations the peak rate of reconnection occurred
     when about 30\% of the original flux through the y-axis had reconnected.
  \fhrule
  \label{fig:recon_time_per_isoperiod_halfscale}
 \end{figure}

\begin{figure}
  \begin{tabular}{c c}
   \begin{tabular}{c}
     \includegraphics[width=0.55\linewidth]{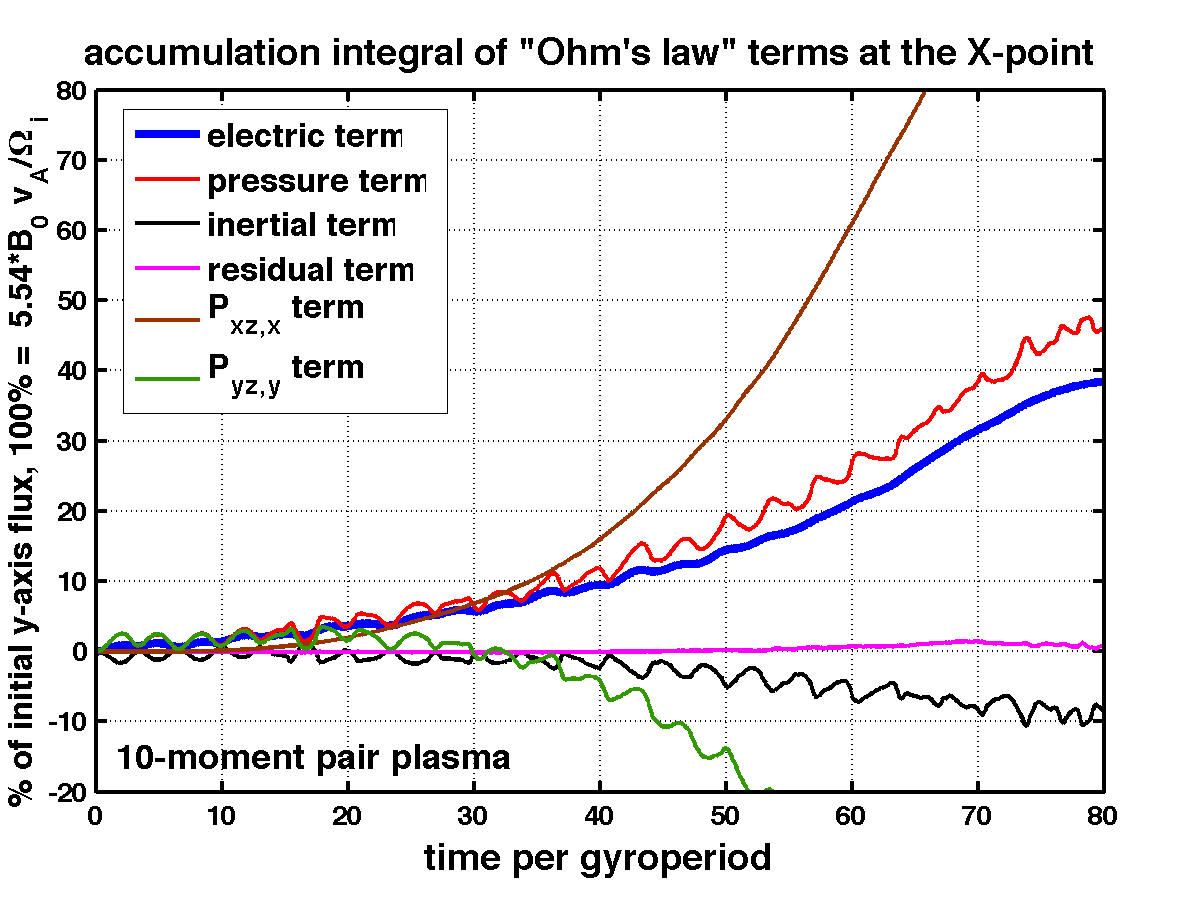}
  \\ \includegraphics[width=0.55\linewidth]{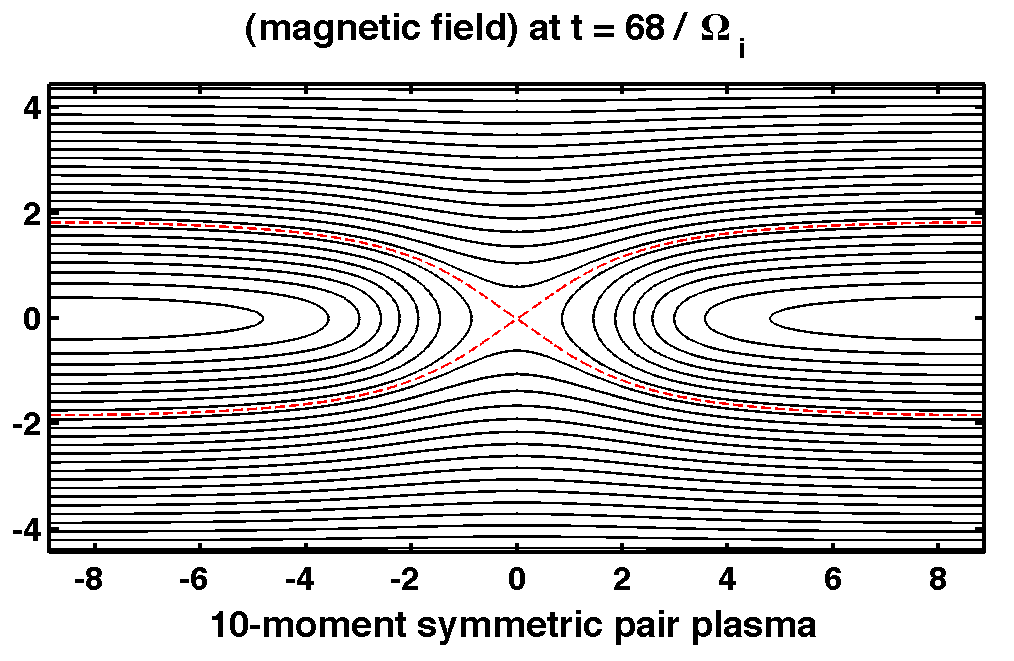}
   \end{tabular}
  &
   \begin{tabular}{c}
      \includegraphics[width=0.45\linewidth]{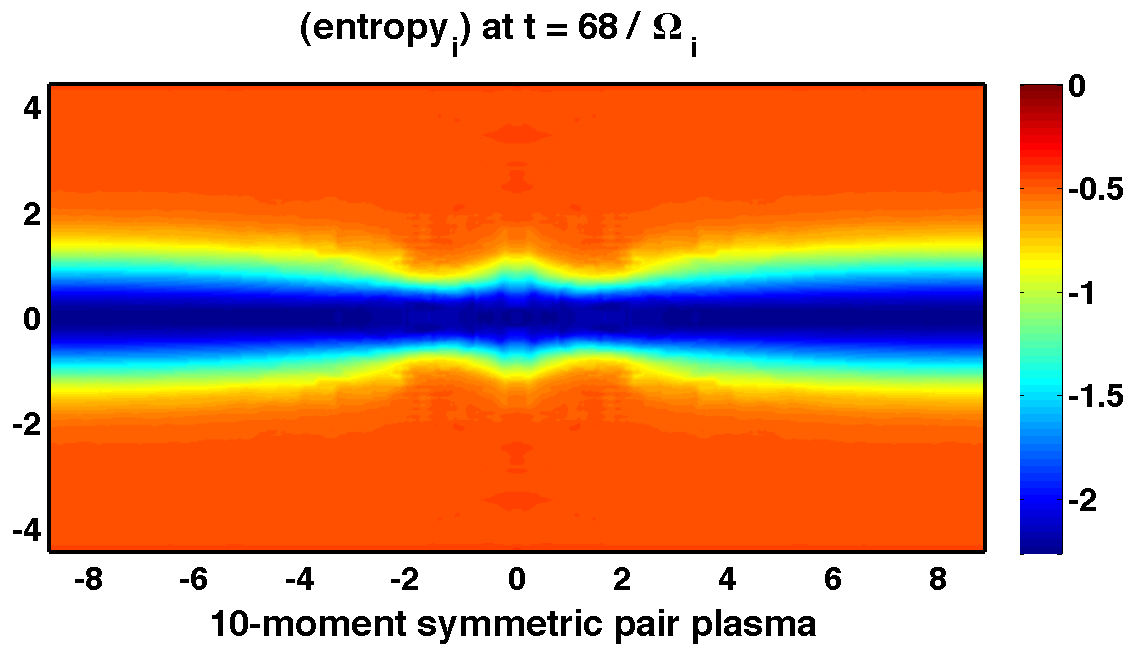}
   \\ \includegraphics[width=0.45\linewidth]{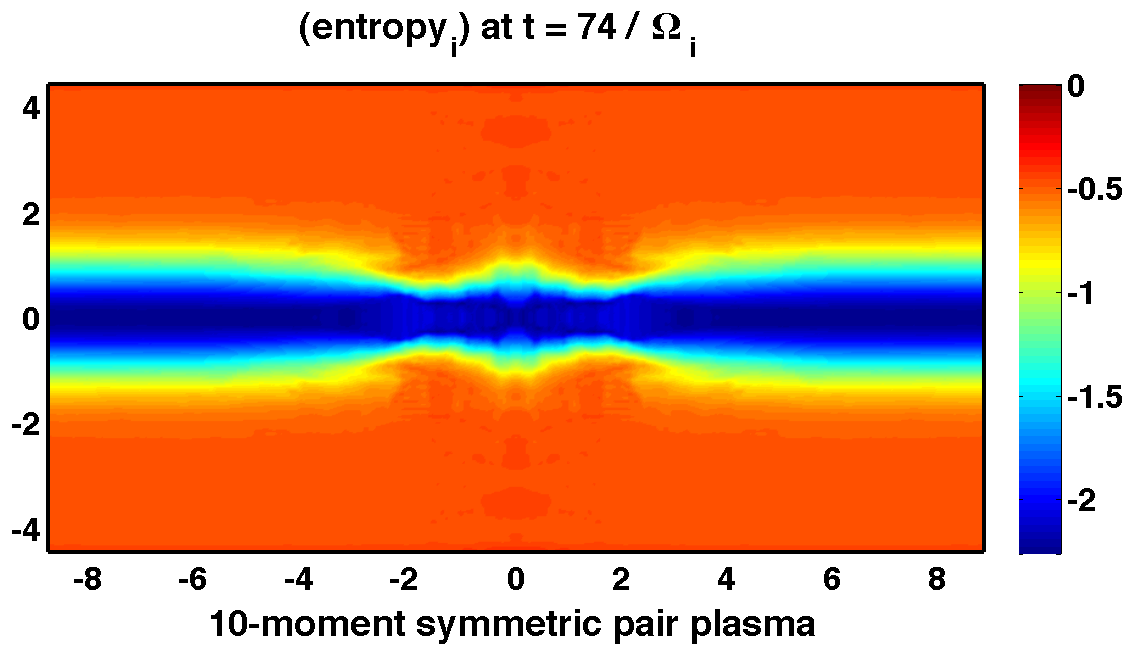}
   \\ \includegraphics[width=0.45\linewidth]{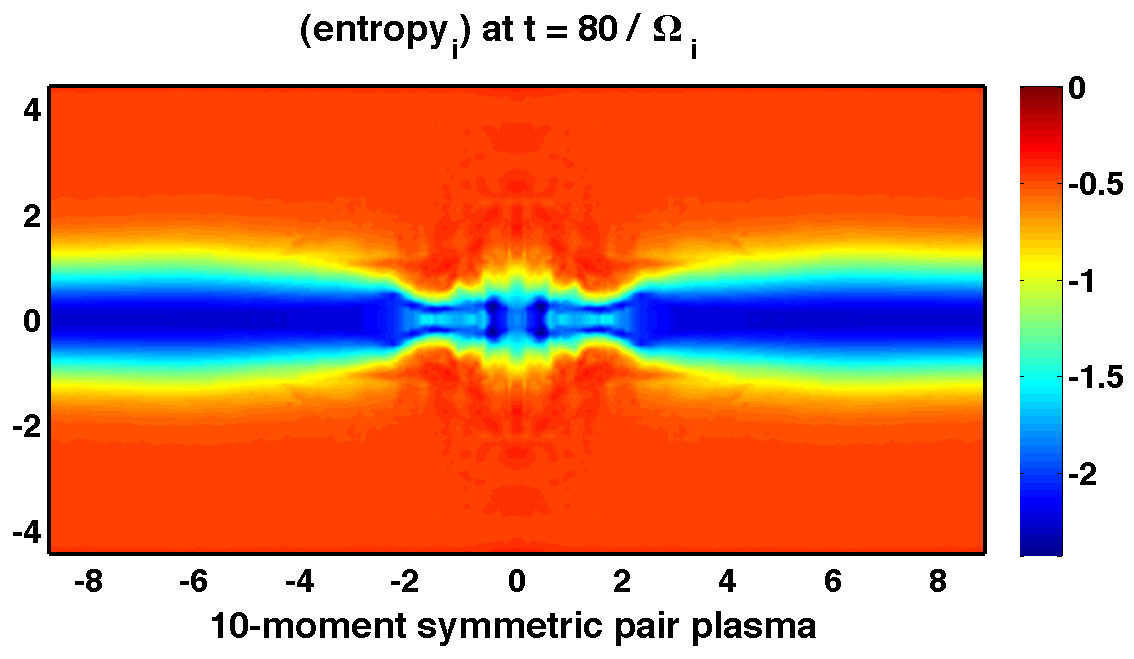}
   \end{tabular}
  \end{tabular}
   \caption{Reconnection for no isotropization (hyperbolic 10-moment model)
     in the half-scale symmetric pair plasma GEM problem.}
   Notice the undamped oscillatory exchange between the pressure 
   and inertial terms. 
   Notice also the U-shaped magnetic field line pattern near the X-point;
   as we increase the isotropization rate in the following figures
   the magnetic field becomes
   more V-shaped or even Y-shaped, evidently due to decreasing viscosity.
   and increasing wont to isotropy.

   This model is hyperbolic
   and requires high resolution for convergence.
  \vspace{3.5ex}
  \fhrule
  \label{fig:infty}
 \end{figure}

\begin{figure}
  \begin{center}
  \end{center}
  \begin{tabular}{c c}
   \begin{tabular}{c}
     \includegraphics[width=0.55\linewidth]{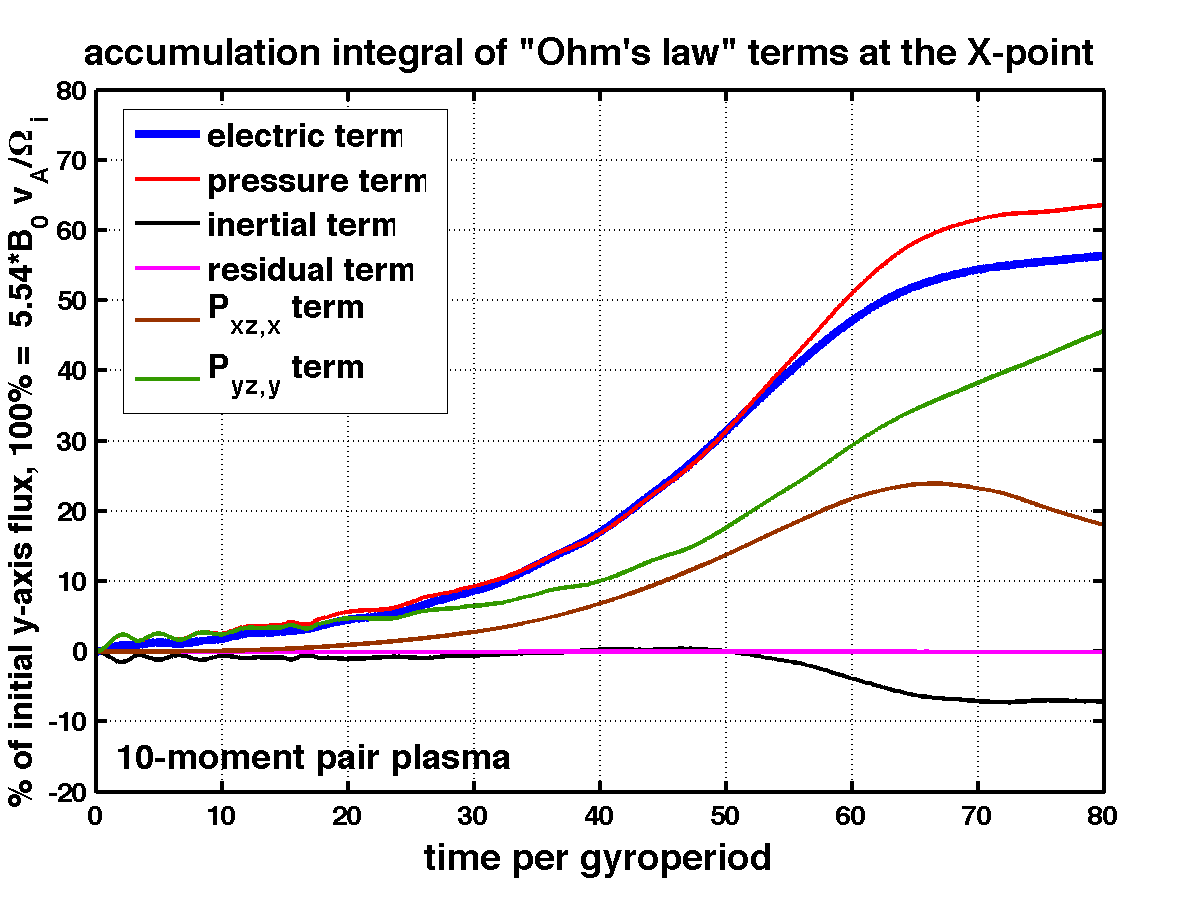}
  \\ \includegraphics[width=0.55\linewidth]{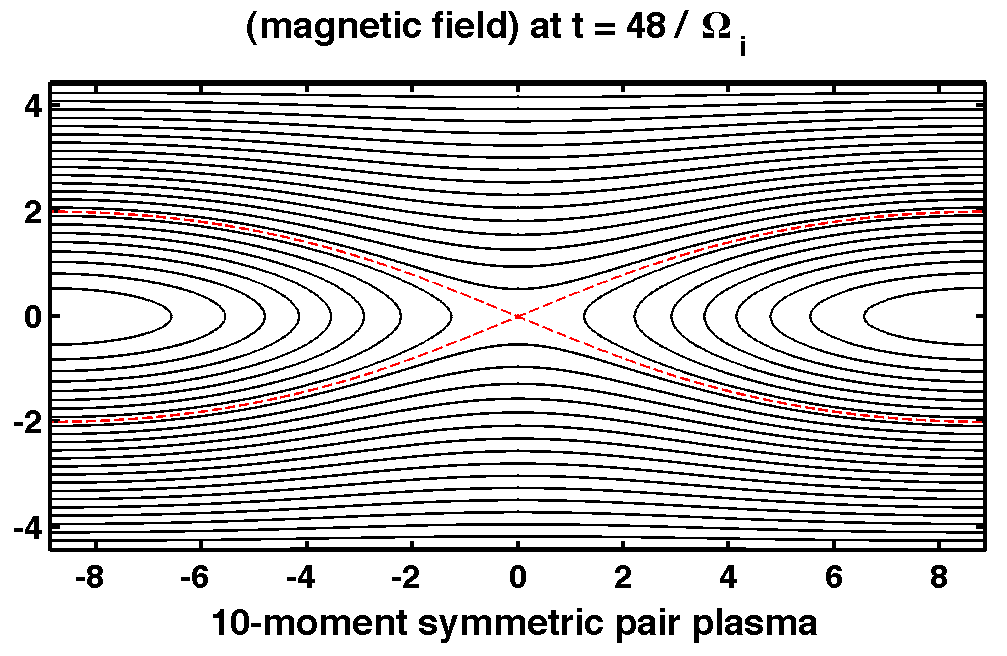}
   \end{tabular}
  &
   \begin{tabular}{c}
      \includegraphics[width=0.45\linewidth]{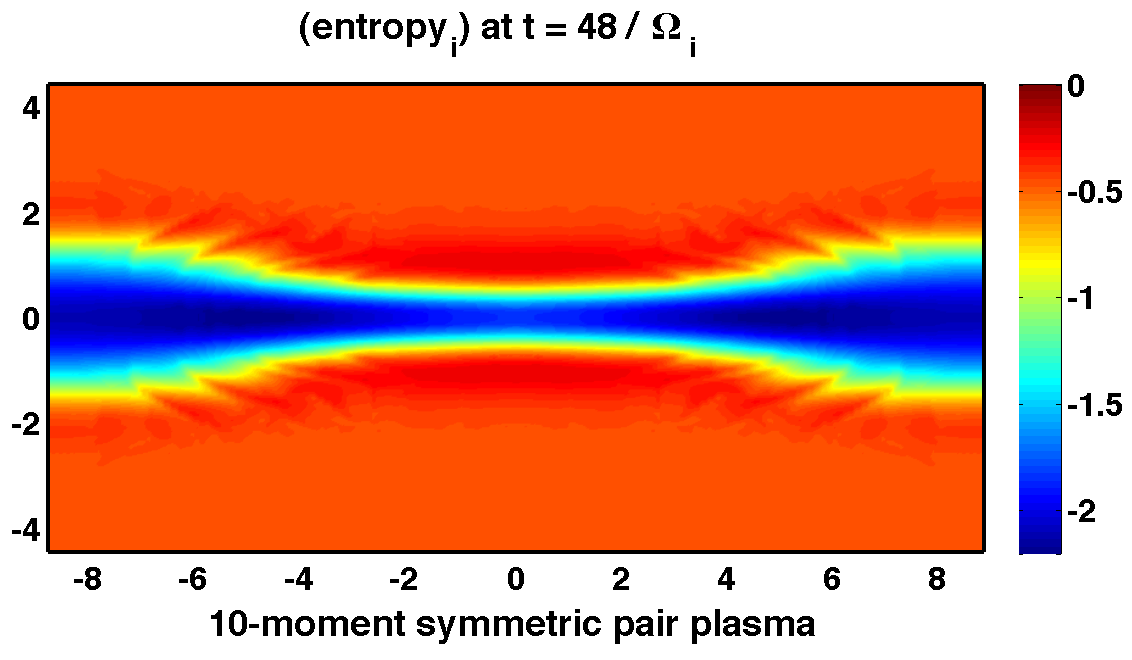}
   \\ \includegraphics[width=0.45\linewidth]{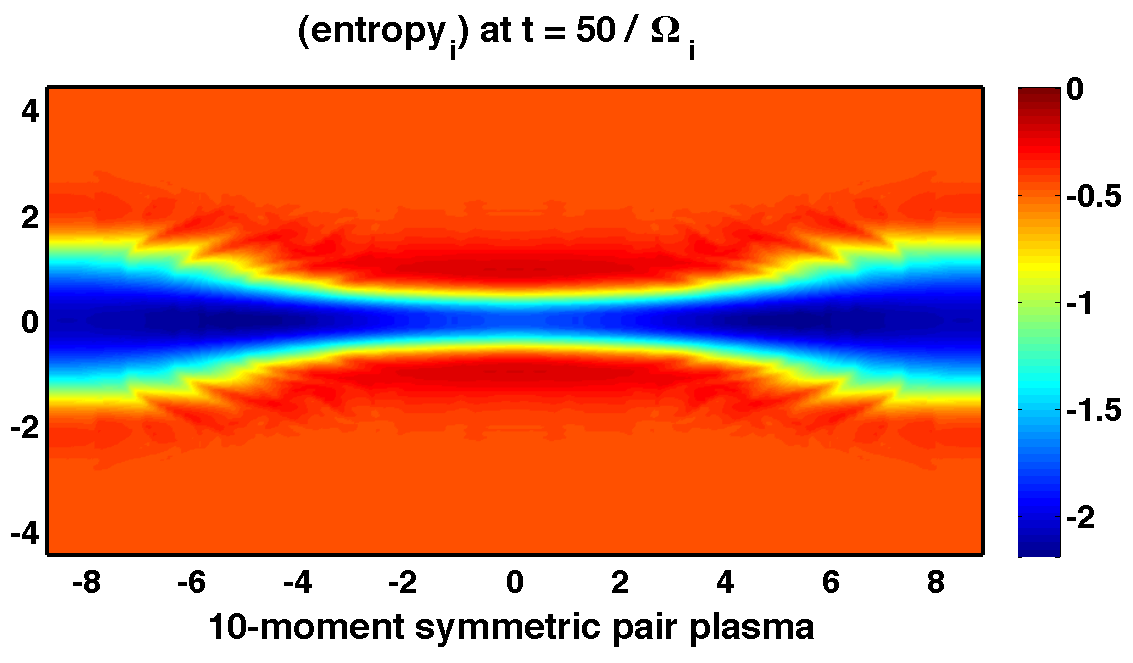}
   \\ \includegraphics[width=0.45\linewidth]{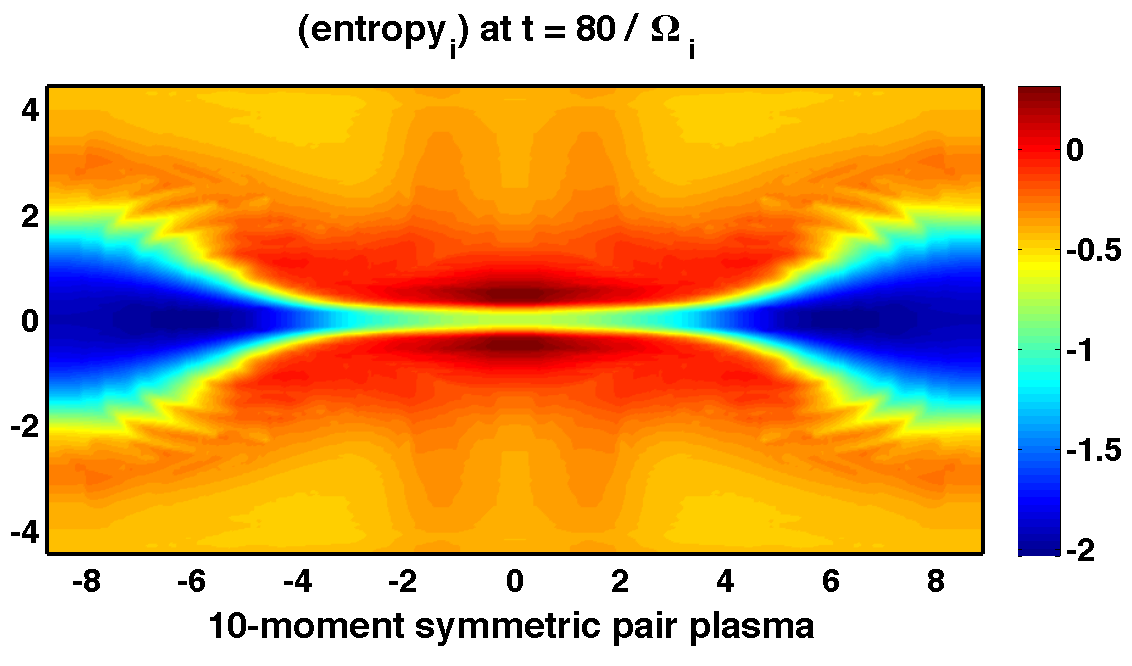}
   \end{tabular}
  \end{tabular}
   \caption{Reconnected flux for moderate isotropization
     in the half-scale symmetric pair plasma GEM problem.}
    Isotropization of the pressure tensor 
    dampens the oscillatory exchange between the pressure
    and inertial terms.  The pressure term supports reconnection,
    in agreement with steady-state theory and PIC simulations.
    The maximum Knudsen number is roughly 6 over the course of this simulation.
    The five-moment closure for the pressure tensor is completely incorrect
    and is an order of magnitude too large.
  \vspace{8ex}
  \fhrule
  \label{fig:6}
 \end{figure}

\begin{figure}
  \begin{tabular}{c c}
   \begin{tabular}{c}
     \includegraphics[width=0.55\linewidth]{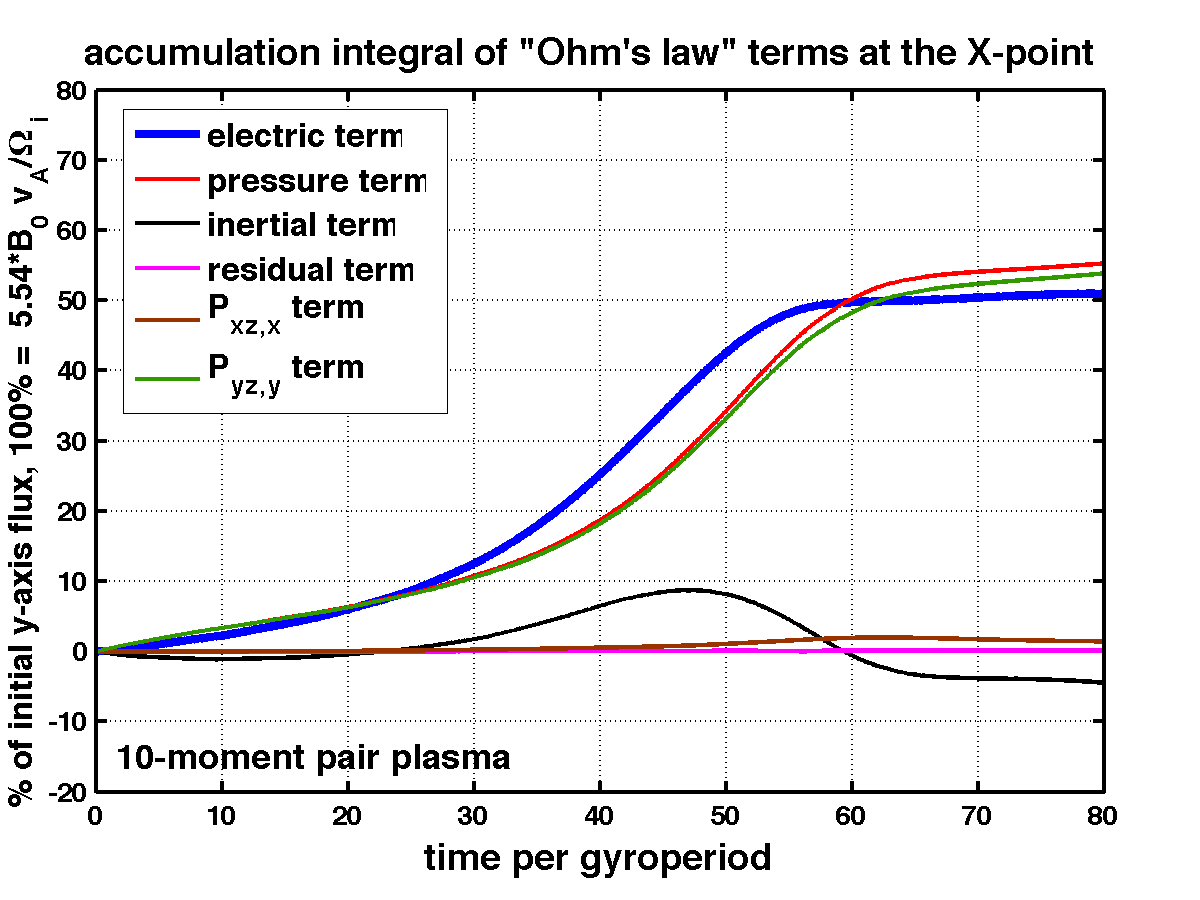}
  \\ \includegraphics[width=0.55\linewidth]{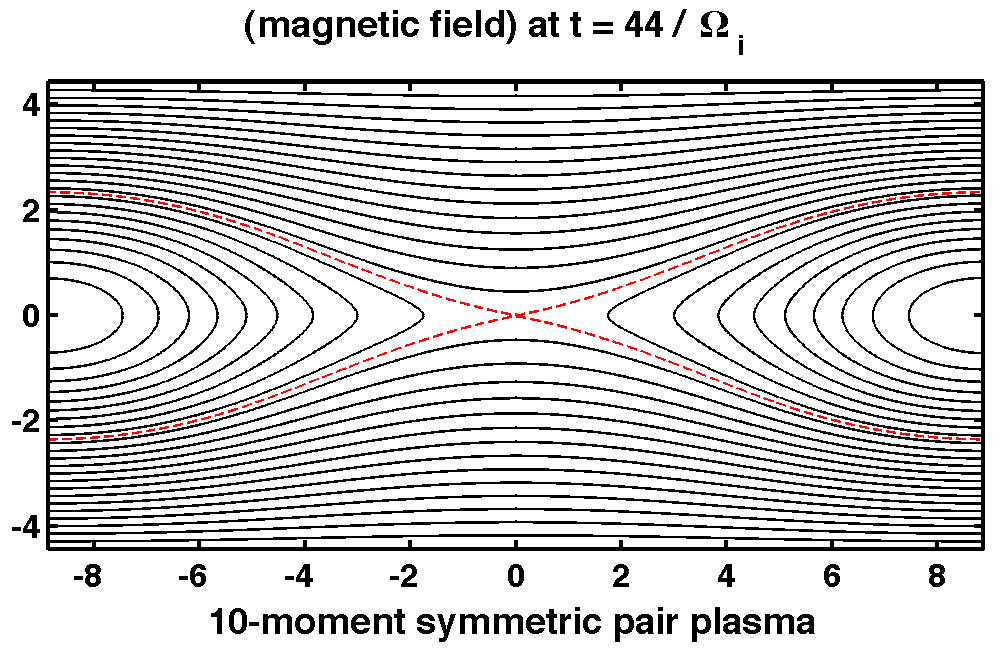}
   \end{tabular}
  &
   \begin{tabular}{c}
      \includegraphics[width=0.45\linewidth]{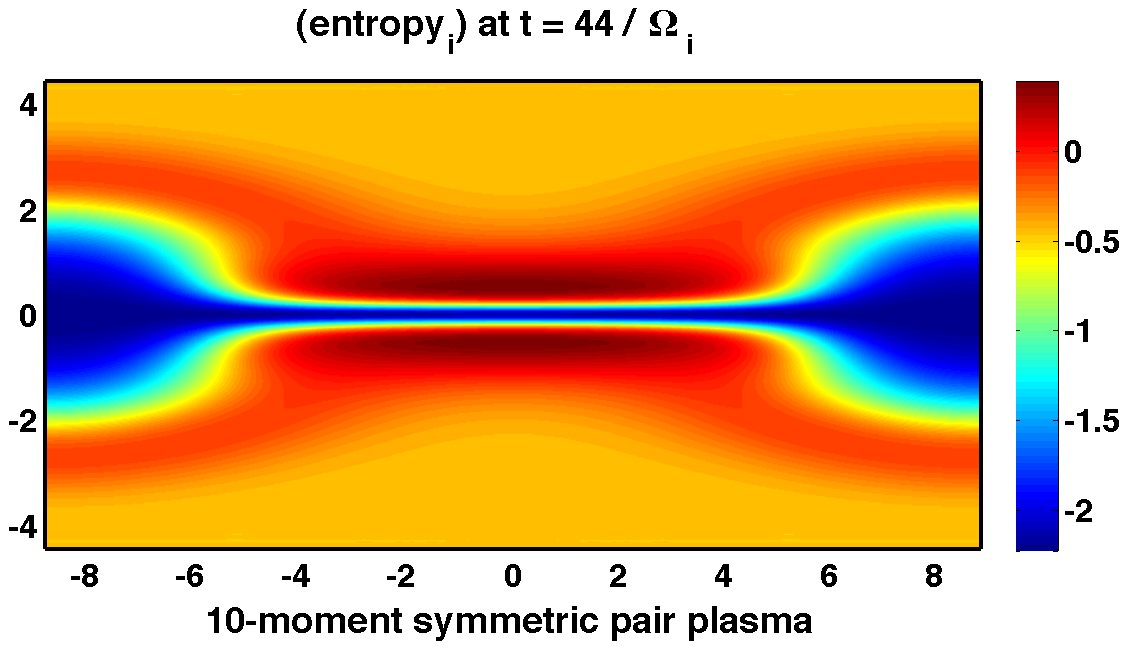}
   \\ \includegraphics[width=0.45\linewidth]{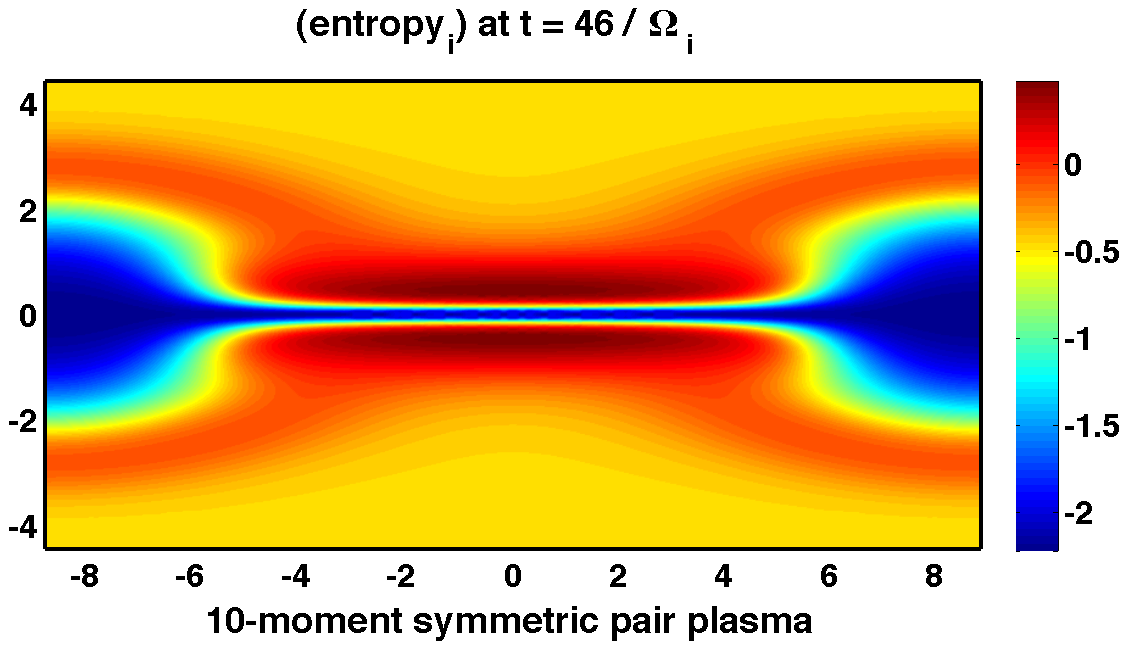}
   \\ \includegraphics[width=0.45\linewidth]{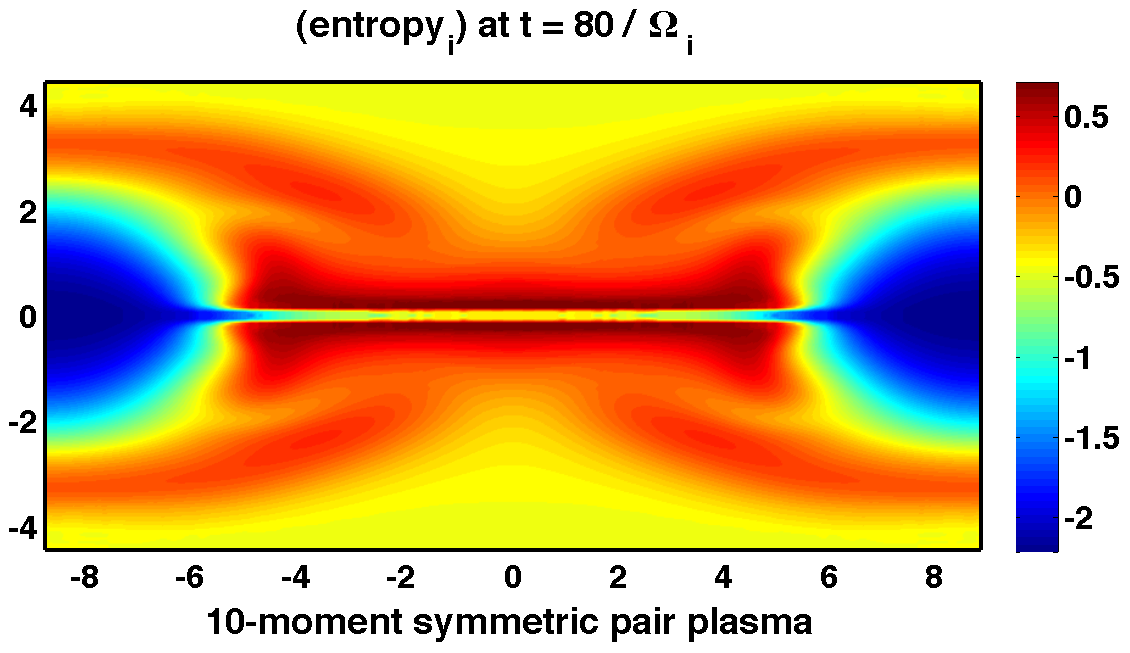}
   \end{tabular}
  \end{tabular}
   \caption{Reconnected flux for fast isotropization
     in the half-scale symmetric pair plasma GEM problem.}
    When the rate of isotropization is fast the
    inertial term initially provides some support for
    reconnection, although the pressure term ultimately 
    provides the support.
    Entropy production occurs in the low-entropy sheet
    along the outflow axis due to pressure isotropization.
    Viscosity prevents turbulence from developing, and
    the solution remains regular over the whole course
    of the simulation.
    The maximum Knusden number is roughly .25 over the course of this simulation.
    The pressure tensor generally shows good agreement with the five-moment closure.
  \fhrule
  \label{fig:0_2}
 \end{figure}

\begin{figure}
  \begin{tabular}{c c}
   \begin{tabular}{c}
     \includegraphics[width=0.55\linewidth]{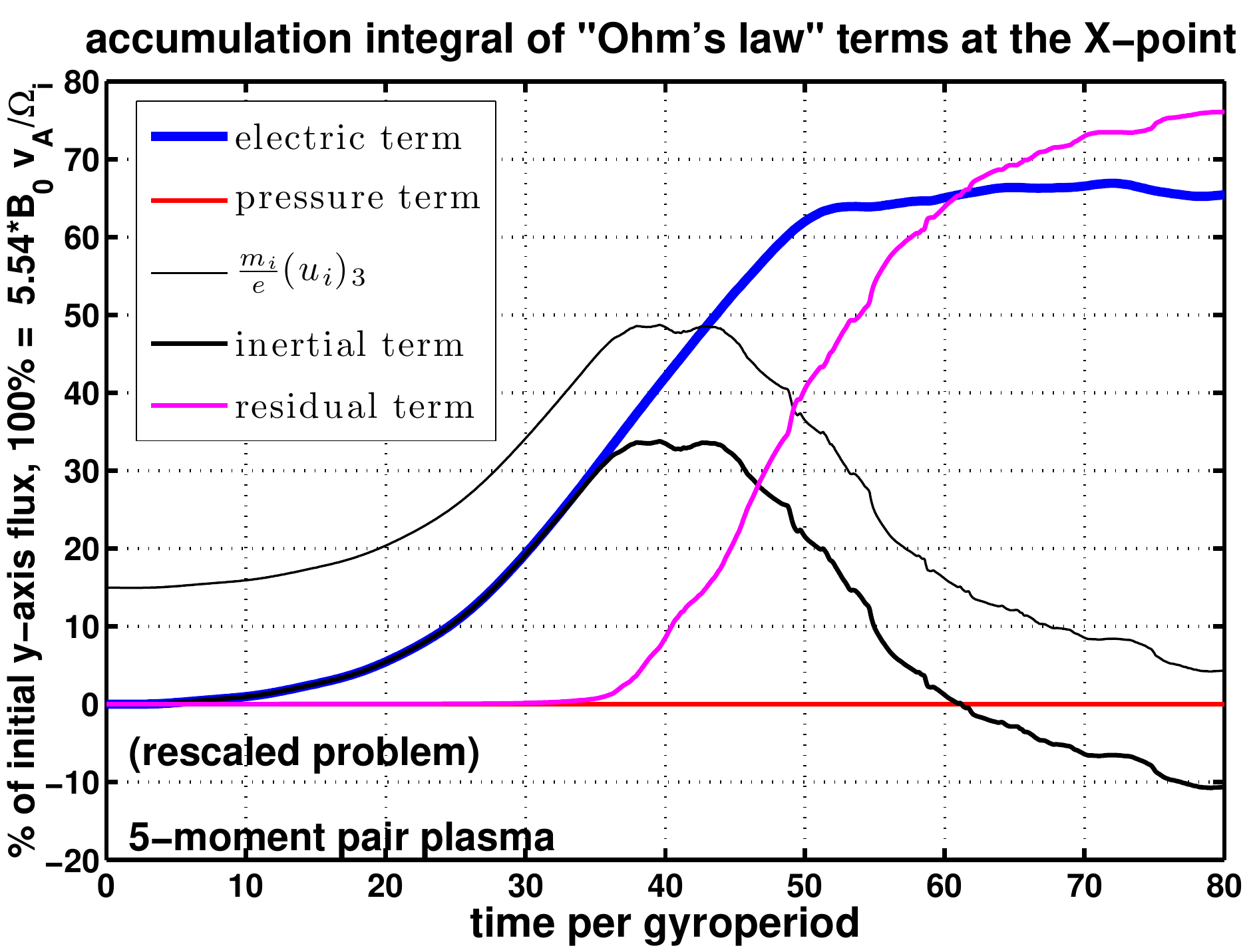}
  \\ \includegraphics[width=0.55\linewidth]{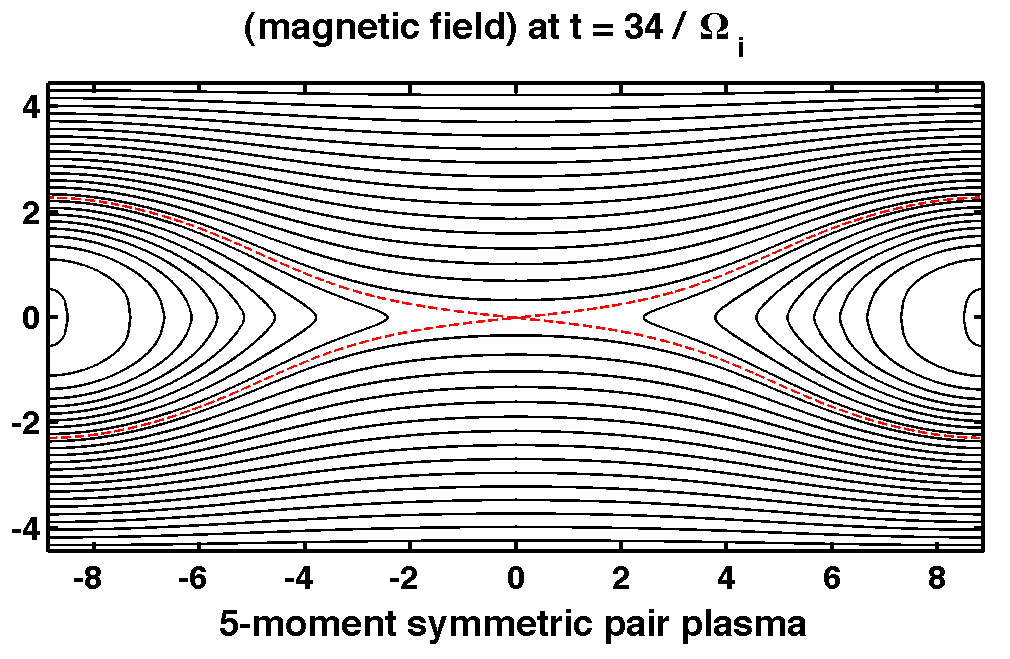}
   \end{tabular}
  &
   \begin{tabular}{c}
      \includegraphics[width=0.45\linewidth]{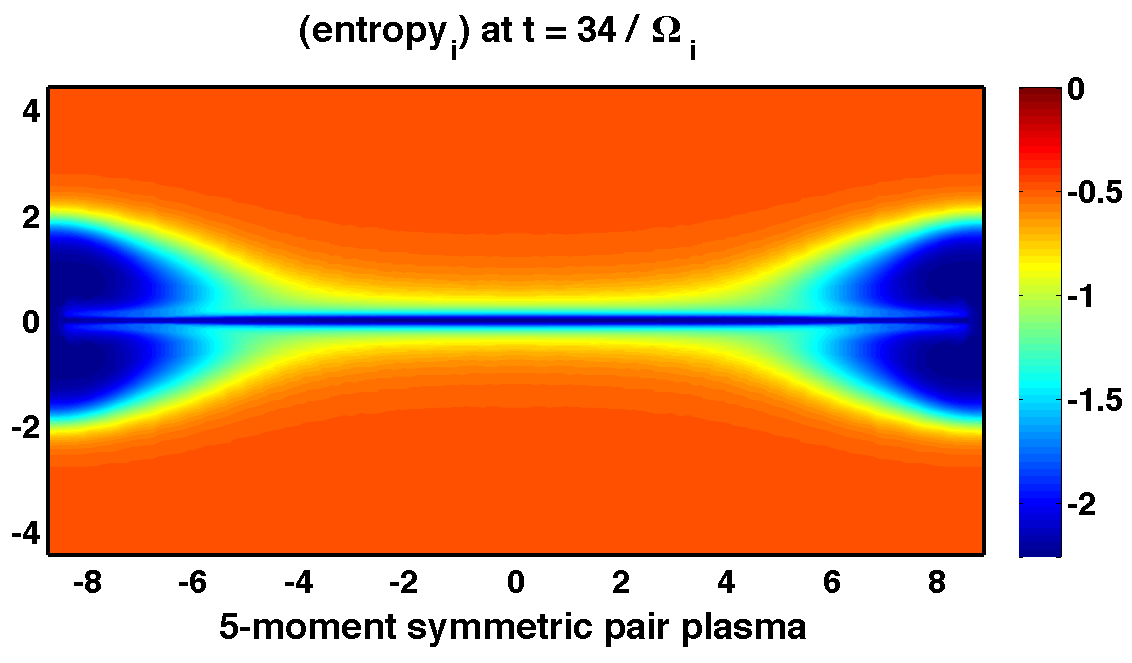}
   \\ \includegraphics[width=0.45\linewidth]{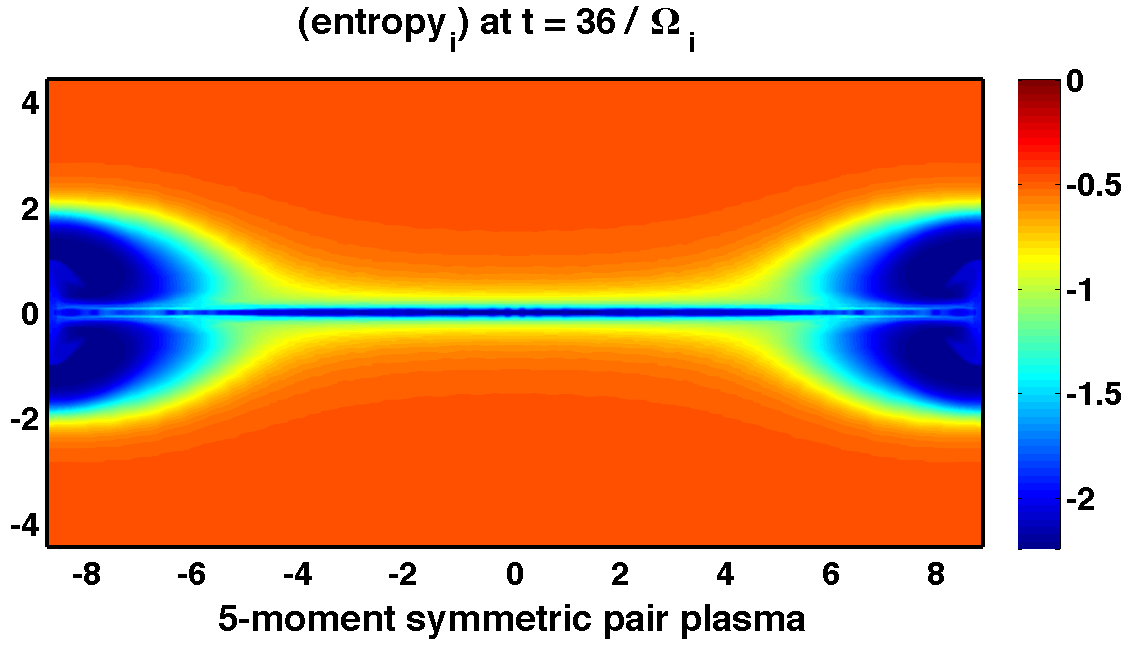}
   \\ \includegraphics[width=0.45\linewidth]{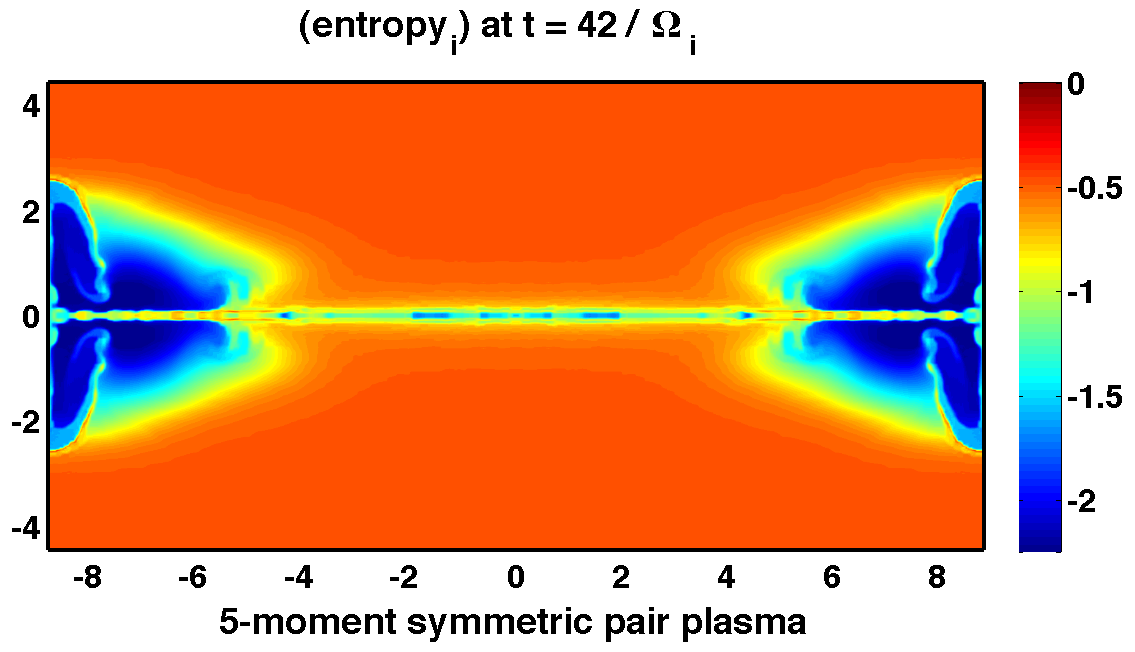}
   \end{tabular}
  \end{tabular}
  \caption{Reconnected flux for the five-moment model (instantaneous isotropization)
     in the half-scale symmetric pair plasma GEM problem.}

   In the five-moment model the inertial term initially tracks with
   reconnected flux until there is a narrow sheet of low entropy along the
   outflow axis, but beginning at a time of about 36 
   numerical residual kicks in.  At almost exactly the
   the same time turbulence begins to form at the walls and
   the sheet of low entropy at the center begins to form beads.
   Sharp outflow jets generate strong turbulence at late times.

   In this simulation we used a five-moment
   model computed on a 256 by 128 mesh.
   For coarser meshes we verified that
   the instantaneously relaxed ten-moment model
   agrees with the five-moment model.
  \fhrule
  \label{fig:0}
 \end{figure}

For all simulations the peak rate of reconnection occurred
when about 30\% of the original flux through the y-axis had reconnected.
The time until 30\% reconnection showed trends similar to
the trends in the peak isotropization rate
(see figure \ref{fig:recon_time_per_isoperiod_halfscale}).
In general we can say that the rate of reconnection is not very
sensitive to the rate of isotropization, which seems to
agree qualitatively with  \cite{article:KuHeWi01}.
We display the accumulation integral of the
terms of \eqref{proxyOhmXpoint} for some of these isotropization periods
($\infty$, $6$, $0.2$, and $0$) in figures (\ref{fig:infty}, \ref{fig:6},
\ref{fig:0_2}, and \ref{fig:0}).

 
 \clearpage
 
\section{Full-scaled pair plasma GEM problem}
  To compare with Bessho and Bhattacharjee
   \cite{article:BeBh05,article:BeBh07}, we also used their
  GEM-like (full) scaling ($r_s=1$). We allowed the isotropization period
  to vary over a smaller range (between $.3$ and $30$), since for
  extreme isotropization periods central magnetic islands formed.
  Over this range our reconnection rate varied between $.13$
  and $.14$, and our time to peak (i.e.\ 30\%) reconnection
  increased from $32$ for a slow isotropization period of $30$ to
  $37$ for a fast isotropization period of $.3$.
 
\begin{figure}
    \begin{center}
      \includegraphics[width=\linewidth,height=.6\linewidth]{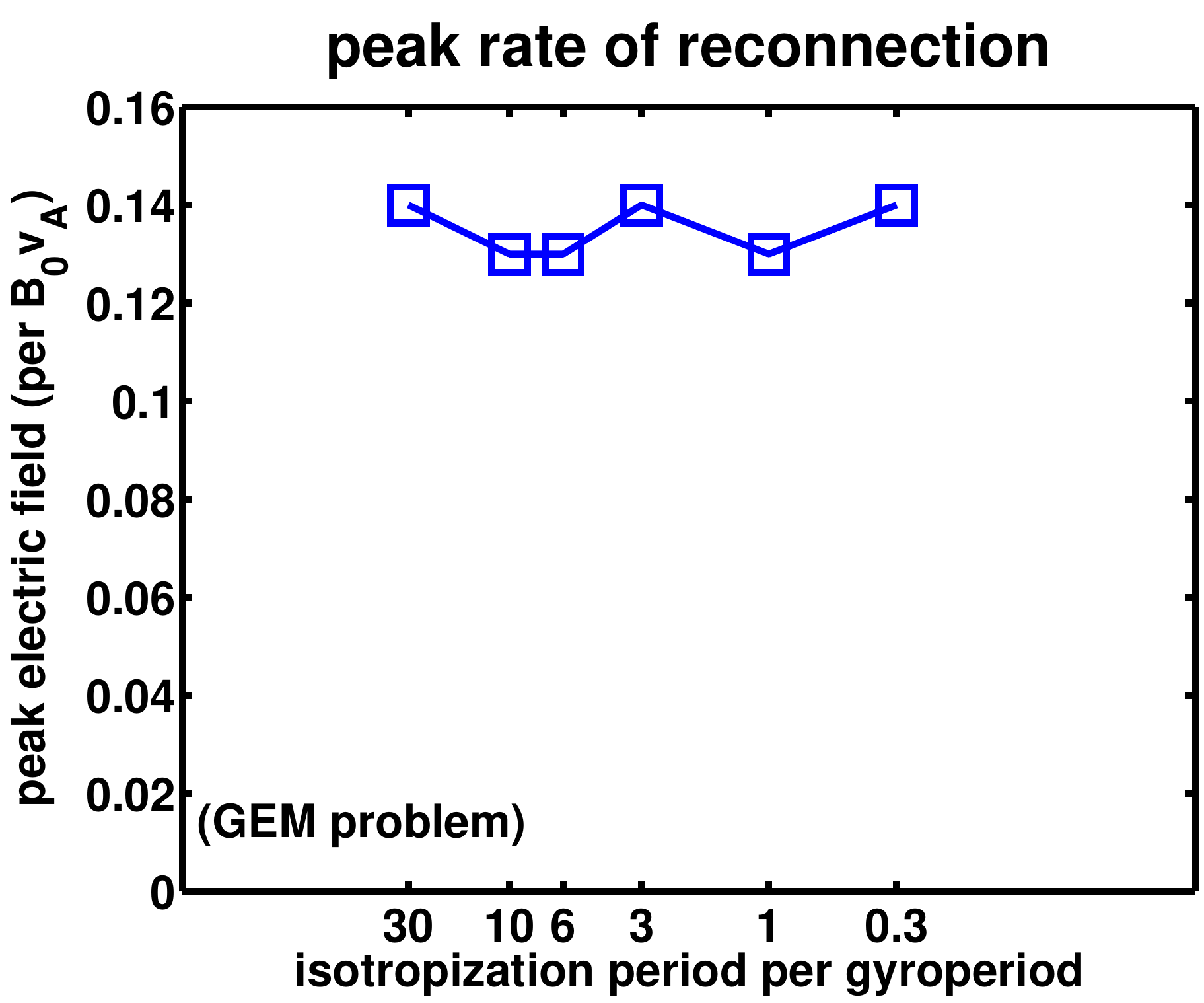}
    \end{center}
   \caption{Peak rate of reconnection versus isotropization period
     in the full-scale symmetric pair plasma GEM problem.}
  \fhrule
  \label{fig:recon_rate_per_isoperiod_fullscale}
 \end{figure}

\begin{figure}
   \begin{center}
     \includegraphics[width=\linewidth,height=.6\linewidth]{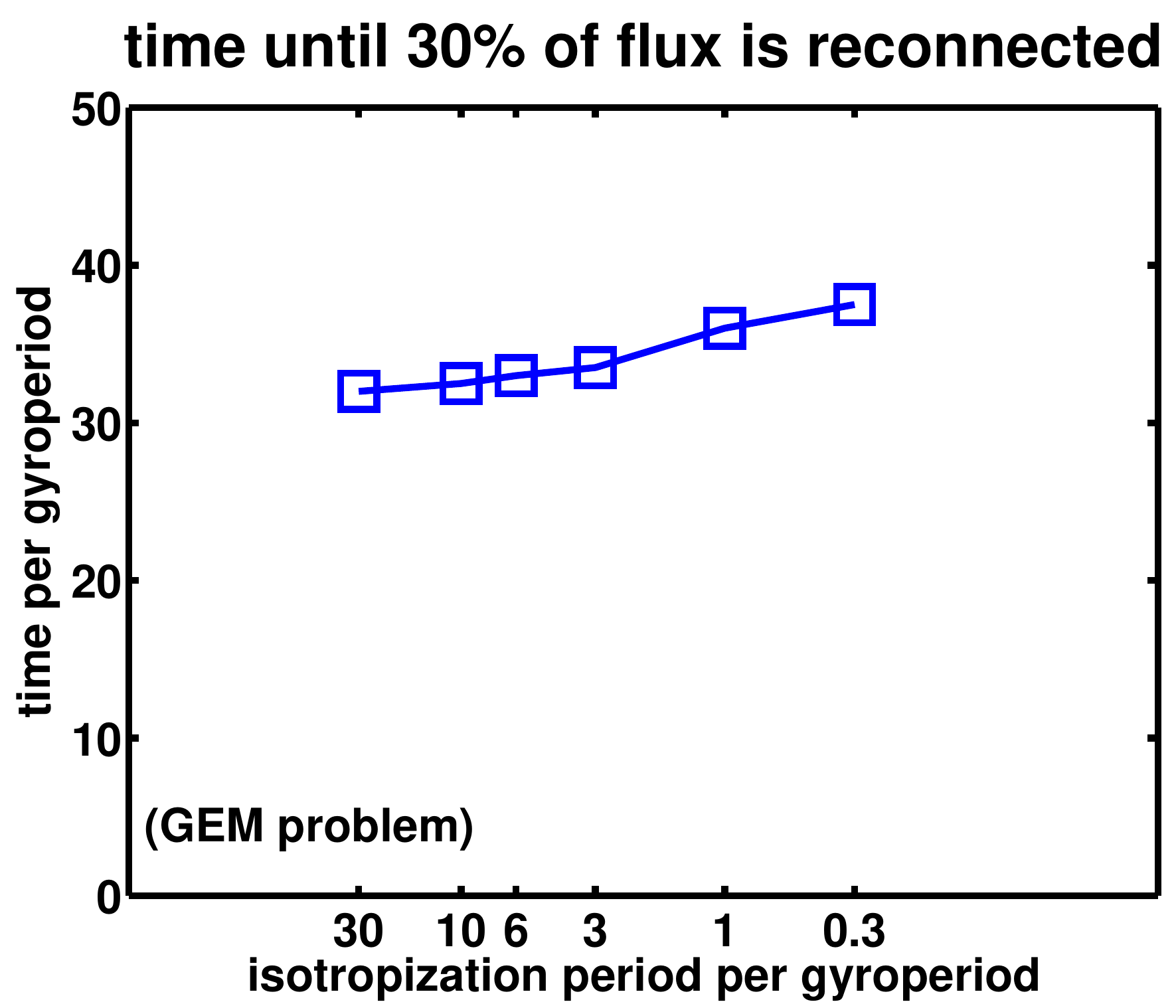}
   \end{center}
   \caption{Time until 30\% of flux is reconnected in 
     the full-scale symmetric pair plasma GEM problem.}
     For all simulations the peak rate of reconnection occurred
     when about 30\% of the original flux through the y-axis had reconnected.
  \fhrule
  \label{fig:recon_time_per_isoperiod_fullscale}
 \end{figure}

  Our nondimensionalized peak rate of reconnection was about 60\%
  of Bessho and Bhattacharjee's peak rate of about .23, and our
  time to peak reconnection was roughly twice their time of about
  18 (angular) gyroperiods
  (as seen in the red curve displayed in Figure $2$
  in either of  \cite{article:BeBh05,article:BeBh07}).
  See figures
  \ref{fig:recon_rate_per_isoperiod_fullscale} and
  \ref{fig:recon_time_per_isoperiod_fullscale}.

\begin{figure}
  \begin{center}
   \includegraphics[width=\linewidth]{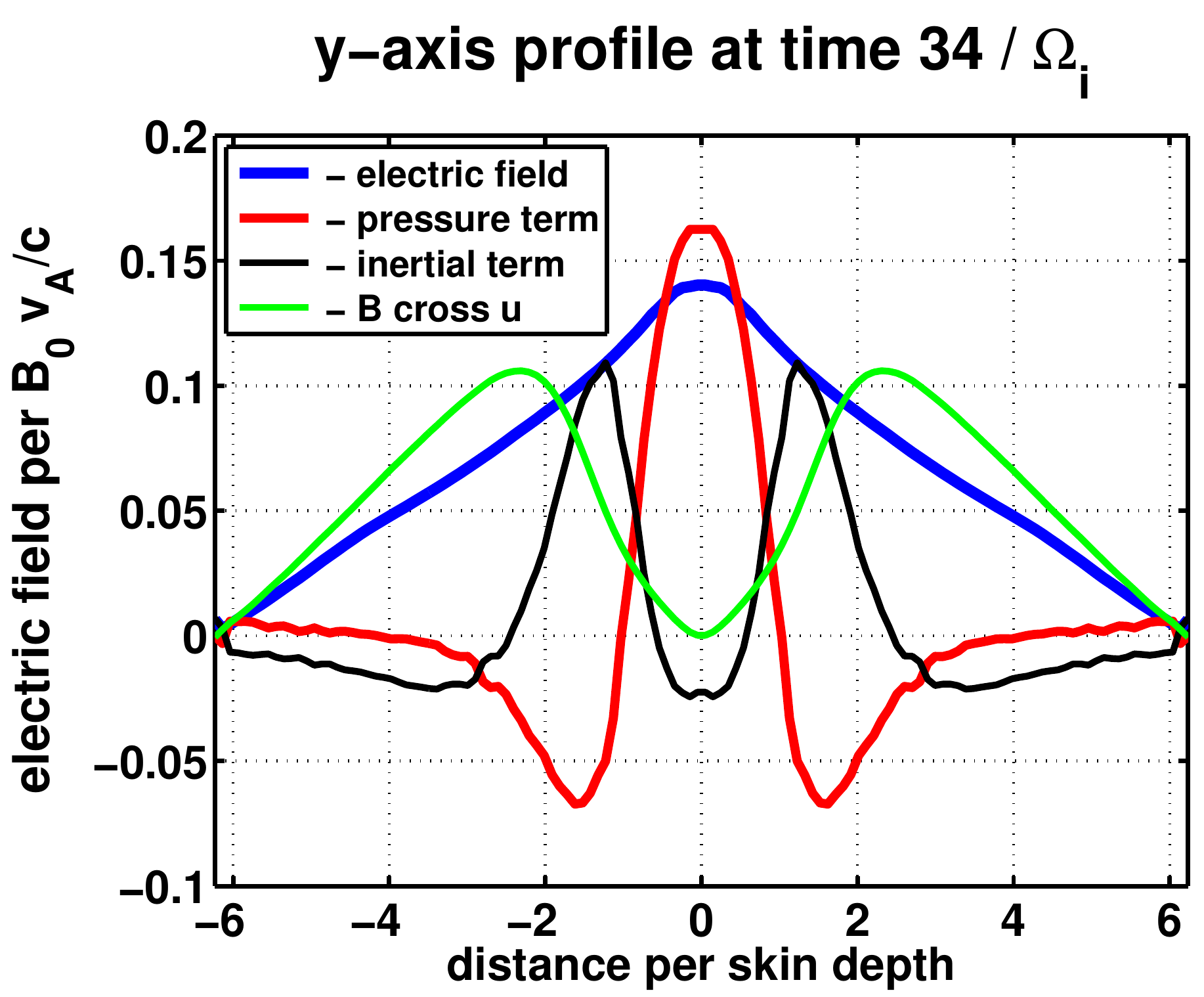}
  \end{center}
  \caption{
    Profile of ``Ohm's law'' terms at the peak reconnection time of
    34 (angular) gyroperiods in the full-scale symmetric pair plasma GEM problem.}
    Compare this figure with FIG. 5 in \cite{article:BeBh07}.
  \fhrule
  \label{fig:ohmslawYaxis}
 \end{figure}

\begin{figure}
  \begin{center}
   \includegraphics[width=\linewidth]{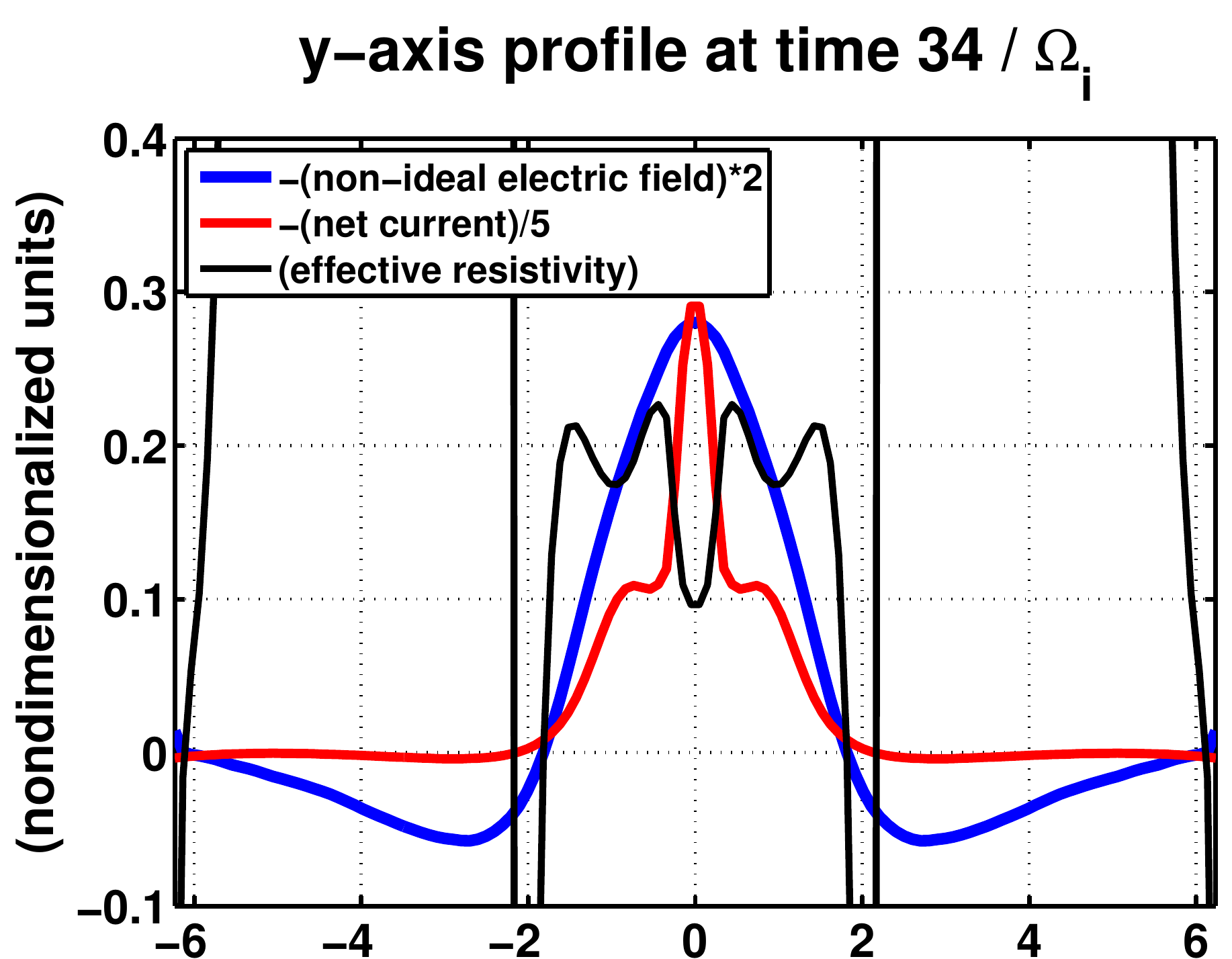}
  \end{center}
  \caption{Anomalous resistivity along the y-axis at the peak reconnection time of
    34 (angular) gyroperiods in the full-scale symmetric pair plasma GEM problem.}
    
    At the peak reconnection rate our effective anomalous
    resistivity at the X-point was about $0.1$ (in units of
    $\frac{B_0}{e n_0}$), in contrast to the value of roughly $0.19$
    reported in \cite{article:BeBh07} (compare this figure
    with FIG. 5 in \cite{article:BeBh07}).
  \fhrule
  \label{fig:anomalousResistivity}
 \end{figure}

  To compare with their results we plot Ohm's law terms on the
  inflow axis in figure \ref{fig:ohmslawYaxis}
  and plot calculated anomalous resistivity in
  \ref{fig:anomalousResistivity}.
  Our plots indicate a value of anomalous resistivity
  roughly half that reported in \cite{article:BeBh07}.

 \begin{figure}
  \begin{center}
   \includegraphics[width=\linewidth]{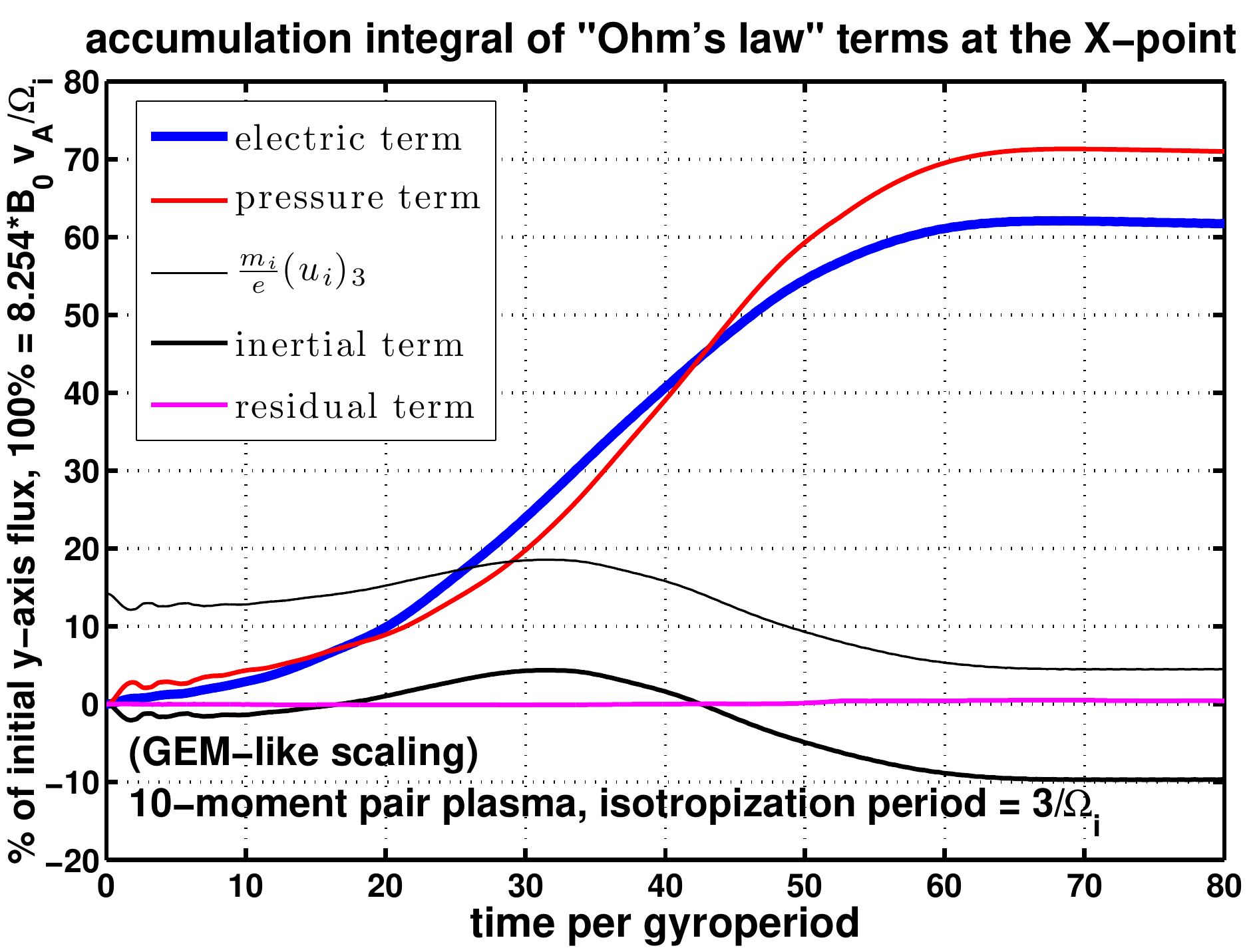}
  \end{center}
   \caption{Reconnected flux for moderate isotropization
     in the full-scale symmetric pair plasma GEM problem.}
    The pressure term is the dominant contribution, in agreement
    with PIC simulations.
    The time to 30\% reconnection is about 34 (angular) gyroperiods.
  \fhrule
 \end{figure}

 \begin{figure}
  \begin{center}
   \includegraphics[width=\linewidth]{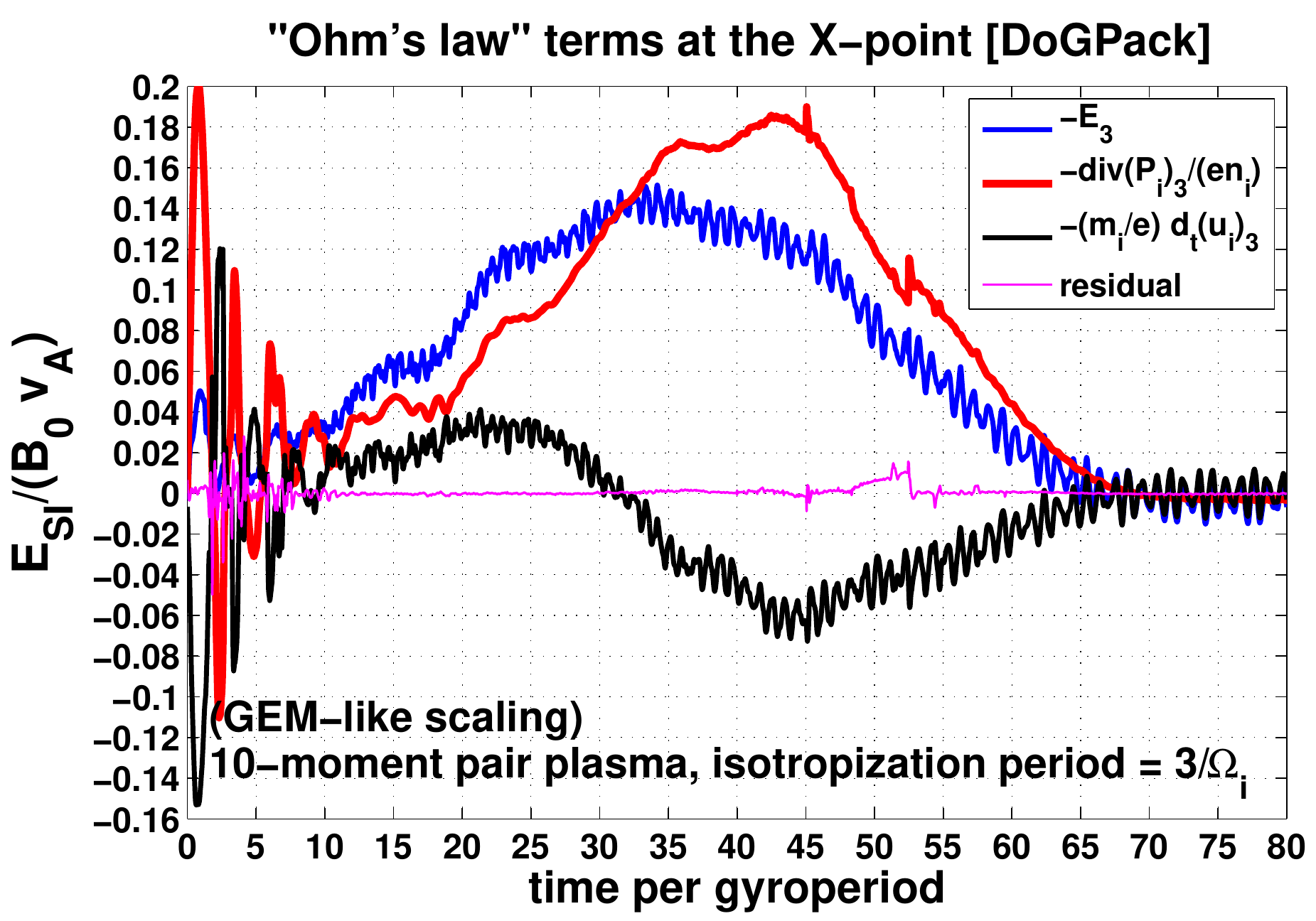}
  \end{center}
   \caption{Reconnection rate for moderate isotropization
     in the full-scale symmetric pair plasma GEM problem.}
    The blue curve in this figure is the rate of reconnection and
    corresponds to the red curve in Figure $2$ of  \cite{article:BeBh07}.
    Our peak reconnection rate is about $.14$.
  \fhrule
 \end{figure}
 
  We display profiles and plots at the peak reconnection time of 34
  gyroperiods.
  The shape of our profiles is similar to those of 
  \cite{article:BeBh07}, though our rates are generally smaller.
 
 \begin{figure}
  \begin{center}
   \includegraphics[width=\linewidth]{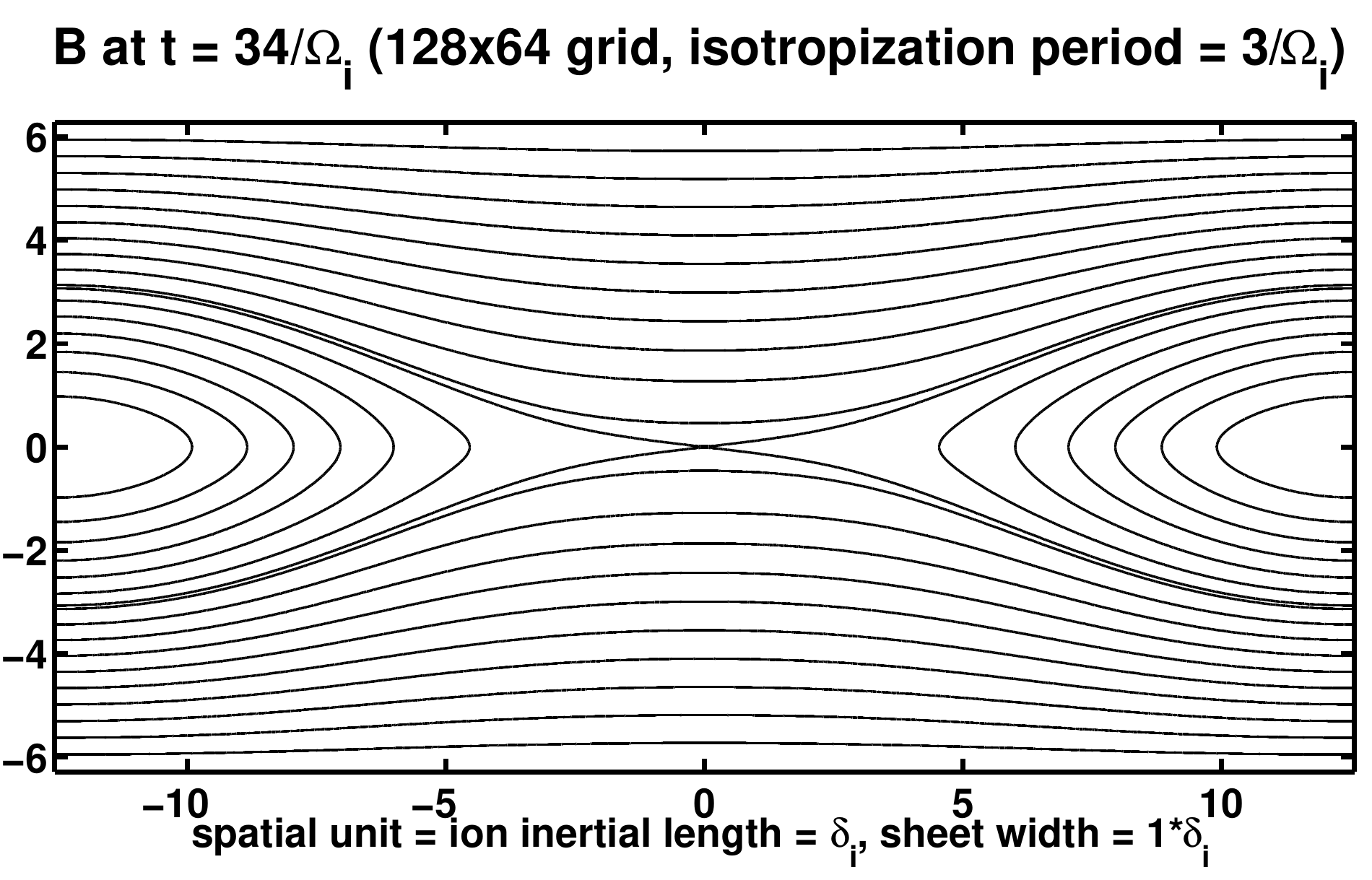}
  \end{center}
  \begin{center}
   \includegraphics[width=\linewidth]{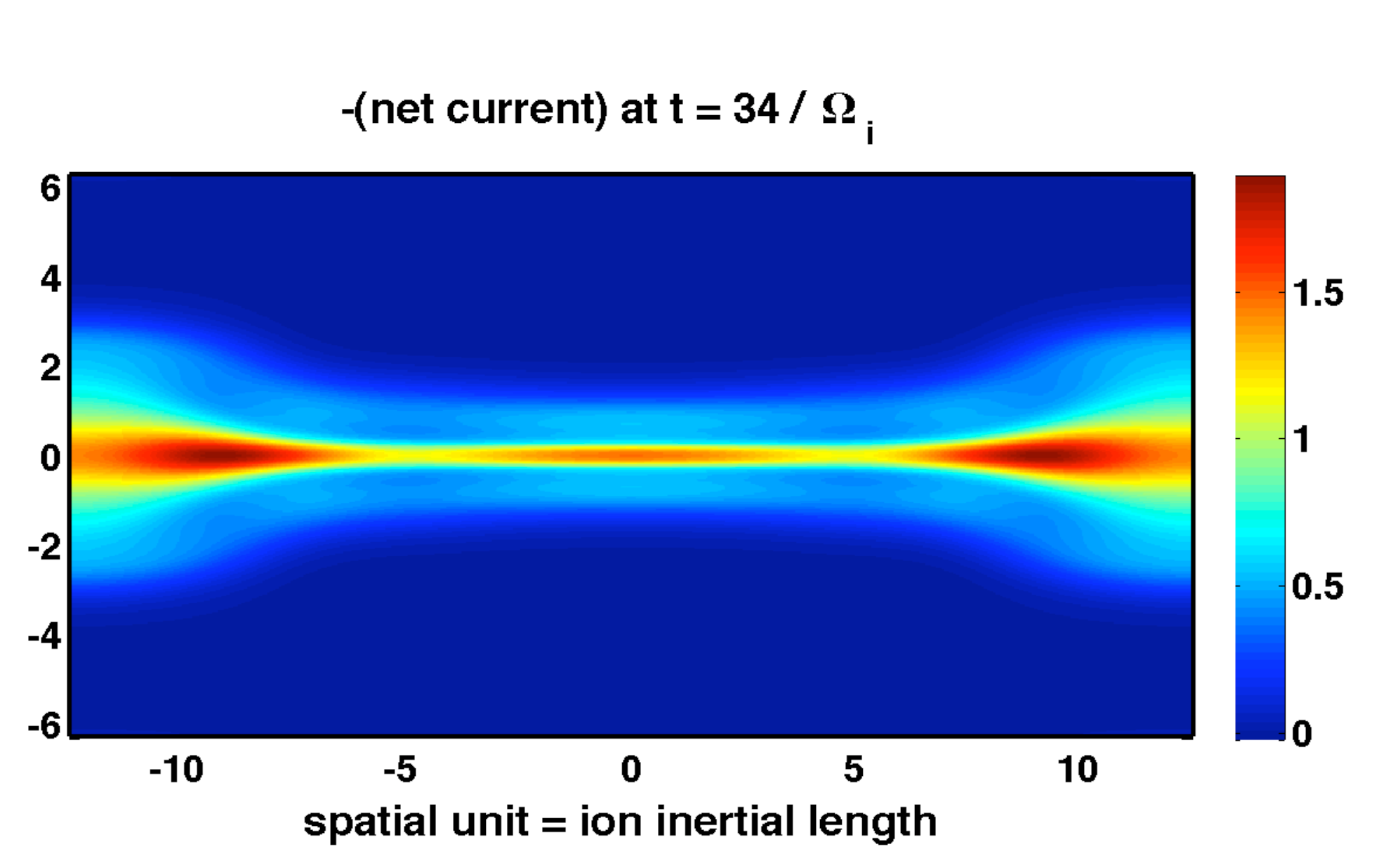}
  \end{center}
  \caption{
    Magnetic field lines and current at the time of
    peak reconnection rate.
  }
  \fhrule
 \end{figure}

\section{Conclusion}


Our simulations of the GEM problem
indicate that for antiparallel reconnection
in pair plasma the ten-moment two-fluid model with
isotropization admits a rate of reconnection 
a bit more than half the rate seen in particle simulations.

This raises the question of how to modify the fluid model we
used so that we can get better agreement.
An obvious way to get agreement is to use an anomalous
resistivity.  In the context of a collisionless 
simulation it would seem more physical to impose an
anomalous viscosity.
Our rate of reconnection was insensitive to the rate of
isotropization. This might suggest use of a different (nonzero)
nonzero heat flux, but heat flux should not affect the rate of
reconnection much until significant temperature
gradients have had time to develop. 

In contrast to these results for pair plasmas,
we find in the next chapter that for the original GEM problem
a two-fluid ten-moment model with
relaxation toward isotropy gives reconnection rates
that agree well with kinetic simulations \cite{article:JoRo10b}.
Perhaps the improved agreement can be attributed to the presence of
the Hall effect as the primary driver of fast reconnection
in hydrogen plasmas.

%% file: chap6.tex
\chapter{GEM hydrogen plasma simulations}


This chapter reports simulations of GEM magnetic reconnection challenge problem
described in chapter \ref{GEMchapter}
(for the case of hydrogen plasma with $\frac{\mi}{\me}=25$).
The main results of this chapter
were also reported in \cite{article:JoRo10b}.

We again used the adiabatic 10-moment two-fluid-Maxwell model
with pressure tensor isotropization to simulate the problem,
as well as the hyperbolic (adiabatic inviscid) 5-moment two-fluid Maxwell model.
We compared our results with published Vlasov and PIC simulations
at the time of peak reconnection.
The magnetic field generally showed good agreement (see figure 
\ref{fig:magneticField}).
We obtained good agreement for the
time until one unit of flux was reconnected for both the 5-moment model
and the 10-moment isotropizing model (see figure \label{fig:xpointRecon}).
For the 10-moment isotropizing model we obtained
qualitatively good agreement for plots of pressure tensor components at the
time of peak reconnection, as shown in figures \ref{fig:offDiagPressure}
and \ref{fig:diagPressure}.

As noted in table \ref{table:recon}, it is standard to compare
results at the time when one unit of flux has been reconnected,
\emph{including the .2 units of reconnected flux due to the
initial perturbation.}  I unfortunately have made my plots at
the point in time when one unit of flux has been reconnected
\emph{not including the .2 units of reconnected flux due to the
the GEM problem}, which puts my plots about one unit of time later
than the time I should be looking at to compare against the
simulations of others.

 \begin{table}[ht]
   \hfrule
    \begin{center}
  \begin{tabular}{l|l|l}
    model & source & time when 16\% flux reconnected
 \\ \hline
    Vlasov & \cite{article:SmGr06}=[ScGr06]
      & $t = 17.7/\Omega_i$:
 \\ PIC & \cite{pritchett01}=[Pritchett01]
      & $t = 15.7/\Omega_i$:
 \\ 10-moment & \cite{article:JoRo10b}=[JoRo10]
      & $t = [16.2,17.2]/\Omega_i$:
 \\ 5-moment & \cite{article:JoRo10b}=[JoRo10]
      & $t = [12.0,12.9]/\Omega_i$:
 \\ 5-moment & \cite{LoHaSh11}=[LoHaSh11]
      & $t = [15,16]/\Omega_i$:
 \\ 10-5-moment & \cite{article:Hakim08}=[Hakim08]
      & $t = 17.6/\Omega_i$:
  \end{tabular}
    \end{center}
\caption{Comparison of time required to reconnect one unit of flux for kinetic
  and two-fluid simulations of the GEM problem.}
  The 10-5-moment simulation modeled ions as a 10-moment gas
  and electrons as a 5-moment gas \cite{article:Hakim08}
  and used hyperbolic (truncation) closures for both species.
  The 5-moment simulations of \cite{article:JoRo10b}
  and of \cite{LoHaSh11} were hyperbolic in both species.

  The rough range given for \cite{LoHaSh11} reflects
  that this paper reports a number of simulations in its figures
  16 and 17.

  The ranges given for $\cite{article:JoRo10b}$
  (the results reported here) reflect that
  \bfit{there is an inconsistency in the time
  when I have chosen to report my results, both
  here and in $\cite{article:JoRo10b}$}.
  (I am unfortunately out of time to redo
  the plots and calculations before the deadline to
  deposit this work.)
  For the GEM problem it is standard to compare
  results at the point in time when one unit of flux
  has been reconnected, \emph{including the initial
  perturbation}.  I neglected to include the initial
  perturbation, and therefore my times were really
  the time until 1.2 units of flux were reconnected
  in the sense understood by the GEM standard.
  Therefore I have listed both times: the time until
  1.0 total units of flux have been reconnected followed
  by the time until 1.2 units have been reconnected.
  \label{table:recon}
  \fhrule
\end{table}
  
\begin{figure}
 \hfrule
 \begin{gather*}
  \begin{matrix}
      \mbox{\parbox{.2\textwidth}{ 
        The ten-moment model attained 16\% flux reconnected at about
        $t = 18/\Omega_i$: }}
     & 
       \begin{matrix}
        \mbox{\includegraphics[width=.75\linewidth]{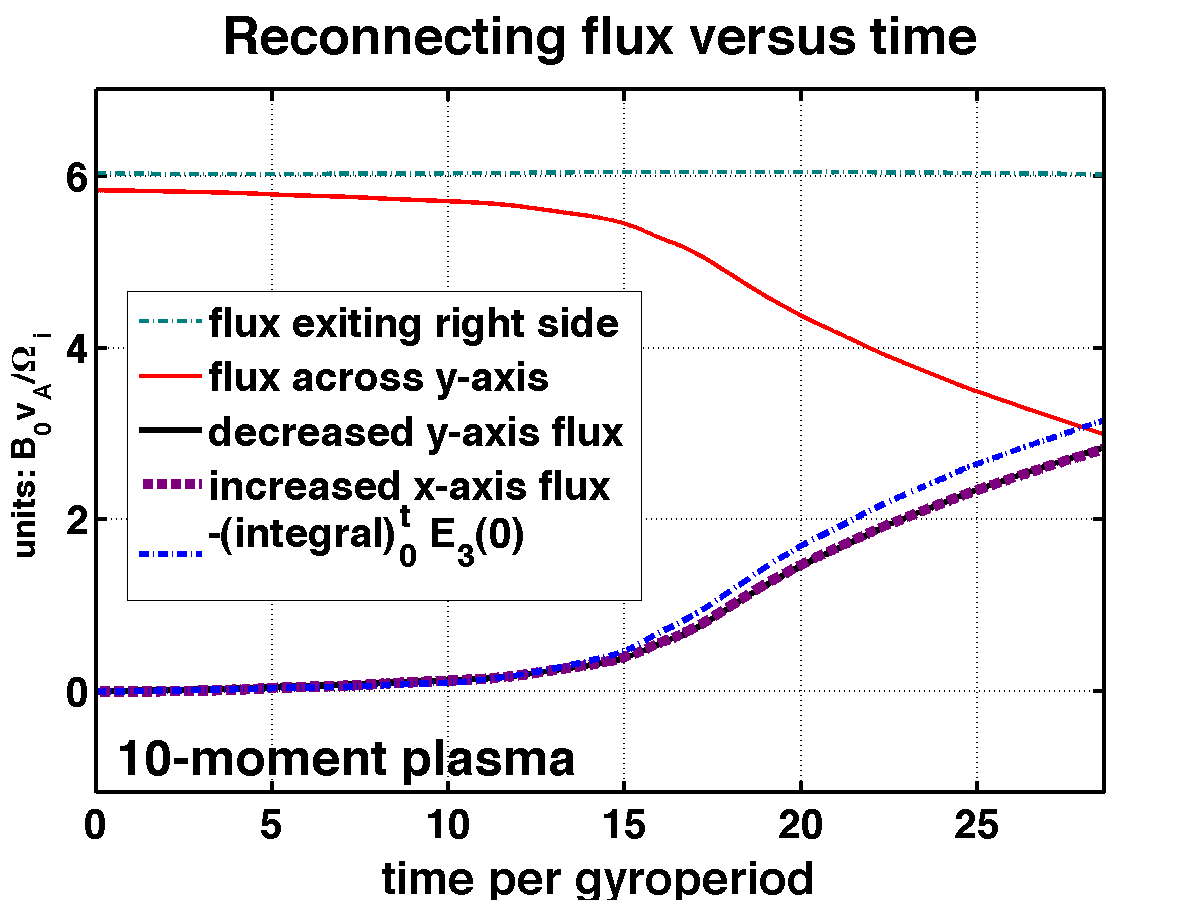}}
       \end{matrix}
     \\
      \mbox{\parbox{.2\textwidth}{ 
         The five-moment model attained 16\% flux reconnected at about
         $t = 13.5/\Omega_i$:}}
     & 
       \begin{matrix}
        \mbox{\includegraphics[width=.75\linewidth]{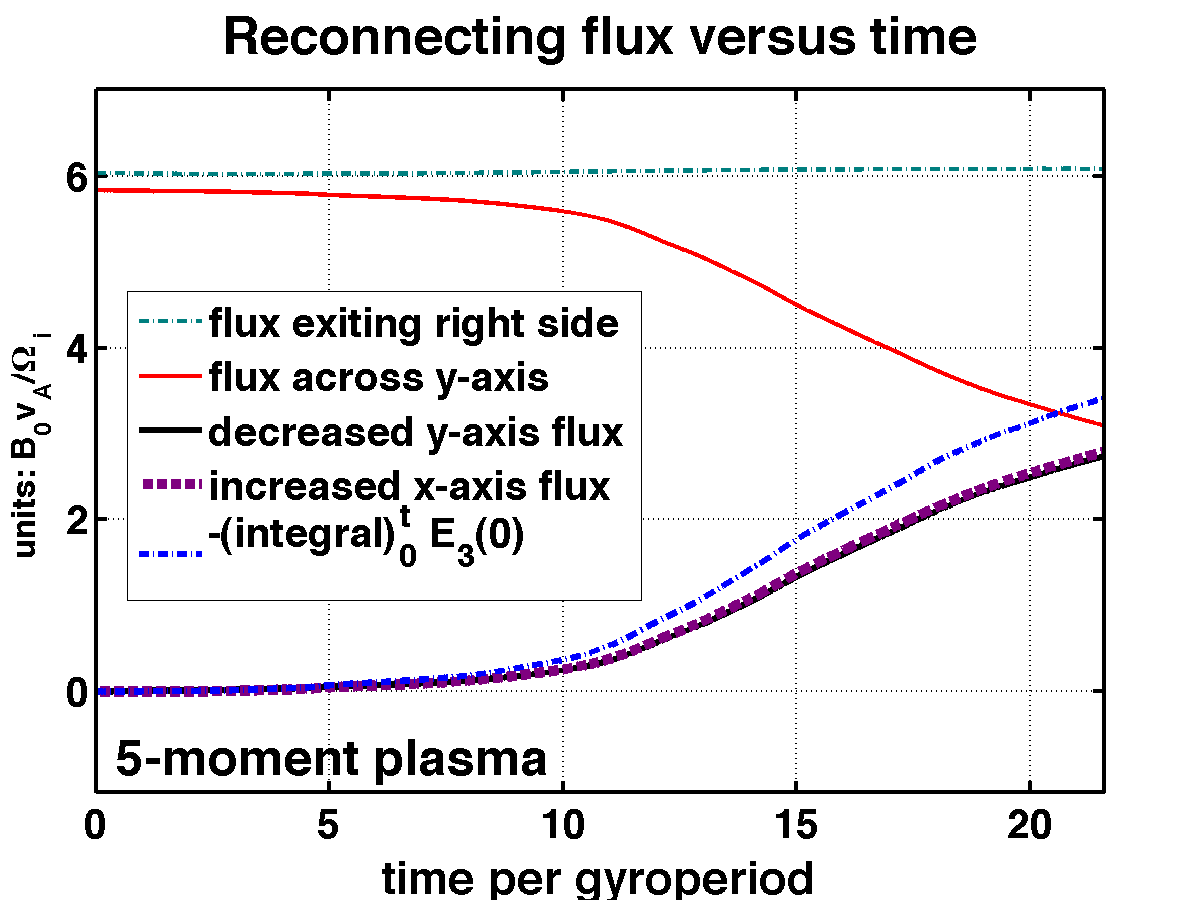}}
       \end{matrix}
  \end{matrix}
 \end{gather*}
 \caption{GEM simulations: reconnecting flux for 10- and 5-moment models.}
 Each simulation crashed at the end time in its plot of reconnecting
 flux versus time.
 \label{fig:xpointRecon}
 \fhrule
\end{figure}

\begin{figure}
 \hfrule
 \begin{tabular}{c c}
  \begin{tabular}{c}
 \\ \includegraphics[width=.5\linewidth]{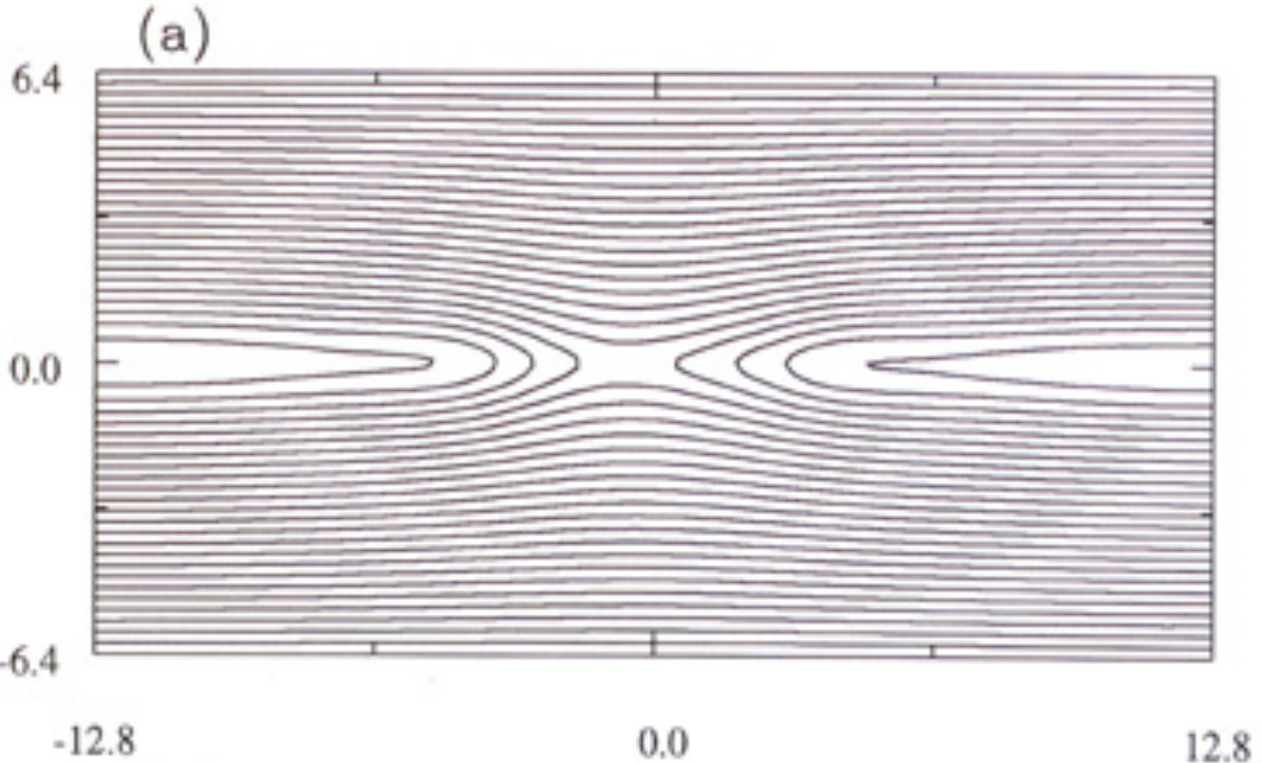}
 \\ Magnetic field lines for PIC
 \\ at $\Omega_i t = 15.7$ [Pritchett01]=\cite{pritchett01}
 \vspace{4ex}
 \\ \includegraphics[width=.5\linewidth]{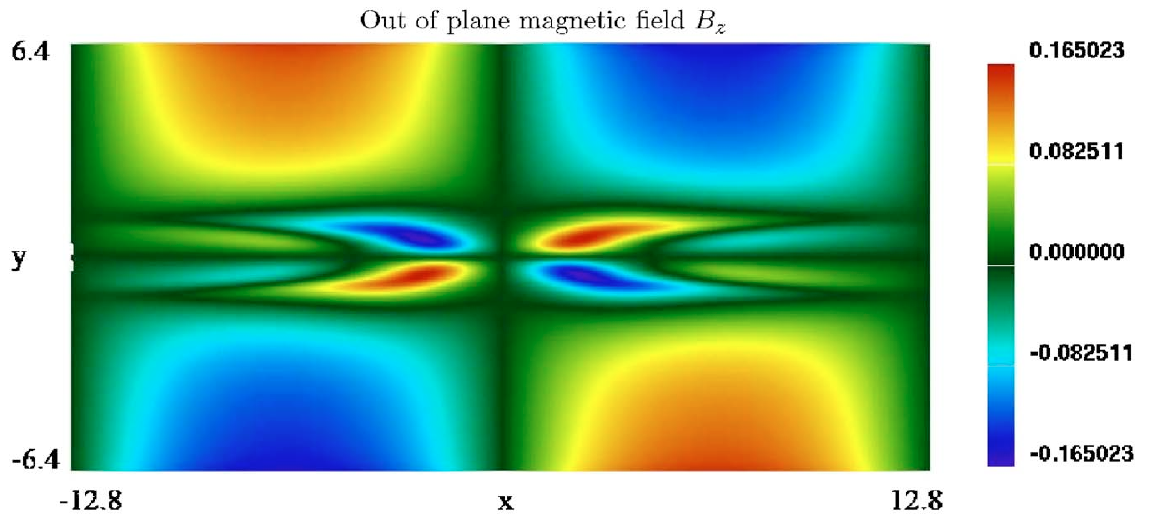}
 \\ Magnetic field for Vlasov
 \\ at $\Omega_i t = 17.7$ [ScGr06]=\cite{article:SmGr06}
  \end{tabular}
  &
  \begin{tabular}{c}
    \includegraphics[width=.5\linewidth]{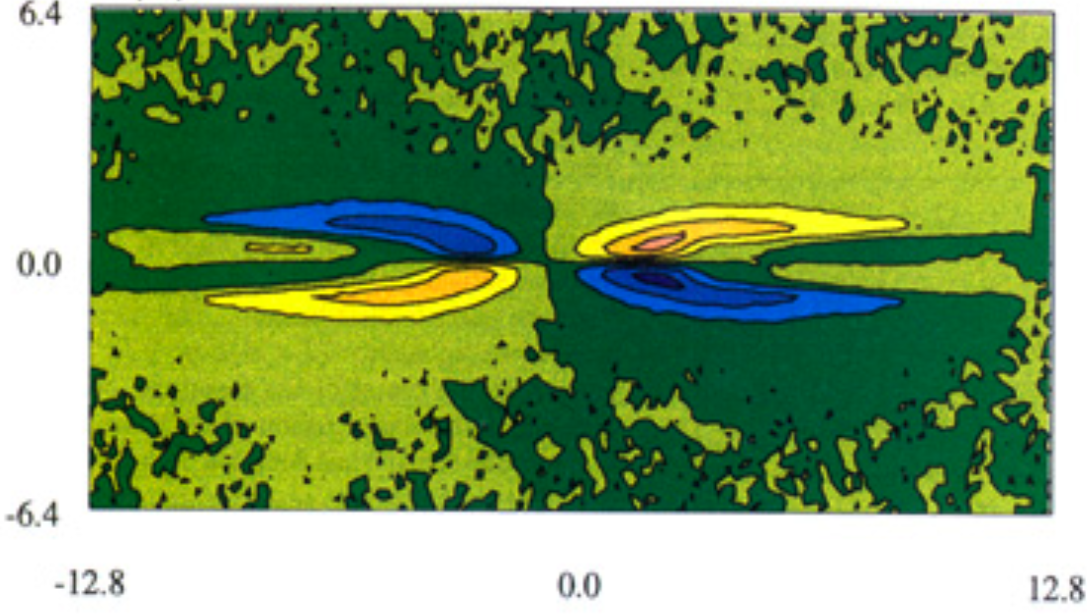}
 \\ magnetic field of [Pritchett01]=\cite{pritchett01}
 \vspace{5ex}
 \\ \includegraphics[width=.5\linewidth]{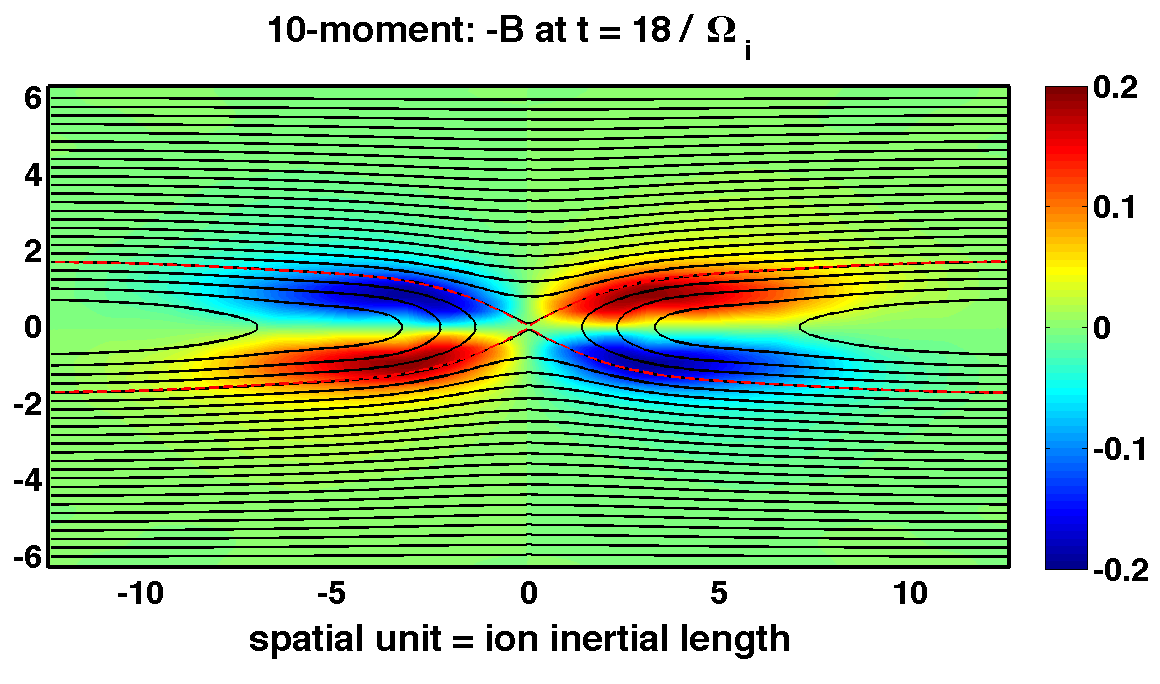}
 \\ \includegraphics[width=.5\linewidth]{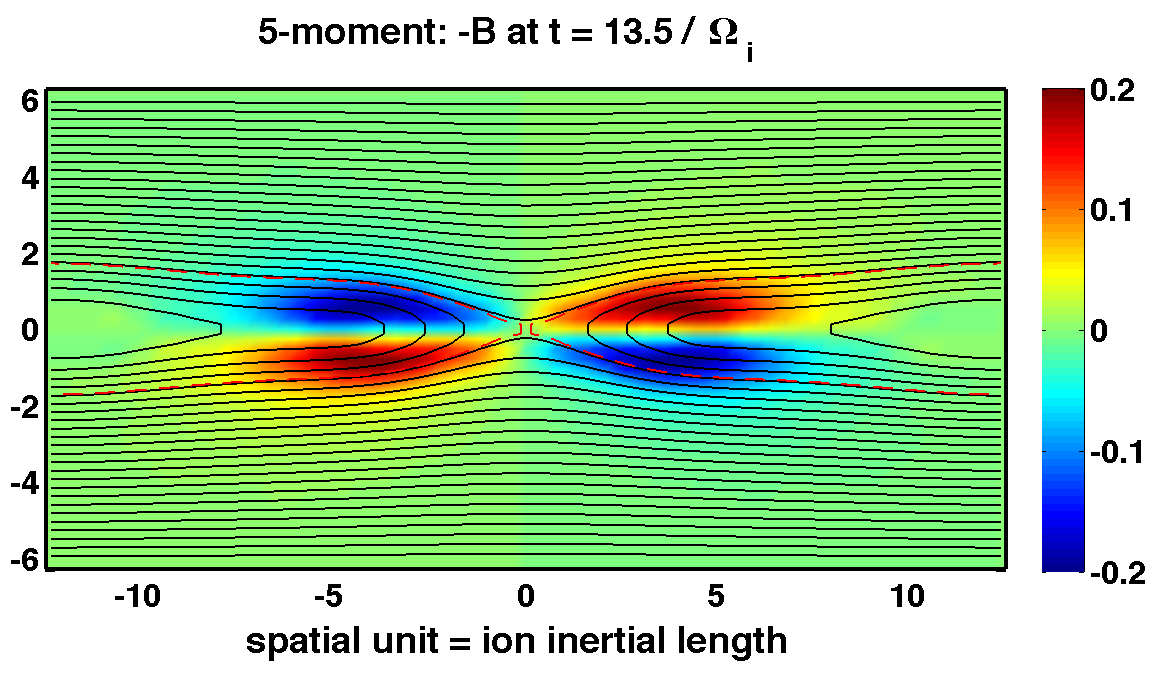}
  \end{tabular}
\end{tabular}
\caption{Magnetic field at 16\% reconnected}
\label{fig:magneticField}
\fhrule
\end{figure}

\begin{figure}
 \hfrule
 \begin{gather*}
  \begin{matrix}
    \begin{matrix}
     \mbox{\includegraphics[width=.47\textwidth]{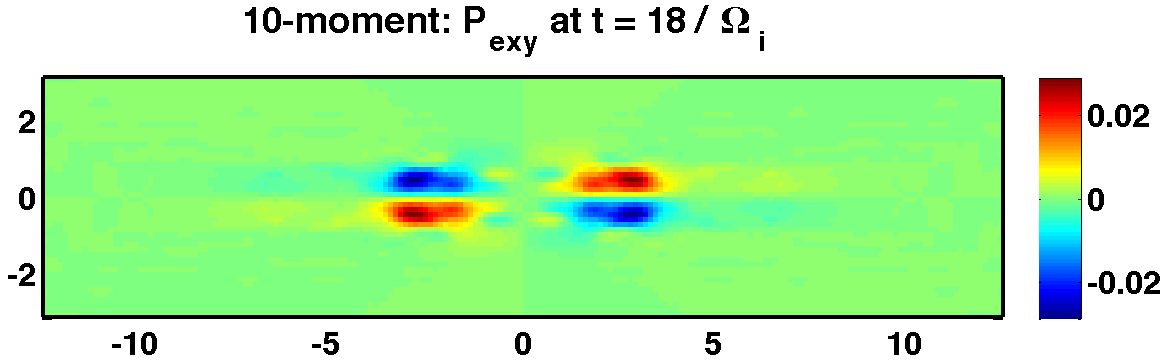}}
  \\ \mbox{\includegraphics[width=.47\textwidth]{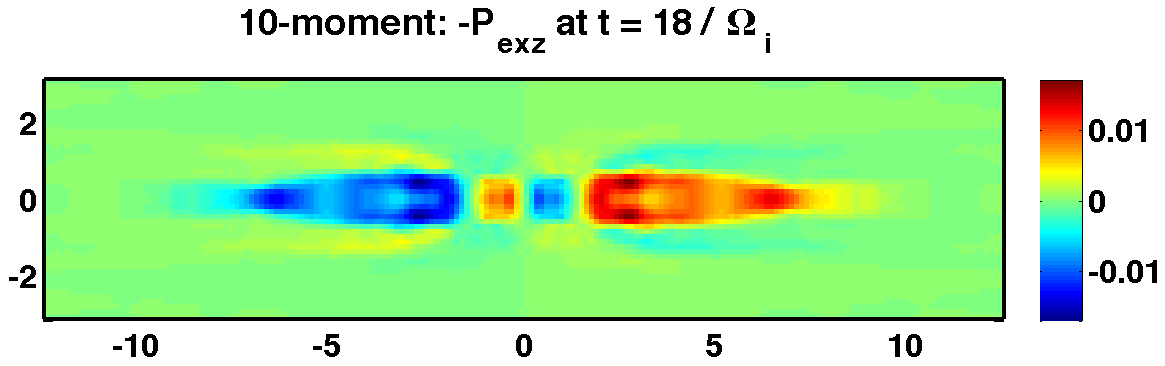}}
  \\ \mbox{\includegraphics[width=.47\textwidth]{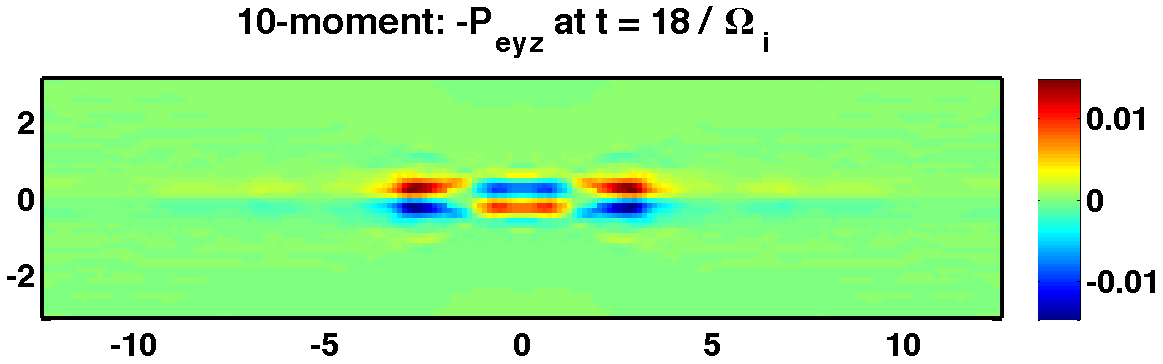}}
    \end{matrix}
   &
    \begin{matrix}
      \mbox{ \includegraphics[width=.53\linewidth]{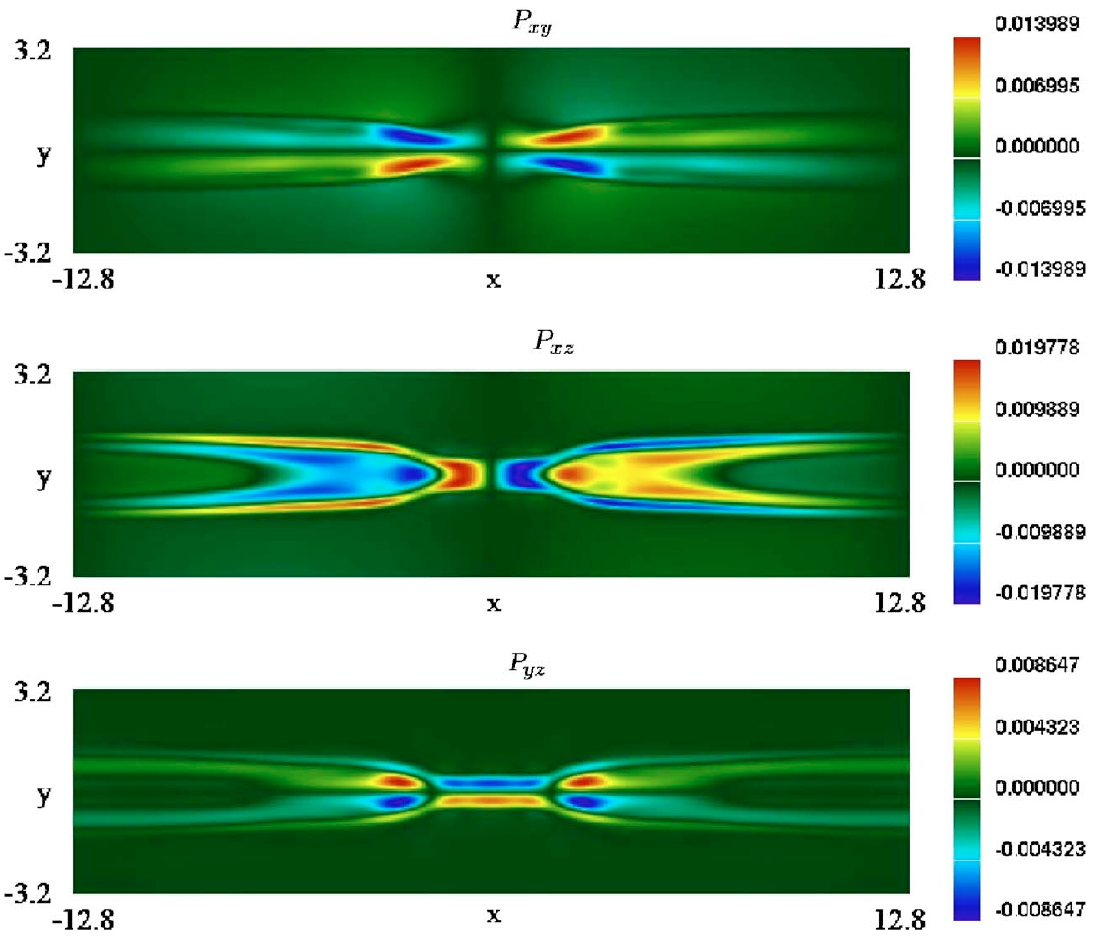} }
    \end{matrix}
   \\
    \begin{matrix}
      \mbox{\parbox{.4\textwidth}{\small
       Off-diagonal components of the electron pressure tensor
         for 10-moment simulation at $\Omega_i t = 18$}}
    \end{matrix}
    &
    \begin{matrix}
      \mbox{\parbox{.4\textwidth}{\small
       Off-diagonal components of the electron pressure tensor for Vlasov
       simulation at $\Omega_i t = 17.7$ [ScGr06]=\cite{article:SmGr06}}}
    \end{matrix}
  \end{matrix}
 \end{gather*}
\caption{Off-diagonal components of electron pressure tensor}
\label{fig:offDiagPressure}
\fhrule
\end{figure}

\begin{figure}
 \hfrule
 \begin{gather*}
  \begin{matrix}
    \begin{matrix}
     \mbox{\includegraphics[width=0.47\textwidth]{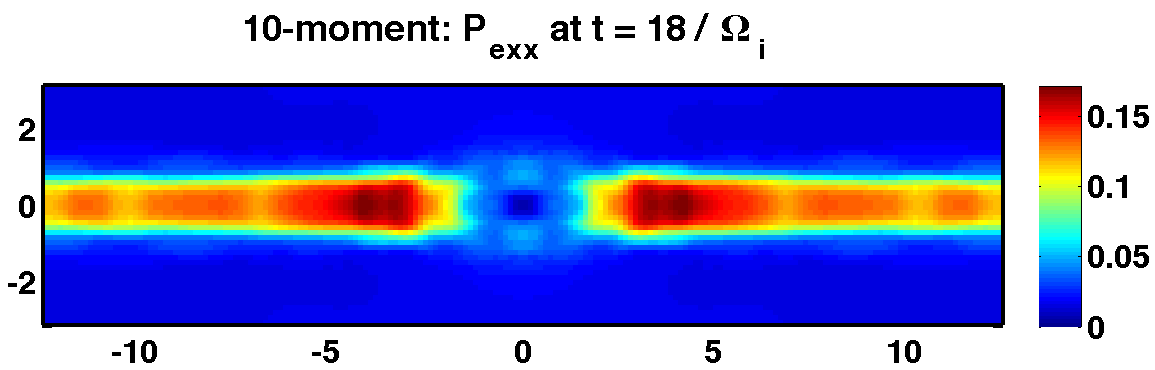}}
  \\ \mbox{\includegraphics[width=0.47\textwidth]{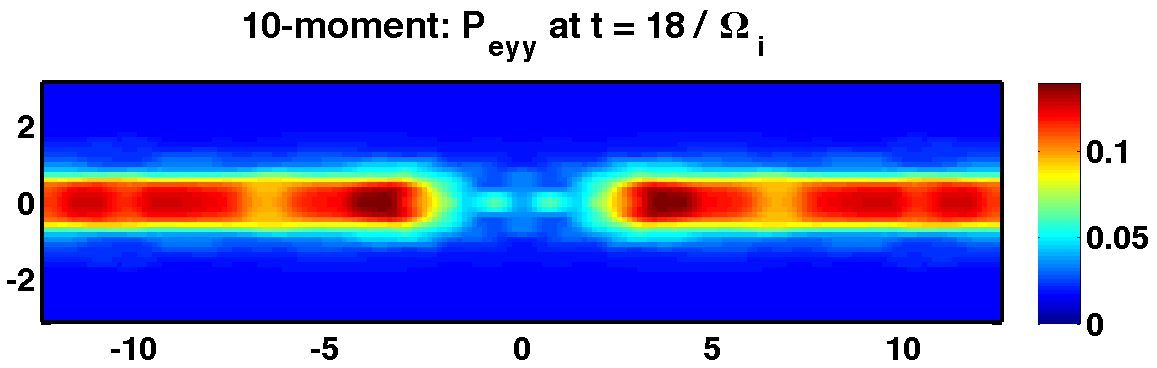}}
  \\ \mbox{\includegraphics[width=0.47\textwidth]{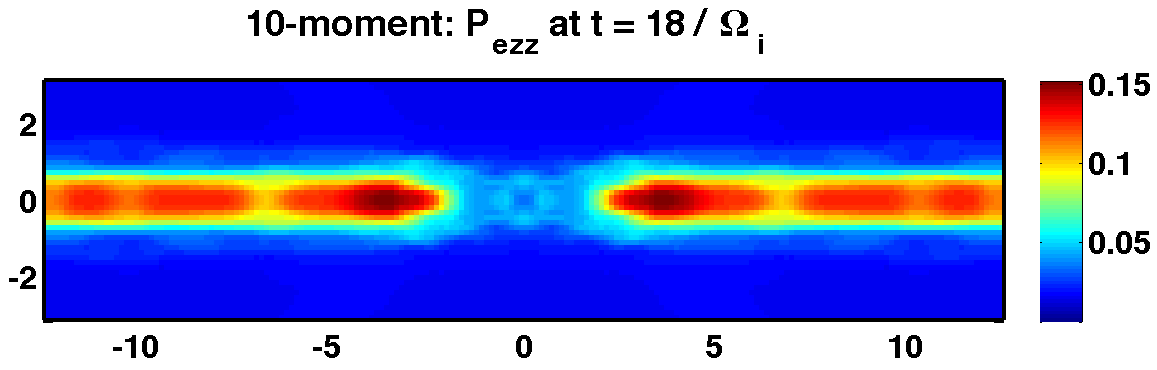}}
    \end{matrix}
   &
    \begin{matrix}
      \mbox{\includegraphics[width=0.53\textwidth]{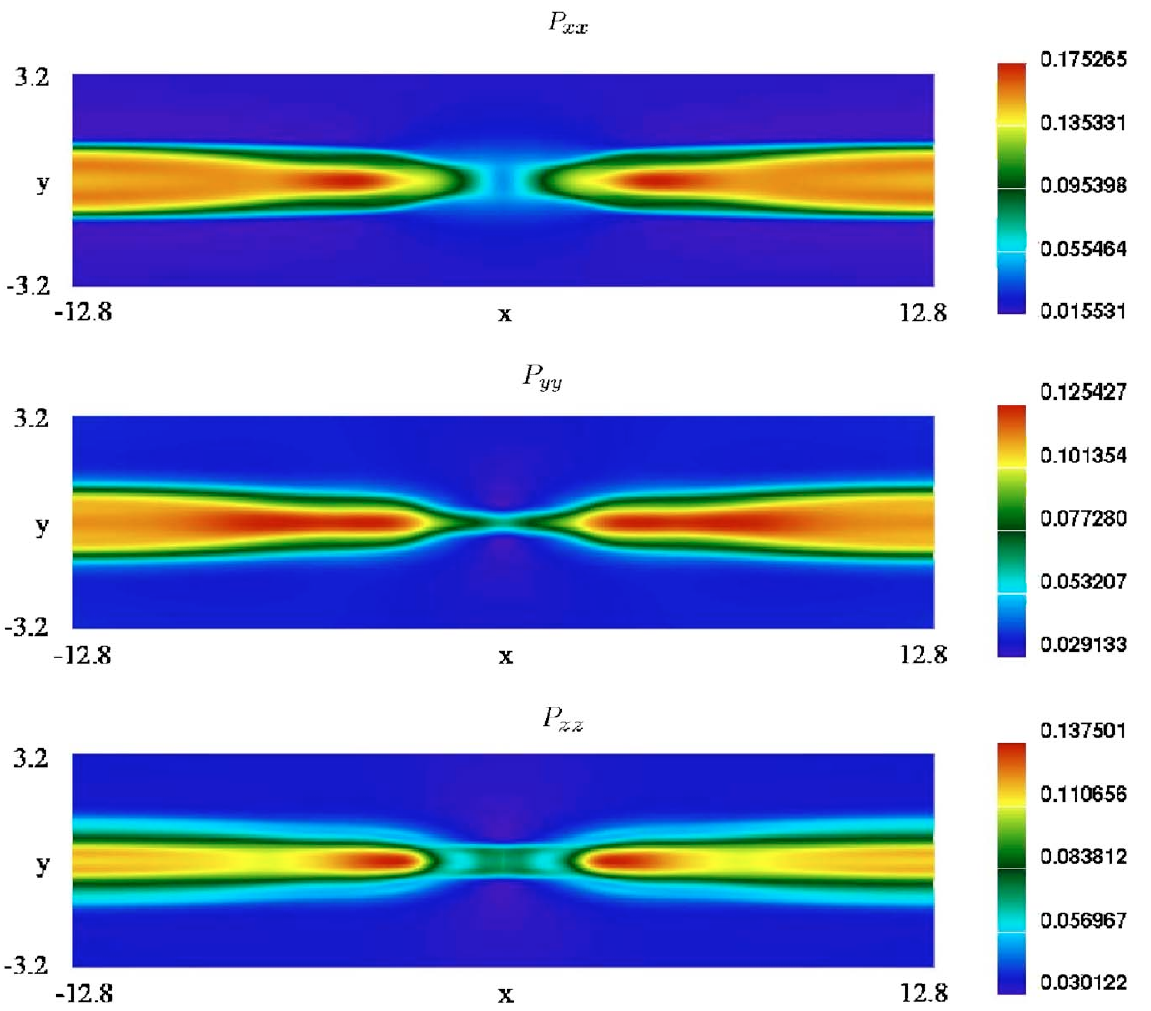}}
    \end{matrix}
   \\
    \begin{matrix}
      \mbox{\parbox{.4\textwidth}{\small
       Diagonal components of the electron pressure tensor
         for 10-moment simulation at $\Omega_i t = 18$}}
    \end{matrix}
    &
    \begin{matrix}
      \mbox{\parbox{.4\textwidth}{\small
       Diagonal components of the electron pressure tensor for Vlasov
        simulation at $\Omega_i t = 17.7$ [ScGr06]=\cite{article:SmGr06}}}
    \end{matrix}
  \end{matrix}
 \end{gather*}
\caption{Diagonal components of electron pressure tensor}
\label{fig:diagPressure}
\fhrule
\end{figure}

\clearpage

For a coarse mesh we can extend the time duration of the
simulation to get a solution that appears quite regular. But for a
refined mesh, after peak reconnection our simulations would
eventually crash due to negative pressures or even densities,
typically near the X-point. This happened chiefly in the zero
guide field case. By implementing positivity limiters we were
able to prevent the simulation from crashing, but the solution exhibits
secondary instabilities past the point where positivity limiting
becomes necessary.

These difficulties prompted a consideration of steady-state
reconnection. The GEM problem is a 2D rotationally symmetric
problem, and as shown in section \ref{steadyReconNeedsHeatflux},
nonsingular 2D steady reconnection is not possible in an
adiabatic model for a rotationally symmetric problem.

The implications are that steadily driven reconnection 
of an adiabatic model cannot go to steady state and
will therefore either develop a singularity or will
exhibit intermittent reconnection.  (On a sufficiently
small scale we expect a singularity, whereas for a larger
scale, especially if symmetry is not enforced, we
expect intermittent reconnection.)

It is reasonable to assume that near the time when the rate of
reconnection peaks the evolution of the solution near the X-point
may be approximated by steadily driven reconnection. If so, then
since we are in the case of antiparallel reconnection, along the
$y$ axis near the origin we expect the solution to approach
a singular steady state of the form
\def\Txx{\TT_{xx}}
\begin{gather*}
  \Txx = C y^\lambda,
\end{gather*}
as seen in \eqref{TxxOrigin} and \eqref{Worigin}.

\begin{figure}
  \begin{tabular}{c c}
    \includegraphics[width=0.50\textwidth]{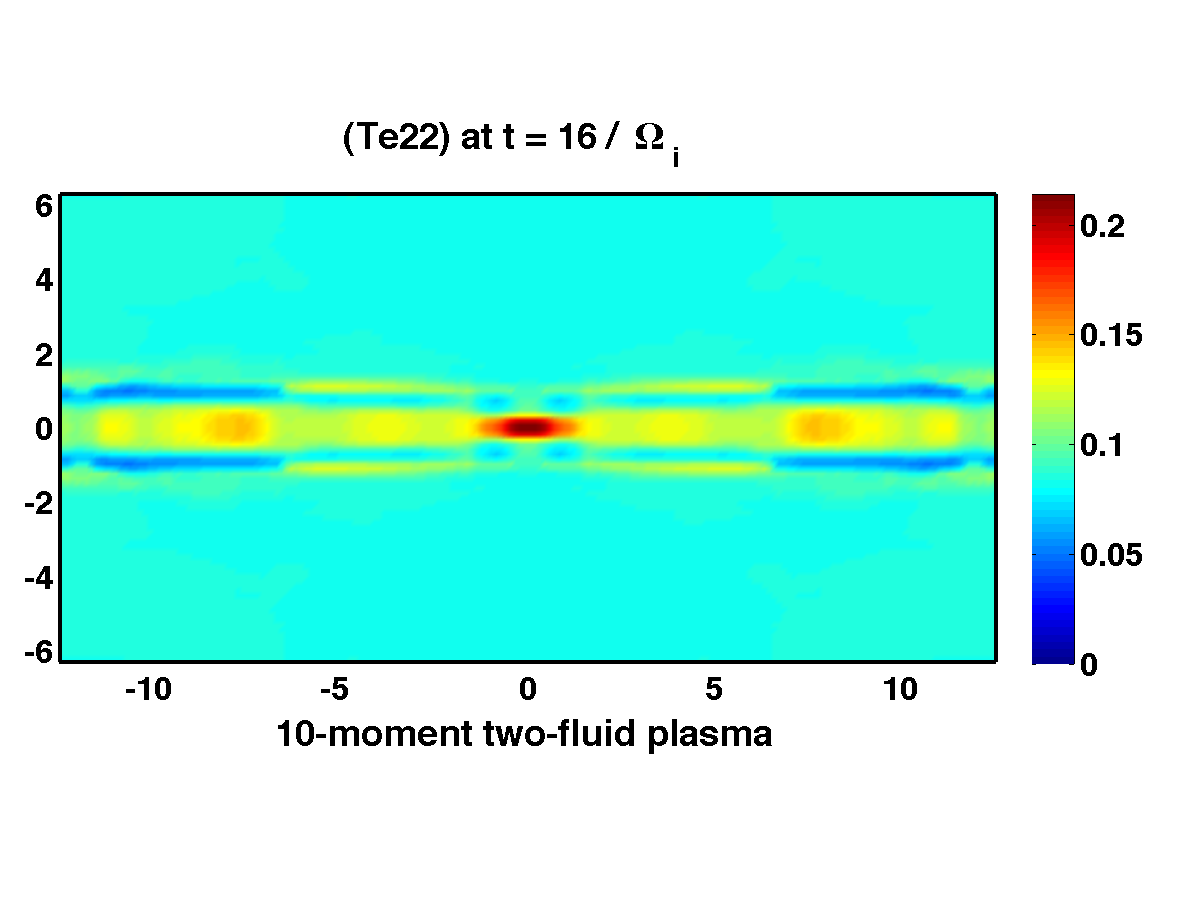}
   &\includegraphics[width=0.50\textwidth]{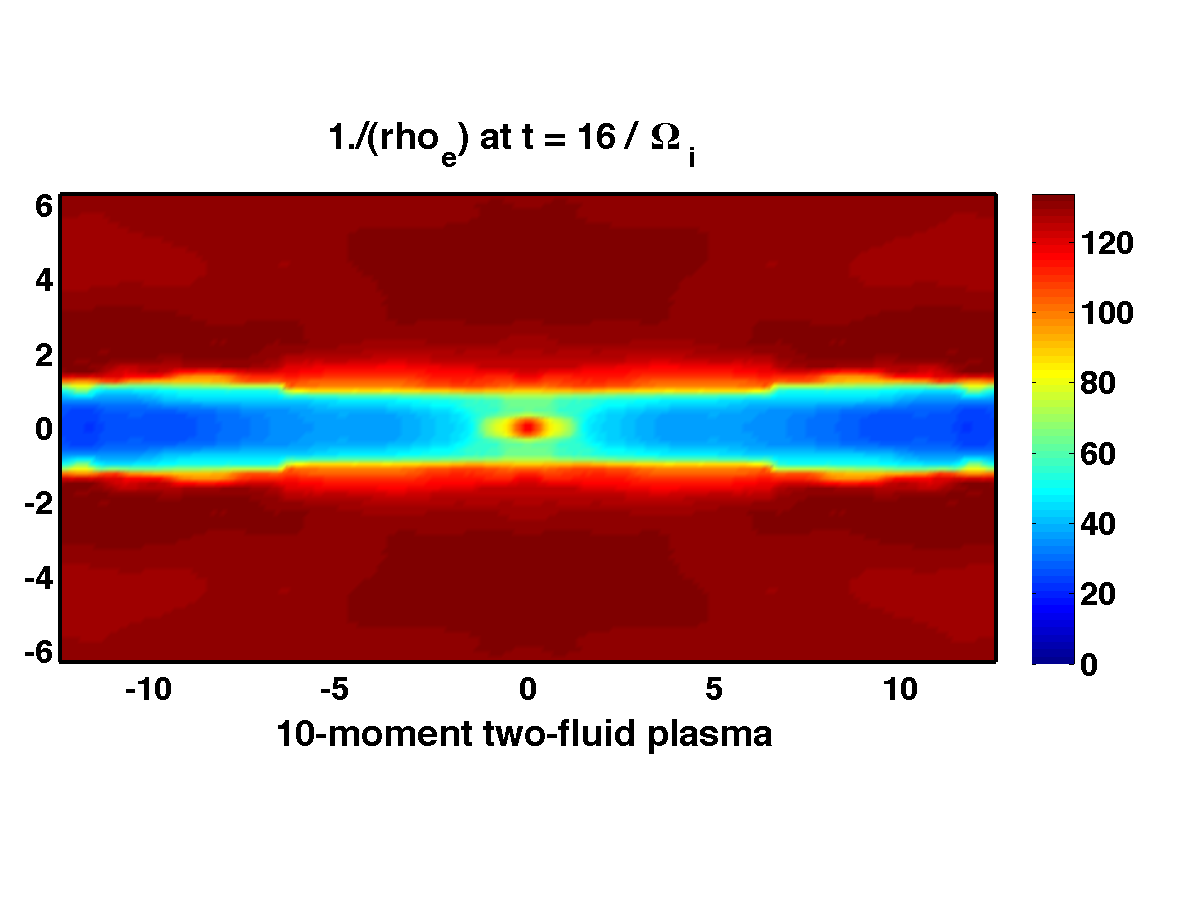}
   \\
    \vspace{-5ex}
    \\
    \includegraphics[width=0.50\textwidth]{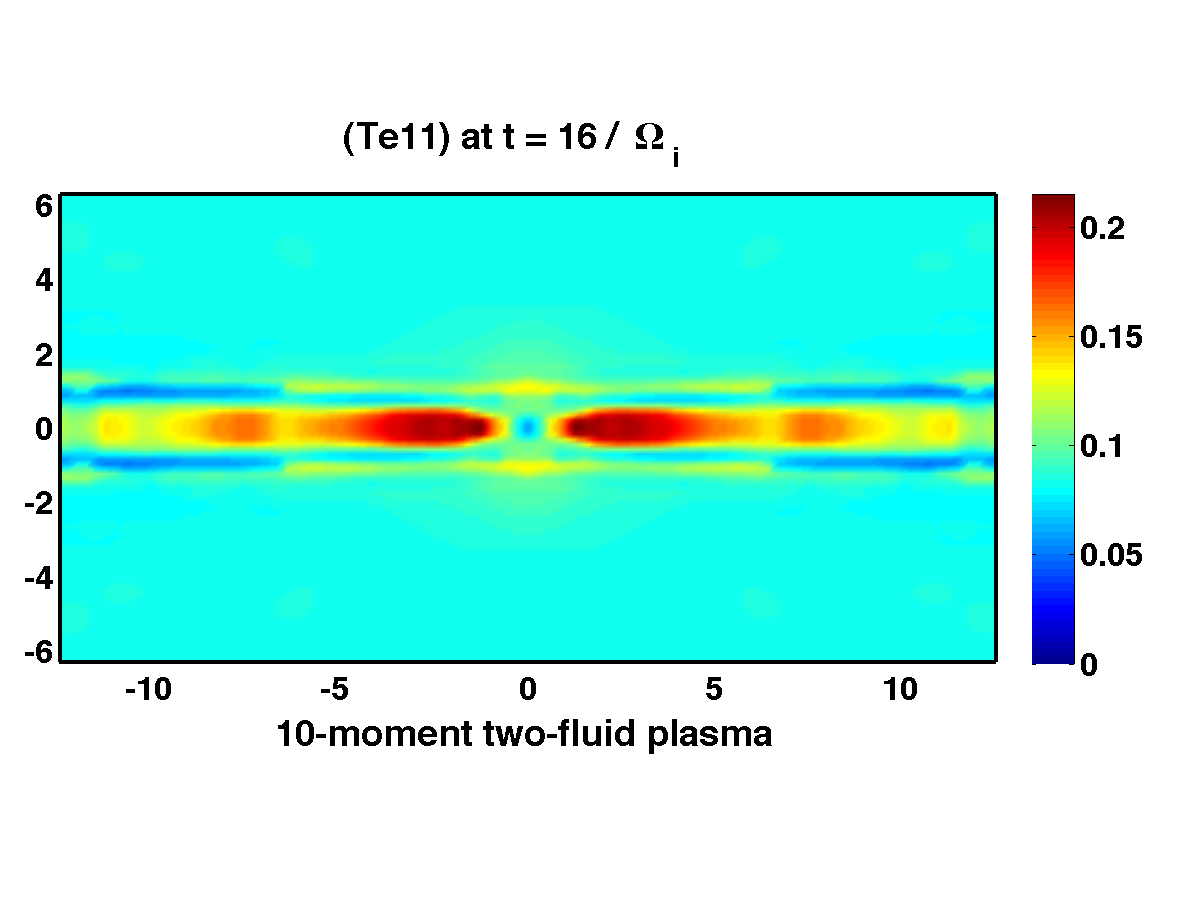}
   &\includegraphics[width=0.50\textwidth]{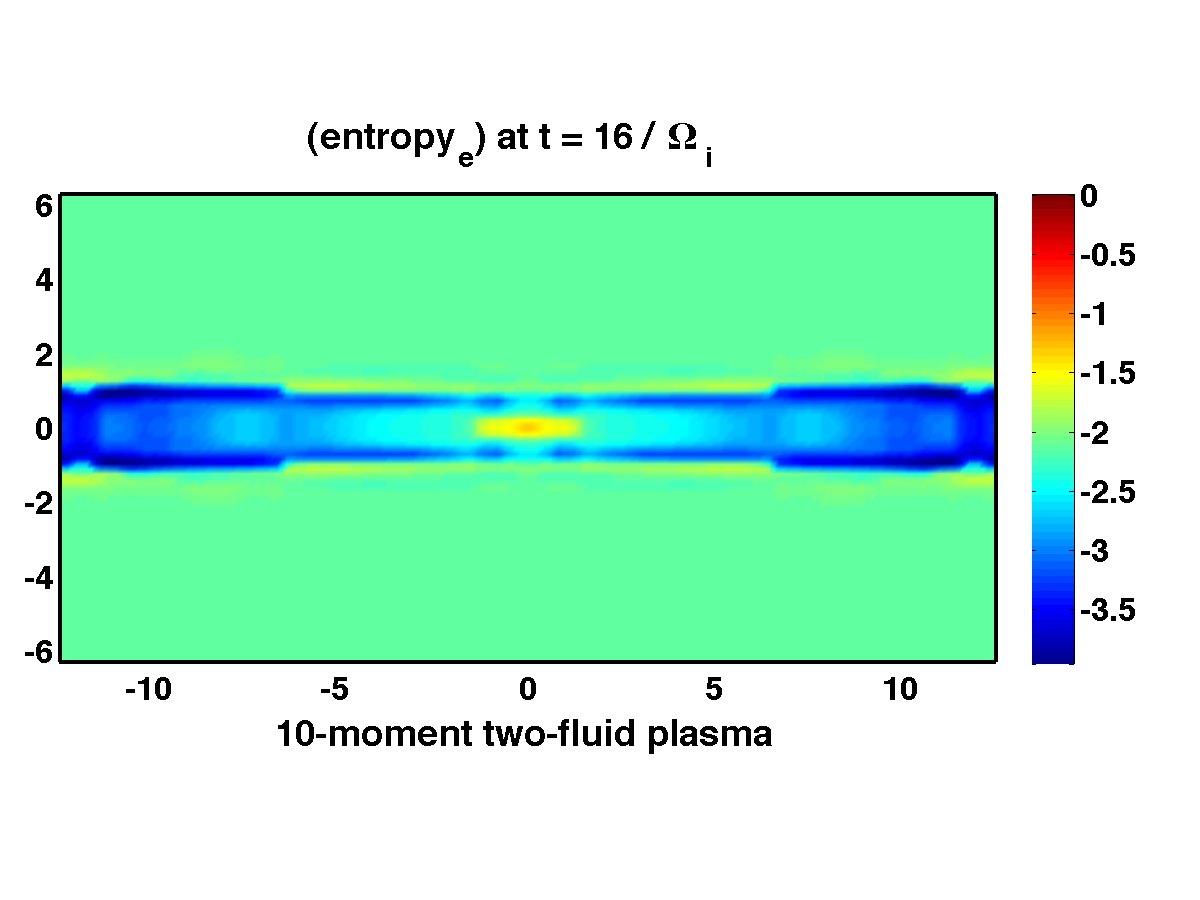}
  \end{tabular}
  \caption{Electron gas at $t=16$ (nascent singularity)}
  \label{fig:elc16}
\end{figure}

\begin{figure}
  \begin{tabular}{c c}
    \includegraphics[width=0.50\textwidth]{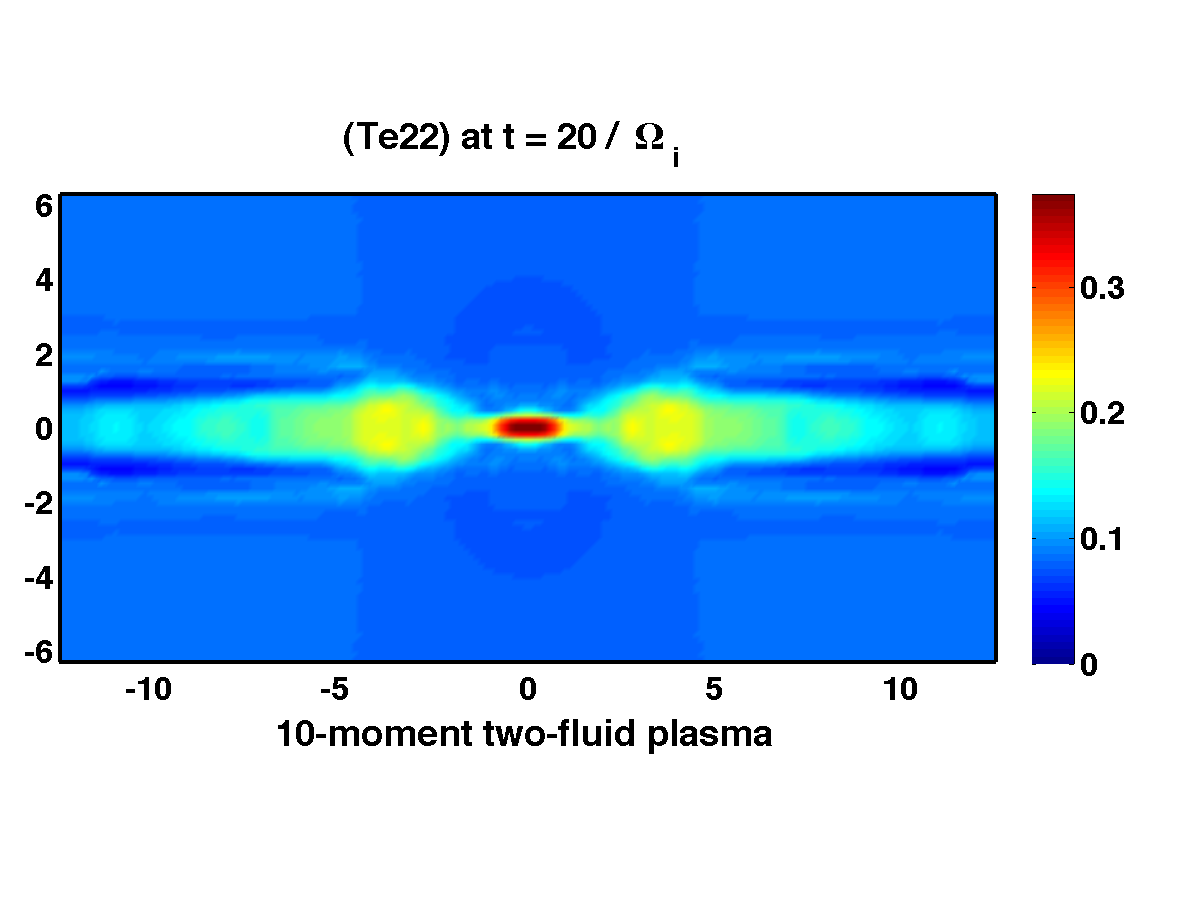}
   &\includegraphics[width=0.50\textwidth]{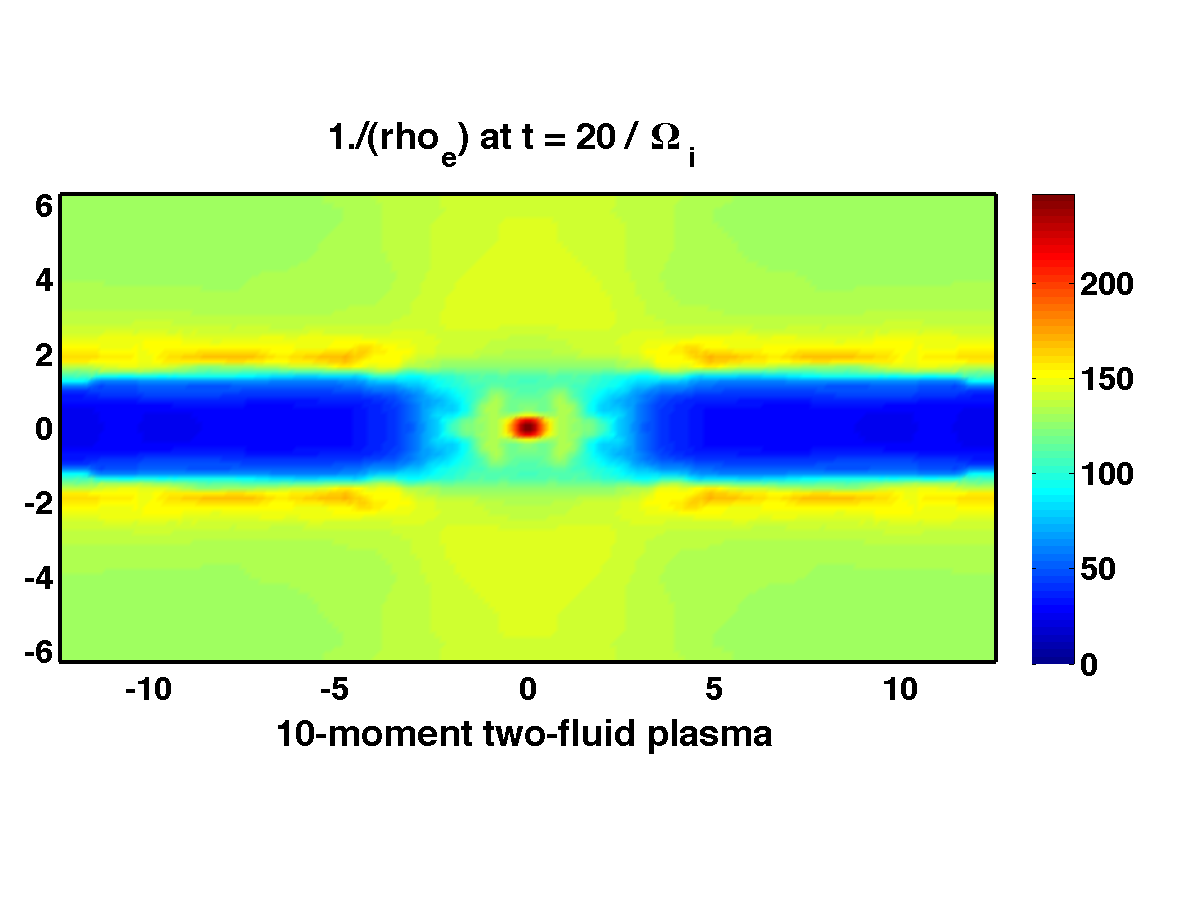}
   \\
    \vspace{-5ex}
    \\
    \includegraphics[width=0.50\textwidth]{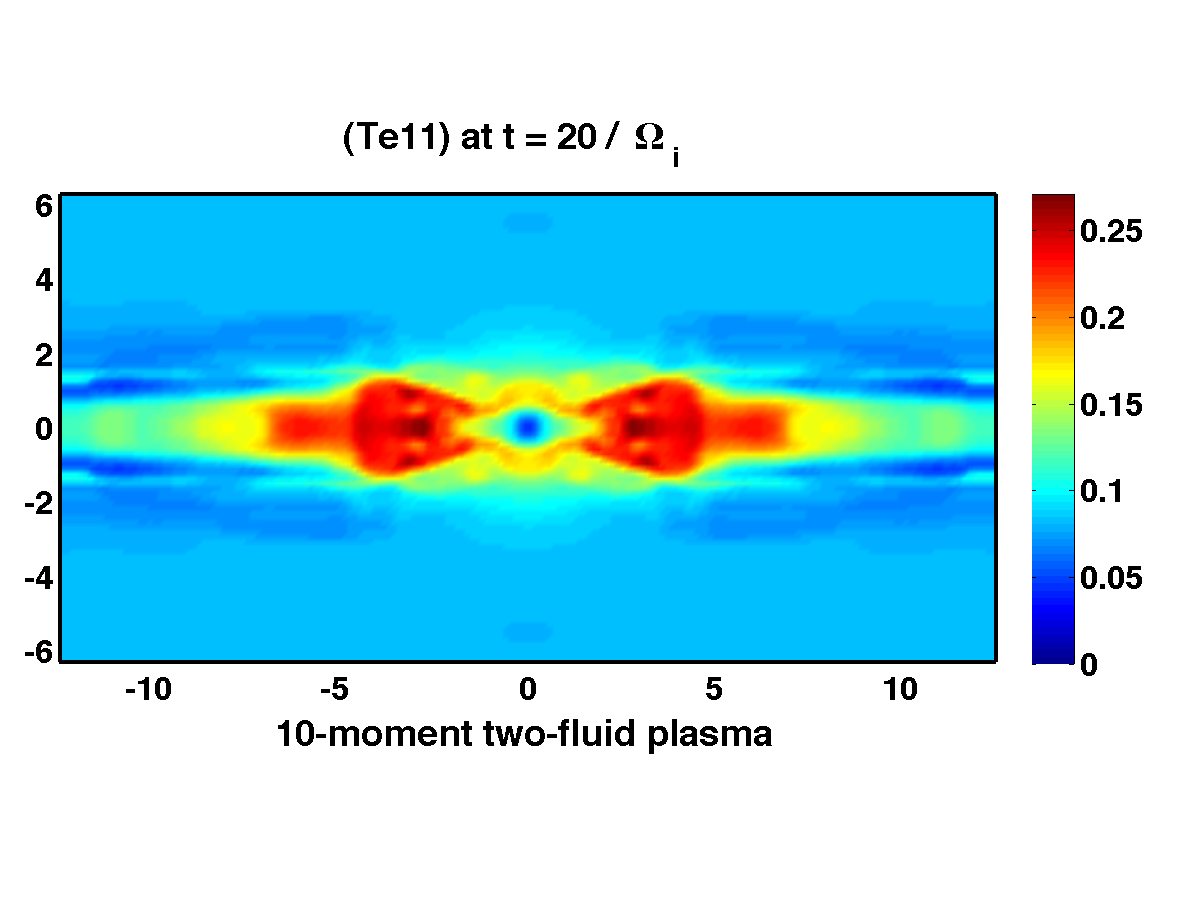}
   &\includegraphics[width=0.50\textwidth]{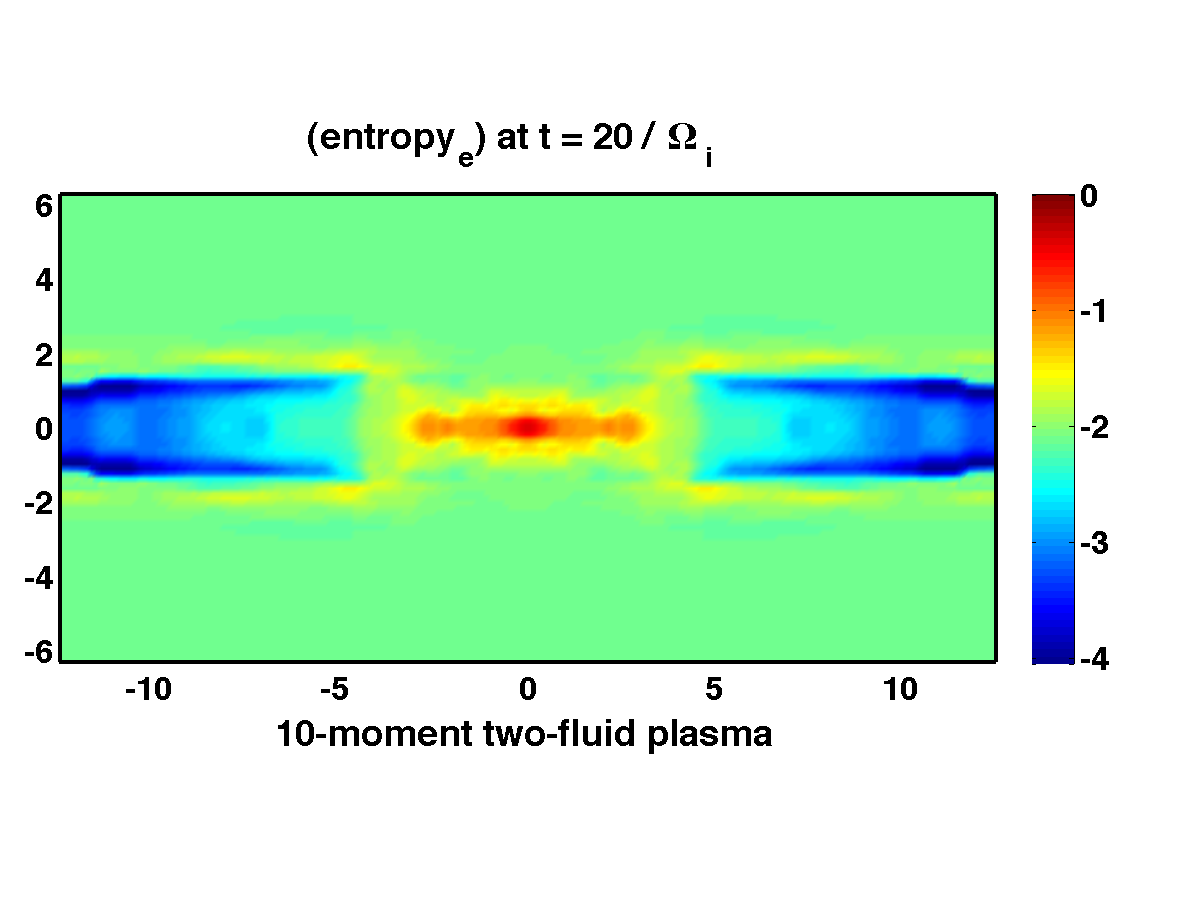}
  \end{tabular}
  \caption{Electron gas at $t=20$ (developing singularity)}
  \label{fig:elc20}
\end{figure}

\begin{figure}
  \begin{tabular}{c c}
    \includegraphics[width=0.50\textwidth]{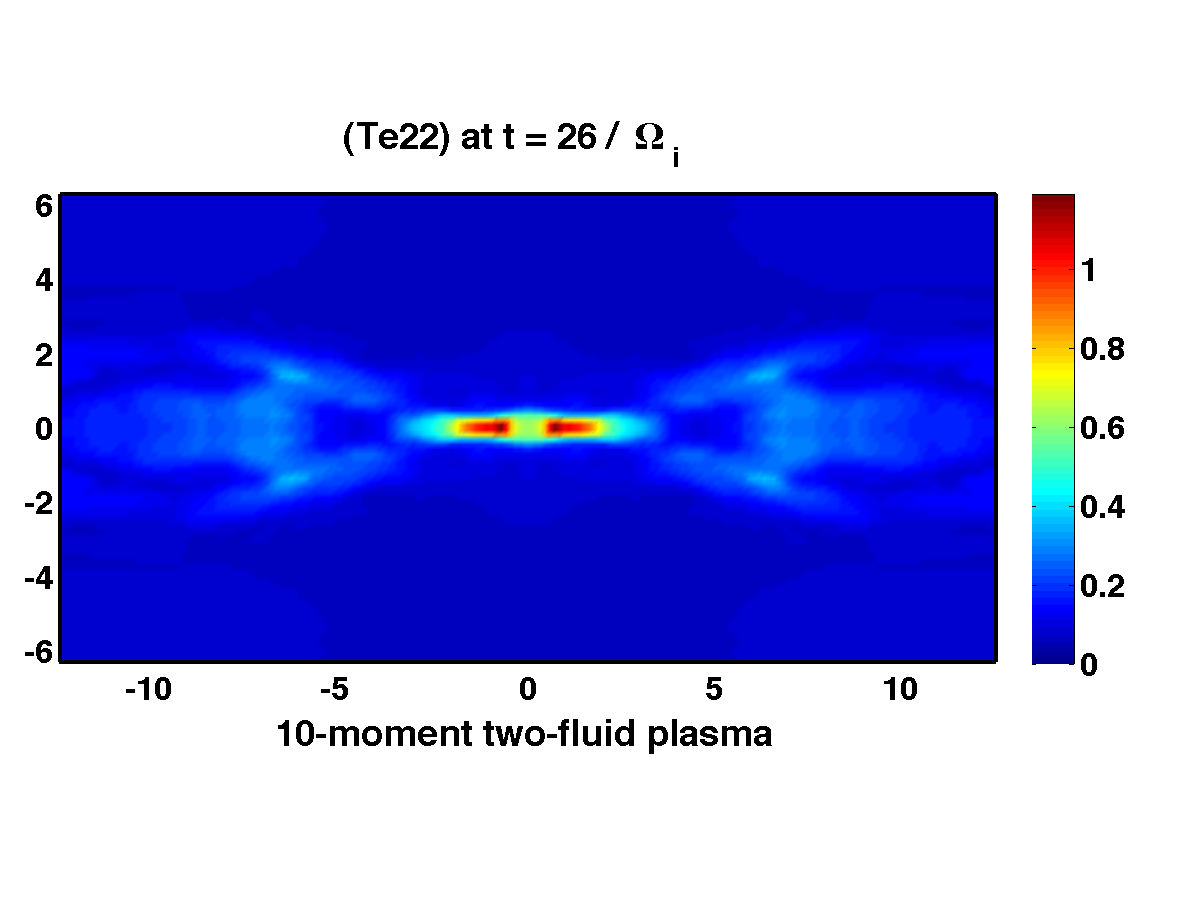}
   &\includegraphics[width=0.50\textwidth]{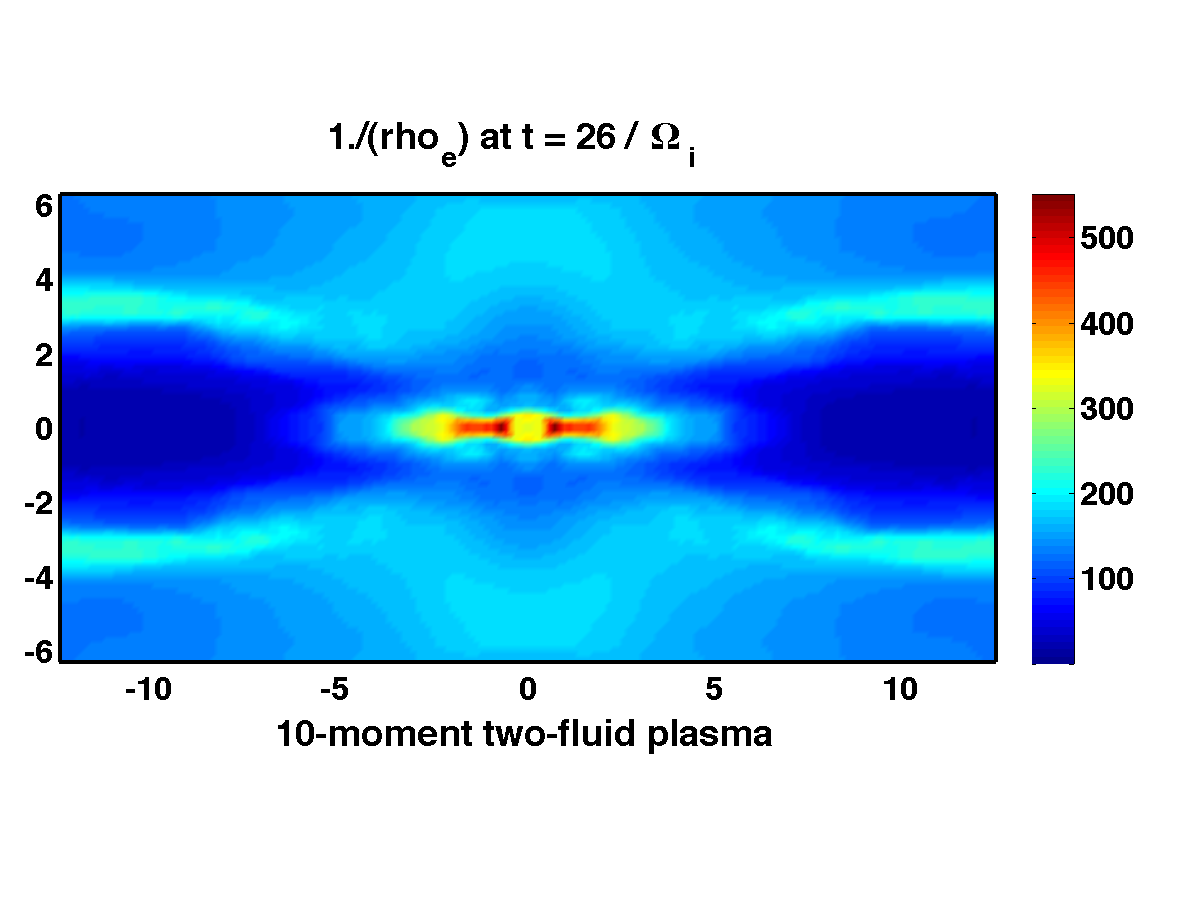}
   \\
    \vspace{-5ex}
    \\
    \includegraphics[width=0.50\textwidth]{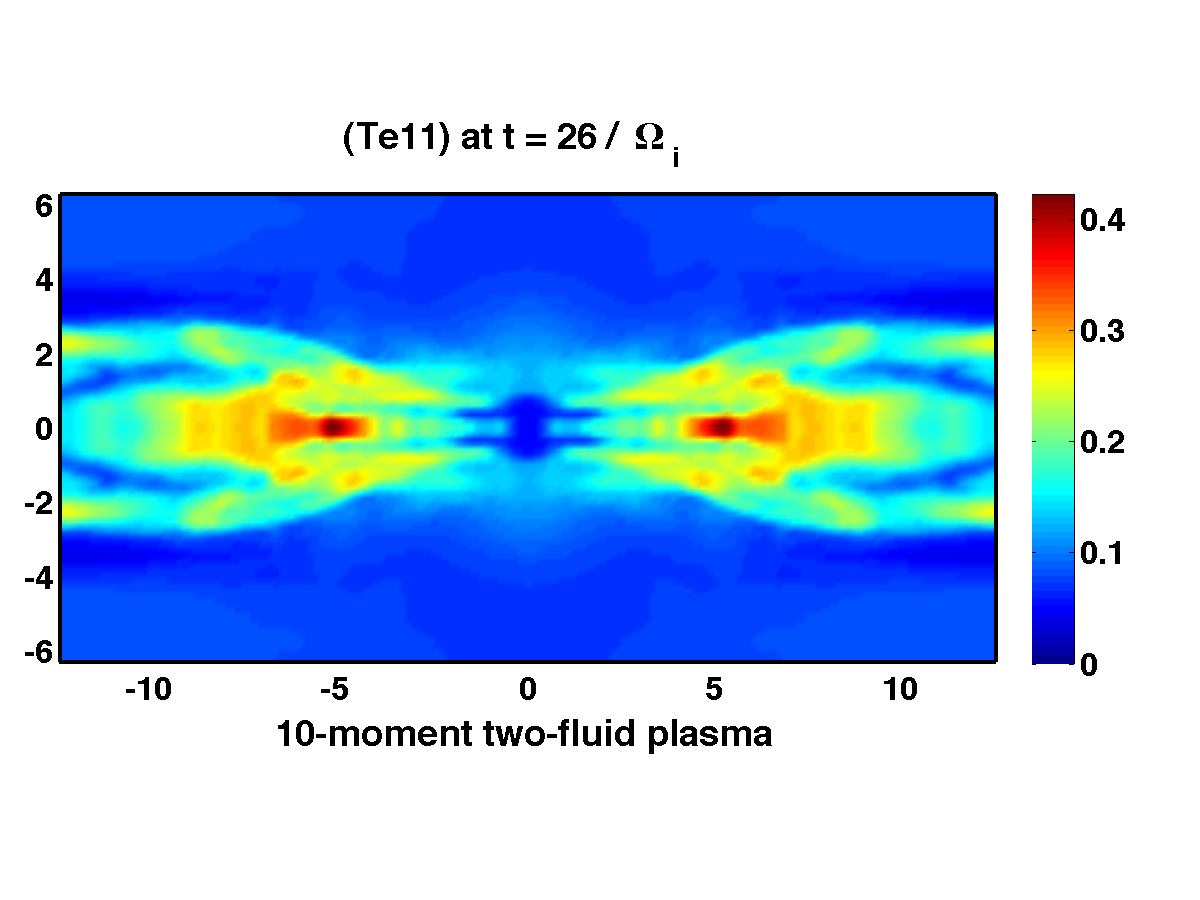}
   &\includegraphics[width=0.50\textwidth]{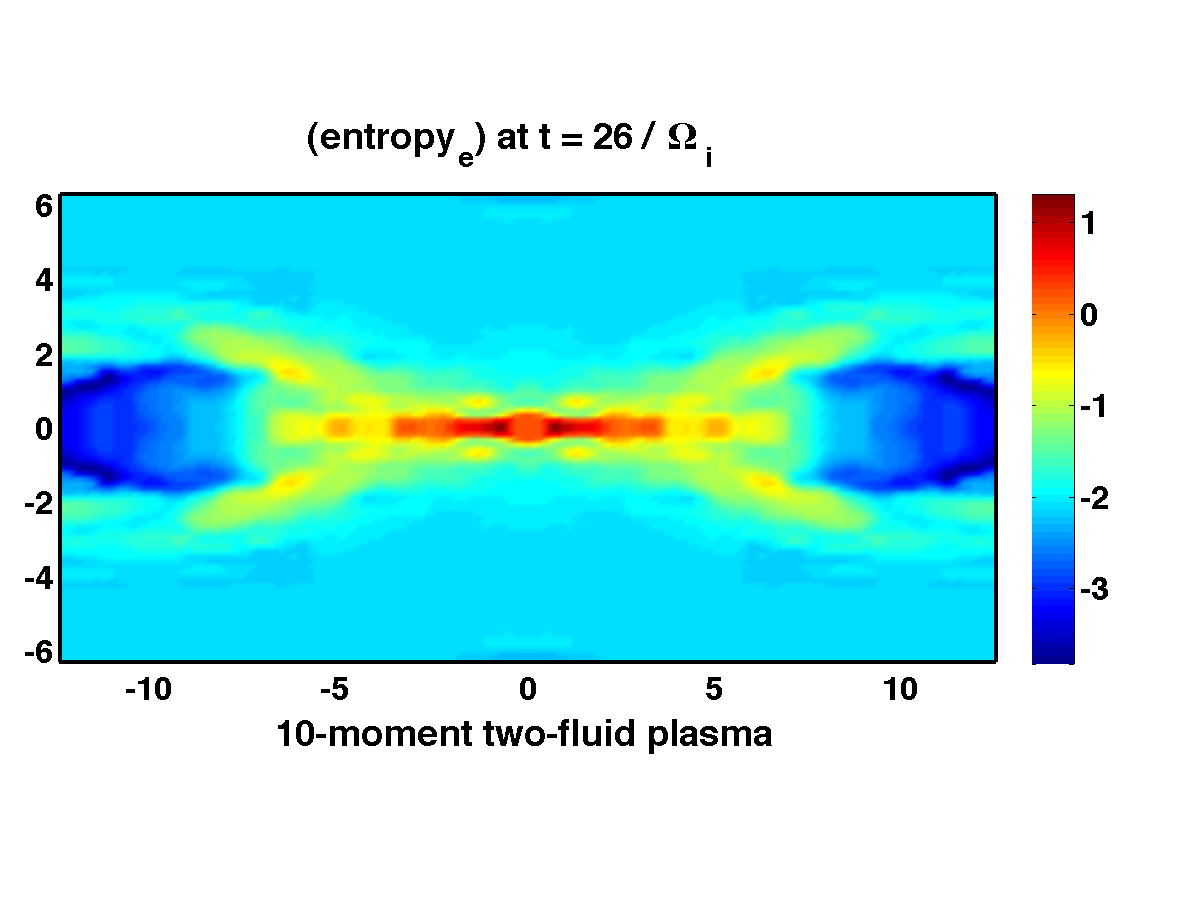}
  \end{tabular}
  \caption{Electron gas at $t=26$ (splitting singularity)}
  \label{fig:elc26}
\end{figure}

\begin{figure}
  \begin{tabular}{c c}
    \includegraphics[width=0.50\textwidth]{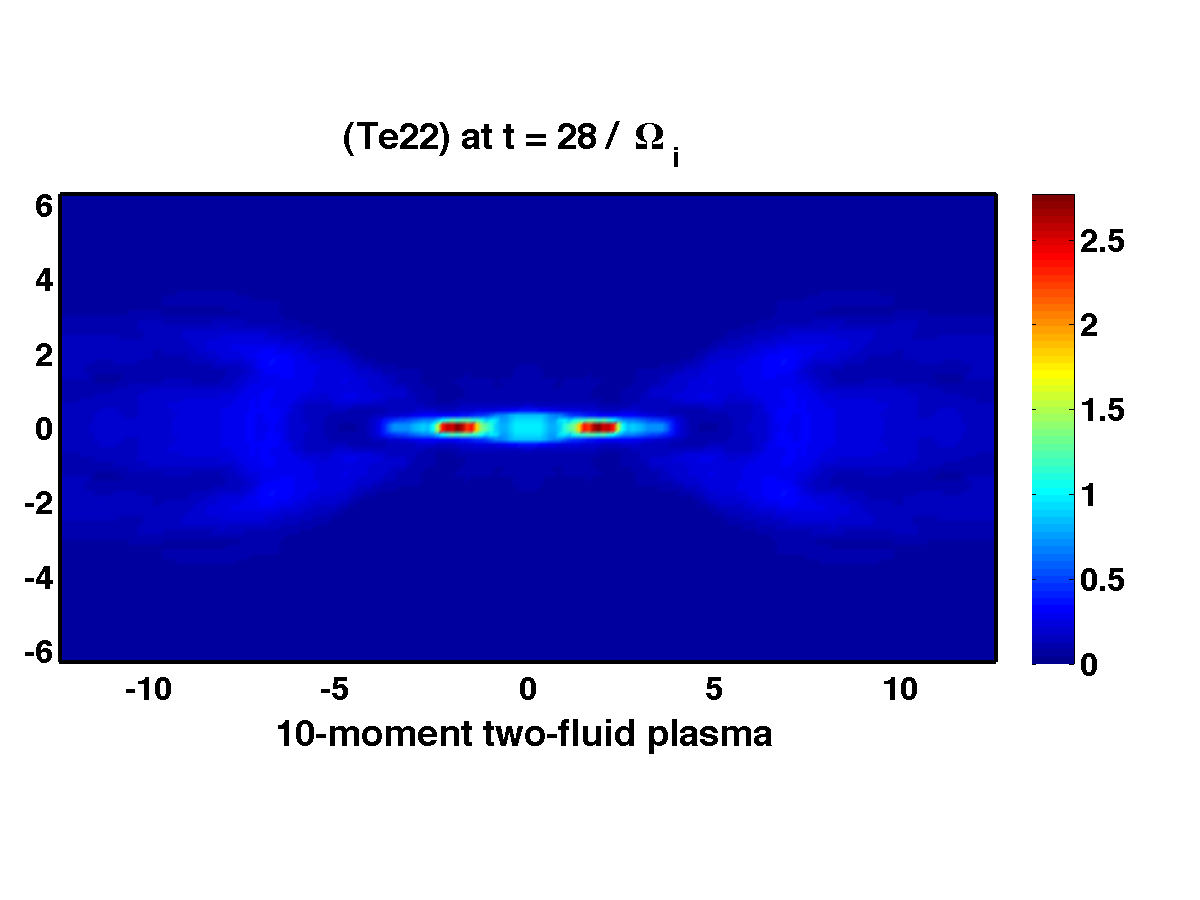}
   &\includegraphics[width=0.50\textwidth]{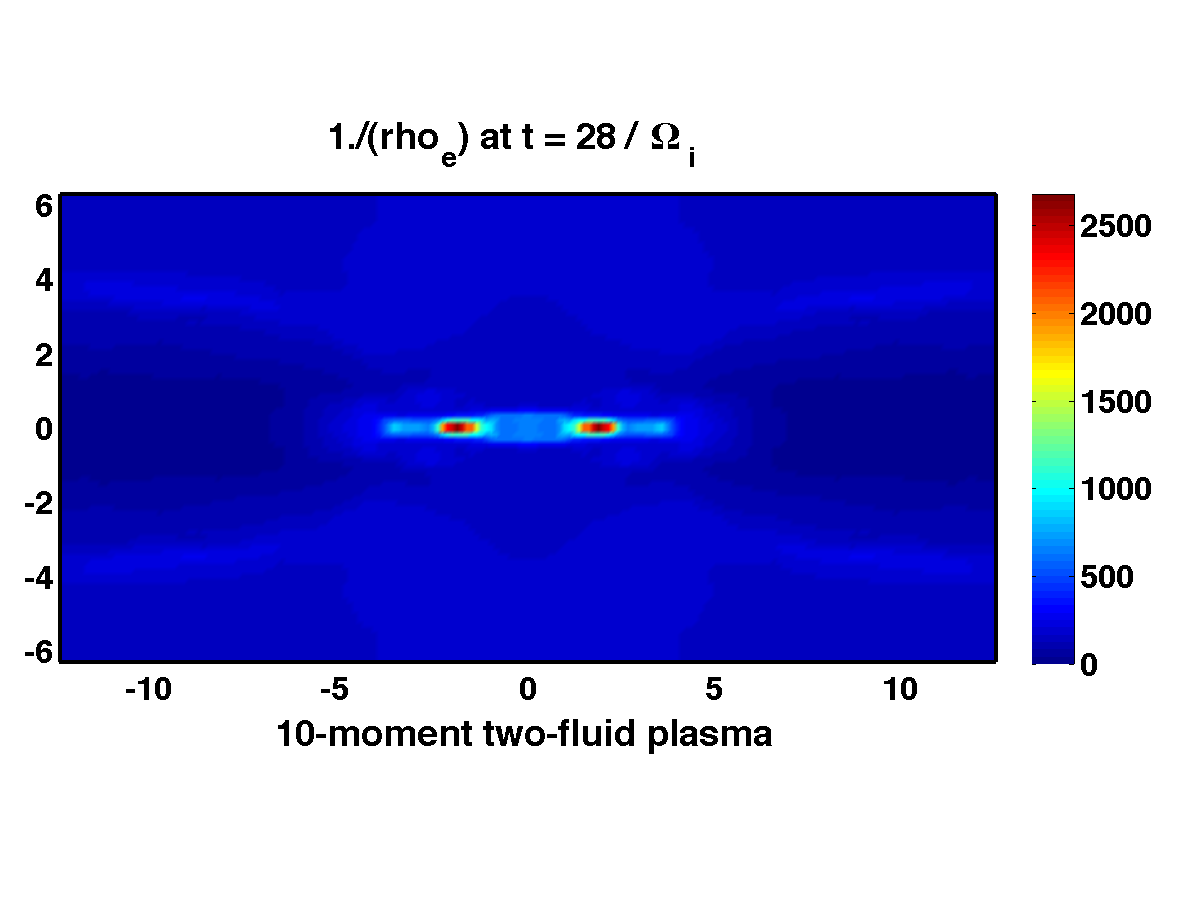}
   \\
    \vspace{-5ex}
    \\
    \includegraphics[width=0.50\textwidth]{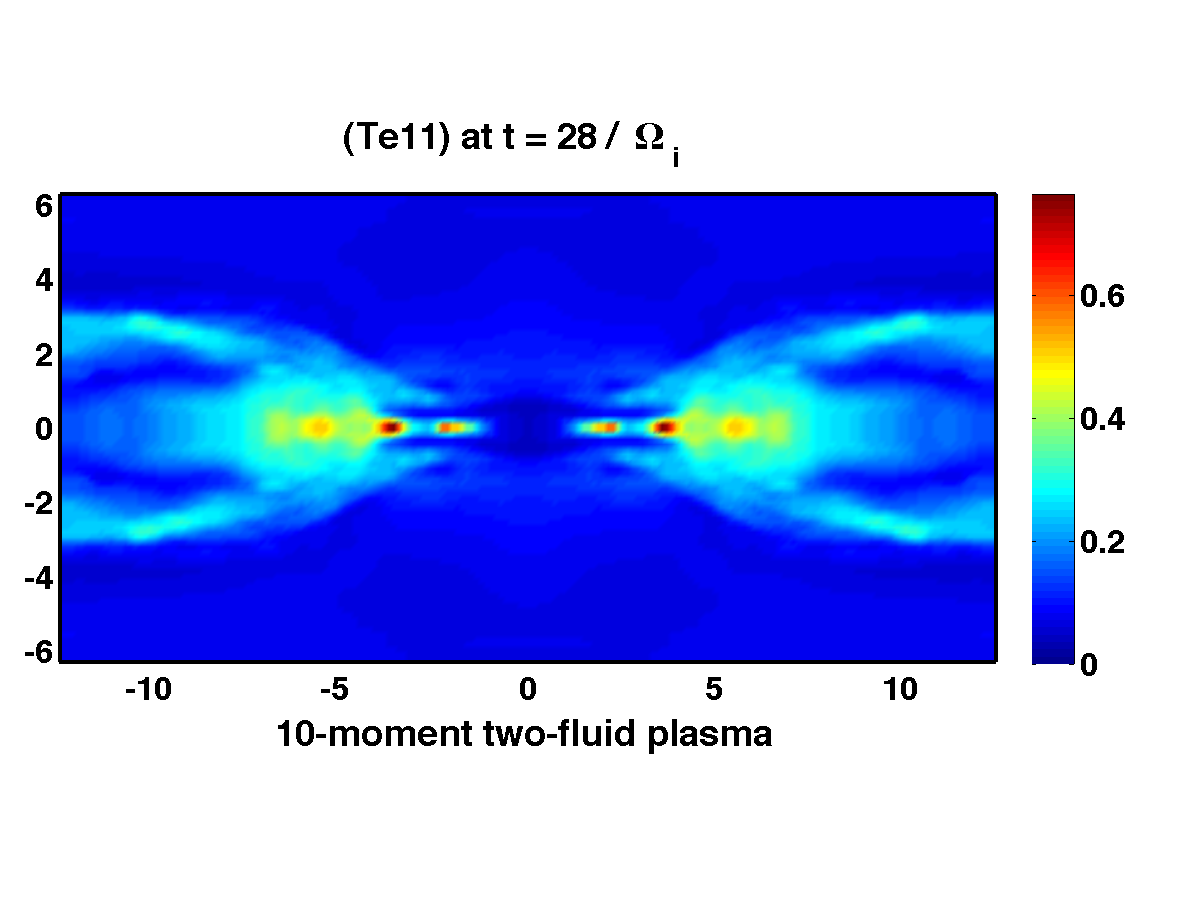}
   &\includegraphics[width=0.50\textwidth]{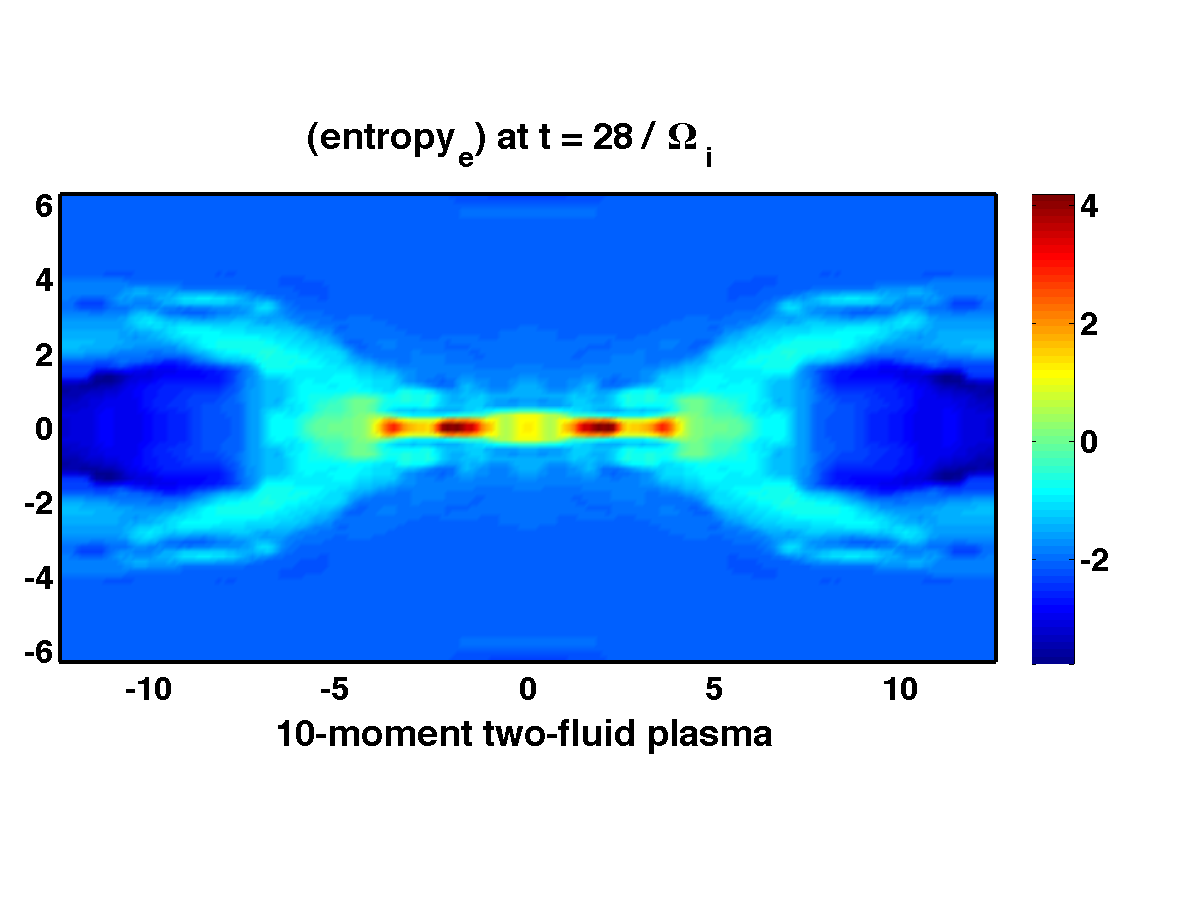}
  \end{tabular}
  \caption{Electron gas at $t=28$ (split singularity, just before crashing)}
  \label{fig:elc28}
\end{figure}

In practice, we typically observe an anisotropy that
increasingly looks like a singularity
at the X-point (see figures \ref{fig:elc16} and \ref{fig:elc20}).
The near-singularity
then splits into a pair of near-singular points
moving outward along the $x$ axis (see figures \ref{fig:elc26} and \ref{fig:elc28}).
If we enforce positivity
to prevent the simulation from crashing,
then the origin becomes an O-point and reconnection ceases.
If we do not enforce symmetry, then when this secondary
island forms symmetry is broken and the island is ejected
to one side or the other.
I remark that in the adiabatic model very strong
temperature gradients develop, especially if a fine
mesh is used.
For a coarse mesh the approximate singularities
fail to become sharp enough to disrupt reconnection
e.g.\ by causing second islands to form.
%

\newpage

\clearpage

%% file: chap7.tex
\chapter{Numerical method}

This chapter describes the numerical method used in this
dissertation to solve the two-fluid-Maxwell system. Recall
that for the two-fluid-Maxwell system we use the coupled
systems \eqref{MaxwellWithCP} for Maxwell's equations,
\eqref{TenMomentIsoSystem} for each ten-moment fluid, and
\eqref{FiveMomentIsoSystem} for each five-moment fluid.
Since we neglect all diffusive terms, the
composite system we solve fits the hyperbolic
conservation form
\def\uu{{\underline{u}}}
\def\FF{\underline{\textbf{f}}}
\def\ss{\underline{s}}
\begin{gather}
  \label{eqn:conslawFramework}
  \partial_t \uu + \Div\FF = \ss,
\end{gather}
where $\uu(t,\xb)$ is the state, $\FF(\uu)$ is the flux function,
and $\ss(\uu)$ is an undifferentiated source term.
In two spatial dimensions we can write
\def\ff{\underline{f}}
\def\gf{\underline{g}}
\begin{gather}
  \label{eqn:conslaw}
  \partial_t \uu + \partial_x\ff + \partial_y\gf = \ss,
\end{gather}
where $\ff:=\ebas_x\dotp\FF$ and $\gf:=\ebas_y\dotp\FF$.
We solve this system using the explicit discontinuous
Galerkin (DG) method.

\section{Discontinuous Galerkin method}


The discontinuous Galerkin (DG) method was developed into a 
modern tool for computing solutions to hyperbolic PDEs in
a series of papers by Cockburn and Shu (see \cite{article:CoShu5} and
references therein).

We define a Cartesian mesh $S_h$ made up of $N$ mesh cells.
Solutions are represented by members of the {\it broken}
finite element space
\def\Tm{{\mathcal T}}
\begin{equation*}
    V^h = \left\{ v^h \in L^{\infty}(S): \,
    v^h |_{\Tm} \in P^k, \, \forall \Tm \in S_h \right\},
\end{equation*}
where $h$ is the grid spacing, $\Tm$ is a mesh cell,
and $P^k$ is the set of polynomials of degree at most $k$.
Each cell can be mapped to the canonical
mesh cell, $[-1,1] \times [-1,1]$, via a simple affine
transformation. On the canonical mesh cell we define the following
normalized Legendre polynomials up to degree two:
\begin{equation}
\label{eqn:Legendre1}
\begin{split}
\{\phi^{(\ell)}\}_{\ell=1}^{6} = \Bigl\{ &1, \  &&\sqrt{3}\,\xi, \ &&\sqrt{3}\,\eta, \ 
&&3\,\xi\eta, \ && \frac{\sqrt{5}}{2} \left(3\,\xi^2 - 1 \right),
\ &&\frac{\sqrt{5}}{2} \left(3\,\eta^2 - 1 \right) \Bigr\}.
\end{split}
\end{equation}
These basis functions are orthonormal with respect
to a cell-average inner product:
\begin{equation}
   \frac{1}{4}  \int_{-1}^{1} \int_{-1}^{1}
    \phi^{(m)}(\xi,\eta) \, \phi^{(n)}(\xi,\eta) \, d\xi \, d\eta = 
     \delta_{mn}.
\end{equation}

We seek approximate solutions of (\ref{eqn:conslaw})
which are a linear combination of basis functions,
\begin{equation}
\label{eqn:q_ansatz}
q^h(\xi, \eta ,t) \bigl|_{\Tm_{ij}}  :=
  \sum_{k=1}^6  U^{(k)}_{ij}(t) \, \phi^{(k)}(\xi, \eta).
\end{equation}
Multiplying (\ref{eqn:conslaw}) by a basis function and
integrating by parts yields
the following semi-discrete evolution equations for
the Legendre coefficients, $U^{(\ell)}_{ij}$:
\begin{equation}
\label{eqn:semidiscrete}
\frac{d}{dt} \,
   U^{(\ell)}_{ij} = N^{(\ell)}_{ij}
- \frac{{\mathcal F}_{ij}^{(\ell)}}{\Delta x}
- \frac{{\mathcal G}_{ij}^{(\ell)}}{\Delta y}
+ S^{(\ell)}_{ij},
\end{equation}
where $\Delta x$ and $\Delta y$ are the dimensions of
the mesh cell and where
\begin{gather}
	\label{eqn:Nvals}
	N^{(\ell)}_{ij} = \frac{1}{4}
   \int_{-1}^{1} \int_{-1}^{1} \left[ \,
   \phi^{(\ell)}_{, x} \, f(q^h) 
   	+  \phi^{(\ell)}_{, y} \, g(q^h) \, \right] \, 
   d\xi  \, 
   d\eta, \\
	\label{eqn:Fl1}
	{\mathcal F}^{(\ell)}_{ij} = 
	\left[ \frac{1}{2} \int_{-1}^{1}
	\phi^{(\ell)} \,
		f(q^h)  \, 
		 d\eta \right]_{\xi=-1}^{\xi=1}, \\
		 \label{eqn:Fl2}
		 {\mathcal G}^{(\ell)}_{ij} = \left[
		 \frac{1}{2} \int_{-1}^{1}
	\phi^{(\ell)} \,
		g(q^h)  \, d\xi
		 \right]_{\eta=-1}^{\eta=1},
   \\
   \label{eqn:Psivals}
   S^{(\ell)}_{ij} = \frac{1}{4}
   \int_{-1}^{1} \int_{-1}^{1} 
   \phi^{(\ell)} \, \psi(q^h) \, 
   d\xi  \, 
   d\eta.
\end{gather}
We approximate these averaging integrals
with Gaussian quadrature rules.  We use
a Riemann solver to determine the flux values used in
the edge integrals.  In particular,
we approximate the integrals in (\ref{eqn:Psivals}) 
via the standard 2D 9-point rule,
the integrals in (\ref{eqn:Nvals}) via the standard
2D 4-point rule, and
the integrals in (\ref{eqn:Fl1})--(\ref{eqn:Fl2}) 
with the standard 1D 3-point rule
and local Lax-Friedrichs Riemann solvers.

We handle time-stepping via a third-order TVD Runge-Kutta
method \cite{article:GoShu98}.
After each time stage we apply a minmod limiter to the coefficients
of the quadratic basis functions and, if limiting
occured, to the coefficients of the mixed and linear terms
(a modification of the method in \cite{article:Kriv07}).
This method was implemented in the C++ code
{\sc DoGPack}, which was developed by James Rossmanith
and collaborators at UW-Madison.

\section{Source term}

For the two-fluid-Maxwell system,
the source term $\ss$ in equation \eqref{eqn:conslawFramework}
is extremely stiff.  Fortunately it is also linear.
To handle the source term we use Strang splitting
and alternately solve the flux equation
(the \mention{hyperbolic part})
\begin{gather}
  \label{fluxeqn}
  \partial_t \uu + \Div\FF = 0
\end{gather}
and the \mention{source ODE}
\begin{gather}
  \label{sourceODE}
  \partial_t \uu = \ss.
\end{gather}
This is a linear ODE with constant coefficients (with the
caveat below) and can be solved exactly.
The source ODE can be further decomposed into
a sum of three commuting operators, for which
splitting incurs no error:
\begin{gather}
  \label{sourceODE}
  \partial_t \uu = \ss_\text{electro-momentum}
                 + \ss_\text{(pressure tensor rotation)}
                 + \ss_\text{(pressure isotropization)}.
\end{gather}
The electro-momentum and pressure tensor rotation terms have
imaginary eigenvalues.  The pressure isotropization term
has a negative real eigenvalue equal in magnitude
to the isotropization rate.

As a caveat, the isotropization ODE has constant coefficients
if the isotropization rate remains constant as isotropization
proceeds, which holds if the isotropization period is defined
in terms of $\tr\TT_\s$ but not if it is defined in terms of
$\det\TT_\s$ (see section \ref{PositivityPreservingClosure}).

Splitting off the source ODE and solving it exactly
(or via the trapezoid rule) ensures that the source
ODE does not cause violation of positivity, energy
conservation, or numerical stability.  This technique
of splitting off the two-fluid source term and solving
it exactly was used for the five-moment two-fluid Maxwell
system in \cite{Kumar10}.

\section{Flux term}
\label{FluxTerm}

The chief challenges in solving the flux term are to ensure
numerical stability and to guarantee that positivity of the
solution is maintained.

To ensure numerical stability we take a time step which
respects a CFL stability condition
\begin{gather*}
  \frac{\Delta t |\lambda|_{\mathrm{max},e}\,\mathrm{d}A_e}{|K|} \le \frac{1}{2m-1},
\end{gather*}
for all edges $e$ of the rectangle, where $m$ is the order of the method
(which in these simulations was $m=3$);
$|\lambda|_{\mathrm{max},e}$ is the maximum wave speed
perpendicular to the edge, $\mathrm{d}A_e$ is the ``area'' (length)
of the edge, and $|K|$ is the ``volume'' (area) of the rectangle.


\section{Variation limiters}

To ensure numerical stability we apply limiters to the
solution after each time step which limit oscillations.
Effective limiting requires that we transform into
the characteristic variables defined by the eigenstructure
of the flux Jacobian of each separable subsystem.
These eigenstructures are calculated in section \ref{tenMomentEigs}.

Note that for a PDE of the form \eqref{eqn:conslawFramework},
\begin{gather}
  \label{eqn:conslawFrameworkCopy}
  \partial_t \uu + \Div\FF = \ss,
\end{gather}
a flux Jacobian such as $\ebas_x\dotp\partial\FF/\partial\uu$
is independent of the source term $\ss$ and therefore
application of limiters to equation \eqref{eqn:conslawFrameworkCopy}
is identical to application of limiters to the conservation law
\begin{gather}
  \partial_t \uu + \Div\FF = 0;
\end{gather}
we call this conseration law the
\defining{hyperbolic part} of the system \eqref{eqn:conslawFrameworkCopy}.
For the hyperbolic part the system is decoupled into three
noninteracting subsystems:
gas dynamics for the ions, 
gas dynamics for the electrons, and
Maxwell's equations.  To limit the solution we
transform each subsystem into the characteristic
variables defined by the eigenstructure of
its flux Jacobian.  We compute the eigenstructure
for MHD (which generalizes a five-moment gas),
a ten-moment gas, and Maxwell's equations in
appendix \ref{Numerics}.

\subsection{Limiting of a 1D scalar problem}
In accordance with the recipe of Krivodonova
\cite{article:Kriv07}, the coefficients
$c_i^\ell$ of the order-$\ell$ Legendre polynomial
basis function in mesh cell $i$
are limited only if limiters are triggered
for all higher-order polynomials; in this case
we define the limited coefficient
\begin{gather*}
  \bar c_i^\ell = \minmod(c_i^\ell,D^{+\ell},D^{-\ell}),
\end{gather*}
where
\begin{gather*}
  \minmod(a,b,c) := \left\{
   \begin{array}{ll}
     \textrm{sgn}(a)\min(|a|,|b|,|c|) & \textrm{if}\ \sgn(a)=\sgn(b)=\sgn(c),
     \\
     0 & \textrm{otherwise}
   \end{array}\right.
\end{gather*}
and where $D^{+\ell}$ is a tunable constant times 
the forward difference
$c_{i+1}^{\ell-1}- c_{i}^{\ell-1}$
and $D^{-\ell}$ is the same constant times 
the backward difference
$c_{i}^{\ell-1}- c_{i-1}^{\ell-1}$;
we say that limiting has been triggered for order $\ell$
if $\bar c_i^\ell \ne c_i^\ell$.
The tunable constant determines the aggressiveness
of the limiters and is bounded by the constant
that defines a first-order-accurate finite-difference
estimate of $c_i^\ell$.

\subsection{Limiting of a 1D system}
\label{1Dlimiting}

To perform limiting on a 1D problem,
we try to reduce it to the scalar case.
A 1D hyperbolic problem is of the form
\begin{gather*}
  \partial_t \uu + \partial_x\ff(\uu) = 0.
\end{gather*}
Applying the chain rule, we have
\begin{gather*}
  \label{quasilinearPDE}
  \partial_t \uu + \ff_\uu\dotp\uu_x = 0,
\end{gather*}
where $\ff_\uu:=\partial\ff/\partial\uu$ is the flux Jacobian.
If $\uu$ is smooth and oscillations are small,
then the quasilinear ODE \eqref{quasilinearPDE}
is approximated by the linear ODE
\def\AA{\tensorb{A}}
\def\RR{\tensorb{R}}
\def\LL{\tensorb{L}}
\def\LT{\LL^{T}}
\begin{gather*}
  \label{quasilinearPDE}
  \partial_t \uu + \AA\dotp\uu_x = 0,
\end{gather*}
where in a given cell we take $\AA$ to
be $\ff_\uu$ evaluated at the cell average.
By assumption this system is hyperbolic
and therefore we can write the diagonalization
\begin{gather*}
  \AA = \RR\dotp\Lambda\dotp\RR\inv,
\end{gather*}
where $\Lambda$ is diagonal and real.
The columns of $\RR$ are right eigenvectors
and the rows of $\LT:=\RR\inv$ are left
eigenvectors.
Multiplying equation \eqref{quasilinearPDE}
by $\LT$ gives
\def\vv{{\underline{v}}}
\begin{gather*}
  \partial_t \vv + \Lambda\dotp\vv_x = 0,
\end{gather*}
where the variables $\vv:=\LT\dotp\uu$ are
called \defining{characteristic variables};
for smooth solution on a fine mesh this
effectively decouples the PDE into scalar
advection equations in the characteristic
variables.  We can apply 1D limiting
in these variables and then transform back.

In detail, suppose the solution representation
\def\UU{\underline{U}}
\def\VV{\underline{V}}
\begin{gather*}
  \UU = \sum_\ell \UU^\ell\phi\sup{\ell}.
\end{gather*}
To get characteristic variables multiply by $\LT$:
\begin{gather*}
  \VV = \LT\dotp\UU = \sum_\ell \left(\LT\dotp\UU^\ell\right)\phi\sup{\ell}.
\end{gather*}
We apply 1D limiters to the components $\LT\dotp\UU^\ell$.


\subsection{2D limiting for a Cartesian mesh}
\label{2DcartesianLimiting}

For a 2D Cartesian mesh we
limit polynomial basis functions $\phi\sup{\ell}$ that depend
only on $x$ in each row of cells aligned with the $x$
axis as if the problem were homogeneous perpendicular
to the $x$ axis.  We do analogously for $y$.

For a third-order method one must also limit the coefficients
of the mixed polynomial basis element $xy$. To limit these
coefficients, we first limit in the $x$ direction to obtain
$x$-limited coefficients, then limit in the $y$ direction to
obtain $y$-limited coefficients, and then define the limited $xy$
coefficient to be the minmod of the $x$-limited coefficient and
the $y$-limited coefficient.

\section{Positivity limiting}

To guarantee positivity we use positivity limiters. Positivity
limiters modify solutions which violate positivity by damping
deviations from the cell average to ensure that positivity of
the cell average is maintained from one time step to the next.
We have applied the method of Zhang and Shu, generalized to
10-moment gas dynamics \cite{article:ZhangShu10rectmesh}. To
ensure that the source term does not cause positivity violations
we use time splitting and solve the source term equation exactly.

%% file: chap8.tex
\chapter{Conclusions and future work}

\section{New results}
This dissertation contains the following new results:
\begin{enumerate}
\item \textit{Simulations of the \textbf{GEM magnetic reconnection
  challenge problem} with an adiabatic two-fluid
  model with pressure tensor evolution and isotropization.}
 \begin{enumerate}
  \item Simulations of a \textbf{pair plasma} version of the GEM problem
    which dial the isotropization rate from zero to infinity.
    \begin{itemize}
     \item The rate of reconnection is insensitive to the rate of
      isotropization and is roughly 60\% of the rate of reconnection
      in published particle simulations.
     \item For isotropic pressure the ramp-up of out-of-plane
       velocity at the origin with reconnecting flux is unphysical
       and ultimately not sustainable.  For feasible convergence
       a small but sufficiently large viscosity must be present.
    \end{itemize}
   \item Simulations of the GEM problem with an isotropizing pressure tensor.
    \begin{itemize}
     \item At the time of peak reconnection rate the
      electron pressure tensor agrees well with published kinetic
      simulations and is strongly agyrotropic.
      
      The isotropization rate which yields this agreement is well below
      the rate for which the five- and ten-moment models agree,
      as would be expected in a regime that allows strong pressure anisotropy.
      Due to the strong pressure anisotropy at the X-point
      one would not expect to get this level of agreement with
      a viscous five-moment model, although I have not confirmed this
      with five-moment viscous simulations.
     \item For late times and fine mesh resolution 
      singularities develop in the adiabatic five- and ten-moment models,
      crashing the simulation code unless positivity limiting is
      enforced.  This prompts the need for a ten-moment heat flux closure.
    \end{itemize}
 \end{enumerate}
\item \textit{A proof that \textbf{steady reconnection must be singular} for
  models which neglect intraspecies collisional terms
  when applied to problems that are invariant under 180-degree
  rotation about and translation along a symmetry axis.}
  \begin{itemize}
    \item \bfit{The Vlasov model does not admit nonsingular
      rotationally symmetric steady reconnection.}
      The implication is that a nonzero collision operator
      for the kinetic equation must be specified to define a standard of truth
      for which converged solutions can be obtained.  
      \textit{\textbf{To obtained converged solutions} it is not enough
      to say that collisions (whether particle-particle or wave-particle
      interactions) are present in a simulation; \textbf{a nonzero collision operator
      must be specified (and known)}.}
      To get converged solutions \emph{you must know what equations you are solving
      and they must have a solution.}
    \item This does not imply that it is necessary to have an explicit
      heat flux.  Specifically, viscous incompressible models typically
      do not evolve an energy equation and can admit steady fast reconnection.
    \item In general we can say that use of truncation closure
      (for deviatoric stress or heat flux or the collision operator)
      results in a singularity.
  \end{itemize}
\item \textbf{A closure for tensor heat flux in the presence of a magnetic field.}
  \begin{itemize}
    \item
      The closure is based on a Chapman-Enskog expansion about a Gaussian
      distribution of particle velocities and assumes a Gaussian-BGK collision model.
    \item
      The trace of this closure gives a closure for vector heat flux
      which Woods has shown agrees well with the more sophisticated closures
      of Braginskii and of Chapman and Cowling \cite{book:Woods04}.
      These more sophisticated closures assume Coulomb collision operators
      and define closure coefficients in terms of rational functions
      with higher order polynomials in the ratio \gls{hgf}
      of gyrofrequency to thermal relaxation rate.
   \item
     A closure for the heat flux tensor could be derived
     which is more accurate for Coulomb collisions.  It is doubtful,
     however, that there is much benefit to be gained.
     In collisional regimes one might as well use the five-moment model.
     In general, extended moment modeling does not aim for high accuracy
     but for reasonable accuracy for the weakest possible collisionality.
     In the context of fast magnetic reconnection, the diffusion region
     is at most weakly collisional.  Outside the diffusion region the
     pressure tensor (and thus the heat flux) does not make a significant
     contribution to the electric field in Ohm's law.  Of course one
     can alway cook up perfect agreement with kinetic simulations by
     using spatially and temporally dependent problem-specific
     anomalous closures, but the prospects for a generic closure
     with substantially improved agreement would seem not to be promising.
  \end{itemize}
\end{enumerate}

\section{Questions raised}
These results of this work raise the following questions for further investigation:
\begin{itemize}
 \item \textit{Is steady reconnection singular for truncation
   closure in non-symmetric 2D reconnection configurations?}

   Symmetry under 180-degree rotation about the X-point
   is physically unlikely and unstable.
   The argument of section \ref{steadyReconNeedsHeatflux} turns on 
   this symmetry and in particular that the stagnation point coincides
   with the X-point.  In asymmetric reconnection these two points do not coincide
   \cite{murphySovinec08,murphy09}.
   In simulations of reconnection, 180-degree rotational symmetry about the
   X-point is unstable to spontaneous symmetry-breaking.
 \item \textit{Is steady reconnection singular for truncation closures for 3D
   reconnection?  Is it possible to make a more general
   argument for the need for entropy production and diffusive entropy
   flux?}
 \item \textit{What insight can be gained from a ten-moment two-fluid linear tearing
   analysis for weakly collisional plasma?}
   
   How does the reconnection electric field and resistive,
   ``viscous,'' and inertial components depend on the
   viscosity and resistivity? Simple isotropization of the
   pressure tensor implies gyroviscosity already ``built-in'',
   potentially facilitating and simplifying the linear tearing
   analysis developed in \cite{mirnov04} and generalized in
   \cite{HoBiVe09}.
 \item \bfit{For the kinetic equation what collision operator should one
   specify as a standard of truth for magnetic reconnection problems?}

   In current practice models used to simulation reconnection
   often rely on numerical diffusion and fail to admit converged
   reconnection.  At a practical level people generally get away with it.
   The physical collision operator is generally unknown, and so
   for the purposes of validation against experiment and observed phenomena
   it is plausible that one might as well rely on numerical diffusivity. 
   One can hope to get reliable results only for stastistics which are
   insensitive to the choice of collision operator.  If a statistic
   of interest does not change significantly as the mesh is refined
   then this is evidence that the statistic is insensitive to the choice
   of collision operator and the simulations are reasonably regarded as
   having physical validity.

   The trouble with relying on numerical diffusivity, however, is that
   verification of solutions against equations is impossible because
   it is unknown what equations one is actually solving.
   It becomes impossible to quantify error.  

   These issues are particularly pertinent when attempting to adjudicate among
   models when simulations that use them disagree.  The natural recourse
   is to refer to a modeling hierarchy and to prefer model A over model B
   if simulations using B differ from simulations using A and if B is a simplification
   of A.  This heuristic can be used as a practical standard to adjudicate
   among any two models if all models are derived from a common ancestor,
   taken to be the \mention{standard of truth}.
   Without a standard of truth it can become difficult or impossible to
   adjudicate among competing models.

   The practical standard of truth in the plasma community has been
   PIC simulations.  The natural question to ask becomes,
   ``\emph{What equation does a given PIC simulation solve?}'' 
   To answer this question one must consider 
   convergence in the limit as the mesh is refined.
   Unless a nonzero collision operator is explicitly
   incorporated into the solution, one is presumably
   solving the Vlasov equation.
   The point here is that in simulating a problem
   one should be able to state the equations that
   one is solving and confirm that they actually have a solution.
   It is a critical result of this work that \emph{the Vlasov
   equation does not admit nonsingular steady reconnection for 2D
   problems rotationally symmetric about the X-point}; a mechanism
   for entropy production is necessary.
   For such problems,
   unless an explicit nonzero collision operator is incorporated
   into a PIC simulation, it evidently fails to satisfy these
   two criteria.

   For such problems, the standard of truth evidently requires
   specification of a nonzero collision operator, and there does
   not seem to be an obvious single standard of truth.

   How sensitive are solutions to the choice of collision operator?
   What happens as one takes the collision operator to zero?
   The fact that steady-state solutions are singular in the
   absence of collisions indicates that, at least by some measure,
   the steady-state solution is \mention{structurally unstable}
   in the limit as the collision operator goes to zero. 
   (We say that a problem is \defining{structurally unstable}
   when it depends discontinuously on the model parameters.)

   Can we define a universal standard of truth using a statistical approach?
   In particular, can we show that for steadily driven reconnection
   a statistical steady state exists when the collision operator is zero?
   Can we study statistical steady state solutions to the Vlasov equation
   by defining a \mention{stochastic collision operator} and taking it to zero?
   Can we define a set of statistics of interest and
   show that these statistics are structurally stable?
   If so, perhaps we can justify the use of PIC simulations
   for these problems without an explicit collision operator.

\item \emph{Can we formulate a fluid analog to the kinetic equation
   with Gaussian-BGK collision operators which
   allows us to study magnetic reconnection with vanishing collisionality?}

  This work shows that for converged simulation of steady reconnection
  it is necessary to have nonzero heat flux and nonzero deviatoric pressure.
  Diffusive closures do not provide a way of exploring
  the zero-collisionality limit. 
  The kinetic equation with a Gaussian-BGK collision operator
  makes it possible to study magnetic reconnection in
  the zero-collisionality limit, but only with great computational
  expense.  What we seek in an extended moment fluid model is not
  detailed agreement but \emph{analogous behavior}
  in comparison to the kinetic equation
  as the collision operator is taken to zero.
  \bfit{We need a fluid model which evolves
  a heat flux tensor as well as a pressure tensor.}
  With such a fluid model we can dial from infinity to zero the rate
  at which the heat flux tensor and deviatoric pressure tensor are relaxed
  to zero and gain insight into the zero-collisionality limit.
  A fluid model with heat flux evolution is, I propound, the Holy Grail of 
  fluid modeling of collisionless magnetic reconnection.
\end{itemize}

%% file: appendix2.tex
\chapter{Modeling}
\def\cyan{\color{cyan}}
\section{Maxwellian and Gaussian distributions}
\label{distributions}

\emph{
In this section entropy decreases, in accordance with the 
convention of mathematicians; that is, the definition
of entropy is minus the definition of entropy elsewhere
in this dissertation.
}

The Maxwellian and Gaussian distributions are
the two working examples of Galilean-invariant
\emph{entropy-minimizing} closures for the
equations of gas dynamics. 
The Maxwellian distribution is the assumed
distribution of hyperbolic five-moment gas-dynamics
(the compressible Euler equations).
The Gaussian distribution is the assumed distribution
of hyperbolic ten-moment gas-dynamics.
A Maxwellian distribution is a normal distribution
that is isotropic in the reference frame of the fluid.
A Gaussian distribution is a distribution
that in the reference frame of the fluid
is a product of normal distributions with
possibly different distribution
widths in three principle orthogonal directions.

An entropy-minimizing closure (for a given set of moments)
requires that particle distributions minimize entropy over all
distributions which share the same given moments. Only variation
in velocity is considered, not variation in space. This is
consonant with the fact that collision operators ignore variation
of particle density in space and only consider variation in
velocity space.  Thus in this document we ignore variation
in space.

\def\mean#1{\langle #1 \rangle}
\def\tr{\,\textrm{tr}\,}
\def\v{\ensuremath{\mathbf{v}}}
\def\u{\mathbf{u}}
\def\c{\mathbf{c}}
\def\M{\mathbf{M}}
\def\energy{\mathcal{E}}
\def\Energy{\mathbb{E}}

\emph{Definitions of conserved moments.}
Let $f(t,\v)$ be the distribution of particle
mass over velocity space.
\begin{align*}
   \mdens &:= \int_\v f &\hbox{is mass (density)},
\\ \M &:= \int_\v f \v &\hbox{is momentum (density)},
\\ \energy &:= \int_\v f \v^2/2 &\hbox{is energy (density)},
\\ \Energy &:= \int_\v f \v \v &\hbox{is energy tensor (density)},
\end{align*}
\emph{Definitions of statistical averages.}
Let $\chi(\v)$ be a ``generic'' moment.
Denote and define its statistical average by
\begin{gather*}
  \mean{\chi} := \frac{\int_\v f\chi}{\mdens}.
\end{gather*}
Primitive variables are naturally
defined in terms of statistical averages:
\begin{align*}
   \u &:= \mean{\v} &\hbox{ is bulk velocity},
\\ \c &:= \v-\u &\hbox{ is thermal velocity},
\\ \glslink{ThetaT}{\ThetaT} &:= \mean{\c\c}
   &\hbox{ is the \defining{pseudo-temperature tensor}, and}
\\ \glslink{theta}{\theta} &:= \mean{c^2/3}
   &\hbox{ is \defining{pseudo-temperature}.}
\end{align*}
Relationships among
primitive and conserved variables are
\begin{align*}
   \mdens\u  &= \M,
\\ \energy &= (\mdens \u^2 + 3\mdens\theta)/2,
\\ \Energy &= \mdens\u\u + \mdens\ThetaT,
\\ \theta  &= \tr\ThetaT/3,
\\ \energy &= \tr\Energy/3.
\end{align*}
Recall that entropy $S$ is defined by
\begin{align*}
  S &:= \int_\v \eta, &&\hbox{where} &
  \eta &:= \beta f\ln f + \alpha f,
\end{align*}
where $\beta>0$ may be chosen arbitrarily
for a single-species gas and
where $\alpha$ is a constant that can be
freely chosen; we will choose $\alpha = 3(\ln(2\pi)+1)/2$
to make the formula for the gas-dynamic
entropy simple.
Note that 
\begin{gather*}
  \eta' = \beta \ln f + (\beta +\alpha).
\end{gather*}
So to minimize the entropy of a single species
we may conveniently choose $\eta=f\ln f-f$,
for which $\eta'=\ln f$.

\subsection{Maxwellian case}

In the Maxwellian case
we minimize $S$ subject to the constraints that
\begin{align*}
  \int_\v f &= \mdens,
  &\int_\v \v f &= \M,
  &\int_\v \v^2 f &= 2\energy.
\end{align*}
We use the technique of Lagrange multipliers.  Define
\def\ft{\tilde f}
\begin{gather*}
  g := \int_\v \eta
    +\lambda\left(\mdens - \int_\v f\right)
    +\mu\cdot\left(\M - \int_\v \v f\right)
    +\nu\left(2\energy - \int_\v \v^2 f\right).
\end{gather*}
Assume $f$ minimizes entropy.
Consider a perturbation $\ft = f + \epsilon f_1$.
Then
\def\tlambda{\tilde\lambda}
\begin{align*}
  0  = d_\epsilon g|_{\epsilon=0} 
     &= \int_\v \eta' f_1
       -\lambda\int_\v f_1
       -\mu\cdot\int_\v \v f_1
       -\nu\int_\v \v^2 f_1
  \\
     &= \int_\v f_1\left(\ln f
       -\tlambda 
       -\mu\cdot\v 
       -\nu\v^2\right),
\end{align*}
where $\tlambda:=\lambda-\alpha-1$.
Since the integral must be zero for arbitrary
perturbation $f_1$ the multiplier of $f_1$ in the integrand
must be zero.  Thus, $f$ must be an exponential
of an ``isotropic'' quadratic polynomial in $\v$:
\begin{align}
\label{expquad}
  f = \exp\left(\tlambda + \mu\cdot\v +\nu\v^2\right).
\end{align}
We impose the finiteness requirement
that $\int_\v f < \infty$; that is, $\nu < 0$.

It remains to compute the moments of such a polynomial
so that we can match them up with the constrained moments.
We will show that
\begin{gather}
  \label{MaxwellianDistribution}
  \boxed{
  f=\fMaxwell:=\frac{\mdens}{(2\pi\theta)^{3/2}}\exp\left(\frac{-|\v-\u|^2}{2\theta}\right)
  }.
\end{gather}
It is evident by completing the square that any 
exponential of a quadratic polynomial in $\v$
of the form \eqref{expquad} can be put in this form.
The issue is whether we indeed have that
$\mdens = \int_\v f$,
$\mdens\u := \int_\v \v f$, and
$\mdens\theta:=\int_\v f |\c|^2/2$.
So it remains to confirm these moments by computation.

Shifting into the reference frame of the fluid,
\begin{gather*}
  \fMaxwell = \frac{\mdens}{(2\pi\theta)^{3/2}}\exp\left(\frac{-\c^2}{2\theta}\right).
\end{gather*}
It will be enough to show that:
\begin{align*}
  \int_\c \fMaxwell &= \mdens, 
  &\int_\c \c \fMaxwell &= 0, 
  &\int_\c c^2 \fMaxwell &= 3\mdens\theta.
\end{align*}
Recall how to integrate a Gaussian normal distribution:
\begin{gather*}
\int_{-\infty}^\infty e^{-x^2/2}\,dx
\int_{-\infty}^\infty e^{-y^2/2}\,dy
   = \int_0^{2\pi}\int_0^\infty e^{-r^2/2}r\,dr\,d\theta
   = 2\pi\left[e^{-r^2/2}\right]_\infty^0
   = 2\pi,
\end{gather*}
so
\begin{gather*}
\int_{-\infty}^\infty e^{-x^2/2}\,dx
   = \sqrt{2\pi},
\end{gather*}
so
\begin{gather*}
\boxed{
\int_{-\infty}^\infty \exp\left(\frac{-x^2}{2 T}\right)\,dx
   = \sqrt{2\pi T}
}.
\end{gather*}
The first moment is zero by even-odd symmetry:
\begin{gather*}
\boxed{
\int_{-\infty}^\infty x \exp\left(\frac{-x^2}{2 T}\right)\,dx
   = 0
}.
\end{gather*}
For the temperature we will need second moments.
Integrating by parts,
\begin{gather}
  \int_{-\infty}^\infty x^2 \exp\left(\frac{-x^2}{2}\right)\,dx
  = \int_{-\infty}^\infty x \left(x \exp\left(\frac{-x^2}{2}\right)\right)\,dx
  = \int_{-\infty}^\infty \exp\left(\frac{-x^2}{2}\right)\,dx
  = \sqrt{2\pi}.
  \label{squareMomentOfNormal}
\end{gather}
So
\begin{gather*}
\boxed{
\int_{-\infty}^\infty x^2 \exp\left(\frac{-x^2}{2 T}\right)\,dx
   = T \sqrt{2\pi T}
}.
\end{gather*}
For density we verify that
\begin{gather*}
     \int_\c \exp\left(\frac{-\c^2}{2\theta}\right) d^3\c
    = (2\pi\theta)^{3/2}.
\end{gather*}
For momentum we compute that
\begin{gather*}
     \int_\c c_1 \exp\left(\frac{-\c^2}{2\theta}\right) d^3\c
    = 0
\end{gather*}
by even/odd symmetry.
For temperature we compute that
\begin{gather*}
  \int_\c c^2 \exp\left(\frac{-c^2}{2\theta}\right)
  = \int_{c_1}\!\!\!
       c_1^2 \exp\left(\frac{-c_1^2}{2 \theta}\right)
    \int_{c_2}\!\!\!
       c_2^2 \exp\left(\frac{-c_2^2}{2 \theta}\right)
    \int_{c_3}\!\!\!
       c_3^2 \exp\left(\frac{-c_3^2}{2 \theta}\right)
  = 3 \theta \sqrt{2\pi \theta}.
\end{gather*}

\def\q{\mathbf{q}}
Maxwellian distributions have the property that the
heat flux $\q:=\int_\c \c c^2f$ is zero.  Indeed,
\begin{gather*}
  \q=\int_\c \c c^2 \exp\left(\frac{-c^2}{2\theta}\right) = 0,
\end{gather*}
because the integrand is odd.

\subsection{Gaussian case}

In the Gaussian case
we minimize $S$ subject to the constraints that
\begin{align*}
  \int_\v f &= \mdens,
  &\int_\v \v f &= \M,
  &\int_\v \v\v f &= \Energy.
\end{align*}
We use the technique of Lagrange multipliers.  Define
\def\ft{\tilde f}
\begin{align*}
  g :=& \int_\v \eta
    +\lambda\left(\mdens - \int_\v f\right)
    +\mu\cdot\left(\M - \int_\v \v f\right)
     +\nu\left(\Energy - \int_\v \v\v f\right).
\end{align*}
Assume $f$ minimizes entropy.
Consider a perturbation $\ft = f + \epsilon f_1$.
Then
\begin{align*}
  0 = d_\epsilon g|_{\epsilon=0} 
     &= \int_\v \eta' f_1
       -\lambda\int_\v f_1
       -\mu\cdot\int_\v \v f_1
       -\nu:\int_\v \v\v f_1
  \\
     &= \int_\v f_1\left(\ln f
       -\tlambda 
       -\mu\cdot\v 
       -\nu:\v\v\right).
\end{align*}
Since the integral must be zero for arbitrary
perturbation $f_1$ the multiplier of $f_1$ in the integrand
must be zero.  Thus, $f$ must be an exponential
of a quadratic polynomial in $\v$:
\begin{align}
\label{expquad2}
  f = \exp\left(\lambda + \mu\cdot\v +\nu:\v\v\right).
\end{align}
We may require that $\nu$ is symmetric.
We impose the finiteness requirement
that $\int_\v f < \infty$; that is, $\nu < 0$,
i.e., $\nu$ is negative definite.

We will show that
\def\u{\mathbf{u}}
\begin{gather*}
  f=\fGauss:=\frac{\mdens}{\sqrt{\det(2\pi\ThetaT)}}\exp\left(-(\v-\u)\cdot\ThetaT^{-1}\cdot(\v-\u)/2\right).
\end{gather*}
That is (shifting into the reference frame of the fluid),
\begin{gather}
  \label{GaussianDistribution}
  \boxed{
  \fGauss = \frac{\mdens}{\sqrt{\det(2\pi\ThetaT)}}\exp\left(-\c\cdot\ThetaT^{-1}\cdot\c/2\right)
  },
\end{gather}
where recall that $\c:=\v-\u$.

By substituting the expansion
$(\v-\u)\cdot\ThetaT^{-1}\cdot(\v-\u)
 = \ThetaT^{-1}:\v\v-2\u\cdot\ThetaT^{-1}\cdot\v+\u\cdot\ThetaT^{-1}\cdot\u$
and matching up with the terms in \eqref{expquad2},
it is evident that we can complete the square to put
any entropy-minimizing closure in this form.

The issue is whether we indeed have that
$\mdens = \int_\v \fGauss$,
$\mdens\u := \int_\v \v \fGauss$, and
$\mdens\ThetaT:=\int_\v \fGauss \c\c$.

It will be enough to show that
\begin{align*}
  \int_\c \fGauss &= \mdens, &
  \int_\c \c \fGauss &= 0, &
  \int_\c \c\c \fGauss &= \mdens\ThetaT.
\end{align*}
Since $\ThetaT$ is positive definite we may choose orthogonal coordinates in which
it is diagonal.  So without loss of generality
$\ThetaT = \mathrm{diag}(\theta_1, \theta_2, \theta_3)$.

For the momentum we compute that
{\small
\begin{gather*}
  \int_\c \c \exp\left(-\c\cdot\ThetaT^{-1}\cdot\c/2\right) = 0
\end{gather*}
}
because the integrand is odd.
For the density we compute that
{\small
\begin{gather*}
  \int_\c \exp\left(-\c\cdot\ThetaT^{-1}\cdot\c/2\right)
  \\
  = \int_{c_1}
       \exp\left(\frac{-c_1^2}{2 \theta_1}\right)
    \int_{c_2}          
       \exp\left(\frac{-c_2^2}{2 \theta_2}\right)
    \int_{c_3}          
       \exp\left(\frac{-c_3^2}{2 \theta_3}\right)
  \\
  = \sqrt{2\pi \theta_1}
    \sqrt{2\pi \theta_2}
    \sqrt{2\pi \theta_3}
  \\
  = \sqrt{\det(2\pi\ThetaT)}.
\end{gather*}
}
For the temperature we compute that
\begin{gather*}
  \int_\c c_1^2 \exp\left(-\c\cdot\ThetaT^{-1}\cdot\c/2\right)
  \\
  = \int_{c_1}
       c_1^2 \exp\left(\frac{-c_1^2}{2 \theta_1}\right)
    \int_{c_2}
       \exp\left(\frac{-c_2^2}{2 \theta_2}\right)
    \int_{c_3}
       \exp\left(\frac{-c_3^2}{2 \theta_3}\right)
  \\
  = \theta_1 \sqrt{2\pi \theta_1}
    \sqrt{2\pi \theta_2}
    \sqrt{2\pi \theta_3}
  \\
  = \theta_1 \sqrt{\det(2\pi\ThetaT)}
\end{gather*}
and that
\begin{gather*}
  \int_\c c_1 c_2 \exp\left(-\c\cdot\ThetaT^{-1}\cdot\c/2\right)
  \\
  = \int_{c_1}
       c_1 \exp\left(\frac{-c_1^2}{2 \theta_1}\right)
    \int_{c_2}
       c_2 \exp\left(\frac{-c_2^2}{2 \theta_2}\right)
    \int_{c_3}
       \exp\left(\frac{-c_3^2}{2 \theta_3}\right)
  \\
  = 0.
\end{gather*}

\subsubsection{Third-order moments}

Gaussian distributions have the property that the
heat flux tensor $\q:=\int_\c \c\c\c \fGauss$ is zero
(because for any component at least one of the
three independent integrals has an odd integrand).

\subsubsection{Fourth-order moments}

Fourth-order moments are needed in computing
heat flux closure with a Chapman-Enskog expansion.

This requires computing the fourth moment of a normal
distribution.  Integrating by parts,
\begin{gather*}
  \int_{-\infty}^\infty x^4 \exp\left(\frac{-x^2}{2}\right)\,dx
  = \int_{-\infty}^\infty x^3 \left(x \exp\left(\frac{-x^2}{2}\right)\right)\,dx
  = \int_{-\infty}^\infty 3 x^2\exp\left(\frac{-x^2}{2}\right)\,dx
  = 3\sqrt{2\pi}
\end{gather*}
by \eqref{squareMomentOfNormal}.  Thus,
\begin{gather*}
\int_{-\infty}^\infty x^4 \exp\left(\frac{-x^2}{2 T}\right)\,dx
   = 3 T^2 \sqrt{2\pi T}.
\end{gather*}
Representative nonvanishing fourth moments of the Gaussian distribution are
\begin{gather*}
  \int_\c c_1^4 \exp\left(-\c\cdot\ThetaT^{-1}\cdot\c/2\right)
  \\
  = \int_{c_1}
       c_1^4 \exp\left(\frac{-c_1^2}{2 \theta_1}\right)
    \int_{c_2}
       \exp\left(\frac{-c_2^2}{2 \theta_2}\right)
    \int_{c_3}
       \exp\left(\frac{-c_3^2}{2 \theta_3}\right)
  \\
  = 3\theta_1^2 \sqrt{2\pi \theta_1}
    \sqrt{2\pi \theta_2}
    \sqrt{2\pi \theta_3}
  \\
  = 3\theta_1^2 \sqrt{\det(2\pi\ThetaT)}
\end{gather*}
and
\begin{gather*}
  \int_\c c_1^2 c_2^2 \exp\left(-\c\cdot\ThetaT^{-1}\cdot\c/2\right)
  \\
  = \int_{c_1}
       c_1^2 \exp\left(\frac{-c_1^2}{2 \theta_1}\right)
    \int_{c_2}
       c_2^2 \exp\left(\frac{-c_2^2}{2 \theta_2}\right)
    \int_{c_3}
       \exp\left(\frac{-c_3^2}{2 \theta_3}\right)
  \\
  = \theta_1 \sqrt{2\pi \theta_1}
    \theta_2 \sqrt{2\pi \theta_2}
    \sqrt{2\pi \theta_3}
  \\
  = \theta_1\theta_2 \sqrt{\det(2\pi\ThetaT)}.
\end{gather*}
So we have shown that if $\Theta$ is diagonalized
along the principal axes then (for a Gaussian distribution)
$\mean{\c_1\c_1\c_1\c_1} = 3\Theta_{11}^2$
and
$\mean{\c_1\c_1\c_2\c_2} = \Theta_{11}\Theta_{22}$.
These two representative moments imply that
\begin{gather*}
  \mean{\c\c\c\c} = \SymC(\Theta\Theta)
\end{gather*}
for a Gaussian distribution.
That is,
\begin{gather}
  \boxed{\intc \c\c\c\c\fGauss = \SymC(\PT\PT)/\mdens}.
  \label{fourthMomentGauss}
\end{gather}

\subsection{Expressions for entropy}

Now that we have found the distribution that minimizes
entropy, what is the entropy?

Recall the Gaussian distribution,
\def\Gauss{\mathcal{G}}
\begin{gather*}
  \Gauss = \frac{\mdens}{\sqrt{\det(2\pi\ThetaT)}}\exp\left(\frac{-\c\cdot\ThetaT^{-1}\cdot\c}{2}\right).
\end{gather*}
By definition the entropy of the Gaussian distribution is
\begin{gather*}
  S = \int_\c \Gauss\ln \Gauss + \alpha \Gauss.
\end{gather*}
By definition,
\begin{gather*}
  \int_\c \Gauss = \mdens.
\end{gather*}
Observe that
\begin{gather*}
  \ln\Gauss =
    \ln\left(\frac{\mdens}{\sqrt{\det(2\pi\ThetaT)}}\right)
    + \frac{-\c\cdot\ThetaT^{-1}\cdot\c}{2}.
\end{gather*}
To compute $\int_\c \Gauss\ln \Gauss$ the main result we need is:
\begin{gather*}
  \int_\c \left(\c\cdot\ThetaT^{-1}\cdot\c\right)\Gauss = 3\mdens.
\end{gather*}
To verify this claim, choose coordinates in which $\ThetaT$ is diagonal.
By definition of $\theta_i$,
\begin{gather*}
  \int_\c (c_i)^2 \Gauss = \theta_i\mdens,\ \ \ \hbox{i.e.}, \ \ \ 
  \int_\c \left(c_i\theta_i^{-1}c_i\right) \Gauss = \mdens.
\end{gather*}
Summing over all three dimensions yields the claim.

We now compute the entropy:
\begin{align*}
  S &= \int_\c \Gauss\ln \Gauss + \alpha \Gauss
 \\ &= \mdens\ln\left(\frac{\mdens}{\sqrt{\det(2\pi\ThetaT)}}\right) -\frac{3}{2}\mdens + \alpha\mdens.
 \\ &= -\mdens\ln\left(\frac{\sqrt{\det(\ThetaT)}}{\mdens}\right)
   +\mdens\left(\alpha - \frac{3}{2} - \frac{3}{2}\ln(2\pi)\right).
\end{align*}
That is,
\begin{gather*}
  \boxed{
  S = -\mdens\ln\left(\frac{\sqrt{\det(\ThetaT)}}{\mdens}\right)
  }
\end{gather*}
if we choose $\alpha = 3\left(1+\ln(2\pi)\right)/2$.

The five-moment formula is a special case:
\begin{gather*}
  \boxed{
  S = -\mdens\ln\left(\frac{\theta^{3/2}}{\mdens}\right)
  }.
\end{gather*}

\def\fn{\widetilde f}
\def\Gn{\widetilde \Gauss}
\subsection{Number density}

Hitherto $f$ has represented mass density.
Let $\fn$ denote particle number density.
Then $\fn=f/m$, where $m$ is particle mass.
We define $n:=\int_\v \fn = \mdens/m$
to be the number density.
So the expression \eqref{MaxwellianDistribution}
for the 5-moment distribution becomes
\begin{gather*}
  \label{MaxwellianParticleDensity}
  \boxed{
  \Gn = \frac{n}{(2\pi\theta)^{3/2}}\exp\left(\frac{-|\v-\u|^2}{2\theta}\right)
  }
\end{gather*}
and the expression \ref{GaussianDistribution}
for the 10-moment distribution becomes
\begin{gather}
  \label{GaussianParticleDensity}
  \boxed{
  \Gn = \frac{n}{\sqrt{\det(2\pi\ThetaT)}}\exp\left(\frac{-\c\cdot\ThetaT^{-1}\cdot\c}{2}\right)
  }.
\end{gather}
\def\Pressure{\PT}

The true temperature $T=m\mean{\c^2}/3$ is related to the
scalar pressure $p = \mdens\mean{\c^2}/3$
and to the
pseudo-temperature $\theta:=\mean{\c^2}/3$ by the
relations
\begin{align*}
  n T = p = \mdens\theta,
  &&\hbox{i.e.,} && \theta = T/m.
\end{align*}

The true temperature tensor $\TT:=m\mean{\c\c}$ is related to the
pressure tensor $\Pressure = \mdens\mean{\c\c}$
and to the pseudo temperature tensor $\ThetaT:=\mean{\c\c}$ by the
relations
\begin{align*}
  n\TT = \Pressure = \mdens\ThetaT,
  &&\hbox{i.e.,} && \ThetaT = \TT/m.
\end{align*}

Note that 
\begin{align*}
  \mean{\chi} = \frac{\int_\v f\chi}{\mdens} = \frac{\int_\v \fn\chi}{n}.
\end{align*}

\subsection{Consistent entropy for interacting species}

For a gas with multiple species we should define
the entropy of each species consistently so that
the total entropy obeys an entropy inequality when
species interact.  For such a consistent entropy
we define the true entropy of each species in terms of the
number density rather than the mass density:
\def\etan{\widetilde\eta}
\def\alphan{\widetilde\alpha}
\def\Sn{\widetilde S}
\def\minv{m^{-1}}
\begin{align*}
  \Sn &:= \int_\v \etan, & \hbox{where} &&
  \etan &:= \fn\ln \fn + {\cyan \alphan} \fn.
\end{align*}
In the remainder of this section the casual reader may regard
factors involving $\cyan \alpha$ or $\cyan\alphan$
(cyan text) as an arbitrary irrelevant constant.
Since $\fn=\minv f$,
\begin{align*}
  \etan &= \fn\ln(\minv f) + {\cyan \alphan} \fn
 \\     &= \minv \eta + {\cyan (\ln \minv + \alphan - \alpha)} \fn.
\end{align*}
So
\begin{align*}
  \Sn &= \int_\v \etan
       = \int_\v \minv \eta + {\cyan (\ln \minv + \alphan - \alpha)}\fn
   \\ &= \minv S + {\cyan (\ln \minv + \alphan - \alpha)} n.
\end{align*}
Recall that for $f=\Gauss$,
\begin{gather*}
  \!\!\!  \!\!\!  \!\!\!
  S = -\mdens\ln\left(\frac{\sqrt{\det(\ThetaT)}}{\mdens}\right)
   +\mdens{\cyan\left(\alpha - \frac{3}{2} - \frac{3}{2}\ln(2\pi)\right)}.
\end{gather*}
So for $\fn=\Gn$,
\begin{align*}
  \Sn = &-n\ln\left(\frac{\sqrt{\det(\ThetaT)}}{\mdens}\right)
         +n{\cyan \left(\alpha - \frac{3}{2} - \frac{3}{2}\ln(2\pi)\right)}
   \\   &+n{\cyan \left(\ln \minv + \alphan - \alpha\right)}
 \\ \phantom{\Sn}
      = &-n\ln\left(\frac{\sqrt{\det(\TT)}}{n}\right)
   \\   &+n{\cyan \left(\alphan
            + \frac{3}{2}\ln m - \frac{3}{2} - \frac{3}{2}\ln(2\pi)\right)}.
\end{align*}
In summary, the \defining{ten-moment gas-dynamic entropy},
defined to be the entropy that the distribution would have if
it were relaxed to minimum entropy subject to the constraint
that all moments of order two or lower are conserved,
may be consistently defined to be
\begin{gather}
  \label{GaussianEntropy}
  \boxed{
  S_\Gauss = -n\ln\left(\frac{\sqrt{\det(\TT)}}{n}\right)
  },
\end{gather}
where we have chosen
${\cyan \alphan = \frac{3}{2}\left(1 + \ln(2\pi/m)\right)}$.

The \defining{five-moment gas-dynamic entropy},
defined to be the entropy that the distribution would have if
it were relaxed to minimum entropy subject to the constraint
that all conserved moments are conserved,
is a special case and may be consistently defined to be
\begin{gather}
  \label{MaxwellianEntropy}
  \boxed{
  S_\Maxwell = -n\ln\left(\frac{T^{3/2}}{n}\right)
  }.
\end{gather}

\section{Collisionless Nondimensionalization}
\label{CollisionlessNondimensionalization}

Physical constants that define an ion-electron plasma are:
\begin{enumerate}
\item $e$, the magnitude of the charge of an electron,
\item $m_i$, $m_e$, the ion and electron mass, and
\item $c$, the speed of light.
\end{enumerate}

Three fundamental parameters that characterize the
state of a plasma are:
\begin{enumerate}
\item $n_0$, a typical particle density,
\item $T_0$, a typical temperature (often per species), and
\item $B_0$, a typical magnetic field strength.
\end{enumerate}
In quasineutral equilibrium
we can take
$n_0 = n_i = n_e$
and
$T_0 = T_i = T_e$.
The thermal pressure is $p_0:=n_0 T_0$
and the magnetic pressure is
$p_B := \frac{B_0^2}{2\mu_0}$.

\def\vA{v_A}
\def\vAs{v_{A,s}}
\def\rgs{r_{g,s}}
\def\Rgs{\tilde r_{g,s}}
\def\omegags{\omega_{g,s}}
\def\omegaps{\omega_{p,s}}
\def\vt{v_t}
\def\vts{v_{t,s}}
\def\Vt{\tilde v_{t}}
\def\Vts{\tilde v_{t,s}}
\def\vAi{v_{A,i}}
\def\rgi{r_{g,i}}
\def\Rgi{\tilde r_{g,i}}
\def\omegagi{\omega_{g,i}}
\def\omegapi{\omega_{p,i}}
\def\vti{v_{t,i}}
\def\Vti{\tilde v_{t,i}}
Subsidiary space, time, and velocity
scale parameters derived from the fundamental parameters are
\begin{gather*}
 \begin{aligned}
  \hbox{\defining{gyrofrequencies}:} &&
  \omegags &:= \frac{e B_0}{m_s},
  \\
  \hbox{\defining{plasma frequencies}:} &&
  \glslink{plasmafreq}{\omegaps}^2 &:= \frac{n_0 e^2}{\epsilon_0 m_s},
  \\
  \hbox{\defining{Alfv\'en speeds}:} &&
  \vAs^2 &:= \frac{B_0^2}{\mu_0 m_s n_0} = \frac{2 p_B}{\mdens_s},
  \\
  \hbox{\defining{thermal velocities}:} &&
  \glslink{vt}{\vts}^2 &:= \frac{T_s}{m_s} = \frac{p_s}{\mdens_s},
     \ \ \ \ \Vts^2:=2\vts^2,
  \\
  \hbox{\defining{gyroradii}:} &&
  \rgs &:= \frac{\vts}{\omegags} = \frac{m_s \vts}{e B_0},
     \ \ \ \ \Rgs  := \frac{\Vts}{\omegags},
  \\
  \hbox{\defining{Debye length}:} &&
  \glslink{DebyeLength}{\ZdebyeLength}^2
  &:= \left(\frac{v_{t,s}}{\omegaps}\right)^2
  = \frac{\epsilon_0 T_0}{n_0 e^2},
  \\
  \hbox{\defining{inertial lengths}:} &&
  \glslink{delta}{\Zdelta_s}^2 &:= \left(\frac{c}{\omegaps}\right)^2
  = \left(\frac{\vAs}{\omegags}\right)^2
  = \frac{m_s}{\mu_0 n_s e^2};
 \end{aligned}
\end{gather*}
the inertial length is also called the \defining{skin depth}.
Nondimensional parameters for the bulk fluid
of a two-species quasi-neutral fluid are
\begin{gather*}
 \begin{aligned}
  \hbox{\defining{plasma frequency}:} &&
  \glslink{plasmafreq}{\omega_p}^2 &:= \omega_{p\e}^2 + \omega_{p\i}^2,
  \\
  \hbox{\defining{Alfv\'en speed}:} &&
  \glslink{vA}{\ZvA}^2 &:= \frac{B_0^2}{\mu_0 \mdens} = \frac{2 p_B}{\mdens},
  \\
  \hbox{\defining{thermal velocity}:} &&
  \vt^2 &:= \v_{t,\i}^2+\v_{t,\e}^2 = \frac{T}{\mred},
     \ \ \ \ \Vt^2:=2\vt^2,
  \\
  \hbox{\defining{inertial length}:} &&
  \Ldelta^2 &:= \Ldelta_\i^2 + \Ldelta_\e^2.
 \end{aligned}
 \label{AlfvenSpeed}
\end{gather*}
Note that \note{(most?)} often in the literature
the thermal velocity is defined as $\Vts$
rather than $\vts$.
We say that two parameters are \emph{equivalent} if 
one is a constant multiple of the other.
For example, the thermal velocities are equivalent to
one another and to
the sound speed $\sqrt{\frac{\gamma p_0}{\mdens_0}}$.
Important nondimensional ratios are
the \defining{plasma beta} $\beta := \frac{p_0}{p_B}$
and the ratio of the speed of light to the Alfv\'en
speed.
Other nondimensional ratios can be
defined in terms of these ratios:
\def\vt{v_t}
\def\Vt{\tilde v_t}
\begin{gather*}
 \begin{aligned}
  \hbox{plasma $\beta$}: && 
       \beta &:= \frac{p_0}{p_B}
    = \left(\frac{\Vts}{\vAs}\right)^2
    = \left(\frac{\Rgs}{\Ldelta_s}\right)^2,
\\
  \hbox{\note{unnamed?}}: && 
       \frac{c}{\vAs} &= \frac{\rgs}{\lambda_D} = \frac{\omegaps}{\omegags}.
 \end{aligned}
\end{gather*}

The subsidiary parameters (except for the temperature-related
parameters $\vts$ and $\lambda_D$) emerge from
a generic nondimensionalization
of the particle (or Vlasov or 2-fluid)
equations.

Choose values for: 
\begin{center}
\begin{tabular}{l l l}
  $t_0$ & (time scale) & (e.g. ion gyroperiod $1/\omegagi$), \\
  $x_0$ & (space scale) & (e.g. ion skin depth $\Ldelta_i$), \\
  $m_0$ & (mass scale) & (e.g. ion mass $m_i$), \\
  $e=q_0$ & (charge scale) & (e.g. ion charge $e$), \\
  $B_0$ & (magnetic field) & (e.g. $\omegagi m_i/e$), and \\
  $n_0$ & (number density) & (e.g. something $\gg 1/x_0^3$).
\end{tabular}
\end{center}
This implies typical values for: 
\begin{center}
\begin{tabular}{l l}
  $v_0 = x_0/t_0$ & (velocity), \\
  $E_0 = B_0 v_0$ & (electric field), \\
  $\sigma_0 = e n_0$ & (charge density), \\
  $J_0 = e n_0 v_0$ & (current density), and \\
  $S_0 = n_0$ & (no. particles per unit number density). \\
\end{tabular}
\end{center}

Making the substitutions
\def\alt{\magenta\widetilde}
\def\tt{{\magenta\tilde t}}
\def\tx{{\alt\Zxb}}
\def\tq{{\alt q}}
\def\tm{{\alt m}}
\def\tB{{\alt\ZB}}
\def\tn{{\alt n}}
\def\tE{{\alt\ZE}}
\def\tJ{{\alt\ZJ}}
\def\tsigma{{\alt\sigma}}
\def\tS{{\alt S}}
\def\tc{{\alt c}}
\def\tv{{\alt \Zv}}
\def\tnabla{\alt\nabla}
\begin{align*}
  t &=\tt t_0, & \E        &=\tE          B_0 v_0,  \\
 \xb&=\tx x_0, & \sigma    &=\tsigma      e n_0,    \\
  q &=\tq e  , & \J        &=\tJ          e n_0 v_0,\\
  m &=\tm m_0, &  S_p(\xb_p) &=\tS_p(\tx_p) n_0, \\
  n &=\tn n_0, &  c        &=\tc          v_0,      \\
 \B &=\tB B_0, & \v        &=\tv          v_0 
\end{align*}
in the fundamental equations
\def\xcurl{\nabla_{\xb}\times}
\def\xDiv{\nabla_{\xb}\cdot}
\def\tcurl{\nabla_{\tx}\times}
\def\tDiv{\nabla_{\tx}\cdot}
\begin{align*}
   &\partial_{t} \B = -\xcurl \E,                    & \xDiv\B&=0,
\\ &\partial_t \E = c^2\xcurl \B - \J/\epsilon_0, & \xDiv\E&=\sigma/\epsilon_0,
\\ &\J = \sum_p S_p(\xb_p) q_p v_p, 
   &\sigma &= \sum_p S_p(\xb_p) q_p,
\end{align*}
and
\begin{align*}
    d_t (\gamma\v_p) &= \frac{q_p}{m_p}
      \Big(\E(\xb_p)+\v_p\times\B(\xb_p)\Big),
 & d_t \xb_p &= \v_p
\end{align*}
gives the almost identical-appearing nondimensionalized system
\def\neps{{\color{cyan} \epsilon}}
\begin{align*}
   &\partial_\tt \tB = -\tcurl \tE,                & \tDiv\tB&=0,
\\ &\partial_\tt \tE = \tc^2\tcurl \tB - \tJ/\neps,   & \tDiv\tE&=\tsigma/\neps,
\\ &\tJ = \sum_p \tS_p(\tx_p) \tq_p \tv_p, 
   &\tsigma &= \sum_p \tS_p(\tx_p) \tq_p,
\end{align*}
and
\begin{align*}
    d_\tt (\gamma\tv_p) &= {\color{cyan} (t_0\omega_g)} \frac{\tq_p}{\tm_p}
      \Big(\tE(\tx_p)+\tv_p\times\tB(\tx_p)\Big),
 & d_\tt \tx_p &= \tv_p;
\end{align*}
here $\cyan{(t_0\omega_g)} = t_0 \frac{q_0 B_0}{m_0}$
is the gyrofrequency nondimensionalized by a choice
of $t_0$ (which can be chosen to be the gyroperiod
in order to set this factor to unity) and
\begin{gather*}
 \frac{1}{\neps} = \frac{x_0 n_0 e}{v_0 B_0 \epsilon_0}
   = t_0 \frac{e B_0}{m_0} \frac{\mu_0 m_0 n_0}{B_0^2} c^2
   = (t_0 \omega_g) \left(\frac{c}{v_A}\right)^2.
\end{gather*}
Note that we can also write
$(t_0 \omega_g) = \frac{x_0}{r_g}$.

It is desirable to make the nondimensionalized system
look exactly like the original system; in this way
we can take any formula derived from the original
system and interpret it as a nondimensional formula.
If the gyrofrequency is not chosen to be the gyroperiod
one can still accomplish this by absorbing the factor
$(t_0 \omega_g)$ into the electromagnetic field and
the definition of $\neps\inv$.

\section{Knudsen Number}

\section{Prandtl Number}
\label{PrandtlNumber}

The \first{Prandtl number} $\glslink{Pr}{\Pr}$ of a gas (unmagnetized)
is the rate of momentum diffusion divided by the rate of temperature diffusion
for small perturbations from a global equilibrium.
To obtain an expression for the Prandtl number we need evolution
equations for the momentum and temperature.

Since perturbations from equilibrium are small, we can
make simplifying approximating assumptions.
In particular, the pressure is assumed constant and
for the momentum equation
the flow is assumed approximately incompressible.
Near a Maxwellian linearized entropy-respecting closures are justified.

\subsection{Rate of momentum diffusion}

To calculate the rate of momentum diffusion we use the Stokes closure,
derived in section \ref{isoViscStressClosure},
for the momentum evolution equation \eqref{momentumEvolution},
which we here write in the form
\begin{gather}\label{momEvolution}
  \mdens d_t \u + \D p = \Div\dstress;
\end{gather}
here $\stress = -\PT$ is the stress tensor
and  $\dstress = -\dPT = -(\PT-p\id)$
is its deviatoric part, also called the viscous stress tensor.
The Stokes closure \eqref{dPTclosure} says that
\begin{gather*}
  \dstress = 2\viscosity\dstrain,
\end{gather*}
where \gls{viscosity} is the viscosity coefficient and
$\glslink{strain}{\strain}=\Sym(\D\u)$ is the strain rate tensor and
$\glslink{dstrain}{\dstrain}=\Sym(\D\u-\Div\u\id/3)$ is its deviatoric part.
Assuming constant pressure ($\D p=0$) and incompressible
flow ($\Div\u=0$), the momentum equation simplifies to
\def\kinematicViscosity{\nu}
\begin{gather}\label{kinematicViscosity}
  d_t \u = \kinematicViscosity\laplacian\u,
\end{gather}
where $\kinematicViscosity:=\frac{\viscosity}{\mdens}$,
the \defining{kinematic viscosity}, is the rate of momentum
diffusion.

\subsection{Rate of temperature diffusion}

To calculate the rate of temperature diffusion
we derive a temperature evolution equation.
Rather than begin with the temperature evolution
equation \eqref{Tevolution},
which assumes that molecules have no internal energy modes,
only translational modes,
we begin with scalar energy balance \eqref{scalarEnergyEvolution}
in the form
\def\tnrg{\nrg_\mathrm{int}}
\begin{gather}\label{tnrgEvolution}
  \mdens d_t(\tnrg+|\u|^2/2)+\Div(\u p) + \Div\q = \Div(\u\cdot\dstress),
\end{gather}
where $\tnrg$ is the internal (nontranslational) energy per mass.
Dotting $\u$ with the momentum evolution equation \eqref{momEvolution}
gives the kinetic energy evolution equation
\begin{gather*}
  \mdens d_t |\u|^2/2 + \u\cdot\D p + \Div\q = \dstress:\D\u;
\end{gather*}
subtracting this from energy evolution \eqref{tnrgEvolution}
gives thermal energy evolution
\begin{gather}\label{tnrgEvolution}
  \mdens d_t \tnrg+p\Div\u + \Div\q = \dstress:\D\u.
\end{gather}

Assume the classical thermodynamic relations
for an ideal gas,
\def\cv{c_v}
\begin{align*}
  p&=\mdens R T, &\hbox{and}&& \tnrg&=\cv T,
\end{align*}
where $R$ is the gas constant (equal to $1/m$ for a gas with a single
species of particle with particle mass $m$)
and $\cv$ is the heat capacity at constant volume.
Assume the linear heat flux closure
$ 
  \q=-\heatConductivity\D T,
$ 
in accordance with equation \eqref{scalarHeatFluxClosure}.
Assume that $k$ is constant.
(In fact $k$ is a function of temperature,
which is approximately constant.)
Neglect viscous heat production $\dstress:\D\u.$
{\em Assume that the pressure $p$ is constant.}
(So we do \emph{not} assume that $\Div\u=0$.)
Then $p=\mdens R T$ says that $\mdens T$ is constant,
so $p\Div\u = -p d_t\ln\mdens = p d_t\ln T = \mdens R d_t T$,
so thermal energy evolution \eqref{tnrgEvolution}
reduces to 
\begin{gather*}
  \mdens \cv d_t T + \mdens R d_t T = \heatConductivity\laplacian T,
\end{gather*}
that is,
\def\thermalConductivity{\kappa}
\def\cp{c_p}
\begin{align}\label{thermalConductivity}
  d_t T &= \thermalConductivity\laplacian T,
\end{align}
where
$ \thermalConductivity := \frac{\heatConductivity}{\mdens\cp} $,
the thermal conductivity, is the rate of temperature diffusion
and $\cp:=\cv+R$ is the heat capacity at constant pressure.

\subsection{Formula for the Prandtl number}

Putting the results of equations
\eqref{kinematicViscosity}
and
\eqref{thermalConductivity}
together,
the \mention{Prandtl number} is
\begin{align*}
  \boxed{
  \Pr = \frac{\kinematicViscosity}{\thermalConductivity}
      = \frac{\cp\viscosity}{\heatConductivity}
      = \frac{\gamma}{\gamma-1}\frac{R\viscosity}{\heatConductivity}
  },
\end{align*}
where $\gamma:=\cp/\cv$ is the \defining{adiabatic index};
that is,
\begin{align*}
    \Pr = \frac{\gamma}{\gamma-1}\frac{\viscosity}{m\heatConductivity}
\end{align*}
for an ideal gas.  For a monatomic gas, $\gamma-1=R/\cv = 2/3$,
$\gamma=5/3$, and 
\begin{align}
  \glslink{Pr}{\Pr} = \frac{5}{2}\frac{\viscosity}{m\heatConductivity},
  \label{PrFormula}
\end{align}
which is the formula assumed in this dissertation.

\section{Explicit closures}
\label{ExplicitClosures}
{

In this section we work out explicit five-moment closures for
the heat flux and viscosity and an explicit ten-moment closure
for the heat flux tensor beginning with the implicit closures
obtained in section \ref{PerturbativeClosure}.

\subsection{Splice tensor operators}
\label{SpliceTensorOperators}

To solve these equations, it will be convenient to define
\defining{splice tensor operations} which operate on tensors
with an even number of indices.
To define these operators, for a given tensor
we partition the indices into the initial half and
the final half and pair corresponding indices in the
initial and final half, like this:
\begin{align*}
 A_{i_1 i_2\cdots i_m|j_1 j_2\cdots j_m};
\end{align*}
as here, for clarity we sometimes insert the symbol $|$ to separate the initial
and final half of the indices of a tensor.
Splice operators do exactly the same thing to the
initial indices and final indices and can be thought of
as operating on pairs of indices.
An identity in terms of splice operators becomes an identity in
terms of standard tensor operators if you replace splice operators with
their corresponding standard tensor operators and
delete the initial half (or the final half) of the indices
of each tensor.  For a given ordinary tensor operator 
(e.g.\ $\otimes$) we denote
its splice operator equivalent with a superior tilde symbol
(e.g.\ $\widetilde{\otimes}$).
Examples of splice products for simple cases are
\begin{align*}
 (A\stimes B)_{ijkl} &= A_{ik} B_{jl}&
 \hbox{and}&&
 (K\stimes L)_{i j_1 j_2 k l_1 l_2} &= K_{i k} L_{j_1 j_2 l_1 l_2}.
\end{align*}
So in general the splice tensor product $\glslink{stimes}{\Zstimes}$
is defined by
\begin{align*}
 (K\stimes L)_{i_1\cdots i_m j_1\cdots j_n|k_1\cdots k_m l_1\cdots l_n}
   &= K_{i_1\cdots i_m|k_1\cdots k_m} L_{j_1\cdots j_n|l_1\cdots l_n}.
\end{align*}
Recall that we define the symmetric tensor product by
$A\veebar B:=\Sym(A\otimes B)$, where $\Sym$ averages over
all permutations of its argument tensor.
Similarly, we define a \defining{splice symmetric tensor product}
$\glslink{sveebar}{\Zsveebar}$ by
\begin{gather*}
 A\sveebar B = \sSym(A\stimes B),
\end{gather*}
where $\glslink{sSym}{\ZsSym}$ averages over all permutations
which permute the initial half of the indices and
the final half of the indices in exactly the same manner.

In this section all tensors will be of even order
and will be built from splice symmetric 
products of second-order tensors exactly as one can build
standard symmetric tensors from symmetric products of vectors.
Elsewhere in this dissertation the default product of tensors
is the simple tensor product $\otimes$, but
\emph{\bf in this section we will take the default tensor
product to be the splice symmetric product $\sveebar$},
and when we wish to denote the splice symmetric product
explictly we will write it simply as $\o$ rather than $\sveebar$.
Let $A$, $B$, and $C$ be second-order tensors.

Examples of splice symmetric products are
$2AB = A\stimes B + B\stimes A$ and
\begin{align*}
  3!ABC = &A\stimes B\stimes C + A\stimes C\stimes B  
      + B\stimes A\stimes C + B\stimes C\stimes A  
      + C\stimes A\stimes B + C\stimes B\stimes A.
\end{align*}

\subsection{Gyrotropic tensors}

\def\i{\bm{\delta}}
\def\w{{\bm{\wedge}}}
\def\l{{\bm{\mypara}}}
\def\p{{\bm{\perp}}}
For ease on the eyes, in this section
we use $\i:=\id$ for the identity matrix.
(Its components are given by the Kronecker delta.)
We will build gyrotropic basis tensors by 
splice symmetric products of the following
fundamental gyrotropic tensors:
\begin{alignat*}{7}
     \w &:= \i_\wedge &&:= \b\times\i = \i\times\b,
  \\ \l &:= \i_\parallel &&:= \b\b,
  \\ \p &:= \i_\perp &&:= \i-\i_\parallel.
\end{alignat*}
I remark that this does not generate a basis that spans all
gyrotropic tensors --- for that we would need to take ordinary tensor
products and splice products --- but we will see that it does
generate exactly the basis needed to solve our implicit closures
for the unknown.

Under the operation of matrix multiplication
the fundamental gyrotropic tensors plus the zero tensor plus $\i$
are a commutative subring with unity with the following multiplication table:
\begin{equation}
 \begin{array}{c|r r r}
    \dotp & \w & \p & \l
    \\ \hline
       \w & -\p & \w & 0
    \\ \p & \w & \p & 0
    \\ \l & 0 & 0 & \l
 \end{array}
 \label{multiplicationTable}
\end{equation}
In particular, we will use the following mappings:
\begin{equation}
 \begin{array}{c|r r}
    M & \w\dotp M & \i\dotp M
    \\ \hline
       \w & -\p & \w
    \\ \p & \w & \p
    \\ \l & 0 & \l
 \end{array}
 \label{basicMappings}
\end{equation}

\subsection{Gyrotropic linear operators on symmetric matrices}

Splice symmetric products allow us to express gyrotropic
linear operators on symmetric matrices.

Recall that matrices are a group under the multiplication $\dotp$.
Fourth-order tensors are a group under the multiplication $\ddotp$
and sixth-order tensors are a group under the multiplication $\dddotp$.

Just as $\i$ is the identity matrix,
the tensor $\i\i = \i\stimes\i$ acts
as the identity tensor on second-order
tensors and $\i\i\i = \i\stimes\i\stimes\i$ acts
as the identity tensor on third-order tensors:
\def\AA{\tensorb{A}}
\def\AAA{\tensorc{A}}
\begin{gather*}
\begin{aligned}
  \i\i\ddotp\AA &= \AA,& &\hbox{and}&
  \i\i\i\dddotp\AAA &= \AAA,
\end{aligned}
\end{gather*}

\def\V{\underline{V}}
\def\AA{\tensorb{M}}
\def\VV{V_{[2]}}
\def\AAAA{M_{[2]}}
\def\VVV{V_{[3]}}
\def\AAAAAA{M_{[3]}}
Recall that any linear transformation on vectors (first-order tensors)
has a matrix (second-order tensor) $\AA$ which represents it by $\V\mapsto \AA\dotp\V$
and a unique inverse matrix $\AA\inv$ satisfying $\AA\inv\dotp\AA=\i=\AA\dotp\AA\inv$.
Likewise, any linear transformation on second-order tensors has a
``second-order matrix'' (fourth-order tensor) $\AAAA$ which represents it by
$\VV\mapsto \AAAA\ddotp\VV$
and a unique inverse $\AAAA\inv$ satisfying
$\AAAA\inv\ddotp\AAAA=\i\i=\AAAA\ddotp\AAAA\inv$,
and any linear transformation on third-order tensors has a
``third-order matrix'' (sixth-order tensor) $\AAAAAA$ which represents it by
$\VVV\mapsto \AAAAAA\dddotp\VVV$.
and a unique inverse $\AAAAAA\inv$ satisfying
$\AAAAAA\inv\dddotp\AAAAAA=\i\i\i=\AAAAAA\dddotp\AAAAAA\inv$.

\subsection{Implicit closures in splice product form}

The implicit closures obtained in section \ref{PerturbativeClosure}
assuming a Gassian-BGK collision operator and using a Chapman-Enskog expansion were
equation \eqref{implicitHeatFluxClosure}
\begin{gather}
  \label{implicitHeatFluxClosureCopy}
   \q + \hgf\b\times\q = -\heatConductivity \D T
\end{gather}
for the heat flux,
equation \eqref{implicitDeviatoricStressClosure}
\begin{gather}
  \dPT + \SymB(\pgf\b\times\dPT) = -\viscosity 2 \dstrain
  \label{implicitDeviatoricStressClosureCopy}
\end{gather}
for the deviatoric stress,
and equation \eqref{implicitHeatFluxTensorClosure}
\begin{gather}
  \label{implicitHeatFluxTensorClosureCopy}
  \qT + \SymC(\hgf\b\times\qT) =
    -\tfrac{2}{5}\heatConductivity \SymC\left(\Pshape\dotp\D\TT\right)
\end{gather}
for the heat flux tensor.

\def\v{\pgf}
Since we will solve each of these equations individually,
we will neglect the distinction between $\hgf$
and $\pgf$ that arises when $\Pr\ne 1$ and will simply
write $\v$ for both.

Using that
$
  \SymB(\b\times\dPT) = \SymB(\w\dotp\dPT) 
                      = \SymB(\w\dotp\dPT\dotp\i)
                      = \SymB(\w\stimes\i \ddotp\dPT)
                      = 2\w\i \ddotp\dPT
$
and that
$
  \SymC(\b\times\qT) = \SymC(\w\dotp\qT) 
                     = \SymC(\w\i\i\dddotp\qT) 
                      = 3\w\i\i \dddotp\qT,
$
we can rewrite the implicit closure equations in the form
\begin{gather}
     (\i+\v\w)\dotp\q = -\heatConductivity \D T,
     \label{impHeatFlux}
  \\ (\i\i+\v 2\w\i)\ddotp\dPT = -\viscosity 2 \dstrain,
     \label{impStress}
  \\ (\i\i\i+\v 3\w\i\i)\dddotp\qT = 
    -\frac{2}{5}\heatConductivity \SymC\left(\Pshape\dotp\D\TT\right).
     \label{impStress}
\end{gather}
To solve these equations we need to invert the matrices in parentheses.

\subsection{Heat flux}
\label{HeatFlux}

We need to solve equation \eqref{impHeatFlux} for $\q$.
Let $M_1$ denote the inverse of the matrix $A_1:=\i+\v \w$.
It must be gyrotropic. The tensors $\l$, $\p$, and $\w$ comprise
a basis for the space of gyrotropic tensors.
So we can expand $M_1$ as
\begin{gather*}
  M_1 = m_n \l + m_0\p + m_1 \w.
\end{gather*}
Substituting this into the required relation $A_1\dotp M_1=\i$
yields
\begin{gather*}
  m_n \l + m_0\p + m_1 \w
  + \v(m_0\w - m_1 \p)
 = \l + \p,
\end{gather*}
where we have used that $\i=\l+\p$ and that
$\w\dotp\p=\w$ and $\w\dotp\w=-\p$.
Matching coefficients reveals that $m_n=1$ and gives
a linear system
\begin{gather}
  \begin{matrix}
    \p: \\
    \w: 
  \end{matrix}
  \left(
   \begin{bmatrix}
     1 & 0 \\
     0 & 1
   \end{bmatrix}
   + \v
   \begin{bmatrix}
     0 & -1 \\
     1 & 0
   \end{bmatrix}
  \right)
  \begin{bmatrix}
    m_0 \\
    m_1
  \end{bmatrix}
  = 
  \begin{bmatrix}
    1 \\
    0
  \end{bmatrix}
  \Longrightarrow
  \begin{bmatrix}
    m_0 \\
    m_1
  \end{bmatrix}
  = \frac{1}{1+\v^2}
  \begin{bmatrix}
    1 \\
    -\v
  \end{bmatrix}.
  \label{systemM}
\end{gather}
That is,
\begin{gather}
  M_1 = \theatConductivity := \l + \frac{1}{1+\v^2}\Big(\p - \v \w\Big).
  \label{theatConductivityClosure}
\end{gather}
Observe that when $\v=0$ then $M_1$ is the identity matrix,
and in the limit of large magnetic field,
\begin{gather*}
  M_1 \approx \l - \v\inv \w,
\end{gather*}
which effectively shuts down heat flux perpendicular
to the magnetic field.
So the heat flux closure is
\begin{gather}
 \q = -\heatConductivity\left(\l + \frac{1}{1+\v^2}\Big(\p - \v \w\Big)\right)\cdot\D T,
 \label{expHeatFluxInternal}
\end{gather}
that is,
\begin{gather}
 \q = -\heatConductivity\left(\b\b + \frac{1}{1+\v^2}\Big((\id-\b\b) - \v \b\times\id\Big)\right)\cdot\D T.
 \label{expHeatFlux}
\end{gather}
This is the closure given by Woods in equations (5.107)--(5.109) in \cite{book:Woods04}.
He identifies this result as that of Chapman and Cowling's kinetic theory.
His figure 5.6 plots $m_0/m_n$ versus $\v$
and shows that it agrees well those of Braginskii's closure
given on pages 249--51 of \cite{article:Braginskii65}.

\def\AB{A_{[2]}}
\def\AC{A_{[3]}}
\def\MB{M_{[2]}}
\def\MC{M_{[3]}}
\def\SpanMi{\Span\{M_i\}}
\def\Linf{L^\infty}
\subsection{Deviatoric pressure tensor}
\label{DeviatoricPressureTensor}
Let $\MB$ denote the inverse of the matrix $\AB:=\i\i+\v 2\w\i$.
We will seek $\MB$ as a linear combination of basis tensors $\{M_i\}$,
$\MB = \sum_i m_i M_i$, where $\SpanMi$ is closed under the map
$M\mapsto \AB\ddotp M$, that is, which is closed under the map
$L:M\mapsto 2\w\i\ddotp M$.  We need that $(\i\i+\v L)\ddotp \MB = \i\i$,
so it is evident that $\i\i$ should be in $\SpanMi$,
so we simply compute what repeated applications of $L$ can generate
starting with $\i\i$.
First, observe that $L$ satisfies
\begin{gather*}
  2\w\i\ddotp XY = \w\dotp XY
                 + \w\dotp YX.
\end{gather*}
Using that $\w\dotp\w=-\p$ and $\w\dotp\p=\w$,
we generate a basis.  Under the mapping $M\mapsto 2\w\i\ddotp M$,
the calculations
\begin{gather*}
 \begin{array}{r r r r r}
    M'^0: &M'^0_0:=&\i\i &\mapsto &2\w\i 
 \\ M'^1: &M'^1_0:=&\p\i &\mapsto &\w\p +\w\i
 \\       &M'^1_1:=&\w\i &\mapsto &\w\w  -\p\i 
 \end{array}
\end{gather*}
and
\begin{gather*}
 \begin{array}{r r r r r}
    M^2:  &M^2_0:=&\p\p &\mapsto & 2\w\p
 \\       &M^2_1:=&\w\p &\mapsto & \w\w - \p\p
 \\       &M^2_2:=&\w\w &\mapsto & -2\w\p
 \end{array}
\end{gather*}
exhibit such a basis.
The span of the $M^2$ subsystem is closed under the map
$M\mapsto \AB\dddotp M$, but the span of the $M^1$
subsystem is not.  This is easily remedied.
Subtract the map of $\p\p$ from the map of $\p\i$
and subtract the map of $\w\i$ from the map of $\w\p$ to get
the closed system
\begin{gather*}
 \begin{array}{r r r r r}
     M^1: &M^1_0:=&\p\l &\mapsto &\w\l
 \\       &M^1_1:=&\w\l &\mapsto &-\p\l 
 \end{array}
\end{gather*}
Similarly, beginning with the map of $\i\i$
we first subtract the map of $\p\i$
to get a map of $\i\l$ and then
subtract the maps of $\p\l$ to get the map
\begin{gather*}
 \begin{array}{r r r r r}
     M^0: &M^0_0:=&\l\l &\mapsto &0.
 \end{array}
\end{gather*}
The basis elements $M^i_j$ generate the same space
as the basis elements $M'^i_j$ and comprise a decoupled
set of cycles.

On the $M^1$ subsystem $L$ is an invertible map,
but on $M^2$ we can separate out a one-dimensional null-space.
Basis vectors for a decoupled system are
\def\sd{\sigma_d}
\def\ss{\sigma_s}
\def\sk{\sigma_k}
\begin{gather*}
 \begin{array}{r r r r r r}
    L^2:  &2\sd :=&\p\p-\w\w &\mapsto & 4\w\p
 \\       &M^2_1:=&\w\p &\mapsto & \w\w - \p\p & = -2\sd
 \\ N^2:  &2\ss :=&\p\p+\w\w &\mapsto & 0
 \end{array}
\end{gather*}
We now decompose $\i\i$ into a sum of the decoupled basis vectors:
\(
  \i\i = (\l + \p)(\l + \p)
       = \l\l + \p\p + 2\p\l
       = \left(\l\l + \frac{\p\p + \l\l}{2}\right) + \frac{\p\p - \l\l}{2} + 2\p\l,
\) that is,
\begin{gather*}
  \i\i = \sk + 2\p\l + \sd,
\end{gather*}
where $\sk:= \left(\l\l + \frac{\p\p + \l\l}{2}\right)$;
note that there is one element from the null space
and one element from each invertible subsystem.
The matrix $M$ is a linear combination of the basis
we have defined:
\begin{gather*}
   M = s_k\sk + 2(m^1_0 \p\l + m^1_1 \w\l) +  (s_d\sd +  m^2_1 \w\p),
\end{gather*}
where we use parentheses for components that correspond to decoupled subsystems.
To find the five unknown viscosity coefficients we substitute this into
the required identity
$
   (\i\i+\v\w\i)\ddotp M = \i\i,
$
that is,
\begin{align*}
    s_k\sk &+ 2(m^1_0 \p\l + m^1_1 \w\l) +  (s_d\sd +  m^2_1 \w\p)
 \\        &+ \v\left[2(m^1_0 \w\l - m^1_1 \p\l) +  (2s_d\w\p - 2 m^2_1 \sd)\right]
 \\      = &\sk + 2\p\l + \sd.
\end{align*}
Matching up coefficients gives $ s_k = 1 $
and two linear systems,
\begin{gather}
  \begin{matrix}
    \p\l: \\
    \w\l: 
  \end{matrix}
  \left(
   \begin{bmatrix}
     1 & 0 \\
     0 & 1
   \end{bmatrix}
   + \v
   \begin{bmatrix}
     0 & -1 \\
     1 & 0
   \end{bmatrix}
  \right)
  \begin{bmatrix}
    m^1_0 \\
    m^1_1
  \end{bmatrix}
  = 
  \begin{bmatrix}
    1 \\
    0
  \end{bmatrix}
  \Longrightarrow
  \begin{bmatrix}
    m^1_0 \\
    m^1_1
  \end{bmatrix}
  = \frac{1}{1+\v^2}
  \begin{bmatrix}
    1 \\
    -\v
  \end{bmatrix}
  \label{system2M1}
\end{gather}
and
\begin{gather}
  \begin{matrix}
    \sd: \\
    \w\p: 
  \end{matrix}
  \left(
   \begin{bmatrix}
     1 & 0 \\
     0 & 1
   \end{bmatrix}
   + \v
   \begin{bmatrix}
     0 & -2 \\
     2 & 0
   \end{bmatrix}
  \right)
  \begin{bmatrix}
    s_d   \\
    m^2_1
  \end{bmatrix}
  = 
  \begin{bmatrix}
    1 \\
    0
  \end{bmatrix}
  \Longrightarrow
  \begin{bmatrix}
    s_d \\
    m^2_1
  \end{bmatrix}
  = \frac{1}{1+4\v^2}
  \begin{bmatrix}
    1 \\
    -2\v
  \end{bmatrix}.
  \label{system2M2}
\end{gather}
Putting it all together,
\begin{gather*}
   M = \sk + \frac{2}{1+\v^2}(\p\l - \v \w\l) + \frac{1}{1+4\v^2}(\sd -  2\v \w\p).
\end{gather*}
Recalling the definitions
$\sk=\l\l+\frac{\p\p+\w\w}{2}$ and $\sd=\frac{\p\p-\w\w}{2}$,
\begin{align}
   \label{tviscosityClosure}
   M = \tviscosity :=& \left(\l\l+\frac{\p\p+\w\w}{2}\right)
     + \frac{2}{1+\v^2}(\p\l - \v \w\l)
     + \frac{1}{1+4\v^2}\left(\frac{\p\p-\w\w}{2} -  2\v \w\p\right).
\end{align}
As a check, observe that when $\v$ is zero
$M$ is the identity $\i\i$.
In the limit of strong magnetic field,
\begin{align*}
   M \approxeq& \left(\l\l+\frac{\p\p+\w\w}{2}\right) -2\v\inv\w\i + \O(\v^{-2}).
\end{align*}
Using the general properties that
\begin{gather*}
  \begin{aligned}
  XY\ddotp 2\dstrain &= \SymB(X\dotp\dstrain\dotp Y^T)&
  &\hbox{and}&
  \w^T &= -\w,
  \end{aligned}
\end{gather*}
the closure for the deviatoric stress,
$\dPT=-\viscosity 2 M\ddotp \dstrain,$ is thus
{
\def\l{\i_\parallel}
\def\p{\i_\perp}
\def\w{\i_\wedge}
\begin{equation}
\begin{aligned}
   \dPT &= -\viscosity \SymB \Bigg(
     &&\l\dotp\dstrain\dotp\l + \frac{\p\dotp\dstrain\dotp\p-\w\dotp\dstrain\dotp\w}{2}
  \\ &&&
         + \frac{2}{1+\v^2}\big(\p\dotp\dstrain\dotp\l+\v\l\dotp\dstrain\dotp\w\Bigg)
 \\ &&&+ \frac{1}{1+4\v^2}\left(\frac{\p\dotp\dstrain\dotp\p + \w\dotp\dstrain\dotp\w}{2}
               + 2\v\l\dotp\dstrain\dotp\w\right)
        \Bigg)
\end{aligned}
\label{gyrodPTclosure}
\end{equation}
}
where we have reverted to the notation $\w=\i_\w=\b\times\i$, $\p=\i_\p=\id-\b\b$,
and $\l=\i_\para=\b\b$.
The result \label{gyrodPTclosure}
agrees with equations (5.125)--(5.128) in \cite{book:Woods04}
if one corrects the typo in his equation (5.127) by replacing his definition
$\mathbf{W}_2:=(\id_\parallel\Ldiamond\id_\parallel+\id_\wedge\Ldiamond\id_\wedge)/2$
(which would result in an $M$ (his $\mathbf{W}$) which fails to be
the identity operator for $\v=0$) with the correct definition
$\mathbf{W}_2:=(\id_\perp\Ldiamond\id_\perp+\id_\wedge\Ldiamond\id_\wedge)/2$.

\subsection{Heat flux tensor}

\def\SpanMi{\Span\{M_i\}}
\def\Linf{L^\infty}
Let $\MC$ denote the inverse of the matrix $\AC:=\i\i\i+\v 3\w\i\i$.
We will seek $\MC$ as a linear combination of basis tensors $\{M_i\}$,
$\MC = \sum_i m_i M_i$, where $\SpanMi$ is closed under the map
$M\mapsto \AC\dddotp M$, that is, which is closed under the map
$L:M\mapsto 3\w\i\i\dddotp M$.  We need that $(\i\i\i+\v L)\dddotp \MC = \i\i\i$,
so it is evident that $\i\i\i$ should be in $\SpanMi$,
so we simply compute what repeated applications of $L$ can generate
starting with $\i\i\i$.
First, observe that $L$ satisfies
\begin{gather*}
  3\w\i\i\dddotp XYZ = \w\dotp X YZ
                     + \w\dotp Y XZ
                     + \w\dotp Z XY.
\end{gather*}
Using that $\w\dotp\w=-\p$ and $\w\dotp\p=\w$,
we generate a basis.  Under the mapping $M\mapsto 3\w\i\i\dddotp M$,
the calculations
\begin{gather*}
 \begin{array}{r r r r r r r}
    M^0:  &M'^0_0:=&\i\i\i &\mapsto &3\w\i\i &         &     
 \\ M^1:  &M'^1_0:=&\p\i\i &\mapsto &2\w\p\i &+\w\i\i
 \\       &M'^1_1:=&\w\i\i &\mapsto &2\w\w\i &-\p\i\i 
 \\ M^2:  &M'^2_0:=&\p\p\i &\mapsto & \w\p\p &+2\w\p\i
 \\       &M'^2_1:=&\w\p\i &\mapsto & \w\w\p &+ \w\w\i &- \p\p\i
 \\       &M'^2_2:=&\w\w\i &\mapsto & \w\w\w &-2\w\p\i
 \end{array}
\end{gather*}
and
\begin{gather}
 \begin{array}{r r r r r r r r}
    M^3:  &M^3_0:=&\p\p\p &\mapsto &3\w\p\p &         &  
 \\       &M^3_1:=&\w\p\p &\mapsto &2\w\w\p & -\p\p\p
 \\       &M^3_2:=&\w\w\p &\mapsto &\w\w\w  &-2\w\p\p
 \\       &M^3_3:=&\w\w\w &\mapsto &-3\w\w\p&         &   
 \end{array}
 \label{systemM3}
\end{gather}
exhibit such a basis.
The span of the $M^3$ subsystem is closed under the map
$M\mapsto \AC\dddotp M$, but the span of the other $M^k$
systems is not.  This is easily remedied.
For each mapping in an $M^k$ subsystem subtract the
corresponding mapping in the $M^{k+1}$ system
and repeat until all $\i$'s have been turned into
$\l$'s.  This gives a new set of basis elements
which satisfy a simpler, decoupled set of mappings and
which still span the same set:
\begin{gather*}
 \begin{array}{r r r r r r r}
    M^0:  &M^0_0:=&\l\l\l &\mapsto & 0
 \\ M^1:  &M^1_0:=&\p\l\l &\mapsto &\w\l\l
 \\       &M^1_1:=&\w\l\l &\mapsto &-\p\l\l 
 \\ M^2:  &M^2_0:=&\p\p\l &\mapsto & 2\w\p\l
 \\       &M^2_1:=&\w\p\l &\mapsto & \w\w\l &- \p\p\l
 \\       &M^2_2:=&\w\w\l &\mapsto & -2\w\p\l
 \end{array}
\end{gather*}
Ignoring the final $\l$ in all these maps, this is
identical to the decoupled system with the same variable
names $M^k_j$ that we obtained in section \ref{DeviatoricPressureTensor}
when closing the deviatoric pressure.
We therefore separate out the null space for this system in the same way,
and we replace the system $M^2$ with
\begin{gather*}
 \begin{array}{r r r r r r}
    L^2:  &2\sd\l :=&(\p\p-\w\w)\l &\mapsto & 4\w\p\l
 \\       &M^2_1:=&\w\p\l &\mapsto & (\w\w - \p\p)\l & = -2\sd\l
 \\ N^2:  &2\ss\l :=&(\p\p+\w\w)\l &\mapsto & 0
 \end{array}
\end{gather*}
Again, the null space for the subsystems $M^k$ where $k$ is odd is trivial.

We now decompose $\i\i\i=\i^3$ into a sum of decoupled basis vectors:
\begin{align*}
  \i^3 &= (\p+\l)^3
    \\ &= \p\p\p + 3\p\l\l + 3\p\p\l+\l\l\l
    \\ &= \p\p\p + 3\p\l\l + 3\sd\l + \left(3\ss\l + \l\l\l\right).
\end{align*}
The matrix $M$ is a linear combination of the basis that we have defined:
\begin{equation}
 \begin{aligned}
   M = &m_n(3\ss\l+\l\l\l) + 2(m^1_0 \p\l\l + m^1_1 \w\l\l) +  (s_d\sd\l +  m^2_1 \w\p\l)
   \\  &+ (m^3_0 \p\p\p + m^3_1 \w\p\p + m^3_2 \w\w\p + m^3_3 \w\w\w),
 \end{aligned}
 \label{3Mexpansion}
\end{equation}
where we use parentheses for components that correspond to decoupled subsystems.
To find the nine heat flux coefficients we substitute this into
the required identity
$
   (\i\i\i+\v\w\i\i)\dddotp M = \i\i\i,
$
that is,
\begin{gather*}
   m_n(3\ss\l+\l\l\l)
   \\ + 2(m^1_0 \p\l\l + m^1_1 \w\l\l)
      + \v2(m^1_0 \w\l\l - m^1_1 \p\l\l)
   \\ +  (s_d\sd\l +  m^2_1 \w\p\l)
      + \v(2s_d\w\p\l -  2m^2_1 \sd\l)
   \\ + (m^3_0 \p\p\p + m^3_1 \w\p\p + m^3_2 \w\w\p + m^3_3 \w\w\w)
   \\ \ \ \ \ + \v(-3m^3_0 \w\p\p + m^3_1 (2\w\w\p-\p\p\p) + m^3_2(\w\w\w-2\w\p\p) -3 m^3_3 \w\w\p)
   \\ = (3\ss\l + \l\l\l) + 3\p\l\l + 3\sd\l + \p\p\p.
\end{gather*}
Matching up coefficients gives the expected equation $m_n=1$,
two linear systems almost identical to those in section
\ref{DeviatoricPressureTensor},
and one new system corresponding to the $M^3$ system \eqref{systemM3}.
For the first two systems, we have
one nearly identical to system \eqref{system2M1},
\begin{gather*}
  \begin{matrix}
    \p\l\l: \\
    \w\l\l: 
  \end{matrix}
  \left(
   \begin{bmatrix}
     1 & 0 \\
     0 & 1
   \end{bmatrix}
   + \v
   \begin{bmatrix}
     0 & -1 \\
     1 & 0
   \end{bmatrix}
  \right)
  \begin{bmatrix}
    m^1_0 \\
    m^1_1
  \end{bmatrix}
  = 
  \begin{bmatrix}
    3/2 \\
    0
  \end{bmatrix}
  \Longrightarrow
  \begin{bmatrix}
    m^1_0 \\
    m^1_1
  \end{bmatrix}
  = \frac{3/2}{1+\v^2}
  \begin{bmatrix}
    1 \\
    -\v
  \end{bmatrix},
\end{gather*}
the only real modification being the replacement $2\to 3/2$ due to the
appearance of $3\p\l\l$ in place of $2\p\l$ on the right hand side,
and one system nearly identical to system \eqref{system2M2},
\begin{gather*}
  \begin{matrix}
    \sd\l: \\
    \w\p\l: 
  \end{matrix}
  \left(
   \begin{bmatrix}
     1 & 0 \\
     0 & 1
   \end{bmatrix}
   + \v
   \begin{bmatrix}
     0 & -2 \\
     2 & 0
   \end{bmatrix}
  \right)
  \begin{bmatrix}
    s_d   \\
    m^2_1
  \end{bmatrix}
  = 
  \begin{bmatrix}
    3 \\
    0
  \end{bmatrix}
  \Longrightarrow
  \begin{bmatrix}
    s_d \\
    m^2_1
  \end{bmatrix}
  = \frac{3}{1+4\v^2}
  \begin{bmatrix}
    1 \\
    -2\v
  \end{bmatrix},
\end{gather*}
the only real modification being the replacement $1\to 3$ due to the
appearance of $3\sd\l$ in place of $\sd$ on the right hand side.
For the new system corresponding to the $M^3$ system \eqref{systemM3},
we have
\begin{gather*}
  \begin{matrix}
    \p\p\p: \\
    \w\p\p: \\
    \w\w\p: \\
    \w\w\w: \\
  \end{matrix}
  \left(
   \begin{bmatrix}
     1 & 0 & 0 & 0 \\
     0 & 1 & 0 & 0 \\
     0 & 0 & 1 & 0 \\
     0 & 0 & 0 & 1 \\
   \end{bmatrix}
   + \v
   \begin{bmatrix}
     0 & -1 & 0 & 0 \\
     3 & 0 & -2 & 0 \\
     0 & 2 & 0 & -3 \\
     0 & 0 & 1 & 0 \\
   \end{bmatrix}
  \right)
  \begin{bmatrix}
    m^3_0 \\
    m^3_1 \\
    m^3_2 \\
    m^3_3 \\
  \end{bmatrix}
  = 
  \begin{bmatrix}
    1 \\
    0 \\
    0 \\
    0 \\
  \end{bmatrix}.
\end{gather*}
That is, we need to solve the linear system
\begin{gather*}
   \begin{bmatrix}
     1   & -\v & 0    & 0    \\
     3\v & 1   & -2\v & 0    \\
     0   & 2\v & 1    & -3\v \\
     0   & 0   &  \v  & 1    \\
   \end{bmatrix}
  \begin{bmatrix}
    m^3_0 \\
    m^3_1 \\
    m^3_2 \\
    m^3_3 \\
  \end{bmatrix}
  = 
  \begin{bmatrix}
    1 \\
    0 \\
    0 \\
    0 \\
  \end{bmatrix}.
\end{gather*}
We solve by row reduction and back substitution.
Subtracting $3\v$ times row 1 from row 2 gives
\begin{gather*}
   \begin{bmatrix}
     1   & -\v       & 0    & 0    \\
     0   & 1+3\v^2   & -2\v & 0    \\
     0   & 2\v       & 1    & -3\v \\
     0   & 0         &  \v  & 1    \\
   \end{bmatrix}
  \begin{bmatrix}
    m^3_0 \\
    m^3_1 \\
    m^3_2 \\
    m^3_3 \\
  \end{bmatrix}
  = 
  \begin{bmatrix}
    1 \\
    -3\v \\
    0 \\
    0 \\
  \end{bmatrix}.
\end{gather*}
Subtracting $2\v$ times row 2 from $1+3\v^2$ times row 3 gives
\begin{gather*}
   \begin{bmatrix}
     1   & -\v       & 0    & 0    \\
     0   & 1+3\v^2   & -2\v & 0    \\
     0   & 0         & 1+7\v^2    & -3\v(1+3\v^2) \\
     0   & 0         &  \v  & 1    \\
   \end{bmatrix}
  \begin{bmatrix}
    m^3_0 \\
    m^3_1 \\
    m^3_2 \\
    m^3_3 \\
  \end{bmatrix}
  = 
  \begin{bmatrix}
    1 \\
    -3\v \\
    6\v^2 \\
    0 \\
  \end{bmatrix}.
\end{gather*}
Subtracting $\v$ times row 3 from $1+7\v^2$ times row 4 gives
\begin{gather*}
   \begin{bmatrix}
     1   & -\v       & 0    & 0    \\
     0   & 1+3\v^2   & -2\v & 0    \\
     0   & 0         & 1+7\v^2    & -3\v(1+3\v^2) \\
     0   & 0         & 0    & 1+10\v^2+9\v^4    \\
   \end{bmatrix}
  \begin{bmatrix}
    m^3_0 \\
    m^3_1 \\
    m^3_2 \\
    m^3_3 \\
  \end{bmatrix}
  = 
  \begin{bmatrix}
    1 \\
    -3\v \\
    6\v^2 \\
    -6\v^3 \\
  \end{bmatrix}.
\end{gather*}
Back substituting yields
\begin{equation}
 \begin{aligned}
   k_3 := m^3_3 &= \frac{-6\v^3}{1+10\v^2+9\v^4} & = -(2/3)\v\inv + \O(\v^{-3}), \\
   k_2 := m^3_2 &= \frac{6\v^2 + 3\v(1+3\v^2) m^3_3}{1+7\v^2} & = \O(\v^{-2}) , \\
   k_1 := m^3_1 &= \frac{-3\v + 2\v m^3_2}{1+3\v^2} & = -\v\inv + \O(\v^{-3}), \\
   k_0 := m^3_0 &= 1+\v m^3_1 & = \O(\v^{-2}).
 \end{aligned}
 \label{heatCoefsM}
\end{equation}
In summary, the heat flux is given by
\begin{gather}
  \qT = -\tfrac{2}{5}\THeatConductivity \dddotp\SymC\left(\Pshape\dotp\D\TT\right),
  \label{qTclosureForm}
\end{gather}
where the heat conductivity tensor
is given by
$\glslink{THeatConductivity}{\THeatConductivity}=\heatConductivity\tHeatConductivity$,
where $\glslink{tHeatConductivity}{\tHeatConductivity}:=M$
is the shape of the heat conductivity tensor.
The heat conductivity is thus a sixth order gyrotropic tensor
and is given by equation \eqref{3Mexpansion},
\begin{equation}
 \begin{aligned}
  \tHeatConductivity
     = &(3\ss\l+\l\l\l) + \frac{3}{1+\v^2}\left(\p\l\l - \v \w\l\l\right)
        + \frac{3}{1+4\v^2}\left(\sd\l - 2\v \w\p\l\right)
   \\  &+ (k_0 \p\p\p + k_1 \w\p\p + k_2 \w\w\p + k_3 \w\w\w),
 \end{aligned}
 \label{3MexpansionCopy}
\end{equation}
that is, recalling that
  $\ss := \frac{\p\p+\w\w}{2}$ and 
  $\sd := \frac{\p\p-\w\w}{2}$,
and using powers to denote splice symmetric product powers,
\begin{equation}
 \begin{aligned}
  \tHeatConductivity
     = &\left(\l^3 + \tfrac{3}{2}\l(\p^2+\w^2)\right)
    \\& + \frac{3}{1+\v^2}\left(\p\l^2 - \v \w\l^2\right)
        + \frac{3}{1+4\v^2}\left(\frac{\p^2-\w^2}{2}\l - 2\v \w\p\l\right)
   \\  &+ (k_0 \p^3 + k_1 \w\p^2 + k_2 \w^2\p + k_3 \w^3),
 \end{aligned}
 \label{3MexpansionCopy2}
\end{equation}
where the coefficients $k_i$ are given in \eqref{heatCoefsM}.
Note that when $\v$ is zero $\tHeatConductivity$ is the identity $\i\i\i$.
In the limit of strong magnetic field,
\begin{equation*}
 \begin{aligned}
   \tHeatConductivity
     = &\left(\l^3 + \tfrac{3}{2}\l(\p^2+\w^2)\right)
        - \v\inv\Big( 3\w\l^2 + 6\w\p\l + \w\p^2 + \tfrac{2}{3} \w^3\Big)
        + \O(\v^{-2}).
 \end{aligned}
\end{equation*}

For computational efficiency, one can make use of the identity
\begin{gather*}
  XYZ\dddotp S = \Sym\left((X\stimes Y\stimes Z)\dddotp S\right),
\end{gather*}
which holds for any symmetric tensor $S$.
So we can write the heat flux closure as
\begin{gather}
  \qT = -\tfrac{2}{5}\heatConductivity \Sym\left(\tHeatConductivity'
     \dddotp\SymC\left(\Pshape\dotp\D\TT\right)\right),
 \label{qTclosure}
\end{gather}
where the formula for $\tHeatConductivity'$
is the same as for $\tHeatConductivity$ except that
because we symmetrize at the end
we can now take the default product of tensors
to be simple splice products (in order to avoid averaging
over $3!=6$ splice products for each triple product).

}

%% file: appendix3.tex
\chapter{Reconnection}

\def\x{x}
\def\y{y}
\def\z{z}
\section{Inflow ODE for the ten-moment equations.}

Jerry Brackbill has studied two-dimensional antiparallel
reconnection using the 10-moment two-fluid equations, assuming
symmetry across the inflow axis.
With mild approximating assumptions
for steady state he reduces the adiabatic ($\qT=0$)
equations on the $\y$-axis to an ODE \cite{article:Br11}.

Recall that symmetry in the $\x$-axis says that all tensor
components are even or odd in $\x$ based on whether
$\x$ appears as a subscript an even or odd number of times.
(For a pseudo-tensor such as the magnetic field $\B$
components are even iff the number of $\x$ subscripts is odd;
for a proper tensor
components are even iff the number of $\x$ subscripts is even.)
\def\ux{\u_\x}
\def\uz{\u_\z}
\def\uy{\u_\y}
\def\Bx{\B_\x}
\def\Bz{\B_\z}
\def\By{\B_\y}
\def\Ex{\E_\x}
\def\Ez{\E_\z}
\def\Ey{\E_\y}
\def\TTxx{\TT_{\x\x}}
\def\TTxz{\TT_{\x\z}}
\def\TTxy{\TT_{\x\y}}
\def\TTzz{\TT_{\z\z}}
\def\TTyz{\TT_{\y\z}}
\def\TTyy{\TT_{\y\y}}
\def\Pxx{\PT_{\x\x}}
\def\Pxz{\PT_{\x\z}}
\def\Pxy{\PT_{\x\y}}
\def\Pzz{\PT_{\z\z}}
\def\Pyz{\PT_{\y\z}}
\def\Pyy{\PT_{\y\y}}
\begin{itemize}
  \item Even:
    $\uz$, $\uy$,
    $\Pxx$,
    $\Pyy$,
    $\Pzz$,
    $\Pyz$,
    $\Bx$,
    $\Ey$,
    $\Ez$,
  \item Odd:
    $\ux$, 
    $\Pxy$,
    $\Pxz$,
    $\Bz$,
    $\By$,
    $\Ex$,
\end{itemize}
Odd functions (including $\x$-derivatives of even functions
but not $\x$-derivatives of odd functions)
vanish on the $\x=0$ axis.

\def\Dx{\partial_\x}
\def\Dy{\partial_\y}
\def\Dz{\partial_\z}

All quantities are independent of $\z$, so $\Dz=0$.
So on the $\x=0$ axis the convective derivative is
simply $\dt:=\partial_t+\u_\y\Dy$.
By using convective derivative we avoid $\x$-derivatives.
Therefore we work with the temperature evolution
equation \eqref{nTTevolution}
\begin{gather*}
   n d_t \TT
     + \SymB(\PT\dotp\nabla\u) + \Div\qT
     = ({q/m})\SymB(\PT\times\B) + \RT,
\end{gather*}
where we assume a default species index.

\subsection{Component evolution}

\def\a{a}
\def\b{b}
For arbitrary indices $\a$ and $\b$ we have
\begin{alignat*}{5}
  2\Sym(\PT\cdot\D\u)_{\a\b}
  =&  &\PT_{\a\x}\Dx\u_\b + \PT_{\b\x}\Dx\u_\a \\
   &+ &\PT_{\a\y}\Dy\u_\b + \PT_{\b\y}\Dy\u_\a.
\end{alignat*}
In particular, on the $\x=0$ axis we have
\begin{align*}
  2\Sym(\PT\cdot\D\u)_{\y\z} &= \PT_{\y\z}\Dy\uy + \PT_{\y\y}\Dy\uz, \\
  2\Sym(\PT\cdot\D\u)_{\x\x} &= 2 \PT_{\x\x}\Dx\ux, \\
  2\Sym(\PT\cdot\D\u)_{\y\y} &= 2 \PT_{\y\y}\Dy\uy, \\
  2\Sym(\PT\cdot\D\u)_{\z\z} &= 2 \PT_{\y\z}\Dy\uz.
\end{align*}
On the $\x$ axis $\Bz=0=\By$.
So for arbitrary indices, on the $\x$ axis we have
\begin{align*}
  2\Sym(\PT\times\B)_{\a\b}
   &= \Bx\left(\PT_{\a j}\epsilon_{\x\b j} + \PT_{\b j}\epsilon_{\x\a j}\right).
\end{align*}
In particular,
\begin{align*}
  2\Sym(\PT\times\B)_{\y\z} &= \Bx(-\Pzz+\Pyy), \\
  2\Sym(\PT\times\B)_{\x\x} &= 0, \\
  2\Sym(\PT\times\B)_{\y\y} &= -2\Bx\Pyz, \\
  2\Sym(\PT\times\B)_{\z\z} &= 2\Bx\Pyz.
\end{align*}
Along the $\y$ axis the isotropization components are 
\begin{gather*}
  \RT_{\y\z} = -\Pyz/\tau, \\
  \RT_{\x\x} = (p-\Pxx)/\tau, \\
  \RT_{\y\y} = (p-\Pyy)/\tau, \\
  \RT_{\z\z} = (p-\Pzz)/\tau.
\end{gather*}
The foregoing component identities continue to hold if we replace the
pressure tensor with the temperature tensor.
Ignoring heat flux, and assuming $\Dz=0$,
the components of the temperature tensor
evolution equation on the inflow axis are thus
 \begin{gather}
\begin{split}
    &\dt \TTyz + \TT_{\y\z}\Dy\uy + \TT_{\y\y}\Dy\uz
    = {q\over m} \Bx(-\TTzz+\TTyy)
       -\TTyz/\tau, \\
    &\dt \TTxx + 2 \TT_{\x\x}\Dx\ux = (T-\TTxx)/\tau, \\
    &\dt \TTyy + 2 \TT_{\y\y}\Dy\uy
    = -2{q\over m}\Bx\TTyz + (T-\TTyy)/\tau, \\
    &\dt \TTzz + 2 \TT_{\y\z}\Dy\uz
    = 2{q\over m} \Bx\TTyz + (T-\TTzz)/\tau.
\end{split}
\label{eqn:TTBinflow}
 \end{gather}

Along the inflow axis, by symmetry
$\dt = \partial_t + \uy\Dy$.  In steady state this becomes
$\dt = \uy\Dy$.  

\section{Adiabatic ten-moment stagnation point flow.}
\label{AdiabaticStagnation}


\subsection{Abstract framework for linearization about the X-point}
\label{abstractFramework}

What can smooth, steady reconnection look like near the X-point?
Linearize about the X-point.  Then spatial derivatives of $\B$ and
$\u$ become constants, so equation
\eqref{eqn:TTBinflow} is a linear ODE in $\y$ of the form
\def\A{A}
\def\BB{B}
\def\xx{\mathbf{U}}
\begin{gather*}
  \y\xx' = \A\cdot\xx + \y\BB\cdot\xx,
\end{gather*}
where $\A$ and $\BB$ are matrices of constant coefficients
(and where $\BB$ involves magnetic field components),
prime denotes differentiation with respect to $\y$,
and $\xx(\y)$ is a vector containing the four
non-vanishing components of $\TT$.

\def\R{R}
\def\Rinv{R^{-1}}
\def\yz{\mathbf{W}}
Near the X-point the magnetic field vanishes
and we get
\begin{gather}
  \label{eqn:genericStagnation}
  \y\xx' \approx \A\cdot\xx
\end{gather}
We will see that we can take $\xx = (\TTxx,\TTzz,\TTyy)^\mathrm{T}$.
To solve such an ODE you can use the eigenstructure of $\A$.
Suppose that it has a full eigenvector decomposition.
Then $\A = \R\cdot\Lambda\cdot\Rinv$, so we have
\begin{gather*}
  \y(\Rinv\cdot\xx)' = \Lambda\cdot(\Rinv\cdot\xx).
\end{gather*}
Let $\yz:=\Rinv\cdot\xx$.
Then this system decouples into scalar equations of the form
\begin{gather*}
  \y W' = \lambda W.
\end{gather*}
Separating and integrating gives
\begin{gather}
  \label{Worigin}
  W = C \y^\lambda.
\end{gather}
So existence of a smooth,
strictly positive steady-state $\TT$
requires that $\lambda=0$ is an eigenvalue
(else any steady solution must be zero or singular
at the origin) for an eigenvector with positive
components.
I remark that existence of a positive definite
solution (whether singular or not) requires existence
of an eigenvector representing a positive-definite
solution.

\def\d{\mathrm{d}}
To confirm physically
that the neglect of the magnetic field
in \eqref{eqn:genericStagnation} is justifiable,
note that the magnetic field simply rotates the
temperature tensor at a rate proportional to the
strength of the magnetic field.  Since $\uy$ and
$\Bx$ are both proportional to distance from $\y=0$,
$
  \frac{\d\theta}{\d \y}
    = \frac{\d\theta}{\d t} \frac{\d t}{\d \y}   
  = \frac{\y\Dy\Bx}{\y\Dy\uy}
  = \frac{\Dy\Bx}{\Dy\uy},
$
so the angle of rotation is proportional
to the distance moved along the $\y$ axis
and rotation can be neglected as $\y$ approaches $0$.

\subsection{Stagnation point flow for an adiabatic isotropizing
ten-moment gas}

We now use the framework outlined in
section \ref{abstractFramework}
to argue that smooth steady-state stagnation
point flow is singular in an
anisotropic ten-moment model without heat flux.

\def\uyy{\u_{\y,\y}}
Ignoring both magnetic field as well heat flux,
the terms involving $\Bx$ drop out of
\eqref{eqn:TTBinflow}.  At the X-point
the vorticity $\Dy\uz$ is zero, so in the linear
analysis these terms also drop out.
At the X-point steady flow must be incompressible,
$\Div\u = -\dt\ln\mdens = 0$,
so $\Dx\ux = -\Dy\uy$.
In the linearization we replace $\Dy\uy$
with its value at the X-point.
Define $\uyy := \Dy\uy|_{\mathbf{\y=0}}$.
So near the X-point the components of the
temperature tensor evolve according to
 \begin{gather}
\begin{split}
    &\dt \TTyz + \TTyz\uyy = -\TTyz/\tau, \\
    &\dt \TTxx - 2 \TTxx\uyy = (T-\TTxx)/\tau, \\
    &\dt \TTyy + 2 \TTyy\uyy = (T-\TTyy)/\tau, \\
    &\dt \TTzz = (T-\TTzz)/\tau,
\end{split}
\label{eqn:TTinflow}
 \end{gather}
where recall that $T=(\TTxx+\TTzz+\TTyy)/3$.
Steady-state-flow implies that
$\dt = \u\cdot\D = \uy\Dy$.
So, in the linearization,
$\dt = \y\uyy\Dy$.
Note that $\TTyz$ decouples from the other components
and is zero if initially zero.
Dividing by $\uyy$ and expressing \eqref{eqn:TTinflow}
in matrix form,
\def\At{\tilde\A}
\def\ttau{\bar\tau}
\begin{gather*}
  \y\Dy
   \underbrace{
   \begin{bmatrix}
     \TTxx \\
     \TTyy \\
     \TTzz \\
   \end{bmatrix}
   }_{\hbox{\large $\xx$}}
  = 
   \frac{1}{3\ttau}
   \underbrace{
   \begin{bmatrix}
     -2+6\ttau & 1 & 1 \\
     1 & -2-6\ttau & 1 \\
     1 & 1 & -2
   \end{bmatrix}
   }_{\hbox{\small Call \large $\At$}}
   \begin{bmatrix}
     \TTxx \\
     \TTyy \\
     \TTzz \\
   \end{bmatrix},
\end{gather*}
where $\ttau:=\tau\uyy$.
(So $\ttau<0$ if there is inflow along $\y$
and $\ttau>0$ if there is outflow along $\y$.)
We have now matched up with the framework of
\eqref{eqn:genericStagnation}.
Note that the matrix $\At$ is symmetric,
so it must have real eigenvalues and a
full set of orthogonal eigenvectors.

\subsubsection{
  Summary of eigenvalue/eigenvector results.}

\label{eigresults}
\def\lambdaA{\lambda^{(1)}}
\def\Z{Y}
In the inflow case ($\ttau<0$) there is one negative
eigenvalue $\lambdaA$ and all components of its
eigenvector can be positive.  The solution
for this component blows up as $\y$ goes to zero.
The other two eigenvalues are positive (representing
decay as $\y$ goes to zero) and have eigenvector
components of mixed sign which cannot be combined
to give all-positive components.
Using these three eigenvectors you
can match any positive data at a boundary
$\y=\pm \Z$ and positivity then holds between
$\pm \Z$ and $0$.  This solution represents
a heating singularity as inflow approaches the
X-point.

In the outflow case ($\ttau>0$) there is one positive
eigenvalue $\lambdaA$ and all components of its
eigenvector can be positive. The solution for this
component blows up as $\y$ goes to $\infty$.
The other two eigenvalues are negative (representing
decay as $\y$ goes to $\infty$ and
have eigenvector components of mixed sign which
cannot be combined to give all-positive components.
Thus a steady outflow solution which is positive at
$\y=\Z_\mathrm{min}$ will be positive for all
$\y$ between $\Z_\mathrm{min}$ and $\infty$.

If the $\y$ axis is an inflow axis then the
$\x$ axis is an outflow axis.
From the previous two paragraphs we conclude
that a steady-state solution must be not only
singular but discontinuous at the X-point
(even in the topology of the real line
which includes a point at infinity).
In practice absence of heat flux will be violated
in the immediate vicinity of the X-point.

Except when crossing the singular point $\ttau^{-1}=0$ all
eigenvalues (and all eigenvector components,
using consistent scaling) change monotonically with $\ttau$
and do not change in sign.
\def\sC{^{(3)}}
\def\sB{^{(2)}}
\def\sA{^{(1)}}
\def\sN{.2113}
\def\mN{.5774}
\def\lN{.7887}

When there is outflow (rather than inflow)
in the $\y$ axis, the $\y$ and $\x$
components effectively swap places and
eigenvalues are negated
($\lambda\sA > \lambda\sB > \lambda\sC$)
See table \ref{eigValsStagPtFlow}.)

\begin{table}[h!]
\begin{displaymath}
  \begin{array}{l|rcccl}
 \hline
   \ttau^{-1}: &-\infty &\le  &   0        &\le   &\infty
 \\                             \hline
   \lambda\sA: &\infty  &\ge  &     2      &\ge   &0
 \\\TTxx\sA:   &\lN     &\le  &    1       &\ge   &\mN
 \\\TTyy\sA:   &-\sN    &\le  & 0\cdot1    &\le   &\mN
 \\\TTzz\sA:   &-\mN    &\le  & 0\cdot2    &\le   &\mN
 \\                             \hline
   \lambda\sB: &\infty  &\ge  &    0       &\ge   &-\infty
 \\\TTxx\sB:   &\sN     &\ge  & 0\cdot(-1) &\ge   &-\lN
 \\\TTyy\sB:   &-\lN    &\le  & 0\cdot1    &\le   &\sN
 \\\TTzz\sB:   &\mN     &\le  & 1          &\ge   &\mN
 \\                             \hline
   \lambda\sC: &0       &\ge  &   -2       &\ge   &-\infty
 \\\TTxx\sC:   &\mN     &\ge  & 0\cdot(-1) &\ge   &-\sN
 \\\TTyy\sC:   &\mN     &\le  &    1       &\ge   &\lN
 \\\TTzz\sC:   &\mN     &\ge  & 0\cdot(-2) &\ge   &-\mN
  \end{array}
\end{displaymath}
\caption{Eigenstructure for solutions to
  adiabatic ten-moment linearized stagnation point flow.}
In this table
multiplication of $0$ by a constant is used to indicate
asymptotic relative scaling of eigenvector
components near $\ttau\inv\approxeq 0$;
$\mN$ is an approximation for $\sqrt{3}^{-1}$,
$\sN$ is an approximation for $\frac{1-\sqrt{3}^{-1}}{2}$,
and $\lN$ is an approximation for $\frac{1+\sqrt{3}^{-1}}{2}$.
Except when crossing the singular point $\ttau\inv=0$,
eigenvalues and eigenvector components change monotonically
with $\ttau$ and do not change in sign.
\fhrule
\label{eigValsStagPtFlow}
\end{table}

In the case of no isotropization,
$\ttau^{-1}=0$, the ODE simplifies to
\begin{gather*}
  \y\Dy
   \begin{bmatrix}
     \TTxx \\
     \TTyy \\
     \TTzz \\
   \end{bmatrix}
  = 
   \begin{bmatrix}
     2 & 0 & 0 \\
     0 &-2 & 0 \\
     0 & 0 & 0
   \end{bmatrix}
   \begin{bmatrix}
     \TTxx \\
     \TTyy \\
     \TTzz \\
   \end{bmatrix}
\end{gather*}
and the solution is simply
\begin{gather*}
  \TTxx(\y) = \TTxx|_{\y=1}\y^{2}, \\
  \TTyy(\y) = \TTyy|_{\y=1}\y^{-2}, \\
  \TTzz(\y) = \TTzz|_{\y=1},
\end{gather*}
which is singular at the origin, contradicting
the smoothness assumption.

\subsubsection{Calculation of eigenvalues and eigenvectors.}

\def\tlambda{\tilde\lambda}
This section incompletely works out
the results tabulated in section \ref{eigresults}.

If $\tlambda$ is an eigenvalue of $\At$
then $\lambda:=\tlambda/(3\ttau)$ is an eigenvalue
of $\A:=\At/(3\ttau)$.
To find the eigenstructure of the coefficient matrix
$\At$, we solve $(\At-\idtens\tlambda)\cdot\xx=\mathbf{0}$.
We represent this system with the matrix
\begin{gather*}
   \begin{bmatrix}
     -2+6\ttau-\tlambda & 1 & 1 \\
     1 & -2-6\ttau - \tlambda & 1 \\
     1 & 1 & -2-\tlambda \\
   \end{bmatrix}.
\end{gather*}
To facilitate finding the eigenvectors, subtract
the third row from the first and the second and get
\begin{gather*}
   \begin{bmatrix}
     6\ttau-(\tlambda+3) & 0                       & \tlambda+3  \\
     0                   & -6\ttau - (\tlambda+3)  & \tlambda+3  \\
     1                   & 1                       & -2-\tlambda \\
   \end{bmatrix}.
\end{gather*}
For convenience define
\begin{alignat*}{3}
  h&:=6\ttau       &&\color{blue} = 6\tau\uyy, \\
  \mu&:=\tlambda+3 &&\color{blue} = 3\ttau\lambda +3 = \tau\uyy\lambda+3.
\end{alignat*}
Then we have
\begin{gather*}
   \begin{bmatrix}
     h-\mu & 0          & \mu   \\
     0     & -(h + \mu) & \mu   \\
     1     & 1          & 1-\mu 
   \end{bmatrix}.
\end{gather*}
Since the null space must be nonzero the rows
are linearly dependent.  So we can ignore the
last row.  The first and second rows reveal
that an eigenvector must be proportional to
something of the form
\begin{gather*}
  \begin{bmatrix}
    \mu(\mu+h) \\
    \mu(\mu-h) \\
    (\mu-h)(\mu+h) \\
  \end{bmatrix}.
\end{gather*}
The last row then reveals the characteristic
equation that $\mu$ must satisfy to be a (shifted)
eigenvalue:
\begin{gather*}
    \mu(\mu+h) +
    (1-\mu)(\mu^2-h^2) +
    \mu(\mu-h) = 0, \hbox{ i.e.,} \\
    -\mu^3+3\mu^2+h^2\mu-h^2 = 0.
\end{gather*}
For nonzero $h$, $\mu=3$ is not a root,
so $\tlambda=\mu-3$ cannot be zero,
so there is no steady-state finite smooth
solution with nonzero temperature at the origin.

In case $h=\infty$ this becomes
\begin{gather*}
    \mu = 1, \hbox{ i.e., } \tlambda = -2.
\end{gather*}
In the limit $h\rightarrow 0$ of instantaneous isotropization
this becomes 
\begin{gather*}
    \mu^2(3-\mu) = 0, \hbox{ i.e., }
    (\tlambda+3)^2\tlambda = 0,
\end{gather*}
and the eigenvectors collapse to a single eigenvector
with equal components representing isotropic temperature,
agreeing with the fact that isotropic pressure can
support stagnation point flow.

%% file: appendix4.tex
\chapter{Numerics}
\label{Numerics}

\section{Source ODE}
{
\subsection{Basic Equations}

If we neglect spatial derivatives, then
the two-fluid-Maxwell equations reduce to an ODE.
The purpose of this appendix section is to solve
this \defining{source term ODE}; we assume throughout
that spatial derivatives are zero.

If we neglect collisional terms then this ODE is
linear with imaginary eigenvalues and can be solved exactly.
If we incorporate entropy-respecting collisional closures
then the ODE is linear if we make the approximating
assumption that the closure coefficients are frozen.
In fact, closure coefficients are functions of temperature,
and collisions increase the temperature, so the frozen-coefficient
assumption incurs time-splitting error.
The collisional terms come from symmetric matrices
and have negative real eigenvalues.

 
Maxwell's equations \eqref{MaxwellsEqns}
assert that the magnetic field is constant
and that the displacement current balances
the net electrical current:
\def\epsN{\varepsilon_0}
\begin{align*}
   \partial_t \B &= 0,
\\ \partial_t \E &= - \J/\epsN = e(n_\e\u_\e-n_\i\u_i).
\end{align*}

The density evolution equations
\eqref{densityEvolution} and \eqref{spcNumDensEvolution}
assert that the densities (whether mass density or particle
density or charge density) remain constant:
 \begin{align*}
   \begin{aligned}
   \partial_t\mdens_\s &= 0,& &\hbox{i.e.,}&
   \partial_t n_\s &= 0,& &\hbox{i.e.,}&
   \partial_t \qdens_\s &= 0.
   \end{aligned}
 \end{align*}

\def\dragForce{\R}
\def\dragCoef{\tilde\eta}
\def\resistivity{\eta}
The momentum equation \eqref{momentumEvolution}
for species $\s$ is
\begin{gather*}
   \partial_t(\mdens_\s\u_\s) 
     = \frac{q_\s}{m_\s}\mdens_\s (\E + \u_\s\times\B) + \dragForce_{\s}.
\end{gather*}
We will neglect the collisional drag force $\dragForce_{\s}$.
Since densities are constant, we can divide by density.

We thus get the electro-momentum system
\begin{align}
   \label{electromomentumSystem}
  \begin{aligned}
   \partial_t \E &= e(n_\e\u_\e-n_\i\u_i),
\\ \partial_t\u_\i &= \frac{e}{m_\i} (\E + \u_\i\times\B),
\\ \partial_t\u_\e &= \frac{-e}{m_\e} (\E + \u_\e\times\B).
  \end{aligned}
 \end{align}


Evolution of energy is implied by evolution of
momentum (which implies evolution of kinetic energy)
and evolution of pressure (which is equivalent to
evolution of thermal energy).
Note that pressure evolution is temperature evolution
times the constant density.

The five-moment pressure evolution equation \eqref{pEvolution}
\begin{gather*}
   (3/2)\partial_t p_\s = Q_\s
\end{gather*}
says that pressure is constant in the absence of interspecies
collisional heating due to resistive drag and thermal
equilibration. If the drag force is non-negligible, and assuming
that the resistive drag coefficient is a function of temperature,
the electro-momentum system coupled to pressure evolution
comprises a minimally closed system. Neglecting collisional
terms, pressure evolution simply asserts that pressure is
constant.

The pressure tensor evolution equation
\eqref{PTevolution} of a ten-moment gas,
\def\Pressure{\mathbb{P}}
\def\Temperature{\mathbb{T}}
\def\Relaxation{\mathbb{R}}
\def\Sym{\mathrm{Sym}}
\def\QQf{\mathbb{Q}^\mathrm{f}}
\def\QQt{\mathbb{Q}^\mathrm{t}}
\begin{align}
  \label{PTsource}
  \partial_t \Pressure_\s 
      &= (q_\s/ m_\s)\SymB(\Pressure_\s\times\B) + \RT_\s
        + \QQf_\s + \QQt_\s,
\end{align}
depends on the solution to the electro-momentum system
if the frictional heating term $\QQf_\s$ is retained,
but is otherwise independent.
We will neglect $\QQf_\s$ and the thermal equilibration
$\QQt_\s$ and will assume the entropy-respecting
isotropization closure \eqref{isoRTclosure},
\begin{gather}
  \label{isoRTclosureCopy}
  \RT_\s = -\tau_\s\inv\dPT_\s,
\end{gather}
which exponentially dampens the deviatoric pressure
$\dPT_\s:= \PT_\s-p_\s\idtens$.
This equation is linear if $\tau_\s$ is defined
in terms of $T_\s$ (rather than in terms of $\det\TT_\s$).


\subsection{The electro-momentum system}
\def\Bidef{\frac{e\B}{m_\i}}
\def\Bedef{\frac{-e\B}{m_\e}}
\def\mBedef{\frac{e\B}{m_\e}}
\def\emi{\frac{e}{m_i}}
\def\eme{\frac{e}{m_e}}
\def\asi{\frac{e n_\i}{\epsN}}
\def\ase{\frac{e n_\e}{\epsN}}
Written in matrix form, the non-resistive electro-momentum system
\eqref{electromomentumSystem} reads
\begin{gather*}
 \partial_t
  \begin{bmatrix}
    \E \\
    \u_\i \\
    \u_\e \\
  \end{bmatrix}
    =
  \begin{bmatrix}
    0 & -\asi & \ase
 \\ \emi & -\Bidef\times\id  & 0
 \\ -\eme & 0 & \mBedef\times\id
  \end{bmatrix}
  \begin{bmatrix}
    \E \\
    \u_\i \\
    \u_\e \\
  \end{bmatrix}.
\end{gather*}

\def\EN{\E_0}
\def\uiN{{\u_\i}_0}
\def\ueN{{\u_\e}_0}
\def\Et{\widetilde\E}
\def\uit{{{\widetilde\u}_\i}}
\def\uet{{{\widetilde\u}_\e}}
\def\uitC{{{\mathbf{\widetilde\u}}_\i}{}_3}
\def\uetC{{{\mathbf{\widetilde\u}}_\e}{}_3}
We can make this ODE antisymmetric by rescaling.
For a generic rescaling, suppose
\begin{align*}
   \E &= \Et \EN,
\\ \u_\i &= \uit \uiN,
\\ \u_\e &= \uet \ueN.
\end{align*}
Making this substitution gives the system
\def\asit{\asi\frac{\uiN}{\EN}}
\def\aset{\ase\frac{\ueN}{\EN}}
\def\emit{\emi\frac{\EN}{\uiN}}
\def\emet{\eme\frac{\EN}{\ueN}}
\begin{gather*}
 \partial_t
  \begin{bmatrix}
    \Et \\
    \uit \\
    \uet \\
  \end{bmatrix}
    =
  \begin{bmatrix}
    0 & -\asit & \aset
 \\  \emit & -\Bidef\times\id  & 0
 \\ -\emet & 0 & \mBedef\times\id
  \end{bmatrix}
  \begin{bmatrix}
    \Et \\
    \uit \\
    \uet \\
  \end{bmatrix}.
\end{gather*}
If we require this system to be antisymmetric then
\begin{align*}
 \begin{aligned}
  \frac{\EN}{\uiN} &= \sqrt{\frac{\mdens_\i}{\epsN}} &
                   & \hbox{and} &
  \frac{\EN}{\ueN} &= \sqrt{\frac{\mdens_\e}{\epsN}},
 \end{aligned}
\end{align*}
(where recall that $\mdens_\i=m_\i n_\i$ and $\mdens_\e=m_\e n_\e$)
and the system becomes
\def\Bi{\B_\i}
\def\Be{\B_\e}
\def\ai{\Omega_\i}
\def\ae{\Omega_\e}
\def\gi{\omega_\i}
\def\ge{\omega_\e}
\begin{gather*}
 \partial_t
  \begin{bmatrix}
    \Et \\
    \uit \\
    \uet \\
  \end{bmatrix}
    =
  \begin{bmatrix}
    0 & -\ai\id & \ae\id
 \\  \ai\id & -\Bi\times\id  & 0
 \\ -\ae\id & 0 & -\Be\times\id
  \end{bmatrix}
  \begin{bmatrix}
    \Et \\
    \uit \\
    \uet \\
  \end{bmatrix},
\end{gather*}
where each entry in the block matrix represents a 3$\times$3 matrix
and where
\begin{align*}
 \begin{aligned}
  \ai &= e\sqrt{\frac{n_\i}{\epsN m_\i}}& & \hbox{and} &&
  \ae = e\sqrt{\frac{n_\e}{\epsN m_\e}}
 \end{aligned}
\end{align*}
denote the ion and electron plasma frequencies and
\begin{align*}
 \begin{aligned}
  \Bi &= \Bidef & & \hbox{and} &&
  \Be = \Bedef 
 \end{aligned}
\end{align*}
are the magnetic field rescaled for ions and electrons.
Their magnitudes are the ion gyrofrequency
$\gi:=|\Bi|$
and the electron gyrofrequency
$\ge:=|\Be|$.

\subsubsection{Solution of perpendicular system}

To solve the system we decompose into parallel and
perpendicular components.  Without loss of generality
assume that $\B$ is in the direction of the third axis.
Then our system decouples into a parallel system

\begin{gather*}
 \partial_t
  \begin{bmatrix}
    \Et_3 \\
    \uitC \\
    \uetC \\
  \end{bmatrix}
    =
  \begin{bmatrix}
    0 & -\ai & \ae
 \\  \ai & 0 & 0
 \\ -\ae & 0 & 0
  \end{bmatrix}
  \begin{bmatrix}
    \Et_3 \\
    \uitC \\
    \uetC \\
  \end{bmatrix},
\end{gather*}

and a perpendicular system

\def\z{0}
\def\m{\!\!\!-}
\begin{gather*}
  \partial_t
  \begin{bmatrix}
    \Et_1 \\
    \Et_2 \\
    {\uit}{}_1 \\
    {\uit}{}_2 \\
    {\uet}{}_1 \\
    {\uet}{}_2 \\
  \end{bmatrix}'
    \!\!
    =
  \begin{bmatrix}
    \z & \z & \m\ai & \z & \ae & \z
 \\ \z & \z & \z & \m\ai & \z & \ae
 \\ \ai & \z & \z & \gi & \z & \z
 \\ \z & \ai &\m\gi & \z & \z & \z
 \\ -\ae & \z & \z & \z & \z & \m\ge
 \\ \z & \m\ae & \z & \z & \ge & \z
  \end{bmatrix}
  \begin{bmatrix}
    \Et_1 \\
    \Et_2 \\
    {\uit}{}_1 \\
    {\uit}{}_2 \\
    {\uet}{}_1 \\
    {\uet}{}_2 \\
  \end{bmatrix}.
\end{gather*}

This is an antisymmetric matrix and therefore
has imaginary eigenvalues and orthogonal eigenvectors.
If we view the first and second components of each
vector as real and imaginary parts, then this
becomes a 3$\times$3 complex linear differential equation
with a skew hermitian coefficient matrix:
\def\Eo{\Et_\perp}
\def\uio{\uit{}_\perp}
\def\ueo{\uet{}_\perp}
\begin{gather}
 \label{perpSystem}
 \partial_t
  \begin{bmatrix}
    \Eo \\
    \uio \\
    \ueo \\
  \end{bmatrix}
    =
  \begin{bmatrix}
    0 & -\ai & \ae
 \\  \ai & -i\gi & 0
 \\ -\ae & 0 & i\ge
  \end{bmatrix}
  \begin{bmatrix}
    \Eo \\
    \uio \\
    \ueo \\
  \end{bmatrix},
\end{gather}
where we have used the natural isomorphism between
$SO(2,\mathbb{R})$ and complex numbers
\begin{gather*}
  a+ib \longleftrightarrow 
  \begin{bmatrix}
    a & -b
 \\ b & a
  \end{bmatrix}.
\end{gather*}

Observe that the parallel system is the special
case of this system when the magnetic field is zero.

To generalize, suppose we want to solve
the constant-coefficient linear ODE
\def\x{\underline{x}}
\def\v{\underline{v}}
\def\AA{\underline{\underline{A}}}
\def\BB{\underline{\underline{B}}}
\begin{gather*}
  \x' = \AA\cdot\x.
\end{gather*}
Seeking a solution $\x(t)=\v\exp(\lambda t)$
(where $\v\ne \underline{0}$)
leads to the eigenvector problem
\begin{align*}
  \begin{aligned}
  \v\lambda &= \AA\cdot\v, && \hbox{i.e.,} &&
  (\AA-\id\lambda)\cdot\v = 0.
  \end{aligned}
\end{align*}

We recall the theory of skew-Hermitian and Hermitian matrices.
Since $\AA$ is skew-Hermitian (i.e.\ $\AA^*=-\AA$, where
$^*$ denotes the conjugate of the transpose),
$\BB:=i\AA$ is Hermitian (i.e.\ $\BB^* = \BB$).

The eigenvalues of a Hermitian matrix are real.
Indeed, assuming without loss of generality that
$\v^*\v=1$,
\begin{align*}
  \lambda &= \v^*\v\lambda
    = \v^*\BB\v
    = \v^*\BB^*\v
    = (\v^*\BB\v)^*
    = (\v^*\v\lambda)^*
\\  &= \v^*\v\lambda^*
    = \lambda^*,
\end{align*}
and eigenvectors for different eigenvalues are orthogonal:
\begin{align*}
  \v_2^*\v_1\lambda_1
    & = \v_2^*\BB\v_1
    = \v_2^*\BB^*\v_1
    = (\v_1^*\BB\v_2)^*
    = (\v_1^*\v_2\lambda_2)^*
\\  &= \v_2^*\v_1\lambda_2,
\end{align*}
which says that either $\v_2^*\v_1=0$ or $\lambda_1=\lambda_2$.

Note that if $(\v,\omega)$ is an eigenvector-eigenvalue pair for
$\BB$ then $(\v,i\omega)$ is an eigenvector-eigenvalue pair for $\AA$.

To find the eigenstructure we solve
\begin{gather}
 0 = (\AA-i\omega)\cdot\v =
  \begin{bmatrix}
    -i\omega & -\ai & \ae
 \\  \ai & -i(\gi+\omega) & 0
 \\ -\ae & 0 & i(\ge-\omega)
  \end{bmatrix} \cdot \v.
  \label{eigvecEq}
\end{gather}
If this has a nontrivial solution
then the first row is a linear combination of
the second two and we can ignore it.
The second two equations then show that
an eigenvector must be a multiple of the form
\def\bi{\beta_\i}
\def\be{\beta_\e}
\begin{align*}
 \begin{aligned}
  \v &=
  \begin{bmatrix}
    i\be\bi
 \\ \ai\be
 \\ \ae\bi
  \end{bmatrix},
  && \hbox{where} && \bi = \gi + \omega && \hbox{and} &&
  \be = \ge - \omega,
 \end{aligned}
\end{align*}
as is confirmed (for the last two rows) by computing
$(\AA-i\omega)\cdot\v$; the relation implied by the first
row reveals the characteristic equation.  Alternatively,
the calculation
\begin{gather*}
 \AA\cdot\v =
  \begin{bmatrix}
     0 & -\ai & \ae
 \\  \ai & -i\gi & 0
 \\ -\ae & 0 & i\ge
  \end{bmatrix} \cdot
  \begin{bmatrix}
    i\be\bi
 \\ \ai\be
 \\ \ae\bi
  \end{bmatrix}
  = 
  \begin{bmatrix}
    -\ai^2\be+\ae^2\bi
 \\ i\bi(\ai\be) - i\gi(\ai\be)
 \\ -i(\ae\bi)\be + i\ge(\ae\bi)
  \end{bmatrix}
  = \v i\omega =
    \begin{bmatrix}
      i\be\bi
   \\ \ai\be
   \\ \ae\bi
    \end{bmatrix}
  i\omega
\end{gather*}
shows that $\omega$ must satisfy
\begin{gather*}
   \be\bi\omega = \ai^2\be-\ae^2\bi,
\\ \omega = \bi-\gi,
\\ \omega = -\be+\ge.
\end{gather*}
The last two equations confirm that
\begin{gather*}
   \bi = \gi + \omega,
\\ \be = \ge - \omega,
\end{gather*}
and substituting these two relationships into
the first equation gives the characteristic
equation that an eigenvalue must satisfy:
\begin{gather*}
   (\ge - \omega)(\gi + \omega)\omega = \ai^2(\ge - \omega)-\ae^2(\gi + \omega).
\end{gather*}
Expanding in $\omega$ and collecting like terms gives
\begin{align}
 \label{characteristicPolynomialEquation}
 0 = \omega^3 + (\gi-\ge)\omega^2 - (\gi\ge+\ai^2+\ae^2)\omega
     + (\ai^2\ge-\ae^2\gi),
\end{align}
which can be solved using the formula for the roots
of a cubic with three real roots.

Note that
the eigenvector 
$ \v =
\begin{bmatrix}
    i\be\bi
 \\ \ai\be
 \\ \ae\bi
\end{bmatrix} $
is never zero; indeed,
$\ai$ and $\ae$ are strictly positive,
and
$\bi=\gi+\omega$ and $\be=\ge-\omega$
cannot both be zero
since otherwise
$\gi=-\omega$ and $\ge = \omega$,
contradicting that $\gi$ and $\ge$ are both strictly positive.

By the theory of Hermitian matrices
a full set of orthogonal eigenvectors must exist.
Since each eigenvector has a one-dimensional eigenspace,
there must be three distinct eigenvalues $\omega$.

\def\a{\underline{a}}
\def\b{\underline{b}}
Decompose $\v$ into real and imaginary parts:
\begin{gather*}
  \v = \a + i\b = 
    \begin{bmatrix}
      0
   \\ \ai\be
   \\ \ae\bi
    \end{bmatrix}
    + i
    \begin{bmatrix}
      \be\bi
   \\ 0
   \\ 0
    \end{bmatrix}.
\end{gather*}
Observe that the real and imaginary parts are orthogonal.
Note that
\begin{gather*}
  -i\v = \b - i\a = 
    \begin{bmatrix}
      \be\bi
   \\ 0
   \\ 0
    \end{bmatrix}
    - i
    \begin{bmatrix}
      0
   \\ \ai\be
   \\ \ae\bi
    \end{bmatrix}
\end{gather*}
is also an eigenvector with eigenvalue $i\omega$
and that these two eigenvectors are orthogonal:
$\v^*(-i\v) = 2\a\cdot\b = 0$.
The eigenvector-eigenvalue pair $(\v,i\omega)$
corresponds to the solution
\begin{align*}
  \v\exp(i\omega t) &= (\a+i\b)(\cos \omega t + i\sin\omega t)
 \\ &= (\a\cos\omega t - \b\sin\omega t) + i(\b\cos\omega t + \a\sin\omega t)
 \\ &=
    \begin{bmatrix}
     -\be\bi\sin\omega t
   \\ \ai\be\cos\omega t
   \\ \ae\bi\cos\omega t
    \end{bmatrix}
    + i
    \begin{bmatrix}
      \be\bi\cos\omega t
   \\ \ai\be\sin\omega t
   \\ \ae\bi\sin\omega t
    \end{bmatrix},
   \label{complexSoln}
\end{align*}
and the eigenvector-eigenvalue pair $(-i\v,i\omega)$
corresponds to the solution
\begin{align*}
  -i\v\exp(i\omega t) &= (\b-i\a)(\cos \omega t + i\sin\omega t)
 \\ &= (\b\cos\omega t + \a\sin\omega t) + i(\b\cos\omega t - \a\sin\omega t)
 \\ &= 
    \begin{bmatrix}
      \be\bi\cos\omega t
   \\ \ai\be\sin\omega t
   \\ \ae\bi\sin\omega t
    \end{bmatrix}
  + i
    \begin{bmatrix}
      \be\bi\sin\omega t
   \\-\ai\be\cos\omega t
   \\-\ae\bi\cos\omega t
    \end{bmatrix}.
\end{align*}
Observe that in each of these solutions
the ion and electron currents are in phase
and the electric field is 90 degrees out of phase
relative to them.

These two solutions are independent when interpreted (in 
$SO(2,\mathbb{R})$) as real solutions:
\begin{gather*}
  \begin{bmatrix}
    \Et_1 \\
    \Et_2 \\
    {\uit}{}_1 \\
    {\uit}{}_2 \\
    {\uet}{}_1 \\
    {\uet}{}_2 \\
  \end{bmatrix}
 = 
    \begin{bmatrix}
     -\be\bi\sin\omega t
   \\ \be\bi\cos\omega t
   \\ \ai\be\cos\omega t
   \\ \ai\be\sin\omega t
   \\ \ae\bi\cos\omega t
   \\ \ae\bi\sin\omega t
    \end{bmatrix}\ \hbox{and}\ 
  \begin{bmatrix}
    \Et_1 \\
    \Et_2 \\
    {\uit}{}_1 \\
    {\uit}{}_2 \\
    {\uet}{}_1 \\
    {\uet}{}_2 \\
  \end{bmatrix}
 = 
    \begin{bmatrix}
      \be\bi\cos\omega t
   \\ \be\bi\sin\omega t
   \\ \ai\be\sin\omega t
   \\-\ai\be\cos\omega t
   \\ \ae\bi\sin\omega t
   \\-\ae\bi\cos\omega t
    \end{bmatrix}.
\end{gather*}
Evaluated at time 0 these solutions are
\begin{gather*}
  \begin{bmatrix}
    \Et_1 \\
    \Et_2 \\
    {\uit}{}_1 \\
    {\uit}{}_2 \\
    {\uet}{}_1 \\
    {\uet}{}_2 \\
  \end{bmatrix}
 = 
    \begin{bmatrix}
      0
   \\ \be\bi
   \\ \ai\be
   \\ 0
   \\ \ae\bi
   \\ 0
    \end{bmatrix}\ \hbox{and}\ 
  \begin{bmatrix}
    \Et_1 \\
    \Et_2 \\
    {\uit}{}_1 \\
    {\uit}{}_2 \\
    {\uet}{}_1 \\
    {\uet}{}_2 \\
  \end{bmatrix}
 = 
    \begin{bmatrix}
      \be\bi
   \\ 0
   \\ 0
   \\-\ai\be
   \\ 0
   \\-\ae\bi
    \end{bmatrix}.
\end{gather*}

Note that orthogonality of complex solutions
is equivalent to orthogonality of real solutions.
So we have found three distinct imaginary eigenvalues
and 6 orthogonal eigenvectors for the original $6\times 6$
antisymmetric matrix.

\subsubsection{Parallel system}

The parallel system
\begin{gather*}
 \partial_t
  \begin{bmatrix}
    \Et_3 \\
    \uitC \\
    \uetC \\
  \end{bmatrix}
    =
  \begin{bmatrix}
    0 & -\ai & \ae
 \\  \ai & 0 & 0
 \\ -\ae & 0 & 0
  \end{bmatrix}
  \begin{bmatrix}
    \Et_3 \\
    \uitC \\
    \uetC \\
  \end{bmatrix}
\end{gather*}
is the special, singular case of
the perpendicular system \eqref{perpSystem}
when the magnetic field is zero.

In this case the system \eqref{eigvecEq} becomes
\begin{gather*}
 0 = (\AA-i\omega)\cdot\v =
  \begin{bmatrix}
    -i\omega & -\ai & \ae
 \\  \ai & -i\omega & 0
 \\ -\ae & 0 & -i\omega
  \end{bmatrix} \cdot \v.
\end{gather*}
So eigenvalue/eigenvector pairs are
\def\Op{\Omega_p}
\begin{align*}
 \begin{aligned}
  \omega&=0, &
  \v &=
  \begin{bmatrix}
    0
 \\ \ae
 \\ \ai
  \end{bmatrix}
  & \hbox{and} &&
  \omega=\pm\Omega_p, &&
  \v &=
  \begin{bmatrix}
    i\omega
 \\ \ai
 \\ -\ae
  \end{bmatrix},
 \end{aligned}
\end{align*}
where $\Op:=\sqrt{\ai^2+\ae^2}$ is the plasma frequency.

To get real solutions we look at the real and imaginary
parts of one of the complex-conjugate pair of solutions.
Choose $\omega=\Op$.  Write
\begin{gather*}
  \a+i\b = 
  \begin{bmatrix}
    0
 \\ \ai
 \\ -\ae
  \end{bmatrix}
  + i
  \begin{bmatrix}
    \Op
 \\ 0
 \\ 0
  \end{bmatrix}.
\end{gather*}
Analogous to \eqref{complexSoln},
the real and imaginary parts are real solutions:
\begin{gather*}
 \begin{aligned}
  \v\exp(i\Op t) &= (\a+i\b)(\cos \Op t + i\sin\Op t)
  \\
  &= (\a\cos\Op t - \b\sin\Op t) + i(\b\cos\Op t + \a\sin\Op t)
  \\
  &= \begin{bmatrix}
     -\Op\sin\Op t
   \\ \ai\cos\Op t
   \\-\ae\cos\Op t
    \end{bmatrix}
    + i
    \begin{bmatrix}
      \Op\cos\Op t
   \\ \ai\sin\Op t
   \\-\ae\sin\Op t
    \end{bmatrix}.
 \end{aligned}
\end{gather*}
So three orthogonal eigensolutions are
\begin{align*}
  \begin{bmatrix}
    0
 \\ \ae
 \\ \ai
  \end{bmatrix},
 &&
  \begin{bmatrix}
   -\Op\sin\Op t
 \\ \ai\cos\Op t
 \\-\ae\cos\Op t
  \end{bmatrix},
 &&
  \begin{bmatrix}
    \Op\cos\Op t
 \\ \ai\sin\Op t
 \\-\ae\sin\Op t
  \end{bmatrix}.
\end{align*}
Evaluated at time $0$ these solutions are
\begin{align*}
  \begin{bmatrix}
    0
 \\ \ae
 \\ \ai
  \end{bmatrix},
 &&
  \begin{bmatrix}
   0
 \\ \ai
 \\-\ae
  \end{bmatrix},
 &&
  \begin{bmatrix}
    \Op
 \\ 0
 \\ 0
  \end{bmatrix}.
\end{align*}

\def\T{\mathrm{T}}
\textbf{Agreement with perpendicular system.}
If we take the limit as $|\B|\rightarrow 0$ in the
perpendicular system we expect the solutions to
decouple into solutions for
$
  \left(
    \Et_1,
    {\uit}{}_1,
    {\uet}{}_1
  \right)^\T
$
and
$
  \left(
    \Et_2,
    {\uit}{}_2,
    {\uet}{}_2
  \right)^\T
$
that agree with the solutions
for the parallel system.
As $|\B|\rightarrow 0$,
$\gi\rightarrow 0$ and
$\ge\rightarrow 0$
and so $\bi \rightarrow \omega$ and
$\be \rightarrow -\omega$.
For the limiting eigenfrequency
$\omega = \Op$
the parallel solutions
\begin{gather*}
  \begin{bmatrix}
    \Et_1 \\
    \Et_2 \\
    {\uit}{}_1 \\
    {\uit}{}_2 \\
    {\uet}{}_1 \\
    {\uet}{}_2 \\
  \end{bmatrix}
 = 
    \begin{bmatrix}
     -\be\bi\sin\omega t
   \\ \be\bi\cos\omega t
   \\ \ai\be\cos\omega t
   \\ \ai\be\sin\omega t
   \\ \ae\bi\cos\omega t
   \\ \ae\bi\sin\omega t
    \end{bmatrix}\ \hbox{and}\ 
  \begin{bmatrix}
    \Et_1 \\
    \Et_2 \\
    {\uit}{}_1 \\
    {\uit}{}_2 \\
    {\uet}{}_1 \\
    {\uet}{}_2 \\
  \end{bmatrix}
 = 
    \begin{bmatrix}
      \be\bi\cos\omega t
   \\ \be\bi\sin\omega t
   \\ \ai\be\sin\omega t
   \\-\ai\be\cos\omega t
   \\ \ae\bi\sin\omega t
   \\-\ae\bi\cos\omega t
    \end{bmatrix}
\end{gather*}
when divided by $\be=-\Op$ become
\begin{gather*}
  \begin{bmatrix}
    \Et_1 \\
    \Et_2 \\
    {\uit}{}_1 \\
    {\uit}{}_2 \\
    {\uet}{}_1 \\
    {\uet}{}_2 \\
  \end{bmatrix}
 = 
    \begin{bmatrix}
      -\Op\sin\Op t
   \\  \Op\cos\Op t
   \\  \ai\cos\Op t
   \\  \ai\sin\Op t
   \\-\ae\cos\Op t
   \\-\ae\sin\Op t
    \end{bmatrix}\ \hbox{and}\ 
  \begin{bmatrix}
    \Et_1 \\
    \Et_2 \\
    {\uit}{}_1 \\
    {\uit}{}_2 \\
    {\uet}{}_1 \\
    {\uet}{}_2 \\
  \end{bmatrix}
 = 
    \begin{bmatrix}
       \Op\cos\Op t
   \\  \Op\sin\Op t
   \\  \ai\sin\Op t
   \\ -\ai\cos\Op t
   \\ -\ae\sin\Op t
   \\  \ae\cos\Op t
    \end{bmatrix},
\end{gather*}
and for the limiting eigenfrequency
$\omega = -\Op$
the parallel solutions
when divided by $\be=\Op$ become
\begin{gather*}
  \begin{bmatrix}
    \Et_1 \\
    \Et_2 \\
    {\uit}{}_1 \\
    {\uit}{}_2 \\
    {\uet}{}_1 \\
    {\uet}{}_2 \\
  \end{bmatrix}
 = 
    \begin{bmatrix}
      -\Op\sin\Op t
   \\ -\Op\cos\Op t
   \\ \ai\cos\Op t
   \\-\ai\sin\Op t
   \\ -\ae\cos\Op t
   \\  \ae\sin\Op t
    \end{bmatrix}\ \hbox{and}\ 
  \begin{bmatrix}
    \Et_1 \\
    \Et_2 \\
    {\uit}{}_1 \\
    {\uit}{}_2 \\
    {\uet}{}_1 \\
    {\uet}{}_2 \\
  \end{bmatrix}
 = 
    \begin{bmatrix}
      -\Op\cos\Op t
   \\  \Op\sin\Op t
   \\-\ai\sin\Op t
   \\-\ai\cos\Op t
   \\ \ae\sin\Op t
   \\ \ae\cos\Op t
    \end{bmatrix}.
\end{gather*}
When projected onto axis $1$ 
these solutions all agree with the
second and third eigensolutions for the parallel
component, and likewise for axis $2$.

For the limiting eigenfrequency $\omega=0$,
$\be:=\ge+\omega$ and
$\bi:=\gi-\omega$ both
go to zero (since $\ge$ and $\gi$ go to zero as
$\B$ goes to zero).  We may infer that $\be\bi$
goes very quickly to zero.  When $\omega$ is small
$\cos\omega t\approx 1$ and $\sin\omega t\approx 0$.
So for small magnetic field we expect 
\begin{gather*}
  \begin{bmatrix}
    \Et_1 \\
    \Et_2 \\
    {\uit}{}_1 \\
    {\uit}{}_2 \\
    {\uet}{}_1 \\
    {\uet}{}_2 \\
  \end{bmatrix}
 = 
    \begin{bmatrix}
     -\be\bi\sin\omega t
   \\ \be\bi\cos\omega t
   \\ \ai\be\cos\omega t
   \\ \ai\be\sin\omega t
   \\ \ae\bi\cos\omega t
   \\ \ae\bi\sin\omega t
    \end{bmatrix}
 \approx
    \begin{bmatrix}
      0
   \\  0
   \\ \ai\be
   \\  0
   \\ \ae\bi
   \\  0
    \end{bmatrix},
\end{gather*}
which agrees with the direction of
the expected limiting eigensolution
\begin{gather*}
  \begin{bmatrix}
    \Et_1 \\
    {\uit}{}_1 \\
    {\uet}{}_1 \\
  \end{bmatrix}
 \approx
    \begin{bmatrix}
      0
   \\ \ae
   \\ \ai
    \end{bmatrix}
\end{gather*}
if $\frac{\be}{\bi} \rightarrow \left(\frac{\ae}{\ai}\right)^2$
as $\B\rightarrow 0$.
I do not see how to show this in general, but in the
neutral case where $n_\i=n_\e$, $\ai^2\ge=\ae^2\gi$, so
the constant term vanishes in
the characteristic polynomial equation
\eqref{characteristicPolynomialEquation},
$\omega=0$ is
always an eigenvalue, and $\bi=\gi$ and $\be=\ge$,
so $\frac{\be}{\bi} = \left(\frac{\ae}{\ai}\right)^2$
as needed.

\subsection{The pressure tensor system}

Ignoring interspecies collisions, the pressure tensor
evolution equation \eqref{PTsource}
with linear closure \eqref{isoRTclosureCopy}
is
\begin{gather*}
  \partial_t \Pressure_\s = (q_\s/ m_\s)\SymB(\Pressure_\s\times\B)
    -\tau_\s\inv\dPT_\s,
\end{gather*}
Observe that temperature isotropization leaves
temperature invariant.
So if one uses the closure \eqref{isoPeriod},
\begin{gather}
  \tau = \tau_0 \sqrt{m}\frac{T^{3/2}}{n},
\end{gather}
then this is a linear ODE with
constant coefficients.
The $\Pressure\times\B$ term
rotates the pressure tensor around the magnetic
field.  The relaxation term relaxes the pressure
toward isotropy.  These two operations commute.
So we can trivially solve this ODE exactly.

\subsection{Rotation of the pressure tensor}

\def\gyrofreq{\omega^\mathrm{g}}
\def\dt{\mathrm{d}t}
The $\Pressure\times\B$ term
rotates the pressure tensor around the magnetic
field vector.  The rate of rotation is the species gyrofrequency,
so the angle of rotation of the ion pressure tensor
in time interval $\dt$ is $\gi\dt$.

\def\e{\mathbf{e}}
\def\ei{\e_i}
\def\ej{\e_j}
\def\em{\e_m}
\def\en{\e_n}
\def\r{\mathbf{r}}
\def\ri{\e_i'}
\def\rj{\e_j'}
\def\Pij{\Pressure_{ij}}
\def\Pmn{\Pressure_{mn}}
\def\Prirj{\Pressure_{i'j'}}
Let $\ei$ denote the $i$th standard basis vector.
Let $\ri(t)$ denote the rotated version of $\ei$.
Let $\Pressure(t)$ denote the rotated pressure tensor.
The pressure tensor components are
$\Pmn(t) = \em\cdot\Pressure\cdot\en$.
The pressure tensor components 
are invariant in a rotating (primed) coordinate frame:
\begin{gather*}
  \Pressure(t) = \Pij(0)\ri\rj.
\end{gather*}
Therefore, the components in the standard basis are:
\def\dt{\,\mathrm{d}t\,}
\def\Bs{\B_\s}
\def\qms{\frac{q_\s}{m_\s}}
\def\R{\vec{R}}
\def\gs{\omega_\s}
\def\u{\mathbf{u}}
\begin{gather*}
  \Pmn(t) = \Pij(0)(\ri\cdot\em)(\rj\cdot\en).
\end{gather*}
Thus, to evolve the pressure tensor $\Pressure_\s$
for species $\s$ over a time interval $\dt$,
we need to apply to the standard basis vectors
a rotation with rotation vector $\R := \frac{q_\s}{m_\s}\B\dt$,
i.e.\ with direction $\bhat:=\B/|\B|$ and angle $\theta:=\gs\dt$,
where $\gs:=|\B|\qms$ is the gyrofrequency of species $\s$.

To rotate a vector $\u$ by the vector $\R=\theta\bhat$,
where $\theta=|\R|$, decompose it into parallel
and perpendicular components and rotate the perpendicular
component:
\def\upara{\u_\parallel}
\def\uperp{\u_\perp}
\def\uxb{\u\times\bhat}
\def\ubb{\u\cdot\bhat\bhat}
\def\r{\mathbf{u}'}
\begin{gather*}
   \upara = \ubb
\\ \uperp = \u-\upara
\\ \u = \upara + \uperp
\end{gather*}
The rotated vector is
\def\tH{\frac{\theta}{2}}
\begin{align*}
 \r&=\upara + (\cos\theta)\uperp + (\sin\theta)\uxb.
\\ &=\u(\cos\theta) + (1-\cos\theta)\upara + (\sin\theta)\uxb
\end{align*}
We can avoid renormalizing the rotation vector
if we use the sine cardinal function.
\def\sinc{\,\mathrm{sinc}}
\begin{alignat*}{6}
  \r &=\u(\cos\theta) &&+ 2\left(\sin\tH\right)^2\ubb
 \\  
     &&&+ 2\left(\cos\tH\right)\left(\sin\tH\right)\uxb
\\
     &=\u(\cos\theta) &&+ \frac{1}{2}\left(\sinc\tH\right)^2\u\cdot\R\R
 \\  
     &&&+ \left(\cos\tH\right)\left(\sinc\tH\right)\u\times\R
\end{alignat*}

\def\T{\mathrm{T}}
\def\RR{\tensorb{R}}
Recall that $\ri$ is the rotated version of the elementary
basis vector $\ei$.  To express the components of the
rotation matrix $\RR_{ij}:=\ei\cdot\rj$, adopt the abbreviations
$c:=\cos\theta$ and $s:=\sin\theta$.
The rotated vector is
\def\idtens{\mathbb{I}}
\def\tH{\frac{\theta}{2}}
\begin{align*}
 \rj &= c \ej + (1-c)\ej\cdot\bhat\bhat + s \ej\times\bhat
 \\  &= \ej\cdot\underbrace{\left(c\idtens+(1-c)\bhat\bhat+s\idtens\times\bhat\right)
          }_{\hbox{$\RR^\T$}}
\end{align*}
To determine the components of $\idtens\times\bhat$,
match up the identity
\begin{gather*}
  \u\times\bhat = \u\cdot\idtens\times\bhat = (\idtens\times\bhat)^\T\cdot\u
\end{gather*}
with the coordinate expansion
\begin{equation*}
  \u\times\bhat = 
    \begin{pmatrix}
      u_2 b_3 - u_3 b_2 \\
      u_3 b_1 - u_1 b_3 \\
      u_1 b_2 - u_2 b_1
    \end{pmatrix}
     = 
    \underbrace{
    \begin{pmatrix}
      0 & b_3 & -b_2 \\
      -b_3 & 0 & b_1 \\
      b_2 & -b_1 & 0
    \end{pmatrix}
    }_{\hbox{$(\idtens\times\bhat)^\T$}}
    \begin{pmatrix}
      u_1 \\
      u_2 \\
      u_3
    \end{pmatrix}
\end{equation*}
So the rotation matrix $\RR$ is:
\def\ph{\phantom{b_3}}
\begin{gather*}
  \begin{bmatrix}
    b_1b_1(1-c) +c \ph & b_1b_2(1-c) +s b_3 & b_1b_3(1-c) -s b_2 \\
    b_2b_1(1-c) -s b_3 & b_2b_2(1-c) +c \ph & b_2b_3(1-c) +s b_1 \\
    b_3b_1(1-c) +s b_2 & b_3b_2(1-c) -s b_1 & b_3b_3(1-c) +c \ph 
  \end{bmatrix},
\end{gather*}
where we are free to make the replacements
\begin{alignat*}{4}
  (1-c)&b_i b_j &&=& (1/2)\sinc^2(\theta/2)&R_i R_j,
 \\
  s &b_i &&=& \sinc(\theta/2)\cos(\theta/2)&R_i.
\end{alignat*}

%
%
%

}

\section{Eigenstructure of the flux Jacobian for a ten-moment gas}
{
\label{tenMomentEigs}.
\newcommand{\pd}[1][]{\ensuremath{\partial_{#1}}}
\newcommand{\qm}[1][]{\ensuremath{\frac{q_{#1}}{m_{#1}}}}
\def\mom{\mathbf{M}}
\def\transpose{\mathrm{T}}

I have verified the ten-moment eigenstructure calculated
in this section with extensive computational checks.

\subsection{Ten-moment system in conservation form}
To perform limiting in shock-capturing methods,
we consider the hyperbolic ten-moment system in balance law form
\eqref{eqn:conslawFramework},
\def\uu{{\underline{u}}}
\def\FF{\underline{\textbf{f}}}
\def\ff{\underline{f}}
\def\ss{\underline{s}}
\begin{gather}
  \label{eqn:conslawFrameworkCopy}
  \partial_t \uu + \Div\FF = \ss,
\end{gather}
and compute the eigenstructure of the flux Jacobian
$\partial\ff/\partial\uu$,
where $\ff:=\ebas_x\dotp\FF$.
This requires that we express the flux $\FF$
in terms of conserved variables.

For a single species, the full 10-moment system
in balance law form
specifies the flux and sources of 
mass density $\mdens$, momentum $\mom=\mdens\u$,
and the energy tensor $\ET=\mdens\u\u+\PT$,
\begin{gather*}
   \pd[t]\mdens+\Div(\mdens\u)=0,
\\ \pd[t](\mdens\u) + \Div(\mdens\u\u+\PT)
  = \frac{q}{m}\mdens(\E+\u\times\B),
\\ \pd[t](\mdens\u\u+\PT)
    + \Div\big(\mdens\u\u\u+\SymC(\u\PT)\big)
   = \qm\SymB\big(\mdens\u\E+(\PT+\mdens\u\u)\times\B\big).
\end{gather*}
To write this entirely in terms of conserved variables,
we note that $\u=\mom/\mdens$ and 
$\PT=\ET-\mom\mom/\mdens$.
So 
\begin{gather}
   \pd[t]\mdens+\Div\mom=0, \nonumber
\\ \pd[t]\mom + \Div\ET
  = \qm(\mdens\E+\mom\times\B), \nonumber
\\ \pd[t]\ET
    +\Div\Big(
      \frac{\SymC(\mom\ET)}{\mdens}
      -\frac{2\mom\mom\mom}{\mdens^2}
     \Big)
   = \qm\SymB\big(\mom\E+\ET\times\B\big).
   \label{eq:tenMomentConserved}
\end{gather}

\subsection{Quasilinear system in primitive quantities}
Shock-capturing limiters need the eigenstructure
of the quasilinearized system.
To calculate the eigenstructure, we put the system in 
quaslinear form.
The eigenstructure is most easily calculated in primitive variables.
In primitive variables and quasilinear form the full system is
\begin{gather*}
   \pd[t]\mdens + \u\cdot\nabla\mdens+\mdens\Div\u = 0
\\ \pd[t]\u + \u\cdot\nabla\u + \frac{\Div\PT}{\mdens}
    = \qm(\E+\u\times\B)
\\ \pd[t]\PT
       + \u\cdot\nabla\PT
       + \PT \Div\u + \SymB(\PT\cdot\nabla\u)
   = \qm\SymB(\PT\times\B)
\end{gather*}

\subsection{Quasilinear 1-D system in primitive variables}

As discussed in sections \ref{1Dlimiting} and \ref{2DcartesianLimiting},
to perform limiting on a 2D Cartesian mesh
we need to be able to compute the eigenstructure
of the flux Jacobian $\ff_\uu$ for the 1D problem
of the form
\begin{gather*}
  \partial_t \uu + \pd[x]\ff(\uu) = \ss
\end{gather*}
which obtains when the problem is homogeneous
perpendicular to $x$.

For a problem homogeneous perpendicular to the $x$ axis
equation \eqref{eqn:conslawFrameworkCopy} simplifies to
\begin{gather}
  \partial_t \uu + \pd[x]\ff = \ss,
\end{gather}
where $\ff:=\ebas_x\dotp\FF$.

Assuming homogeneity in all space dimensions except the first
($x$), the 1-D system in primitive variables becomes
\begin{gather*}
   \pd[t]\mdens + u_1\pd[x]\mdens+\mdens\pd[x]u_1 = 0,
\\ \pd[t]\u + u_1\pd[x]\u + \frac{\pd[x]\PT_{1\cdot}}{\mdens}
    = \qm (\E+\u\times\B),
\\ \pd[t]\PT
       + u_1\pd[x]\PT_{1\cdot}
       + \PT \pd[x]u_1 + \SymB(\PT_{\cdot 1}\pd[x]\u)
   = \qm\SymB(\PT\times\B).
\end{gather*}
We align derivatives to prepare to put this quasilinear system in
matrix form.
\begin{alignat*}{5}
   0 &= 
   \pd[t]\mdens &&+ u_1\pd[x]\mdens 
            &&+\mdens\pd[x]u_1,
\\
   \qm (\E+\u\times\B) &=
   \pd[t]\u &&&&+ u_1\pd[x]\u
               &&+ \frac{\pd[x]\PT_{1\cdot}}{\mdens},
\\
   \qm\SymB(\PT\times\B) &=
   \pd[t]\PT
         &&&&+\PT \pd[x]u_1+\SymB(\PT_{\cdot 1}\pd[x]\u)
               &&+u_1\pd[x]\PT_{1\cdot}.
\end{alignat*}
In matrix form this reads
{\scriptsize
\begin{gather*}
 \!\!\!\!\!\!\!\!
   \begin{bmatrix}
     \mdens
  \\ u_1
  \\ u_2
  \\ u_3
  \\ \PT_{11}
  \\ \PT_{12}
  \\ \PT_{13}
  \\ \PT_{23}
  \\ \PT_{22}
  \\ \PT_{33}
   \end{bmatrix}_t
 \!\!
 +
   \begin{bmatrix}
     u_1&\mdens&0&0&0&0&0&0&0&0
  \\ 0&u_1&0&0&1/\mdens&0&0&0&0&0
  \\ 0&0&u_1&0&0&1/\mdens&0&0&0&0
  \\ 0&0&0&u_1&0&0&1/\mdens&0&0&0
  \\ 0&3\PT_{11}&0&0&u_1&0&0&0&0&0
  \\ 0&2\PT_{12}&\PT_{11}&0&0&u_1&0&0&0&0
  \\ 0&2\PT_{13}&0&\PT_{11}&0&0&u_1&0&0&0
  \\ 0&\PT_{23}&\PT_{31}&\PT_{21}&0&0&0&u_1&0&0
  \\ 0&\PT_{22}&2\PT_{12}&0&0&0&0&0&u_1&0
  \\ 0&\PT_{33}&0&2\PT_{13}&0&0&0&0&0&u_1
   \end{bmatrix}
  \cdot
   \begin{bmatrix}
     \mdens
  \\ u_1
  \\ u_2
  \\ u_3
  \\ \PT_{11}
  \\ \PT_{12}
  \\ \PT_{13}
  \\ \PT_{23}
  \\ \PT_{22}
  \\ \PT_{33}
   \end{bmatrix}_x
 \!\!\!
 = \qm
   \begin{bmatrix}
     0
  \\ E_1+(u_2 B_3 - u_3 B_2)
  \\ E_2+(u_3 B_1 - u_1 B_3)
  \\ E_3+(u_1 B_2 - u_2 B_1)
  \\ 2(\PT_{12} B_3 - \PT_{13} B_2)
  \\ \PT_{13}B_1-\PT_{23}B_2+(\PT_{22}-\PT_{11})B_3
  \\ \PT_{32}B_3-\PT_{12}B_1+(\PT_{11}-\PT_{33})B_2
  \\ \PT_{21}B_2-\PT_{31}B_3+(\PT_{33}-\PT_{22})B_1
  \\ 2(\PT_{23}B_1-\PT_{21}B_3)
  \\ 2(\PT_{31}B_2-\PT_{32}B_1)
   \end{bmatrix}.
\end{gather*}
}
\subsection{Eigenstructure for primitive variables}
If we neglect the source term, the eigenvalues of the matrix represent wave speeds,
and the corresponding eigenvectors represent the corresponding waves.
Let $u:=:u_1+c$ denote wave speed (i.e., $c$ is wave
speed in the reference frame moving with the fluid, where $u_1=0$).
To find the eigenstructure of the quasilinearized
system we put the matrix
{\small
\begin{gather*}
   \begin{bmatrix}
     -c&\mdens&0&0&0&0&0&0&0&0
  \\ 0&-c&0&0&1/\mdens&0&0&0&0&0
  \\ 0&0&-c&0&0&1/\mdens&0&0&0&0
  \\ 0&0&0&-c&0&0&1/\mdens&0&0&0
  \\ 0&3\PT_{11}&0&0&-c&0&0&0&0&0
  \\ 0&2\PT_{12}&\PT_{11}&0&0&-c&0&0&0&0
  \\ 0&2\PT_{13}&0&\PT_{11}&0&0&-c&0&0&0
  \\ 0&\PT_{23}&\PT_{31}&\PT_{21}&0&0&0&-c&0&0
  \\ 0&\PT_{22}&2\PT_{12}&0&0&0&0&0&-c&0
  \\ 0&\PT_{33}&0&2\PT_{13}&0&0&0&0&0&-c
   \end{bmatrix}
\end{gather*}
}
in upper triangular form.

\subsubsection{Right primitive eigenstructure}
For the right eigenvectors we combine rows to do so:
{\small
\begin{gather*}
   \begin{bmatrix}
     -c&\mdens&0&0&0&0&0&0&0&0
  \\ 0&-\mdens c^2&0&0&c&0&0&0&0&0
  \\ 0&0&-\mdens c^2&0&0&c&0&0&0&0
  \\ 0&0&0&-\mdens c^2&0&0&c&0&0&0
  \\ 0&3\PT_{11}&0&0&-c&0&0&0&0&0
  \\ 0&2\PT_{12}&\PT_{11}&0&0&-c&0&0&0&0
  \\ 0&2\PT_{13}&0&\PT_{11}&0&0&-c&0&0&0
  \\ 0&\PT_{23}&\PT_{31}&\PT_{21}&0&0&0&-c&0&0
  \\ 0&\PT_{22}&2\PT_{12}&0&0&0&0&0&-c&0
  \\ 0&\PT_{33}&0&2\PT_{13}&0&0&0&0&0&-c
   \end{bmatrix}
  \cdot
   \begin{bmatrix}
     \mdens
  \\ u_1
  \\ u_2
  \\ u_3
  \\ \PT_{11}
  \\ \PT_{12}
  \\ \PT_{13}
  \\ \PT_{23}
  \\ \PT_{22}
  \\ \PT_{33}
   \end{bmatrix}'
 = 0.
\end{gather*}
}
So
{\small
\begin{gather*}
   \begin{bmatrix}
     -c&\mdens&0&0&0&0&0&0&0&0
  \\ 0&3\PT_{11}-\mdens c^2&0&0&0&0&0&0&0&0
  \\ 0&2\PT_{12}&\PT_{11}-\mdens c^2&0&0&0&0&0&0&0
  \\ 0&2\PT_{13}&0&\PT_{11}-\mdens c^2&0&0&0&0&0&0
  \\ 0&3\PT_{11}&0&0&-c&0&0&0&0&0
  \\ 0&2\PT_{12}&\PT_{11}&0&0&-c&0&0&0&0
  \\ 0&2\PT_{13}&0&\PT_{11}&0&0&-c&0&0&0
  \\ 0&\PT_{23}&\PT_{31}&\PT_{21}&0&0&0&-c&0&0
  \\ 0&\PT_{22}&2\PT_{12}&0&0&0&0&0&-c&0
  \\ 0&\PT_{33}&0&2\PT_{13}&0&0&0&0&0&-c
   \end{bmatrix}
  \cdot
   \begin{bmatrix}
     \mdens
  \\ u_1
  \\ u_2
  \\ u_3
  \\ \PT_{11}
  \\ \PT_{12}
  \\ \PT_{13}
  \\ \PT_{23}
  \\ \PT_{22}
  \\ \PT_{33}
   \end{bmatrix}'
 = 0.
\end{gather*}
}
Denote the fast and slow speeds by
\begin{equation*}
     c_f:=\sqrt{\frac{3\PT_{11}}{\mdens}},
  \\ c_s:=\sqrt{\frac{\PT_{11}}{\mdens}}.
\end{equation*}
In the most difficult case, where $c=\pm c_f$,
we have that $3\PT_{11}-\mdens c_f^2=0$, so
$\PT_{11}-\mdens c_f^2=-2\PT_{11}$.

So right eigenvectors are:
\newline
{\small
\begin{tabular}{|c|c|c|c|}
   \hline
     $c$:
   & $\pm \sqrt{\frac{3\PT_{11}}{\mdens}}$
   & $\pm \sqrt{\frac{\PT_{11}}{\mdens}}$
   & $0$
\\ \hline \hline
     $
       \begin{bmatrix}
          \mdens
       \\ u_1
       \\ u_2
       \\ u_3
       \\ \PT_{11}
       \\ \PT_{12}
       \\ \PT_{13}
       \\ \PT_{23}
       \\ \PT_{22}
       \\ \PT_{33}
       \end{bmatrix}' \propto
     $
   &
     $
       \begin{bmatrix}
          \mdens\PT_{11}
       \\ \pm c_f\PT_{11}
       \\ \pm c_f\PT_{12}
       \\ \pm c_f\PT_{13}
       \\ 3\PT_{11}\PT_{11}
       \\ 3\PT_{12}\PT_{11}
       \\ 3\PT_{13}\PT_{11}
       \\ \PT_{23}\PT_{11}+2\PT_{13}\PT_{12}
       \\ \PT_{22}\PT_{11}+2\PT_{12}\PT_{12}
       \\ \PT_{33}\PT_{11}+2\PT_{13}\PT_{13}
       \end{bmatrix}
     $
   &
     $
       \begin{bmatrix}
          0
       \\ 0
       \\ \pm c_s
       \\ 0
       \\ 0
       \\ \PT_{11}
       \\ 0
       \\ \PT_{13}
       \\ 2\PT_{12}
       \\ 0
       \end{bmatrix},
       \begin{bmatrix}
          0
       \\ 0
       \\ 0
       \\ \pm c_s
       \\ 0
       \\ 0
       \\ \PT_{11}
       \\ \PT_{12}
       \\ 0
       \\ 2\PT_{13}
       \end{bmatrix}
     $
   &
     $
       \begin{bmatrix}
          1
       \\ 0
       \\ 0
       \\ 0
       \\ 0
       \\ 0
       \\ 0
       \\ 0
       \\ 0
       \\ 0
       \end{bmatrix},
       \begin{bmatrix}
          0
       \\ 0
       \\ 0
       \\ 0
       \\ 0
       \\ 0
       \\ 0
       \\ 1
       \\ 0
       \\ 0
       \end{bmatrix},
       \begin{bmatrix}
          0
       \\ 0
       \\ 0
       \\ 0
       \\ 0
       \\ 0
       \\ 0
       \\ 0
       \\ 1
       \\ 0
       \end{bmatrix},
       \begin{bmatrix}
          0
       \\ 0
       \\ 0
       \\ 0
       \\ 0
       \\ 0
       \\ 0
       \\ 0
       \\ 0
       \\ 1
       \end{bmatrix}
     $
 \\
   \hline
\end{tabular}
}

\subsubsection{Left primitive eigenstructure}
To find the left eigenstructure, we combine columns to reduce to upper
triangular form.
{\small
\begin{gather*}
   \begin{bmatrix}
     \mdens
  \\ u_1
  \\ u_2
  \\ u_3
  \\ \PT_{11}
  \\ \PT_{12}
  \\ \PT_{13}
  \\ \PT_{23}
  \\ \PT_{22}
  \\ \PT_{33}
   \end{bmatrix}'^\transpose
  \cdot
   \begin{bmatrix}
     -c&\mdens&0&0&0&0&0&0&0&0
  \\ 0&-c&0&0&1/\mdens&0&0&0&0&0
  \\ 0&0&-c&0&0&1/\mdens&0&0&0&0
  \\ 0&0&0&-c&0&0&1/\mdens&0&0&0
  \\ 0&3\PT_{11}&0&0&-c&0&0&0&0&0
  \\ 0&2\PT_{12}&\PT_{11}&0&0&-c&0&0&0&0
  \\ 0&2\PT_{13}&0&\PT_{11}&0&0&-c&0&0&0
  \\ 0&\PT_{23}&\PT_{31}&\PT_{21}&0&0&0&-c&0&0
  \\ 0&\PT_{22}&2\PT_{12}&0&0&0&0&0&-c&0
  \\ 0&\PT_{33}&0&2\PT_{13}&0&0&0&0&0&-c
   \end{bmatrix}
 = 0.
\end{gather*}
}
In case $c\ne 0$ this reduces to
{\small
\begin{gather*}
   \begin{bmatrix}
     \mdens
  \\ u_1
  \\ u_2
  \\ u_3
  \\ \PT_{11}
  \\ \PT_{12}
  \\ \PT_{13}
  \\ \PT_{23}
  \\ \PT_{22}
  \\ \PT_{33}
   \end{bmatrix}'^\transpose
  \cdot
   \begin{bmatrix}
     -c&0&0&0&0&0&0&0&0&0
  \\ 0&-c&0&0&0&0&0&0&0&0
  \\ 0&0&-c&0&0&0&0&0&0&0
  \\ 0&0&0&-c&0&0&0&0&0&0
  \\ 0&3\PT_{11}&0&0&
       3\PT_{11}-\mdens c^2&0&0&0&0&0
  \\ 0&2\PT_{12}&\PT_{11}&0&
       2\PT_{12}&\PT_{11}-\mdens c^2&0&0&0&0
  \\ 0&2\PT_{13}&0&\PT_{11}&
       2\PT_{13}&0&\PT_{11}-\mdens c^2&0&0&0
  \\ 0&\PT_{23}&\PT_{31}&\PT_{21}&
       \PT_{23}&\PT_{31}&\PT_{21}&-c&0&0
  \\ 0&\PT_{22}&2\PT_{12}&0&
       \PT_{22}&2\PT_{12}&0&0&-c&0
  \\ 0&\PT_{33}&0&2\PT_{13}&
       \PT_{33}&0&2\PT_{13}&0&0&-c
   \end{bmatrix}
 = 0.
\end{gather*}
}
Here the most difficult case is $c=\pm c_s$,
when $\PT_{11}=\mdens c^2$,
so  $3\PT_{11}-\mdens c^2=2\PT_{11}$.
The left eigenvectors for nonzero speeds are thus:

{\small
\begin{tabular}{|c|c|c|}
   \hline
     $c$:
   & $\pm \sqrt{\frac{3\PT_{11}}{\mdens}}$
   & $\pm \sqrt{\frac{\PT_{11}}{\mdens}}$
\\ \hline \hline
     $
       \begin{bmatrix}
          \mdens
       \\ u_1
       \\ u_2
       \\ u_3
       \\ \PT_{11}
       \\ \PT_{12}
       \\ \PT_{13}
       \\ \PT_{23}
       \\ \PT_{22}
       \\ \PT_{33}
       \end{bmatrix}' \propto
     $
   &
     $
       \begin{bmatrix}
          0
       \\ 3\PT_{11}
       \\ 0
       \\ 0
       \\ \pm c_f
       \\ 0
       \\ 0
       \\ 0
       \\ 0
       \\ 0
       \end{bmatrix}
     $
   &
     $
       \begin{bmatrix}
          0
       \\ -\PT_{12}\PT_{11}
       \\ \PT_{11}\PT_{11}
       \\ 0
       \\ \mp c_s\PT_{12}
       \\ \pm c_s\PT_{11}
       \\ 0
       \\ 0
       \\ 0
       \\ 0
       \end{bmatrix},
       \begin{bmatrix}
          0
       \\ -\PT_{13}\PT_{11}
       \\ 0
       \\ \PT_{11}\PT_{11}
       \\ \mp c_s\PT_{13}
       \\ 0
       \\ \pm c_s\PT_{11}
       \\ 0
       \\ 0
       \\ 0
       \end{bmatrix}
     $
 \\
   \hline
\end{tabular}
}

\vspace{8pt}
The left eigenvectors for $c=0$, are readily obtained from the
original system:
\begin{gather*}
       \begin{bmatrix}
          \mdens
       \\ u_1
       \\ u_2
       \\ u_3
       \\ \PT_{11}
       \\ \PT_{12}
       \\ \PT_{13}
       \\ \PT_{23}
       \\ \PT_{22}
       \\ \PT_{33}
       \end{bmatrix}' =
     \alpha_1
       \begin{bmatrix}
          3\PT_{11}
       \\ 0
       \\ 0
       \\ 0
       \\ -\mdens
       \\ 0
       \\ 0
       \\ 0
       \\ 0
       \\ 0
       \end{bmatrix}
   + \alpha_2
       \begin{bmatrix}
          0
       \\ 0
       \\ 0
       \\ 0
       \\ 4\PT_{12}\PT_{13}-\PT_{23}\PT_{11}
       \\ -3\PT_{11}\PT_{13}
       \\ -3\PT_{11}\PT_{12}
       \\ 3\PT_{11}^2
       \\ 0
       \\ 0
       \end{bmatrix}
   + \alpha_3
       \begin{bmatrix}
          0
       \\ 0
       \\ 0
       \\ 0
       \\ 4\PT_{12}^2-\PT_{11}\PT_{22}
       \\ -6\PT_{12}\PT_{11}
       \\ 0
       \\ 0
       \\ 3\PT_{11}^2
       \\ 0
       \end{bmatrix}
   + \alpha_4
       \begin{bmatrix}
          0
       \\ 0
       \\ 0
       \\ 0
       \\ 4\PT_{13}^2-\PT_{11}\PT_{33}
       \\ 0
       \\ -6\PT_{13}\PT_{11}
       \\ 0
       \\ 0
       \\ 3\PT_{11}^2
       \end{bmatrix}
\end{gather*}

\subsection{Eigenstructure for ``conserved'' variables}

\def\cons{{\underline{q}}}
\def\prim{{\underline{p}}}
\def\energyVec{\widetilde\ET}
\def\pressureVec{\widetilde\PT}
Let
\begin{gather*}
  \cons:=(\mdens,\mom,\energyVec)^\transpose
       =(\mdens,\mdens\u,\widetilde{\mdens\u\u}+\pressureVec)^\transpose
\end{gather*}
denote conserved quantities,
and let
\begin{gather*}
  \prim:=(\mdens,\u,\pressureVec)^\transpose
        =(
          \mdens, \mom/\mdens, \energyVec-\widetilde{\mom\mom}/\mdens
        )^\transpose
\end{gather*}
denote primitive quantities,
where for an arbitrary symmetric tensor $\PT$ we define the tuple
$\widetilde\PT$ to be the six distinct components listed in the following
order:
\begin{gather*}
       \pressureVec :=
       \begin{bmatrix}
          \PT_{11}
       ,  \PT_{12}
       ,  \PT_{13}
       ,  \PT_{23}
       ,  \PT_{22}
       ,  \PT_{33}
       \end{bmatrix}^\transpose.
\end{gather*}
In the one-dimensional case, the balance law states
\begin{equation*}
  \cons_t+f(\cons)_x=s(\cons).
\end{equation*}
\subsubsection{General theory of state variable conversion}
In quasilinear form this reads
\begin{equation*}
  \cons_t+ f_\cons\cdot\cons_x=s.
\end{equation*}
Converting to primitive variables,
  $\cons_\prim\cdot\prim_t+ f_\cons\cdot\cons_\prim\cdot\prim_x=s,$
i.e., 
\begin{equation*}
  \prim_t+ (\prim_\cons\cdot f_\cons\cdot\cons_\prim)\cdot\prim_x=s.
\end{equation*}
\def\leftcons{\cons^\mathrm{L}}
\def\leftprim{\prim^\mathrm{L}}
\def\rghtcons{\cons^\mathrm{R}}
\def\rghtprim{\prim^\mathrm{R}}
Let $\leftcons$ and $\rghtcons$
denote conservative left and
right eigenvectors with eigenvalue $\lambda$:
\begin{equation*}
  f_\cons\cdot\rghtcons = \lambda\rghtcons,
 \\
  \leftcons\cdot f_\cons = \lambda\leftcons.
\end{equation*}
Then
\begin{equation*}
  \rghtprim:=\prim_\cons\cdot\rghtcons,
 \\
  \leftprim:=\leftcons\cdot\cons_\prim
\end{equation*}
are the corresponding primitve left and right eigenvectors
of the primitive-variable
wave propagation matrix $\prim_\cons\cdot f_\cons\cdot\cons_\prim$.
So we can calculate conservative eigenvectors from primitive
eigenvectors using the relations
\begin{equation*}
  \rghtcons=\cons_\prim\cdot\rghtprim,
 \\
  \leftcons=\leftprim\cdot\prim_\cons.
\end{equation*}

Observe that inner product is preserved under transformation
of state variables:
\begin{equation*}
 \leftprim\cdot\rghtprim
= 
 \leftcons\cdot\cons_\prim\cdot\prim_\cons\cdot\rghtcons
=
 \leftcons\cdot\rghtcons
\end{equation*}
\subsubsection{Derivative of state variable transformation}
To convert our primitive eigenvectors to conservative
eigenvectors, we need to calculate the derivative of 
the conserved variables with respect to the primitive
variables.   The conversion matrix for right eigenvectors is
\def\zerovec{\underline{0}}
\def\zeromat{\underline{\underline{0}}}
\def\zeromatvec{\widetilde{\underline{\underline{0}}}}
\begin{gather*}
   \cons_\prim =
      \partial
        \begin{bmatrix}
          \mdens\\ \mdens\u\\ \widetilde{\mdens\u\u}+\pressureVec
        \end{bmatrix}
      /\partial
        \begin{bmatrix}
          \mdens\\ \u\\ \pressureVec
        \end{bmatrix}
   = \begin{bmatrix}
       1 & \zerovec^\transpose & \zeromatvec^\transpose
    \\ \u & \mdens\idtens & \zerovec\zeromatvec^\transpose
    \\ \widetilde{\u\u} & \mdens(\widetilde{\u\u})_{\u} & \idtens
     \end{bmatrix},
\end{gather*}
(where $\idtens$ is the identity tensor), i.e.,
\begin{gather*}
   \cons_\prim =
    \partial
      \begin{bmatrix}
        \mdens
     \\ \mdens u_1
     \\ \mdens u_2
     \\ \mdens u_3
     \\ \mdens u_1 u_1 + \PT_{11}
     \\ \mdens u_1 u_2 + \PT_{12}
     \\ \mdens u_1 u_3 + \PT_{13}
     \\ \mdens u_2 u_3 + \PT_{23}
     \\ \mdens u_2 u_2 + \PT_{22}
     \\ \mdens u_3 u_3 + \PT_{33}
      \end{bmatrix}
    /
      \partial
      \begin{bmatrix}
        \mdens
     \\ u_1
     \\ u_2
     \\ u_3
     \\ \PT_{11}
     \\ \PT_{12}
     \\ \PT_{13}
     \\ \PT_{23}
     \\ \PT_{22}
     \\ \PT_{33}
      \end{bmatrix}
 =\begin{bmatrix}
     1&0&0&0&0&0&0&0&0&0
  \\ u_1&\mdens&0&0&0&0&0&0&0&0
  \\ u_2&0&\mdens&0&0&0&0&0&0&0
  \\ u_3&0&0&\mdens&0&0&0&0&0&0
  \\ u_1 u_1 & 2\mdens u_1&0&0&1&0&0&0&0&0
  \\ u_1 u_2 & \mdens u_2&\mdens u_1&0&0&1&0&0&0&0
  \\ u_1 u_3 & \mdens u_3&0&\mdens u_1&0&0&1&0&0&0
  \\ u_2 u_3 & 0&\mdens u_3&\mdens u_2&0&0&0&1&0&0
  \\ u_2 u_2 & 0&2\mdens u_2&0&0&0&0&0&1&0
  \\ u_3 u_3 & 0&0&2\mdens u_3&0&0&0&0&0&1
  \end{bmatrix}
\end{gather*}
The conversion matrix for left eigenvectors is slightly different:
\begin{gather*}
   \prim_\cons =
      \partial
        \begin{bmatrix}
          \mdens\\ \mom/\mdens\\ \energyVec-\widetilde{\mom\mom}/\mdens
        \end{bmatrix}
      /\partial
        \begin{bmatrix}
          \mdens\\ \mom\\ \energyVec
        \end{bmatrix}
   = \begin{bmatrix}
       1 & \zerovec^\transpose & \zeromatvec^\transpose
    \\ -\mom/\mdens^2 & \idtens/\mdens & \zerovec\zeromatvec^\transpose
    \\ \widetilde{\mom\mom}/\mdens^2 & -(\widetilde{\mom\mom})_{\mom}/\mdens & \idtens
     \end{bmatrix}
   = \begin{bmatrix}
       1 & \zerovec^\transpose & \zeromatvec^\transpose
    \\ -\u/\mdens & \idtens\mdens & \zerovec\zeromatvec^\transpose
    \\ \widetilde{\u\u} & -(\widetilde{\u\u})_{\u} & \idtens
     \end{bmatrix},
\end{gather*}
i.e.,
\begin{gather*}
 \prim_\cons = 
   \begin{bmatrix}
      1&0&0&0&0&0&0&0&0&0
   \\ -u_1/\mdens&1/\mdens&0&0&0&0&0&0&0&0
   \\ -u_2/\mdens&0&1/\mdens&0&0&0&0&0&0&0
   \\ -u_3/\mdens&0&0&1/\mdens&0&0&0&0&0&0
   \\ u_1 u_1 & -2u_1&0&0&1&0&0&0&0&0
   \\ u_1 u_2 & -u_2&-u_1&0&0&1&0&0&0&0
   \\ u_1 u_3 & -u_3&0&-u_1&0&0&1&0&0&0
   \\ u_2 u_3 & 0&-u_3&-u_2&0&0&0&1&0&0
   \\ u_2 u_2 & 0&-2u_2&0&0&0&0&0&1&0
   \\ u_3 u_3 & 0&0&-2u_3&0&0&0&0&0&1
   \end{bmatrix}
\end{gather*}


\subsection{Eigenvectors for ``conserved'' variables}

\subsubsection{Right eigenvectors for ``conserved'' variables}
\def\PrimRghtEigenvectors{P^\mathrm{R}}
\def\PrimLeftEigenvectors{P^\mathrm{L}}
\def\ConsRghtEigenvectors{Q^\mathrm{R}}
\def\ConsLeftEigenvectors{Q^\mathrm{L}}
Let $\PrimRghtEigenvectors$ denote a matrix of primitive
right eigenvectors.
We can compute conservative right eigenvectors from the
primitive right eigenvectors (without ever having to write
down the quasilinearized conservative system) by the relationship
\begin{gather*}
   \ConsRghtEigenvectors = \cons_\prim\cdot\PrimRghtEigenvectors,
\end{gather*}
where $\ConsRghtEigenvectors$ denotes a matrix of conservative
right eigenvectors.  So for the non-fast eigenvectors we have
\begin{align*}
   \ConsRghtEigenvectors &=
  \begin{bmatrix}
     1&0&0&0&0&0&0&0&0&0
  \\ u_1&\mdens&0&0&0&0&0&0&0&0
  \\ u_2&0&\mdens&0&0&0&0&0&0&0
  \\ u_3&0&0&\mdens&0&0&0&0&0&0
  \\ u_1 u_1 & 2\mdens u_1&0&0&1&0&0&0&0&0
  \\ u_1 u_2 & \mdens u_2&\mdens u_1&0&0&1&0&0&0&0
  \\ u_1 u_3 & \mdens u_3&0&\mdens u_1&0&0&1&0&0&0
  \\ u_2 u_3 & 0&\mdens u_3&\mdens u_2&0&0&0&1&0&0
  \\ u_2 u_2 & 0&2\mdens u_2&0&0&0&0&0&1&0
  \\ u_3 u_3 & 0&0&2\mdens u_3&0&0&0&0&0&1
  \end{bmatrix}
  \cdot
     \begin{bmatrix}
       \begin{matrix}
          0
       \\ 0
       \\ \pm c_s
       \\ 0
       \\ 0
       \\ \PT_{11}
       \\ 0
       \\ \PT_{13}
       \\ 2\PT_{12}
       \\ 0
       \end{matrix} \ \ 
       \begin{matrix}
          0
       \\ 0
       \\ 0
       \\ \pm c_s
       \\ 0
       \\ 0
       \\ \PT_{11}
       \\ \PT_{12}
       \\ 0
       \\ 2\PT_{13}
       \end{matrix} \ \ 
       \begin{matrix}
          1
       \\ 0
       \\ 0
       \\ 0
       \\ 0
       \\ 0
       \\ 0
       \\ 0
       \\ 0
       \\ 0
       \end{matrix} \ \ 
       \begin{matrix}
          0
       \\ 0
       \\ 0
       \\ 0
       \\ 0
       \\ 0
       \\ 0
       \\ 1
       \\ 0
       \\ 0
       \end{matrix} \ \ 
       \begin{matrix}
          0
       \\ 0
       \\ 0
       \\ 0
       \\ 0
       \\ 0
       \\ 0
       \\ 0
       \\ 1
       \\ 0
       \end{matrix} \ \ 
       \begin{matrix}
          0
       \\ 0
       \\ 0
       \\ 0
       \\ 0
       \\ 0
       \\ 0
       \\ 0
       \\ 0
       \\ 1
       \end{matrix}
     \end{bmatrix}
\end{align*}
That is,
\begin{align*}
  \ConsRghtEigenvectors &=
 \begin{bmatrix}
  \underbrace{
   \begin{bmatrix}
     \begin{matrix}
       0
    \\ 0
    \\ \mdens (\pm c_s)
    \\ 0
    \\ 0
    \\ \mdens u_1 (\pm c_s)
    \\ 0
    \\ \mdens u_3 (\pm c_s)
    \\ 2\mdens u_2 (\pm c_s)
    \\ 0
     \end{matrix}
     \begin{matrix}
       \phantom{+ 0}
    \\ \phantom{+ 0}
    \\ \phantom{+ 0}
    \\ \phantom{+ 0}
    \\ \phantom{+ 0}
    \\ + \PT_{11}
    \\ \phantom{0}
    \\ + \PT_{13}
    \\ + 2\PT_{12}
    \\ \phantom{+ 0}
     \end{matrix}
   \end{bmatrix},
   \begin{bmatrix}
     \begin{matrix}
       0
    \\ 0
    \\ 0
    \\ \mdens (\pm c_s)
    \\ 0
    \\ 0
    \\ \mdens u_1 (\pm c_s) 
    \\ \mdens u_2 (\pm c_s)
    \\ 0
    \\ 2\mdens u_3 (\pm c_s)
     \end{matrix}
     \begin{matrix}
       \phantom{0}
    \\ \phantom{0}
    \\ \phantom{0}
    \\ \phantom{0}
    \\ \phantom{0}
    \\ \phantom{0}
    \\ + \PT_{11}
    \\ + \PT_{12}
    \\ \phantom{0}
    \\ + 2\PT_{13}
     \end{matrix}
   \end{bmatrix}
  }_{\hbox{$c=\pm c_s:=\sqrt{\frac{\PT_{11}}{\mdens}}$}},
  \underbrace{
   \begin{bmatrix}
       1
    \\ u_1
    \\ u_2
    \\ u_3
    \\ u_1 u_1
    \\ u_1 u_2
    \\ u_1 u_3
    \\ u_2 u_3
    \\ u_2 u_2
    \\ u_3 u_3
   \end{bmatrix},
   \begin{bmatrix}
       0
    \\ 0
    \\ 0
    \\ 0
    \\ 0
    \\ 0
    \\ 0
    \\ 1
    \\ 0
    \\ 0
   \end{bmatrix},
   \begin{bmatrix}
       0
    \\ 0
    \\ 0
    \\ 0
    \\ 0
    \\ 0
    \\ 0
    \\ 0
    \\ 1
    \\ 0
   \end{bmatrix},
   \begin{bmatrix}
       0
    \\ 0
    \\ 0
    \\ 0
    \\ 0
    \\ 0
    \\ 0
    \\ 0
    \\ 0
    \\ 1
   \end{bmatrix}
  }_{\hbox{$c=0$}}
 \end{bmatrix}
\end{align*}

The fast right eigenvectors are more involved.
So we just give names to the primitive components
and multiply:
\def\dmdens{\mdens'}
\def\duA{u'_1}
\def\duB{u'_2}
\def\duC{u'_3}
\def\dPAA{\PT_{11}'}
\def\dPAB{\PT_{12}'}
\def\dPAC{\PT_{13}'}
\def\dPBC{\PT_{23}'}
\def\dPBB{\PT_{22}'}
\def\dPCC{\PT_{33}'}
\begin{gather*}
   \ConsRghtEigenvectors =
  \begin{bmatrix}
     1&0&0&0&0&0&0&0&0&0
  \\ u_1&\mdens&0&0&0&0&0&0&0&0
  \\ u_2&0&\mdens&0&0&0&0&0&0&0
  \\ u_3&0&0&\mdens&0&0&0&0&0&0
  \\ u_1 u_1 & 2\mdens u_1&0&0&1&0&0&0&0&0
  \\ u_1 u_2 & \mdens u_2&\mdens u_1&0&0&1&0&0&0&0
  \\ u_1 u_3 & \mdens u_3&0&\mdens u_1&0&0&1&0&0&0
  \\ u_2 u_3 & 0&\mdens u_3&\mdens u_2&0&0&0&1&0&0
  \\ u_2 u_2 & 0&2\mdens u_2&0&0&0&0&0&1&0
  \\ u_3 u_3 & 0&0&2\mdens u_3&0&0&0&0&0&1
  \end{bmatrix}
  \cdot
     \begin{bmatrix}
       \begin{matrix}
            \dmdens
         \\ \pm (\duA
         \\ \pm (\duB
         \\ \pm (\duC
         \\ \dPAA
         \\ \dPAB
         \\ \dPAC
         \\ \dPBC
         \\ \dPBB
         \\ \dPCC
       \end{matrix}
       \begin{matrix}
            :=
         \\ :=
         \\ :=
         \\ :=
         \\ :=
         \\ :=
         \\ :=
         \\ :=
         \\ :=
         \\ :=
       \end{matrix}
       \begin{matrix}
          \mdens\PT_{11}
       \\ c_f\PT_{11})
       \\ c_f\PT_{12})
       \\ c_f\PT_{13})
       \\ 3\PT_{11}\PT_{11}
       \\ 3\PT_{12}\PT_{11}
       \\ 3\PT_{13}\PT_{11}
       \\ \PT_{23}\PT_{11}+2\PT_{13}\PT_{12}
       \\ \PT_{22}\PT_{11}+2\PT_{12}\PT_{12}
       \\ \PT_{33}\PT_{11}+2\PT_{13}\PT_{13}
       \end{matrix}
     \end{bmatrix}.
\end{gather*}
So
{\small
\begin{gather*}
   \ConsRghtEigenvectors =
 \begin{bmatrix}
   \begin{matrix}
     0
  \\ \dmdens u_1
  \\ \dmdens u_2
  \\ \dmdens u_3
  \\ \dmdens u_1u_1&\pm\duA(2\mdens u_1)
  \\ \dmdens u_1u_2&\pm\duA\mdens u_2\pm\duB\mdens u_1
  \\ \dmdens u_1u_3&\pm\duA\mdens u_3\pm\duC\mdens u_1
  \\ \dmdens u_2u_3&\pm\duB\mdens u_3\pm\duC\mdens u_2
  \\ \dmdens u_2u_2&\pm\duB(2\mdens u_2)
  \\ \dmdens u_3u_3&\pm\duC(2\mdens u_3)
   \end{matrix}
   \begin{matrix}
     +
  \\ \pm
  \\ \pm
  \\ \pm
  \\ +
  \\ +
  \\ +
  \\ +
  \\ +
  \\ +
   \end{matrix}
   \begin{matrix}
        \dmdens
     \\ \duA\mdens
     \\ \duB\mdens
     \\ \duC\mdens
     \\ \dPAA
     \\ \dPAB
     \\ \dPAC
     \\ \dPBC
     \\ \dPBB
     \\ \dPCC
   \end{matrix}
 \end{bmatrix}
   =
 \mdens
 \begin{bmatrix}
       0
    \\ \PT_{11}u_1
    \\ \PT_{11}u_2
    \\ \PT_{11}u_3 
    \\ \PT_{11}u_1(u_1\pm2c_f)
    \\ \PT_{11}u_2(u_1\pm c_f) \pm c_f u_1\PT_{12}
    \\ \PT_{11}u_3(u_1\pm c_f) \pm c_f u_1\PT_{13}
    \\ \PT_{11}u_2 u_3\pm c_f u_3\PT_{12}\pm c_f u_2\PT_{13}
    \\ \PT_{11}u_2 u_2 \pm c_f u_22\PT_{12}
    \\ \PT_{11}u_3 u_3 \pm c_f u_32\PT_{13}
  \end{bmatrix}
  +
  \begin{bmatrix}
          \mdens\PT_{11}
       \\ \pm c_f\PT_{11}\mdens
       \\ \pm c_f\PT_{12}\mdens
       \\ \pm c_f\PT_{13}\mdens
       \\ 3\PT_{11}\PT_{11}
       \\ 3\PT_{12}\PT_{11}
       \\ 3\PT_{13}\PT_{11}
       \\ \PT_{23}\PT_{11}+2\PT_{13}\PT_{12}
       \\ \PT_{22}\PT_{11}+2\PT_{12}\PT_{12}
       \\ \PT_{33}\PT_{11}+2\PT_{13}\PT_{13}
  \end{bmatrix}.
\end{gather*}
}

\subsubsection{Left eigenvectors for conserved variables}
Similarly,
let $\PrimLeftEigenvectors$ denote a matrix of primitive
left eigenvectors.
We can compute conservative left eigenvectors from the
primitive left eigenvectors by the relationship
\begin{gather*}
   \ConsLeftEigenvectors = \PrimLeftEigenvectors\cdot\prim_\cons,
\end{gather*}
where $\ConsLeftEigenvectors$ denotes a matrix of conservative
left eigenvectors.

To avoid big expressions, we
give a simple name to each nonzero matrix component before
multiplying the matrices.

For the slow eigenvector pair for $\PT_{12}$ and $\PT_{13}$
define
\def\du{u'}
\def\dP{\PT'}
\begin{align*}
   \du &:= \PT_{11}\PT_{11}
\\ \dP &:= c_s\PT_{11},
\end{align*}
and define
\def\duAb{{u_1'}_b}
\def\duBb{\du}
\def\duBbFull{{u_2'}_b}
\def\dPAAb{{\PT_{11}'}_b}
\def\dPABb{\dP}
\def\dPABbFull{ {\PT_{12}'}_b}
\def\duAc{{u_1'}_c}
\def\duCc{\du}
\def\duCcFull{ {u_3'}_c}
\def\dPAAc{{\PT_{11}'}_c}
\def\dPACc{\dP}
\def\dPACcFull{ {\PT_{13}'}_c}
\begin{align*}
   \duAb &:= -\PT_{12}\PT_{11},
  &\duAc &:= -\PT_{13}\PT_{11},
\\ \duBbFull &:= \duBb,
  &\duCcFull &:= \duCc,
\\ \dPAAb &:= -c_s\PT_{12},
  &\dPAAc &:= -c_s\PT_{13},
\\ \dPABbFull &:= \dPABb,
  &\dPACcFull &:= \dPACc.
\end{align*}

So in terms of these quantities the left eigenvectors are
\begin{align*}
\ConsLeftEigenvectors &=
     \begin{bmatrix}
       \begin{bmatrix}
          0
       \\ 3\PT_{11}
       \\ 0
       \\ 0
       \\ \pm c_f
       \\ 0
       \\ 0
       \\ 0
       \\ 0
       \\ 0
       \end{bmatrix}
       \begin{bmatrix}
          0
       \\ \duAb
       \\ \duBb
       \\ 0
       \\ \pm\dPAAb
       \\ \pm\dPABb
       \\ 0
       \\ 0
       \\ 0
       \\ 0
       \end{bmatrix},
       \begin{bmatrix}
          0
       \\ \duAc
       \\ 0
       \\ \duCc
       \\ \pm\dPAAc
       \\ 0
       \\ \pm\dPACc
       \\ 0
       \\ 0
       \\ 0
       \end{bmatrix}
     \end{bmatrix}^\transpose
   \cdot
   \begin{bmatrix}
      1&0&0&0&0&0&0&0&0&0
   \\ -u_1/\mdens&1/\mdens&0&0&0&0&0&0&0&0
   \\ -u_2/\mdens&0&1/\mdens&0&0&0&0&0&0&0
   \\ -u_3/\mdens&0&0&1/\mdens&0&0&0&0&0&0
   \\ u_1 u_1 & -2u_1&0&0&1&0&0&0&0&0
   \\ u_1 u_2 & -u_2&-u_1&0&0&1&0&0&0&0
   \\ u_1 u_3 & -u_3&0&-u_1&0&0&1&0&0&0
   \\ u_2 u_3 & 0&-u_3&-u_2&0&0&0&1&0&0
   \\ u_2 u_2 & 0&-2u_2&0&0&0&0&0&1&0
   \\ u_3 u_3 & 0&0&-2u_3&0&0&0&0&0&1
   \end{bmatrix}.
\end{align*}
So
\begin{equation*}
\!\!\!\!\!\!
\ConsLeftEigenvectors =
     \begin{bmatrix}
       \begin{bmatrix}
          \hbox{\small $-3\frac{\PT_{11}}{\mdens}u_1 \pm c_f u_1 u_1$}
       \\ 3\frac{\PT_{11}}{\mdens} \mp 2c_f u_1
       \\ 0
       \\ 0
       \\ \pm c_f
       \\ 0
       \\ 0
       \\ 0
       \\ 0
       \\ 0
       \end{bmatrix},
       \begin{bmatrix}
         \hbox{\small $
          \frac{-\duAb u_1 - \duBb u_2}{\mdens}
             \pm u_1 u_1 \dPAAb
             \pm u_1 u_2 \dPABb $}
       \\ \duAb/\mdens \mp 2 u_1\dPAAb \mp u_2 \dPABb
       \\ \duBb/\mdens \mp u_1 \dPABb
       \\ 0
       \\ \pm \dPAAb
       \\ \pm \dPABb
       \\ 0
       \\ 0
       \\ 0
       \\ 0
       \end{bmatrix},
       \begin{bmatrix}
         \hbox{\small $
          \frac{-\duAc u_1 - \duCc u_3}{\mdens}
             \pm u_1 u_1 \dPAAc
             \pm u_1 u_3 \dPACc $}
       \\ \duAc/\mdens \mp 2 u_1\dPAAc \mp u_3\dPACc
       \\ 0
       \\ \duCc/\mdens \mp u_1\dPACc
       \\ \pm\dPAAc
       \\ 0
       \\ \pm\dPACc
       \\ 0
       \\ 0
       \\ 0
       \end{bmatrix}
     \end{bmatrix}^\transpose
\end{equation*}

For the $c=0$ eigenvectors define
\def\dPAAe{{\PT_{11}'}_e}
\def\dPABe{{\PT_{12}'}_e}
\def\dPACe{{\PT_{13}'}_e}
\def\dPBCe{3\PT_{11}^2} 
\begin{align*}
   \dPAAe &= {4\PT_{12}\PT_{13}
              -\PT_{32}\PT_{11}},
\\ \dPABe &:= {-3\PT_{11}\PT_{13}},
\\ \dPACe &:= {-3\PT_{11}\PT_{12}},
\end{align*}
and define
\def\dPAAf{{\PT_{11}'}_f}
\def\dPABf{{\PT_{12}'}_f}
\def\dPBBf{3\PT_{11}^2} 
\def\dPAAg{{\PT_{11}'}_g}
\def\dPACg{{\PT_{13}'}_g}
\def\dPCCg{3\PT_{11}^2} 
\begin{align*}
   \dPAAf &:= 4\PT_{12}^2-\PT_{11}\PT_{22},
  &\dPAAg &:= 4\PT_{13}^2-\PT_{11}\PT_{33},
\\ \dPABf &:= -6\PT_{12}\PT_{11},
  &\dPACg &:= -6\PT_{13}\PT_{11}.
\end{align*}

The left eigenvectors for $c=0$ are given by
{\small
\begin{align*}
\ConsLeftEigenvectors &=
     \begin{bmatrix}
       \begin{bmatrix}
          3\PT_{11}
       \\ 0
       \\ 0
       \\ 0
       \\ -\mdens
       \\ 0
       \\ 0
       \\ 0
       \\ 0
       \\ 0
       \end{bmatrix},
       \begin{bmatrix}
          0
       \\ 0
       \\ 0
       \\ 0
       \\ \dPAAe
       \\ \dPABe
       \\ \dPACe
       \\ \dPBCe
       \\ 0
       \\ 0
       \end{bmatrix},
       \begin{bmatrix}
          0
       \\ 0
       \\ 0
       \\ 0
       \\ \dPAAf
       \\ \dPABf
       \\ 0
       \\ 0
       \\ \dPBBf
       \\ 0
       \end{bmatrix},
       \begin{bmatrix}
          0
       \\ 0
       \\ 0
       \\ 0
       \\ \dPAAg
       \\ 0
       \\ \dPACg
       \\ 0
       \\ 0
       \\ \dPCCg
       \end{bmatrix}
     \end{bmatrix}^\transpose
   \cdot
   \begin{bmatrix}
      1&0&0&0&0&0&0&0&0&0
   \\ -u_1/\mdens&1/\mdens&0&0&0&0&0&0&0&0
   \\ -u_2/\mdens&0&1/\mdens&0&0&0&0&0&0&0
   \\ -u_3/\mdens&0&0&1/\mdens&0&0&0&0&0&0
   \\ u_1 u_1 & -2u_1&0&0&1&0&0&0&0&0
   \\ u_1 u_2 & -u_2&-u_1&0&0&1&0&0&0&0
   \\ u_1 u_3 & -u_3&0&-u_1&0&0&1&0&0&0
   \\ u_2 u_3 & 0&-u_3&-u_2&0&0&0&1&0&0
   \\ u_2 u_2 & 0&-2u_2&0&0&0&0&0&1&0
   \\ u_3 u_3 & 0&0&-2u_3&0&0&0&0&0&1
   \end{bmatrix}.
\end{align*}
}
So the $c=0$ left eigenvectors are
\begin{align*}
  \leftcons &=
       \begin{bmatrix}
          3\PT_{11} - \mdens u_1 u_1
       \\ 2 \mdens u_1
       \\ 0
       \\ 0
       \\ -\mdens
       \\ 0
       \\ 0
       \\ 0
       \\ 0
       \\ 0
       \end{bmatrix},
       \begin{bmatrix}
          u_1 u_1 \dPAAe + u_1 u_2\dPABe+u_1u_3\dPACe+3u_2u_3\PT_{11}^2
       \\ -(2 u_1 \dPAAe + u_2 \dPABe + u_3 \dPACe)
       \\ -(u_1 \dPABe + 3 u_3 \PT_{11}^2)
       \\ -(u_1 \dPACe + 3 u_2 \PT_{11}^2)
       \\ \dPAAe
       \\ \dPABe
       \\ \dPACe
       \\ \dPBCe
       \\ 0
       \\ 0
       \end{bmatrix},
  \\ &
       \begin{bmatrix}
          u_1 u_1 \dPAAf + u_1 u_2 \dPABf + u_2 u_2 \dPBBf
       \\ -(2u_1\dPAAf + u_2\dPABf)
       \\ -(u_1\dPABf + 6u_2\PT_{11}^2)
       \\ 0
       \\ \dPAAf
       \\ \dPABf
       \\ 0
       \\ 0
       \\ \dPBBf
       \\ 0
       \end{bmatrix},
       \begin{bmatrix}
          u_1 u_1 \dPAAg + u_1 u_3 \dPACg + u_3 u_3 \dPCCg
       \\ -(2u_1\dPAAg + u_3\dPACg)
       \\ 0
       \\ -(u_1\dPACg + 6 u_3 \PT_{11}^2)
       \\ \dPAAg
       \\ 0
       \\ \dPACg
       \\ 0
       \\ 0
       \\ \dPCCg
       \end{bmatrix}.
\end{align*}

}